\documentclass[letterpaper,11pt]{amsart}
\usepackage[margin=1.2in]{geometry}
\usepackage{eucal}
\usepackage{tikz}
\usepackage{xypic}
\usepackage{amsfonts,amssymb}
\usepackage{epsfig}
\usepackage{enumerate}
\usepackage{enumitem}
\usetikzlibrary{decorations.pathreplacing}

\newtheorem{theorem}{Theorem}[section]
\newtheorem{lemma}[theorem]{Lemma}
\newtheorem{prop}[theorem]{Proposition}

\newtheorem{corollary}[theorem]{Corollary}
\newtheorem{defn}[theorem]{Definition}

\newtheorem{conjecture}[theorem]{Conjecture}

\newcommand{\Q}{\mathbb Q}
\newcommand{\C}{\mathbb C}
\newcommand{\R}{\mathbb R}

\renewcommand{\P}{\mathbb P}
\renewcommand{\O}{\EuScript O}

\newcommand{\D}{\mathbb D}

\newcommand{\Z}{\mathbb Z}
\newcommand{\N}{\mathbb N}

\newlist{HF}{enumerate}{1}
\setlist[HF]{label=(HF\arabic*)}

\newlist{PB}{enumerate}{1}
\setlist[PB]{label=(PB\arabic*)}

\newlist{RS}{enumerate}{1}
\setlist[RS]{label=(PB\arabic*)}

\newlist{CZ}{enumerate}{1}
\setlist[CZ]{label=(CZ\arabic*)}

\title{Reeb orbits and the minimal discrepancy of an isolated singularity}
\author{Mark McLean}

\begin{document}

\begin{abstract}
Let $A$ be an affine variety inside a complex $N$ dimensional
vector space which has an isolated singularity at the origin.
The intersection of $A$ with a very small
sphere turns out to be a contact manifold called the link of $A$.
Any contact manifold contactomorphic to the link of $A$ is said to be
Milnor fillable by $A$.
If the first Chern class of our link is torsion then we can assign an
invariant of our singularity called the minimal discrepancy, which
is an important invariant in birational geometry.
We define an invariant of the link up to contactomorphism using
Conley-Zehnder indices of Reeb orbits and then we relate this invariant with
the minimal discrepancy.
As a result we show that the standard contact $5$ dimensional sphere has
a unique Milnor filling up to normalization proving a conjecture by Seidel.
\end{abstract}

\maketitle

\bibliographystyle{halpha}

%-----------------------------------------------------------------
%-----------------------------------------------------------------
%-----------------------------------------------------------------

\setcounter{tocdepth}{1}
\tableofcontents

\section{Introduction}

Suppose we have an irreducible affine variety $A \subset \C^N$ of complex dimension $n$
which has an isolated singularity at $0$ (here we include the case that $A$ might be smooth at $0$).
For any $\epsilon>0$ small enough we have that $L_A := A \cap \{ \sum_{i=1}^N |z_i|^2 = \epsilon^2\}$ is a differentiable
manifold of real dimension $2n-1$ and such a manifold is an invariant of the germ of $A$
at $0$.
We call $L_A$ {\bf the link of $A$}.
Near $0$, we have that $A$ is homeomorphic to the cone over $L_A$.
The simplest example is when $A$ is smooth at $0$ in which case $L_A$ is diffeomorphic to a sphere.
Many people have studied the relationship between the algebraic properties of $A$ at $0$
and the topology of $L_A$.
Such results go back to \cite{Heegaard:thesis}.
% and \cite{Wirtinger:1895}.
There have been particularly powerful results when $\text{dim}_\C A = 2$ but
there are far less powerful results in higher dimensions.
For instance, let's start with the following definition inspired by \cite{Heegaard:thesis}[Page 236 (French Translation)]:
A singularity is {\bf topologically smooth} if its link is diffeomorphic to a sphere.
Mumford in \cite{Mumford:topologyofnormalsingularities} showed that
every normal topologically smooth singularity of complex dimension $2$ is in fact smooth.
%There are examples of non-normal isolated singularities whose link is $S^3$ but if one normalizes them,
%then the result by Mumford says that such a normalization is smooth.
But in complex dimension $3$ or higher there are many examples of isolated
normal singularities which are topologically smooth, but not smooth at $0$
such as $\{x^2 + y^2 + z^2 + w^3 = 0\}$
(see 
\cite{Brieskorn:Examples}, \cite{Brieskorn:sphere},
\cite{Hirzebruch:Singularitiesexotic}
 and
\cite{Brieskorn:SingularitiesHirzebruch}).
Normality of this singularity follows from
Serre's criterion for normality \cite[Theorem 18.15]{EisenbudCommutativeAlgebraViewTowards}.

Having said this, one can put additional structure on the link.
Let $J_0 : T\C^N \to T\C^N$ be the standard complex structure on $\C^N$
viewed as an automorphism of the real tangent bundle whose square is $-\text{id}$.
Then $\xi_A := TL_A \cap J_0(TL_A) \subset TL_A$ is a contact structure for $\epsilon$ small enough
and $(L_A,\xi_A)$ is an invariant of the germ of $A$ at $0$ up to contactomorphism
(see \cite{Varchenko:isolatedsingularities}).
One can also view $\xi_A$ as the kernel of $\sum_{j=1}^N x_j dy_j - y_j dx_j|_{L_A}$ where $z_j = x_j + i y_j$ are coordinates for $\C^N$.
Following \cite{CCNA:Milnorfillable3mdlfs},
a contact manifold $(C,\xi)$ is said to be {\bf Milnor fillable}
if it is contactomorphic to $(L_A,\xi_A)$ for some $A$.
Here $A$ is called a {\bf Milnor filling} of $(C,\xi)$.
%It is {\it Milnor fillable} if it is Milnor fillable by $A$ for some $A$.
An example of a Milnor fillable contact structure is the standard contact sphere
$(S^{2n-1},\xi_{\text{std}})$ which is defined to be the link of $\C^n$
(i.e. $S^{2n-1}$ is the unit sphere in $\C^n$ and $\xi_{\text{std}}$ is the unique
hyperplane distribution which is $J_0$ invariant).

In \cite{Ustilovsky:infinitecontact} it was shown, for each $m > 0$, that there are infinitely many examples of isolated singularities
whose links are diffeomorphic to $S^{4m+1}$, but not contactomorphic to each other.
Hence $(L_A,\xi_A)$ is a stronger invariant than $L_A$ on its own.
Building on the work of \cite{Ustilovsky:infinitecontact}, \cite{kwonkoert} systematically investigated the links of
weighted homogenous hypersurface singularities $\{\sum_j z_j^{k_j} = 0\}$.
In particular using Conley-Zehnder indices of Reeb orbits, \cite[Theorem 6.3]{kwonkoert} (along with its proof) tells us whether $\sum_j 1/k_j>1$ just from $(L_A,\xi_A)$.
Such a result is significant because Reid in \cite[Proposition 4.3]{Reid:canonical3folds} showed that such a singularity
is canonical at $0$ if and only if $\sum_j 1/k_j>1$
(see \cite[Section 1]{Reid:canonical3folds} for a definition of canonical singularity).

For certain singularities called {\it $\Q$-Gorenstein} singularities, one can define an invariant taking values in $\Q$
called the  {\bf minimal discrepancy}
(see \cite{Ambro:minimaldiscrepancysurvey}).
We write $\text{md}(A,0)$ for the minimal discrepancy of $A$ at $0$.
All isolated complete intersection singularities of complex dimension $2$ or higher are $\Q$-Gorenstein,
and it turns out that canonical singularities are defined as $\Q$-Gorenstein singularities with
non-negative minimal discrepancy. Hence there is a direct relationship between the result in
\cite{kwonkoert} mentioned earlier and minimal discrepancy.
Minimal discrepancy can be defined for a larger class of singularities called numerically $\Q$-Gorenstein singularities
in \cite{BFFU:valuation}.
An isolated singularity is {\bf numerically $\Q$-Gorenstein} if $c_1(TA|_{L_A}) = c_1(\xi_A;\Q)$ is torsion in $H^2(L_A,\Z)$
(for related definitions, including the proof that
$c_1(TA|_{L_A}) = c_1(\xi_A;\Q)$, see Section \ref{section:minimaldiscrepancydefinition}).
%In Section \ref{section:minimaldiscrepancydefinition}, we will give both an algebraic and topological definition of minimal discrepancy
%when $H^1(L_A;\R) = 0$.
Singularities with positive minimal discrepancy are called {\it terminal singularities}.
These have special importance in the minimal model program (\cite{KollarMori:birational}).
See \cite[Corollary 5.17]{BFFU:valuation} for a proof that positive minimal discrepancy is equivalent to being terminal.
In fact minimal discrepancy itself has a special importance in the minimal model program
(\cite{Shokurov:ProblemsAboutFanoVarieties}).

Now let $(C,\xi)$ be a cooriented contact manifold of dimension $2n-1$ with $H^1(C;\Q) = 0$ and $c_1(\xi;\Q) = 0 \in H^2(C,\Q)$.
Suppose that $\alpha$ is a contact form with $\text{ker}(\alpha) = \xi$ respecting this coorientation.
To any Reeb orbit $\gamma : \R / L\Z \to C$ of $\alpha$, we have an associated index
$\text{CZ}(\gamma) \in \Q$ called the {\it Conley-Zehnder index}.
One should think of the Conley-Zehnder index as a measure of how much the
nearby Reeb flow lines `wrap' around our Reeb orbit $\gamma$.
This index will be defined in Section \ref{section:conleyzehnderindex}
(we also define it in some cases in the introductory section
\ref{subsection:basicnotions}).
Let $\phi_t : C \to C$, $t \in \R$ be the Reeb flow of $\alpha$.
The differential $D\phi_t : TC \to TC$ of the Reeb flow preserves $\xi$
and so for $p \in C$ let $D_p\phi_t|_\xi : \xi|_{p} \to \xi|_{\phi_t(p)}$
be the restriction of this differential to the contact distribution.
Because the Reeb orbit $\gamma$ has period $L$,
our linearized map $D_{\gamma(0)}\phi_L|_\xi$ sends
$\xi|_{\gamma(0)}$ to itself, because $\xi|_{\gamma(L)} = \xi|_{\gamma(0)}$.
Such a map is called the {\bf linearized return map of $\gamma$}.
We define the {\bf lower SFT index} $\text{lSFT}(\gamma)$
to be 
\[\text{CZ}(\gamma) - \frac{1}{2}\text{dim ker}(D_{\gamma(0)}\phi_L|_\xi - \text{id}) + (n-3).\]
For any contact form $\alpha$
we define the {\bf minimal SFT index of }$\alpha$ to be $\text{mi}(\alpha) := \text{inf}_\gamma \text{lSFT}(\gamma)$
where the infimum is taken over all Reeb orbits $\gamma$ of $\alpha$.
The {\bf highest minimal SFT index} of $(C,\xi)$ is defined as
$\text{hmi}(C,\xi) := \text{sup}_{\alpha} \text{mi}(\gamma)$
where the supremum is taken over all contact forms $\alpha$  with $\text{ker}(\alpha) = \xi$ respecting the coorientation of $\xi$.
%%%%
This is an invariant up to coorientation preserving contactomorphism.
%In practice we will only be able to calculate $\text{sup}_{\alpha} \text{inf}_{\gamma} (\text{lSFT}(\gamma))$
%in the case when the contact forms $\alpha$ respect a chosen coorientation of $\xi$.
Having said that if $(C,\xi)$ admits a coorientation reversing contactomorphism then
$\text{hmi}(C,\xi)$ is in fact a (not necessarily coorientation preserving) contactomorphism invariant.
This is the case for links of singularities $(L_A,\xi_A)$ because $z \to \overline{z}$
is a coorientation reversing contactomorphism.

Our main theorem is:
\begin{theorem} \label{theorem:minimaldiscrepancy}
Let $A$ have a normal isolated singularity at $0$ that is numerically
$\Q$-Gorenstein with $H^1(L_A;\Q) = 0$ then:
\begin{itemize}
\item If $\text{md}(A,0) \geq 0$ then $\text{hmi}(L_A,\xi_A) = 2\text{md}(A,0)$.
\item If $\text{md}(A,0) < 0$ then $\text{hmi}(L_A,\xi_A) < 0$.
\end{itemize}
\end{theorem}
This Theorem will follow from Theorems \ref{theorem:boundingdiscrepancyfromabove} and \ref{theorem:lowerbound}.
%As described earlier, for weighted homogenous hypersurface singularities, it was shown in \cite[Theorem 1.2]{kwonkoert} and \cite[Section 1]{Reid:canonical3folds} that the link remembers whether minimal discrepancy is negative or not.
Our main theorem works for any normal isolated singularity, even if it cannot be smoothed.
We have the following corollary, proving a conjecture by Seidel \cite[Lecture 6]{seidel:symplecticcohomologylecture}.
\begin{corollary} \label{corollary:fillingissmooth}
Suppose that $A$ is normal and that $(L_A,\xi)$ is contactomorphic to the link of $\C^3$
(i.e. the standard contact sphere $(S^5,\xi_{\text{std}})$),
then $A$ is smooth at $0$.
\end{corollary}
The above corollary says that the standard contact $5$ dimensional sphere has a unique Milnor filling up to normalization.
This generalizes the previously stated theorem by Mumford because
the three sphere has a unique strongly fillable contact structure
(see \cite{Eliashberg:fillingbyholomorphicdisks}, \cite{Gromov:Pseudoholomoprhic})
and because Milnor fillable contact structures are strongly fillable
by resolving the singularity (Lemma \ref{lemma:symplecticformonresolution}).
In fact every oriented 3-manifold admits at most one Milnor fillable
contact structure up to orientation preserving diffeomorphism
\cite{CCNA:Milnorfillable3mdlfs}.

The above corollary is a direct consequence of the following conjecture by Shokurov proven in complex dimension $3$
by
\cite[Main Theorem (I)]{Reid:minimalcanonical3fold}
combined with minimal discrepancy calculations from \cite{Markushevich:cDV} and \cite{Kawamata:logflipsappendix}.
\begin{conjecture} \label{conjectuer:shokurov}
(Shokurov \cite[Conjecture 2]{Shokurov:birationalist4}).

Suppose $A$ is normal and numerically $\Q$-Gorenstein with
$\text{md}(A,0) = n-1$ then $A$ is smooth at $0$.
\end{conjecture}

Shokurov has the stronger condition that $A$ is $\Q$-Gorenstein, but
\cite[Corollary 5.17]{BFFU:valuation} ensures that any numerically $\Q$-Gorenstein singularity
with minimal discrepancy $>-1$ is in fact $\Q$-Gorenstein.

As a result we have the following corollary:
\begin{corollary}
Assuming that Conjecture \ref{conjectuer:shokurov} is true,
$A$ is normal and $(L_A,\xi)$ is contactomorphic to the standard contact sphere of any dimension greater than $1$
then $A$ is smooth at $0$.
\end{corollary}
In other words, Shokurov's conjecture combined with Theorem \ref{theorem:minimaldiscrepancy}
implies that the standard contact sphere has a unique Milnor filling up to normalization.

We also have the following corollary:
\begin{corollary}
Assume that
$A$ has an isolated canonical singularity with $H^1(L_A;Q) = 0$, and let $B$
be any other normal isolated singularity whose link is contactomorphic to
that of $A$. Then $B$ also has canonical singularities, and $\text{md}(A) = \text{md}(B)$.
In particular, $A$ is terminal iff $B$ is.
\end{corollary}

We will now wildly speculate on the relationship between the main result of this paper and other results concerning the arc space.
The arc space was introduced by Nash in \cite{Nash:arcpaper}.
Let $\text{Arc}(A)$ be the space of formal disks
$\text{Hom}(\text{Spec}~\C[[z]],A)$
and $\text{Arc}(A,0) \subset \text{Arc}(A)$ the subspace of such disks passing through our singularity.
Very roughly, \cite[Theorem 2.6]{einmustatayasuda:discrepancies} relates the codimension of
$\text{Arc}(A,0)$ inside $\text{Arc}(A)$ with the minimal discrepancy.
Imagine that, as a disk in $\text{Arc}(A,0)$ approaches the origin, it converges
to some Reeb orbit, and that the component of $\text{Arc}(A,0)$ of highest codimension finds the lowest
index Reeb orbit.
%In some sense we bound $\text{md}(A,0)$ from below by using a certain space of holomorphic curves
%which `intersects' this lowest codimension component but not any other component as it has too high codimension.
Hence one can ask, what is the relationship between the space of pseudo holomorphic curves on the symplectization
of this link (such as those curves encoded by symplectic field theory) and the arc space? The arc space
may not quite be the right space to study, instead it might be the space of {\it short arcs}
defined in \cite{KollarNemethi:holomorphicarcs}.

We can also ask other questions. For instance minimal discrepancy is also defined for non-isolated singularities
and more generally for log pairs, and it would be interesting to see if there is some way of characterizing the minimal discrepancy
of such objects using contact geometry.
%Note that \cite[Theorem 2.6]{einmustatayasuda:discrepancies} applies to any $\Q$-Gorenstein variety
%and so this suggests that there might be some way of studying more general log pairs.

This paper is structured as follows:
Section \ref{section:basicnotionsandsketch} is split into two subsections,
Subsection \ref{subsection:basicnotions} and
Subsection \ref{subsection:sketchofproof}.
In Subsection \ref{subsection:basicnotions}
we give some basic notions from symplectic and contact topology.
We also define the Conley-Zehnder index in a restricted situation
for illustrative purposes.
In Subsection \ref{subsection:sketchofproof}, we give a sketch of the proof
in the specific case when $A$ is a cone singularity.
We hope that this section will explain the tight connection
between $\text{md}(A,0)$ and $\text{hmi}(L_A,\xi_A)$.
We also introduce some of the tools used in the proof
such as Gromov-Witten theory (in a very specific case)
and also neck stretching which is used to find Reeb orbits.

In Section \ref{section:minimaldiscrepancydefinition} we give two definitions of numerically $\Q$-Gorenstein
singularities, an algebraic and a topological one, then we prove their equivalence.
We then define the minimal discrepancy.
In Section \ref{section:conleyzehnderindexsection} we define the Conley-Zehnder index
of a Reeb orbit
and then we relate the $\text{lSFT}$ indices
of degenerate orbits with $\text{lSFT}$ indices of non-degenerate orbits coming from perturbations of these degenerate orbits.
In Section \ref{section:neighbourhoods} we show that a resolution of $A$
admits a nice symplectic structure and the boundary of a neighborhood of
the exceptional divisors is a contact manifold contactomorphic to $(L_A,\xi_A)$
and admitting a contact form with nice families of Reeb orbits.
As a result we prove the inequality $\text{hmi}(L_A,\xi_A) \geq 2\text{md}(A,0)$.
In Section \ref{section:gromovwittenopen} we show how to define genus $0$ Gromov-Witten invariants for certain open symplectic manifolds
(i.e. we count compact curves in some compact subset of such manifolds).
We also prove some important technical Lemmas involving these open manifolds.
In Section \ref{section:discrepancyfrombelow}  we use results from the previous section to show
$\text{hmi}(L_A,\xi_A) \leq 2\text{md}(A,0)$ if $\text{md}(A,0) \geq 0$
and $\text{hmi}(L_A,\xi_A) < 0$ if $\text{md}(A,0) < 0$.
This is done by partially compactifying some resolution of $A$ and then using Gromov-Witten invariants
along with a neck stretching argument to find Reeb orbits of the appropriate Conley-Zehnder index.
Appendix A reviews neck stretching and proves a compactness result when the contact structure is degenerate.
Appendix B proves a maximum principle for stable Hamiltonian structures which is a key argument
enabling us to show that we can define Gromov-Witten invariants for some partial compactification of a resolution of $A$.

\bigskip

{\it Acknowledgements}: I would like to thank Chris Wendl, Ivan Smith, Cheuk Yu Mak and Paul Seidel
(who enabled me to generalize the result from numerically Gorenstein to numerically $\Q$-Gorenstein singularities).
I would also like to thank the referees for numerous helpful comments which have improved this
paper.
The main work for this paper was done while I was at the University of Aberdeen and therefore I
would like to thank the people at the mathematics department for providing a great working environment.

\section{Basic Notions and a Sketch of the Proof} \label{section:basicnotionsandsketch}

\subsection{Basic Notions} \label{subsection:basicnotions}

The purpose of this section is to give some basic definitions
from contact and symplectic topology for non-experts,
and also to set up notational conventions in this paper.
A good introduction of the basics of contact/symplectic topology is found in \cite{McduffSalamon:sympbook}.
%Here we give some of the main definitions used in the statement
%and proof of the main theorem in this paper.
%Let $\Omega$ be a linear $2$-form on a finite dimensional $\R$ vector space $V$.
%The {\bf kernal} of $\Omega$ is the set of $v \in V$ satisfying $\Omega(v,w) = 0 \quad \forall w \in V$.
%A $2$-form $\Omega$ on a vector space $V$ is {\bf non-degenerate}
%if the map $V \to V^*$ given by $v \to \omega(v,\cdot)$ is invertible
%(or equivalently, the kernal is trivial).
%One can show that $\Omega$ is non-degenerate if and only if $\omega^{\text{dim}(V)} \neq %0$.
%Let $C$ be a smooth manifold.
%Recall the the Frobenius  integrability theorem tells us that a hyperplane distribution $\xi \subset TC$
%inside a manifold is the tangent space to a foliation if and only if for any pair of vector fields
%tangent to $\xi$, their Lie bracket is also tangent to $\xi$.
%If $\xi = \text{ker}(\alpha)$ for some $1$-form $\alpha$ on $C$,
%then this theorem tells us that
%$\xi$ is the tangent space to a foliation if and only if
%$d\alpha|_\xi = 0$.
%A contact structure is the opposite of this:
\begin{defn}
Suppose that $C$ is a manifold of dimension $2n - 1$.
A {\bf contact structure} is a hyperplane distribution $\xi \subset TC$
which is locally equal to $\text{ker}(\alpha)$ for some $1$-form $\alpha$
where $\alpha \wedge (d\alpha)^{n-1} \neq 0$.
The $1$-form $\alpha$ is called a
{\bf contact form associated to $\xi$}.
The pair $(C,\xi)$ is called a {\bf contact manifold}
if $\xi$ is a contact structure.
A {\bf cooriented contact manifold} $(C,\xi)$ is a contact manifold so that
$\xi = \text{ker}(\alpha)$ for some global $1$-form $\alpha$.
The $1$-form $\alpha$ induced an orientation on the bundle $TX / \xi$
which we call a {\bf coorientation of $\xi$}.
Two contact manifolds $(C_1,\xi_1),(C_2,\xi_2)$ are {\bf contactomorphic}
if there is a diffeomorphism $\phi : C_1 \to C_2$ sending $\xi_1$ to $\xi_2$.
The diffeomorphism $\phi$ is called a {\bf contactomorphism}.
\end{defn}

For any two contact forms $\alpha$ and $\beta$ associated to $\xi$,
there is a smooth function $f : C \to \R \setminus \{0\}$
where $\alpha = f \beta$.
From now on we will assume that every contact manifold is cooriented,
and we will study such manifolds up to coorientation preserving contactomorphism unless stated otherwise.
Having said that our main invariant $\text{hmi}(\xi)$ will be invariant up to general (not necessarily
orientation preserving) contactomorphism for links of singularities.
Unless stated otherwise, we will assume that any contact form associated to $\xi$
respects the chosen coorientation.
%One can also think of a contact structure as an odd dimensional version of a symplectic manifold.
%Recall that a symplectic manifold $M$ is a manifold with a closed non-degenerate $2$-form $\omega$.
%In contact case, we have a $2$-form $d\alpha$ which cannot be non-degenerate since
%$C$ is odd dimensional, but its restriction to $\xi$ is non-degenerate.
We will also assume that all contact manifolds in this paper are compact unless stated otherwise.

The main example of a contact manifold in this paper is the link $L_A$ of an isolated singularity $A \subset \C^N$ at $0$.
The curve selection Lemma \cite[Section 3]{milnorsingular} tells us that the function $\phi := |z|^2\big|_A$
has no singularities on $A \setminus \{0\}$ near $0$.
Hence $L_A := \phi^{-1}(\epsilon^2)$ is a manifold for all $\epsilon>0$ small enough whose diffeomorphism type is independent of $\epsilon$.
For a function $f \in C^\infty(\C^N)$, define $d^c f$ by $d^cf(X) = df(i X)$ for all vectors $X$ on $\C^N$
and where $iX$ is the vector $X$ (viewed as a point in $\C^N$) multiplied by $i$. 
%We have that $- d d^c \phi|_{A \setminus \{0\}}$ is a symplectic structure on $A \setminus \{0\}$ near $0$ and that
%$L_A$ is a contact hypersurface with associated contact form
%$-d^c \phi|_{\phi^{-1}(\epsilon^2)}$ for all small enough $\epsilon$
We define the {\bf link} of $A$ at $0$ to be $(L_A,\xi_A) := (\phi^{-1}(\epsilon^2),\text{ker}(-d^c \phi|_{\phi^{-1}(\epsilon^2)}))$
for $\epsilon$ sufficiently small.
Here $(L_A,\xi_A)$ is a contact structure which only depends on the
analytic germ of $A$ at $0$
(see \cite{Varchenko:isolatedsingularities}).
One can show that the contact structure $\xi_A$ on $L_A$ is equal to
$TL_A \cap i TL_A$.

We have the following theorem:
\begin{theorem} (Gray's stability theorem)
Given a smooth family of contact structures $\xi_t$ on a compact manifold $C$,
there is a smooth family of contactomorphisms $\phi_t$ from
$(C,\xi_0)$ to $(C,\xi_t)$ starting from the identity.
\end{theorem}

%One can show, using an appropriate open version of Gray's stability theorem, that every contact %manifold
%is locally contactomorphic to an open subset of $\R^{2n-1}$
%with contact structure $\text{ker}(dz - \sum_i x_i dy_i)$
%where $(x_1,y_1,\cdots,x_{n-1},y_{n-1},z)$ is the natural coordinate system on
%$\R^{2n-1}$.

\begin{defn}
Let $\alpha$ be a contact form on $C$.
The {\bf Reeb vector field of $\alpha$} is the unique vector field $R$ on $C$ satisfying
$i_R d\alpha = 0, i_R \alpha = 1$.
The {\bf Reeb flow} of $\alpha$ is the flow $(\phi_t : C \to C)_{t \in \R}$ of our vector field $R$.
A {\bf Reeb orbit of $\alpha$} is a smooth map $\gamma : \R / L\Z \to C$
satisfying $\frac{d \gamma(t)}{dt} = R$ for some $L \in (0,\infty)$.
Here, $L$ is called the {\bf period of $\gamma$}.
The {\bf linearized return map} of a Reeb orbit $\gamma : \R / L \Z \to C$
is the linear map $D\phi_L|_{T_{\gamma(0)}C} : T_{\gamma(0)} C \to T_{\gamma(0)} C$.
A Reeb orbit is {\bf non-degenerate} if the kernal of the linearized return
map is $1$-dimensional (here the kernal is spanned by $R$ at $\gamma(0)$).
If every Reeb orbit of $\alpha$ is non-degenerate then we say that
{\bf $\alpha$ is non-degenerate} or {\bf the Reeb flow of $\alpha$ is non-degenerate}.
\end{defn}

The dynamics of the Reeb flow can change a lot if we change $\alpha$ to another contact form associated to $\xi$.
An important fact in contact geometry is the fact that a $C^\infty$ generic
contact form is non-degenerate.

Quite often, to solve a problem in contact geometry it is good to embed the contact manifold
as a hypersurface in a symplectic manifold.
Before we do this we will give a definition of a symplectic manifold
and related objects used in this paper.
\begin{defn}
A {\bf symplectic form} is a closed non-degenerate $2$-form on a manifold.
A {\bf symplectic manifold} is a pair $(M,\omega)$
where $M$ is a smooth manifold and $\omega$ is a symplectic form.
If we have any $1$-form $\nu$ on $M$ then we define
it's {\bf $\omega$-dual}
$X^\omega_\nu$ to be the unique
vector field satisfying $\omega(X^\omega_\nu,\cdot) = \nu$.
If the context is clear we just write $X_\nu$.
Let $H : M \to \R$ be smooth, then the {\bf Hamiltonian vector field} $X_H$
is defined to be $X_{dH}$.
%uniquely defined by the equation $i_{X_H} \omega = dH$.
Let $(H_t : M \to \R)_{t \in [0,a]}$ be a smooth family of Hamiltonians,
then the flow $\phi_t : M \to M$ of $(X_t)_{t \in [0,a]}$ is called
the {\bf Hamiltonian flow} of $H_t$.
Such a flow is {\bf time independent} if $X_t$ does not depend on $t$.
A {\bf $\tau$-periodic orbit} of $H_t$ is a map $\gamma : \R / \tau \Z \to M$
satisfying $\frac{d\gamma(t)}{dt} = X_{H_t}$.
An {\bf orbit} of $H$ is a $\tau$-periodic orbit for some $\tau>0$.
The {\bf linearized return map} of a $\tau$-periodic orbit of $H$
is the map $\phi_\tau : T_{\gamma(0)}(M) \to T_{\gamma(0)}M$.
A $\tau$-periodic orbit is {\bf non-degenerate} if the linearized return map has trivial kernal.
A submanifold $L \subset M$ is called {\bf Lagrangian} if $\omega|_L=0$
and $\text{dim}_\R(L) = \frac{1}{2}\text{dim}_\R(M)$.
A submanifold $S \subset M$ is a {\bf symplectic submanifold} if $\omega|_S$ is a symplectic form on $S$.
The {\bf symplectic normal bundle} of a symplectic manifold $S$ is the bundle $TM/TS$ over $S$
with the induced non-degenerate $2$-form. This is sometimes
identified with the set of vectors $v \in TM|_S$ satisfying
$\omega(TS,v) = 0$.
\end{defn}

The main example of a symplectic manifold is $\C^n$ with the {\bf standard symplectic structure}
$\sum_j dx_j \wedge dy_j$. 
If, on some symplectic manifold  $(M,\omega)$, we have some coordinate chart
$x_1,y_1,\cdots,x_n,y_n$ so that $\omega = \sum_j dx_j \wedge dy_j$,
then $x_1,y_1,\cdots,x_n,y_n$ is called a {\bf symplectic coordinate chart}.
A symplectic analogue of Gray's stability
theorem tells us that every point on every symplectic manifold admits a symplectic coordinate
chart centered at that point.
In this paper, many symplectic manifolds will be non-compact and so we do not assume compactness.

\begin{defn}
Let $(M,\omega)$ be a symplectic manifold.
We say that a hypersurface $C \subset M$ is a {\bf contact hypersurface}
if there is a contact form $\alpha$ on $C$ so that $d\alpha = \omega|_C$.
We will call $\alpha$ a {\bf contact form associated to $C \subset M$}.
We call the natural inclusion map $\iota : C \hookrightarrow M$ a {\bf contact embedding}.
Sometimes we will write $\iota : (C,d\alpha) \rightarrow (M,\omega)$ as our contact embedding
with associated contact form $\alpha$.
By Gray's stability theorem the induced contact structure $\xi$ on $C$ is independent of choice
of $\alpha$ up to contactomorphism due to the fact that the
space of such contact forms of the same coorientation is convex.
We will call $\xi$ a {\bf contact structure associated to $C \subset M$}.
\end{defn}

Even though there are many choices of $\alpha$ associated to a given a contact hypersurface
$C \subset M$, all of these choices have Reeb vector fields
that are non-zero multiples of each other at each point.
This is because the Reeb vector field $R$ of a contact form $\alpha$
satisfying $d\alpha = \omega|_C$ satisfies $i_R(\omega|_C) = 0$
and this relation fixes the direction of $R$.
In particular all of these contact forms have exactly the same Reeb orbits
up to reparameterization and hence we can talk about the Reeb orbits
of a contact hypersurface $C \subset M$ even when a contact form is not specified.
When $\text{dim}_\R M \geq 4$,
 the contact embedding gives us a fixed choice of coorientation
(see \cite[Exercise 3.60]{McduffSalamon:sympbook}).

If a contact hypersurface is a regular level set of a Hamiltonian $H$,
then there is a one to one correspondence between its Reeb orbits
and orbits of $H$ inside this hypersurface.
This means that we can translate questions about Reeb orbits into
to questions about orbits of $H$.
This is done in the proof of Theorem \ref{label:nicecontactneighbourhoodexistence}.

The following construction tells us that every contact manifold can be embedded as a contact hypersurface in a symplectic manifold.
This construction will be used later on to find Reeb orbits (see Subsection \ref{subsection:sketchofproof}).

\begin{defn}
Let $(C,\xi)$ be a contact manifold and let $\alpha$ be a contact form associated to $\xi$.
The {\bf symplectization} of $(C,\xi)$ is the symplectic manifold
$(C \times (0,\infty), d(r \alpha))$ where $r$ parameterizes $(0,\infty)$
and by abuse of notation, the pullback of $\alpha$ via the projection map
$C \times (0,\infty) \to C$ is also written as $\alpha$.
If $R$ is the Reeb vector field of $\alpha$, then by abuse of notation
we also define
$R$ to be the unique vector field on
$C \times (0,\infty)$ tangent to $C \times \{x\}$ for all $x$
and equal to $R$ inside $C = C \times \{x\}$ for all $x \in (0,\infty)$.
We will call $R$ the {\bf Reeb vector field of $\alpha$ on $C \times (0,\infty)$}.

Sometimes the symplectization will be written as $(C \times \R, d(e^\rho \alpha))$ where $\rho$ parameterizes $\R$
(here $r = e^\rho$).
\end{defn}

Such a definition is independent of choice of $\alpha$ associated to $\xi$ because if $\beta = f \alpha$
for some smooth $f : C \to \R \setminus \{0\}$,
then the map $(x,r) \to (x, r f(x))$ is a symplectomorphism from
$(C \times (0,\infty), d(r \beta))$ to $(C \times (0,\infty), d(r \alpha))$.
%The symplectization can be useful as it translates questions about Reeb flows
%into questions involving Hamiltonian flows on the symplectization.
%This is because a Hamiltonian $H = f(r)$ for any $f : (0,\infty) \to \R$
%satisfies $X_H = f'(r) R$.
%Hence the Hamiltonian flow of $H$ on each level set of $H$
%is a scalar multiple of the Reeb flow and hence questions about Reeb orbits
%get translated into questions involving orbits of $H$ for appropriate $f$.

We have the following Lemma:
\begin{lemma} \label{lemma:contacthypersurfaceneighborhood}
Suppose that $C \subset M$ is a contact hypersurface in $M$ with an induced contact structure $\xi$.
Then a neighborhood of $C$ in $M$ is contactomorphic to a neighborhood of $C \times \{1\}$
in the symplectization of $(C,\xi)$.
\end{lemma}
The proof of this Lemma is contained in the last paragraph of the proof
of \cite[Proposition 3.58]{McduffSalamon:sympbook}.
The above Lemma implies the following fact:
If $\iota : (C,d\alpha) \hookrightarrow (M,\omega)$ is a contact embedding,
then for any sufficiently $C^\infty$ close contact form $\beta = f \alpha$,
there is a contact embedding $\iota' : (C,d\beta) \hookrightarrow (M,\omega)$
$C^\infty$ close to $\iota$.
In particular, given any contact hypersurface,
we can perturb it in a $C^\infty$ generic way so that all of
its Reeb orbits are non-degenerate.

\bigskip
{\bf Complex structures and the Conley-Zehnder index.}
The Conley-Zehnder index was originally defined in \cite{arnold:characteristic}.
The Conley-Zehnder index of a Reeb orbit $\gamma$
tells us how many times the Reeb flow near such an orbit
`wraps' around $\gamma$.
Its definition is fairly technical,
and so we leave the details to Section \ref{section:conleyzehnderindexsection}.
Having said that, we will state a property  of the
Conley-Zehnder index (Lemma \ref{CZ:illustrativeproperty})
which will help explain the above intuitive meaning of
this index and which will
also help us explain the relationship between these indices and the
minimal discrepancy.

\begin{defn}
Let $E \to X$ be a vector bundle.
A {\bf symplectic structure} on $E$ is a $2$-form $\Omega$ on $E$ which is fiberwise linear and non-degenerate.
A {\bf complex structure} on $E$ is a fiberwise linear automorphism $J : E \to E$ satisfying $J^2 = -\text{id}$.
We say that $J$ is {\bf compatible with $\Omega$} if $\Omega(\cdot,J(\cdot))$ is a Riemannian metric on $E$.
The triple $(E,\Omega,J)$ is called a {\bf Hermitian vector bundle}.
A {\bf trivialization} of any of these bundles is a bundle 
isomorphism $\tau : E \to X \times \C^k$ so that the symplectic (resp. complex)
structure is the pullback of the standard one on $\C^k$ via
the composition of $\tau$ with the projection map $X \times \C^k \twoheadrightarrow \C^k$.
An {\bf almost complex structure} $J$ on a manifold $M$ is a complex structure on its tangent bundle.
\end{defn}

If $\Omega$ is a fixed symplectic structure on a vector bundle then the space of
complex structures compatible with $\Omega$ is contractible
(see \cite[Proposition 2.50]{McduffSalamon:sympbook}).

\begin{defn}
Let $(C,\xi)$ be a contact manifold with a choice of coorientation for $\xi$.
Let $\alpha$ be a contact form associated to $\xi$
respecting this choice of coorientation,
then $d\alpha|_\xi$ is a symplectic structure on $\xi$.
Let $J$ be a choice of complex structure compatible with $d\alpha|_\xi$.
We define the {\bf first Chern class}
$c_1(\xi) \in H^2(C;\Z)$ of $\xi$ to be the first Chern class of the complex vector bundle $(\xi,J)$.
The {\bf anticanonical bundle} $\kappa^*_\xi$ of $\xi$
is defined to be the highest exterior power of the complex bundle $(\xi,J)$.
The {\bf canonical bundle} $\kappa_\xi$ is the dual of the anticanonical bundle.

If $(M,\omega)$ is a symplectic manifold, then its first Chern class $c_1(M,\omega)$
is defined to be $c_1(M,J)$ for some almost complex $J$ compatible with $\omega$.
Sometimes we will write $c_1(TM)$ if it is clear what symplectic structure we are using.
The {\bf canonical bundle $\kappa_\omega$} of $(M,\omega)$
is the highest exterior power of $(TM,J)$ and the {\bf anticanonical bundle}
is its dual.
\end{defn}

The class $c_1(\xi)$ only depends on $\xi$ and its choice of coorientation and not on $\alpha$.
Also the anticanonical bundle does not depend on $J$ up to complex bundle isomorphism.
This is because the space of $d\alpha|_\xi$ compatible almost complex structures is contractible
and also because the space of contact forms associated to $\xi$ respecting the given choice of
coorientation is contractible.
Similarly, the Chern class is independent of $J$ and the (anti)canonical bundle of $(M,\omega)$ do not depend on $J$ up to isomorphism.

Let $(C,\xi)$ be a contact manifold with a choice of coorientation for $\xi$ and suppose that
$N c_1(\xi) = 0$ for some $N > 0$ and that $H^1(C;\Q)= 0$.
For any Reeb orbit $\gamma : \R / L \Z \to C$ of some contact form $\alpha$ associated to $\xi$,
we can define an index $\text{CZ}(\gamma) \in \frac{1}{N}\Z \subset \Q$ called the
{\bf Conley-Zehnder index}.
The Conley-Zehnder index tells us how much the Reeb flow `wraps' around our Reeb orbit $\gamma$.
The following Lemma illustrates this in a particular case.

\begin{lemma} \label{CZ:illustrativeproperty}
Let $J$ be a complex structure on $\xi$ compatible with $d\alpha|_\xi$,
let $\pi_\C : C \times \C \to \C$ be the projection map
and choose a trivialization $\tau : {\kappa^*_\xi}^{\otimes N} \to C \times \C$.
Let $\phi_t$ be the Reeb flow of $\alpha$.
The restriction $D\phi_t|_\xi$ of the linearization $D\phi_t$ preserves $\xi$ and so we view this as a map from $\xi$ to $\xi$.
Suppose that $D\phi_t|_\xi : \xi \to \xi$ is a $J$ complex linear map (i.e. it commutes with $J$).
This induces a fiberwise complex linear map $\widetilde{\phi}_t : {\kappa^*_\xi}^{\otimes N} \to {\kappa^*_\xi}^{\otimes N}$. Suppose also that the linearized return map of $\gamma$ is the identity map.

Then the Conley-Zehnder index of $\gamma$ is given by $2/N$ multiplied by the degree of the map
$q : \R / L \Z \to U(1)$,
defined by:
\[q(t) := (\pi_\C \circ \tau \circ \widetilde{\phi}(t) \circ (\pi_\C \circ \tau|_{\gamma(0)})^{-1}).\]
\end{lemma}

Here $\tau|_{\gamma(0)}$ means the restriction of $\tau$ to the fiber
of ${\kappa^*}^{\otimes N}$ over $\gamma(0)$.
This Lemma follows directly from property \ref{item:determinantproperty} in Section \ref{section:conleyzehnderindexsection}.
The map $q$ describes how the Reeb flow `wraps' around our orbit $\gamma$.
In the more general situation, the map $q$ would be replaced
by a path $[0,1] \to \text{Sp}(\R^{2n-2})$
where $\text{Sp}(\R^{2n-2})$ is the group of linear symplectomorphisms
and then one uses a recipe in \cite{RobbinSalamon:maslov}
to give us an index from such a path.

\begin{defn}
Let $\gamma$ be a Reeb orbit of $\alpha$. Then
the {\bf lower SFT index} is defined to be
\[\text{lSFT}(\alpha) := \text{CZ}(\gamma) - \frac{1}{2}\text{dim ker}(D_{\gamma(0)}\phi_L|_\xi - \text{id}) + (n-3).\]
For any contact form $\alpha$ associated to $\xi$
we define the {\bf minimal SFT index of }$\alpha$ to be $\text{mi}(\alpha) := \text{inf}_\gamma \text{lSFT}(\gamma)$
where the infimum is taken over all Reeb orbits $\gamma$ of $\alpha$.
The {\bf highest minimal SFT index} of $(C,\xi)$ is defined as
$\text{hmi}(C,\xi) := \text{sup}_{\alpha} \text{mi}(\gamma)$
where the supremum is taken over all nowhere zero $1$-forms $\alpha$  with $\text{ker}(\alpha) = \xi$ respecting the coorientation of $\xi$.
\end{defn}

The above definitions are original, although they are very similar
to the definition of the SFT index (originally from \cite[Proposition 1.7.1]{EliashbergGiventalHoferSFT}, although not stated as an index).
SFT stands for symplectic field theory.
%The meaning of the lower SFT index is slightly more complicated.
Subtracting the term $- \frac{1}{2}\text{dim ker}(D_{\gamma(0)}\phi_L|_\xi - \text{id})$ makes our index have the nice property
that this index is lower semi-continuous with respect to $\alpha$ in some sense
(see Lemma \ref{lemma:perturbedindexcalculation}).
This is why it is called {\it lower} SFT index.
Later on this enables us to perturb our contact form,
do some $\text{lSFT}$ index calculations on this perturbed
contact form so that $\text{hmi}(\alpha)$ can be bounded
(this is done in the proof of Lemma \ref{lemma:reeborbits}
using Corollary \ref{corollary:lsftlowersemicontinuous}).
This term also appears in calculations from \cite{Bourgeois:morsebott}.
The $n-3$ term is used here for two reasons:
\begin{enumerate}
\item It appears naturally in Symplectic Field Theory
\cite[Proposition 1.7.1]{EliashbergGiventalHoferSFT}
as a dimension of a moduli space of curves, suggesting
a deeper relationship between the singularity and holomorphic curves.
\item Because it makes the formulas in Theorem
\ref{theorem:minimaldiscrepancy} look less complicated.
\end{enumerate}

It is clear that $\text{hmi}(C,\xi)$ is a coorientation preserving contactomorphism invariant.
In our case we can show more.

\begin{lemma} \label{lemma:generalcontactomorphisminvariant}
Suppose that $(C,\xi)$ admits a coorientation reversing contactomorphism.
Then $\text{hmi}(C,\xi)$ is a (not necessarily coorientation preserving) contact invariant.
\end{lemma}
\proof
Let $\Phi$ be this coorientation reversing contactomorphism.
Then for any contact form $\alpha$ associated to $\xi$, we have that $\text{mi}(\alpha) = \text{mi}(\Phi^*(\alpha))$.
Also the map $\Phi$ induces a bijection between contact forms associated to $\xi$ which respect the coorientation
and ones which give the opposite coorientation.
Hence $\text{hmi}(C,\xi)$ is a (not necessarily coorientation preserving) contactomorphism invariant.
\qed

\bigskip
The map $z \to \overline{z}$ restricted to the link $L_A$ of our singularity $A$ at $0$
gives us a coorientation reversing contactomorphism of $(L_A,\xi_A)$.
This is because it sends $-d^c \phi$ to $d^c\phi$ where $\phi = \sum_j |z_j|^2 \big|_A$.
Hence if $N c_1(\xi_A) = 0$ for some $N$, we have that
$\text{hmi}(L_A,\xi_A)$ is a (not necessarily coorientation preserving) contact invariant by the above Lemma.

%\subsection{A Brief Description of the Conley-Zehnder index and the Highest Minimal Index}

\subsection{Sketch of the proof} \label{subsection:sketchofproof}

Here we give a sketch of the proof of
Theorem \ref{theorem:minimaldiscrepancy}.
We will do this just for cone singularities
as they are easier to manage, and then at the end
we will briefly explain how to adjust this proof in the case of a general
isolated singularity.

Throughout this section
let $A$ have a normal isolated singularity at
$0$ with $H^1(L_A;\Q) = 0$ and $c_1(\xi_A;\Q) = 0$.
Our main theorem follows from the following two statements:
\begin{enumerate}
\item {\bf The easier statement} (Theorem \ref{theorem:boundingdiscrepancyfromabove}).
\[\text{hmi}(L_A,\xi_A) \geq 2\text{md}(A,0).\]
\item
{\bf The harder statement} (Theorem \ref{theorem:lowerbound}).
\begin{itemize}
\item If $\text{md}(A,0) \geq 0$ then $\text{hmi}(L_A,\xi_A) \leq 2\text{md}(A,0)$.
\item If $\text{md}(A,0) < 0$ then $\text{hmi}(L_A,\xi_A) < 0$.
\end{itemize}
\end{enumerate}

To make the proof easier to explain, we will
assume that $A \subset \C^N$ is the cone over a smooth connected projective
variety $X \subset \P^{N-1}$.
In this case, it can be resolved by blowing up once at the origin.
This resolution is the variety
$\widetilde{A}$ which is the total space of the line
bundle ${\mathcal O}(-1)$ over $X$.
Let $\widetilde{\pi} : \widetilde{A} \to X$ be the natural projection map.
Here we identify $X$ with the zero section of this bundle.
%The numerically $\Q$-Gorenstein condition in this case
%is equivalent to the fact that
%\begin{equation} \label{equation:numericallygorensteincone}
%c_1({\mathcal K}_X;\Q) = (a+1) c_1({\mathcal N}_X;\Q)
%\end{equation}
%for some $a \in \Q$.
%The number $a$ here is the discrepancy of the exceptional divisor
%$X$.
%Hence the minimal discrepancy is $a$ if $a \geq -1$
%and $-\infty$ otherwise.

So far we have not given a definition of the minimal discrepancy $\text{md}(A,0)$.
We will give a general definition in Section \ref{section:minimaldiscrepancydefinition}.
In our particular case, it is calculated as follows:
Let ${\mathcal K}_X$ be canonical bundle of $X$.
The numerical $\Q$-Gorenstein condition is equivalent to the fact that
$c_1({\mathcal K}_{\widetilde{A}} \big|_{L_A};\Q) = 0$
(see Lemma \ref{lemma:topologicalnumericallygorenstein}).
Because $L_A \hookrightarrow \widetilde{A} \setminus X$
is a homotopy equivalence, we get that
$c_1({\mathcal K}_{\widetilde{A}} \big|_{\widetilde{A} \setminus X};\Q) = 0$.
This is equivalent to the existence of a $C^\infty$ section
$s$ of ${\mathcal K}_{\widetilde{A}}^{\otimes N}$
which is non-zero outside a compact set
for some $N > 0$.
By Thom transverality, we can assume that $s$ is transverse to $0$.
The {\bf discrepancy} of $X$ is
equal to $a$ where $a$ satisfies
\[[s^{-1}(0)] = aN [X] \in H_{2n-2}(\widetilde{A};\Q) = H_{2n-2}(X;\Q).\]
The {\bf minimal discrepancy} $\text{md}(A,0)$ is defined to be $a$ if $a \geq -1$
and $-\infty$ otherwise.

{\bf Proof of the easier statement.}
For this it is sufficient to find a contact form $\alpha_A$
associated to $\xi_A$ so that
$\text{mi}(\alpha_A) = 2\text{md}(A,0)$.
We have that $L_A$
is the circle bundle of radius $\epsilon$ on the Hermitian
line bundle ${\mathcal O}(-1)$ over $X$
(the Hermitian metric here is the restriction of the natural
Euclidean metric on $\C^N$ to the fibers projected to $A$).
We have that the map $\pi := \widetilde{\pi}|_{L_A}$
makes $L_A$ in to a circle bundle over $X$.
Let $B : \R / 2\pi \Z \times \widetilde{A} \to \widetilde{A}$ be the $S^1$
action rotating the fibers of $\widetilde{\pi}$.
Define $\alpha_A := -\frac{1}{4\pi\epsilon^2}d^c(\sum_j |z_j|^2)\big|_{L_A}$.
Here $\alpha_A$ is a $\xi_A$ admissible contact form
whose Reeb flow is the natural
circle action $B|_{L_A} : \R /  \Z \times L_A \to L_A$ rotating the fibers
of our circle bundle $\pi$.
Hence one can describe all the Reeb orbits.
For each $k \in \N$ and each point $p$ in $L_A$
there is a Reeb orbit $\gamma : \R /  k \Z \to L_A$ defined by
$\gamma(t) = B(t,p)$.
Every Reeb orbit is of this form.
The linearized return map for any of these orbits is the identity map.
We now have to calculate the Conley-Zehnder index for each of these Reeb orbits
using Lemma
\ref{CZ:illustrativeproperty}.
This calculation really enables us to see the relationship between the Reeb flow and the minimal discrepancy
and so we will now do it in detail.

Let $\gamma : \R /  k \Z \to L_A$  be one of these Reeb orbits defined
by $\gamma(t) = B(t,p)$ as above.
We can assume our section $s$ of ${\mathcal K}_{\widetilde{A}}$
described earlier is non-zero along $L_A$.
Let ${\mathcal N}_X$ the normal bundle of the zero section $X$ in $\widetilde{A}$.
The bundle $\pi^* {\mathcal N}_X$ has a canonical non-zero section
$s_{\mathcal N}$ sending
$x \in L_A$ to its respective point in ${\mathcal N}_X$ and then pulling it
back to $\pi^* {\mathcal N}_X$.
Because $\xi_A$ is transverse to the fibers of $\widetilde{\pi}$, it is canonically isomorphic to $\pi^* TX$
as complex vector bundles.
Therefore
${\kappa^*_{\xi_A}} \otimes \pi^* {\mathcal N}_X \otimes {\mathcal K}_{\widetilde{A}}|_{L_A}\cong \Lambda^nT\widetilde{A}|_{L_A} \otimes {\mathcal K}_{\widetilde{A}}|_{L_A} \cong L_A \times \C$
is canonically trivial and hence has a canonical non-zero section $S$.
Hence we have a unique non-zero section
$s_{\kappa^*}$ of ${\kappa^*_{\xi_A}}^{\otimes N}$ satisfying
$s_{\kappa^*} \otimes s_{\mathcal N}^{\otimes N} \otimes s = S^{\otimes N}$.

Let  $B_{\kappa^*} : \R / \Z \times {\kappa^*}_{\xi_A}^{\otimes N} \to
 {\kappa^*}_{\xi_A}^{\otimes N}$,
$B_{\mathcal N} : \R / \Z \times \pi^* {\mathcal N}_X^{\otimes N} \to
 \pi^* {\mathcal N}_X^{\otimes N}$ and
$B_{{\mathcal K}_{\widetilde{A}}} : \R / \Z \times {\mathcal K}_{\widetilde{A}}^{\otimes N}
\to {\mathcal K}_{\widetilde{A}}^{\otimes N}$
be the respective liftings of the action $B$ to $ {\kappa^*}_{\xi_A}^{\otimes N}$,
$ \pi^* {\mathcal N}_{X}^{\otimes N}$
and ${\mathcal K}_{\widetilde{A}}$
respectively.
Let $P : \C^* \to S^1 = U(1)$ be the map sending $z$ to $z / |z|$.
Then by Lemma \ref{CZ:illustrativeproperty},
we get that the Conley-Zehnder index of $\gamma$ is $2/N$ multiplied by the degree of the map:
\[Q : \R /  k \Z \to U(1), \quad Q(t) = \left[ z \to P( B_{\kappa^*}(t,s_{\kappa^*}(\gamma(0))) / s_{\kappa^*}(\gamma(t)) )\right].\]
Define
\[Q_{\mathcal K} : \R /  k \Z \to U(1), \quad Q_{\mathcal K}(t) = \left[ z \to P( B_{\mathcal K}(t,s(\gamma(0))) / s(\gamma(t)) )\right].\]
The actions $B_{\mathcal N}$, $B_{\kappa^*}$ and $B_{\mathcal K}$
induce an action on ${\kappa^*}_{\xi_A}^{\otimes N} \otimes \pi^*{\mathcal N}_X^{\otimes N} \otimes {\mathcal K}_{\widetilde{A}}^{\otimes N}|_{L_A}$
which sends $S(\gamma(0))$ to $S(\gamma(t))$ for all $t \in \R$.
Also $B_{\mathcal N}(t,s_{\mathcal N}(\gamma(0))) = s_{\mathcal N}(\gamma(t))$ for all $t \in \R$.
Hence the degree of $Q$ is equal to the degree of $Q_{\mathcal K}$ taken with negative sign.

Let $F$ be the fiber containing $\gamma$.
Let $s_F$ be a non-zero section of ${\mathcal K}_{\widetilde{A}}^{\otimes N}\big|_F$.
Define
$Q_F : \R /  k \Z \to U(1)$ by
$Q_F(t) = \left[ z \to P( B_{\mathcal K}(t,s_F(\gamma(0))) / s_F(\gamma(t)) )\right]$.
Then because the action $B$ rotates the fiber $F$,
we have that the degree of $Q_F$
is $- k N$.
Perturb $s$ slightly so that $s^{-1}(0)$ is transverse to $F$.
Then $[s^{-1}(0)\big|_F] = a N \in H_0(F;\Q) \cong \Q$.
Hence $\text{degree}(Q_{\mathcal K}) = - a N + \text{degree}(Q_F)$
and so $\text{degree}(Q) = a N - \text{degree}(Q_F)$.
Hence the Conley-Zehnder index of $\gamma$ is $2/N(a N - \text{degree}(Q_F)) = 2(a+1)k$.

Because the linearized return map of $\gamma$ is the identity,
this implies that the $\text{lSFT}$ index of $\gamma$
is $2(a+1)k - \frac{1}{2}(2n-2) + (n-3) = 2(a+1)k - 2$.
This implies that $\text{mi}(\alpha_A) = 2 \text{md}(A)$ which gives us our result.
\bigskip

{\bf Sketch of the proof of the harder argument.}
Here we will show that any contact form $\beta$
associated to $\xi_A$ admits a Reeb orbit
either with negative $\text{lSFT}$ index
or with $\text{lSFT}$ index bounded above by
$2\text{md}(A,0)$.
There are three main tools we use in the proof. The first is
{\it neck stretching}, the second is {\it Gromov-Witten invariants}
and the third is {\it symplectic dilation}.
Before we talk about the proof directly, we will talk about these
three tools. Some of the definitions here are slightly more detailed than needed
as they are used later on in the paper.
Both neck stretching and Gromov-Witten invariants involve
holomorphic curves which we now define:

\begin{defn} \label{defn:pseudoholomorphiccurve}
A {\bf nodal Riemann surface} is a one dimensional complex analytic variety
with only nodal singularities
(i.e. it is locally analytically isomorphic to $\C$ or $\{xy=0\} \subset \C^2$).
It has {\bf arithmetic genus $0$} if it is a subvariety of a simply connected
nodal Riemann surface.
Suppose we have an almost complex manifold $(S,J)$.
A {\bf nodal $J$-holomorphic curve} is a continuous map
$u : \Sigma \to S$ so that the restriction to
the smooth part of $\Sigma$ is $J$-holomorphic
(i.e. $du \circ j = J \circ du$ where $j : T\Sigma \to T\Sigma$ is the complex structure on $\Sigma$).
The {\bf fundamental class} $[\Sigma] \in H_2(S;\Z)$ is the image of the fundamental class
$[\widetilde{\Sigma}]$ of the normalization of $\Sigma$
(if $\widetilde{\Sigma}$ is non-compact then we use Borel-Moore homology).
Such a curve is ${\bf compact}$ if its domain is compact, {\bf connected} if its domain is connected,
{\bf smooth} if its domain is smooth
and it is of {\bf genus $0$} if $\Sigma$ has {\bf arithmetic genus $0$}.
A compact $J$-holomorphic curve $u : \Sigma \to S$ {\bf represents a class } $[A] \in H_2(S;\Z)$
if $u_*([\Sigma]) = [A]$.
One has a similar definition in the case when $\Sigma$ is non-compact
and $u$ is a proper map,
in which case we use Borel-Moore homology.
A $J$-holomorphic curve $u : \Sigma \to M$ is {\bf somewhere injective},
if each irreducible component of $\Sigma$ has a point $p$
where $u^{-1}(u(p)) = \{p\}$.
An irreducible component $\Sigma_1$ of $\Sigma$ is said to be {\bf multiply covered}
if $u|_{\Sigma_1}$ is not somewhere injective
(this is because $u|_{\Sigma_1}$ factors through a degree $\geq 2$ branched covering
of Riemann surfaces).
\end{defn}

From now on we will assume that all $J$-holomorphic curves are connected,
but they may not be compact or irreducible.
Sometimes in this paper we will talk about $J$-holomorphic curves being {\bf regular}.
A $J$-holomorphic curve is regular if a certain linear operator is surjective.
We will not give a definition in this paper of regularity as it is not needed.
All we need to know is that for a $C^\infty$ generic perturbation of $J$
inside some open region $U$,
we have that all smooth somewhere injective $J$-holomorphic curves
with image intersecting $U$ are regular.
We only need such a statement in order to apply
a result in \cite{Dragnev:transversality}
inside the proof of Lemma \ref{lemma:reeborbits}.

We will first talk about {\it neck stretching}
(described in detail in \cite{BEHWZ:compactnessfieldtheory}
and also in Appendix A).
This will be used to find Reeb orbits. Our neck stretching
construction is contained in the proof of the following Lemma.
In the more general case when $A$ is not a cone singularity, we need a more technical Lemma
(see Lemma \ref{lemma:reeborbits}).

\begin{lemma} \label{label:neckstretchingdemonstration}
Let $(M,\omega)$ be a compact symplectic manifold,
and suppose that it has a contact hypersurface
$C \subset M$ so that:
\begin{enumerate}
\item  $M \setminus C$ has two connected components
$M_-$ and $M_+$.
\item There are two codimension $2$ submanifolds
$Q_- \subset M_-,Q_+ \subset M_+$
and a homology class $[A] \in H_2(M;\Z)$ such that $[A]\cdot[Q_\pm] \neq 0$.
\item
For every almost complex structure $J$ compatible with $\omega$,
there exists a compact arithmetic genus $0$ $J$-holomorphic curve
$u : \Sigma \to M$ representing $[A]$.
\end{enumerate}
Then $C$ has at least one Reeb orbit.
\end{lemma}

At the moment this Lemma is not good enough because it does not
give us a bound on the $\text{lSFT}$ index of the Reeb orbit of $C$.
We will explain how to find such bounds later inside the proof
of the main result of this section.
\smallskip

{\it Sketch of proof.}
First we choose a neighborhood of $C$ equal to
$(1-\epsilon,1+\epsilon) \times C$
with symplectic form $d(r \alpha)$ as in Lemma
\ref{lemma:contacthypersurfaceneighborhood}.
We will call this region a {\bf neck}.
We now wish to `stretch' this neck.
%Define $r_C := \ln(r)$.
Define $\phi_\infty : (1-\epsilon,1+\epsilon) \setminus \{0\} \to (0,\infty)$,
$\phi_\infty := e^{\frac{1}{1-x}}$.
Let
$\phi_i : (1-\epsilon,1+\epsilon) \to (0,\infty)$ be a sequence of $C^\infty$ functions with $\phi_i' > 0$
so that $\phi'_i|_{(1-\epsilon,1+\epsilon) \setminus \{1\}}$ converges in $C^\infty_{\text{loc}}$
to $\phi'_\infty$ as $i \to \infty$.
Let $\widehat{J}$ be an almost complex structure on the symplectization
$C \times (0,\infty)$ compatible with the symplectic form,
invariant under the map $(x,r) \to (x,\kappa r)$ for any $\kappa > 0$,
 so that $\widehat{J}(\text{ker}(\alpha)) = \text{ker}(\alpha)$
and so that $\widehat{J}\big(r\frac{\partial}{\partial r}\big) = R$ where $R$ is the Reeb vector field.
Choose a sequence of almost complex structures $(J_i)_{i \in \N}$
compatible with $\omega$ so that $J_i|_{C \times (1-\epsilon,1+\epsilon)} = (\text{id}_C,\phi_i)^* \widehat{J}$.
%where $(\phi_i,\text{id}_C)$ is a product morphism from the open subset
%$C \times (1-\epsilon,1+\epsilon)$ to $C \times (0,\infty)$.
The maps $(\phi_i,\text{id}_C)$ are the maps which `stretch the neck'.
Let $J_\infty$ be a compatible almost complex structure on $M \setminus C$
equal to $(\phi_\infty,\text{id}_C)^*\widehat{J}$ near $C$.
We can ensure that $J_i$ converges in $C^\infty_\text{loc}$ to $J_\infty$.
The sequence $J_i$ is called a {\bf a sequence of almost complex structures
stretching the neck along $C \subset M$}.
One should imagine the sequence $J_i$ `stretching' $M$ along $C$
until it `breaks' $M$ into two pieces $M_-$ and $M_+$.
Here is a schematic picture:

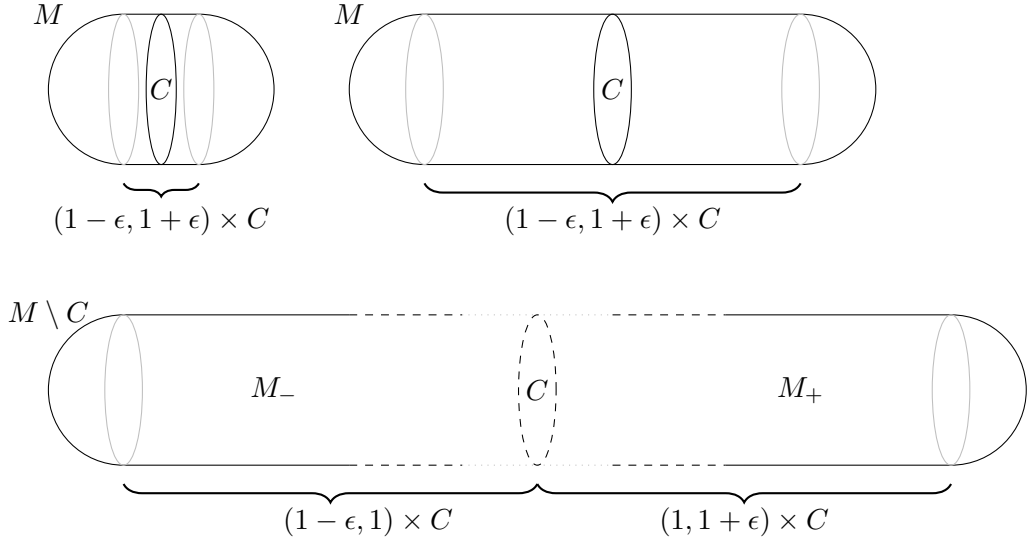
\begin{tikzpicture} [scale=1.0]
% Unstretched
\draw[shift={(5,5)},rotate=0,shift={(0,0)}] (0,1) arc (90:270:1);
\draw[shift={(6,5)},rotate=180,shift={(0,0)}] (0,1) arc (90:270:1);
\draw[shift={(0,-1)}] (5,5) -- (6,5);
\draw[shift={(0,1)}] (5,5) -- (6,5);
\draw[color=lightgray] (5,5) ellipse (0.2 and 1);
\draw[color=lightgray] (6,5) ellipse (0.2 and 1);
\draw[] (5.5,5) ellipse (0.2 and 1);
\node  at (5.5,5) {$C$ };

\draw [
    thick,
    decoration={
        brace,
        mirror,
        raise=0.0cm,
	amplitude=5pt
    },
    decorate
] (5,5-1.25) -- (6,5-1.25);
\node  at (5.5,5-1.75) {$(1-\epsilon,1+\epsilon) \times C$ };
\node  at (5-1,5+1) {$M$ };

%Partially stretched
\draw[shift={(9,5)},rotate=0,shift={(0,0)}] (0,1) arc (90:270:1);
\draw[shift={(14,5)},rotate=180,shift={(0,0)}] (0,1) arc (90:270:1);
\draw[shift={(0,-1)}] (9,5) -- (14,5);
\draw[shift={(0,1)}] (9,5) -- (14,5);
\draw[color=lightgray] (9,5) ellipse (0.25 and 1);
\draw[color=lightgray] (14,5) ellipse (0.25 and 1);
\draw[] (11.5,5) ellipse (0.25 and 1);
\node  at (11.5,5) {$C$ };
\draw [
    thick,
    decoration={
        brace,
        mirror,
        raise=0.0cm,
	amplitude=7pt
    },
    decorate
] (9,5-1.25) -- (14,5-1.25);
\node  at (11.5,5-1.75) {$(1-\epsilon,1+\epsilon) \times C$ };
\node  at (9-1,5+1) {$M$ };

%broken
\draw[shift={(5,1)},rotate=0,shift={(0,0)}] (0,1) arc (90:270:1);
\draw[shift={(16,1)},rotate=180,shift={(0,0)}] (0,1) arc (90:270:1);
\draw[shift={(0,-1)}] (5,1) -- (8,1);
\draw[shift={(0,1)}] (5,1) -- (8,1);
\draw[shift={(0,-1)},dashed] (8,1) -- (9.5,1);
\draw[shift={(0,1)},dashed] (8,1) -- (9.5,1);
\draw[shift={(0,1)},dotted,color=lightgray] (9.5,1) -- (10.5,1);
\draw[shift={(0,-1)},dotted,color=lightgray] (9.5,1) -- (10.5,1);
\draw[shift={(0,-1)}] (13,1) -- (16,1);
\draw[shift={(0,1)}] (13,1) -- (16,1);
\draw[shift={(0,-1)},dashed] (11.5,1) -- (13,1);
\draw[shift={(0,1)},dashed] (11.5,1) -- (13,1);
\draw[shift={(0,1)},dotted,color=lightgray] (10.5,1) -- (11.5,1);
\draw[shift={(0,-1)},dotted,color=lightgray] (10.5,1) -- (11.5,1);
\draw[color=lightgray] (5,1) ellipse (0.25 and 1);
\draw[color=lightgray] (16,1) ellipse (0.25 and 1);
\draw[dashed] (10.5,1) ellipse (0.25 and 1);

\draw [
    thick,
    decoration={
        brace,
        mirror,
        raise=0.0cm,
	amplitude=7pt
    },
    decorate
] (10.5,1-1.25) -- (16,1-1.25);
\node  at (13.25,1-1.75) {$(1,1+\epsilon) \times C$ };

\draw [
    thick,
    decoration={
        brace,
        mirror,
        raise=0.0cm,
	amplitude=7pt
    },
    decorate
] (5,1-1.25) -- (10.5,1-1.25);
\node  at (8.25,1-1.75) {$(1-\epsilon,1) \times C$ };

\node  at (5-1,1+1) {$M \setminus C$ };
\node  at (10.5,1) {$C$ };
\node  at (5 + 2,1) {$M_-$ };
\node  at (16 - 2,1) {$M_+$ };

%\draw (4,5) .. controls (4,6) and (6,6) .. (6,5);
%\draw[color=lightgray] (5,5) ellipse (1 and 0.25);
%\draw (5,5) node {\tiny A} circle (1);
%(7,6) is the START point of the arc
%\draw[shift={(5,5)},rotate=35,shift={(-5,-5)}] (7,5)+(0,1) arc (90:270:1);
%\draw[shift={(5,5)},rotate=130,shift={(-5,-5)}] (7,5)+(0,1) arc (90:270:1);
%\draw[shift={(5,5)},rotate=35,shift={(-5,-5)}] (7,5) node {\tiny B} ellipse (0.25 and 1);
%\draw[shift={(5,5)},rotate=130,shift={(-5,-5)}] (7,5) node {\tiny C} ellipse (0.25 and 1);
%\draw[shift={(5+5,5)},rotate=130,shift={(-5,-5)}]  (5,5) -- (7,5) node[below,left] {\tiny C};
%\draw[shift={(5+5,5)},rotate=35,shift={(-5,-5)}] (5,5) node[below] {\tiny A} -- (7,5) node[below,right] (7,5) {\tiny B};
%\draw[fill=black] (5+5,5) circle (0.05);
%\draw[shift={(5+5,5)},rotate=35,shift={(-5,-5)},fill=black] (7,5) circle (0.05);
%\draw[shift={(5+5,5)},rotate=130,shift={(-5,-5)},fill=black] (7,5) circle (0.05);
%\node  at (4,3.5) {  Riemann surface.};
%\node  at (4+5,4.5) { Dual graph. };

%\node[shift={(5,5)},rotate=35,shift={(-5,-5)}] at (7,5);
%\node[shift={(5,5)},rotate=130,shift={(-5,-5)}] at (7,5);
\end{tikzpicture}

Let us look at the following example:
Here $M = \C \P^1$ with symplectic structure coming from the
Fubini-Study K\"{a}hler form and $C = \R P^1$ is the equator.
We choose $J_i$ so that $(M,J_i)$ is biholomorphic to
$\{xy = z/i\} \subset \P^2$ for each $i \geq 1$
where $[x,y,z] \in \P^2$ are projective coordinates.
We do this in such a way so that $M \setminus C$ with $J_\infty$
is identified holomorphically with $\{(xy = 0)\} \setminus \{(0,0,1)\}$.
This example (and other similar examples) suggest that neck stretching
is a certain contact analogue of a sequence of smooth varieties degenerating
to a smooth normal crossing variety with two irreducible components.

We now return to the general situation.
By assumption, there are arithmetic genus $0$ $J_i$-holomorphic
curves $u_i : \Sigma_i \to M$ representing $[A]$.
After passing to a subsequence, one can show
via a compactness argument
(see \cite{BEHWZ:compactnessfieldtheory} or
Proposition \ref{proposition:compactnessresult})
that $u_i|_{u_i^{-1}(M_+)}$ $C^0$ converges to
some $J_\infty$-holomorphic curve
$u_\infty : \Sigma_\infty \to M_+$ where $u_\infty$ is a proper map
(after identifying the domains $u_i^{-1}(M_+)$ with $\Sigma_\infty$ in a particular way).
Also
\cite[Proposition 5.6]{BEHWZ:compactnessfieldtheory}
and Lemma \ref{lemma:convergence}
tell us that $\Sigma_\infty$
compactifies to a Riemann surface with boundary $\overline{\Sigma}$,
and $u_\infty$ extends to
a continuous map $\overline{u}_\infty : \overline{\Sigma} \to M_+ \cup C$,
where $\overline{u}_\infty$ maps the boundary of
$\Sigma$ to a union of Reeb orbits $\gamma_1,\cdots,\gamma_k$
(here we have to assume that the Reeb flow is non-degenerate,
but this is true if we perturb $C$ slightly, which is OK).
Here we say that $\overline{u}_\infty$
{\bf has negative ends converging to $\gamma_1,\cdots,\gamma_k$}.
Hence $C$ has at least one Reeb orbit. Here is a schematic picture:

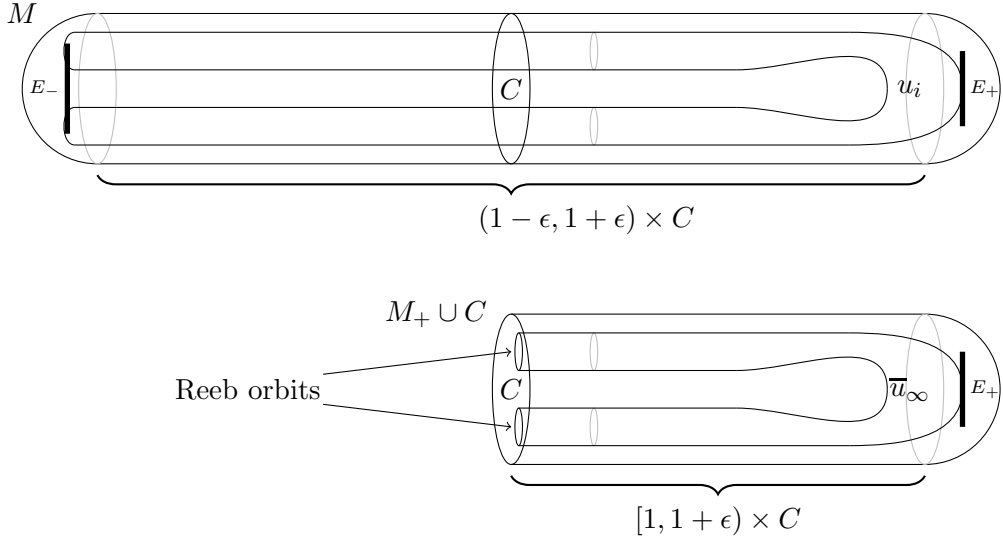
\begin{tikzpicture} [scale=1.0]

%Partially stretched
\draw[shift={(5,5)},rotate=0,shift={(0,0)}] (0,1) arc (90:270:1);
\draw[shift={(16,5)},rotate=180,shift={(0,0)}] (0,1) arc (90:270:1);
\draw[shift={(0,-1)}] (5,5) -- (16,5);
\draw[shift={(0,1)}] (5,5) -- (16,5);
\draw[color=lightgray] (5,5) ellipse (0.25 and 1);
\draw[color=lightgray] (16,5) ellipse (0.25 and 1);
\draw[] (10.5,5) ellipse (0.25 and 1);
\node  at (10.5,5) {$C$ };
\draw [
    thick,
    decoration={
        brace,
        mirror,
        raise=-0.1cm,
	amplitude=7pt
    },
    decorate
] (5,5-1.25) -- (16,5-1.25);
\node  at (11.5,5-1.75) {$(1-\epsilon,1+\epsilon) \times C$ };
\node  at (5-1,5+1) {$M$ };

\begin{scope}[shift={(0,4)}]

%Divisor E_+
\draw[line width=2pt] (16+0.5,1+0.5) -- (16+0.5,1-0.5);
\node at (16+0.8,1) {\tiny $E_+$};

%Divisor E_-
\draw[line width=2pt] (5-0.4,1+0.6) -- (5-0.4,1-0.6);
\node at (5-0.7,1) {\tiny $E_-$};

%\draw[] (10.5+0.1,1-0.5) ellipse (0.25*0.2 and 0.25);
%\draw[] (10.5+0.1,1+0.5) ellipse (0.25*0.2 and 0.25);
\draw[shift={(1,0)},color=lightgray] (10.5+0.1,1-0.5) ellipse (0.25*0.2 and 0.25);
\draw[shift={(1,0)},color=lightgray] (10.5+0.1,1+0.5) ellipse (0.25*0.2 and 0.25);

%upper leg
\draw (5-0.3,1+0.5-0.25) -- (5+8.5,1+0.5-0.25);
\draw (5-0.3,1+0.5+0.25) -- (5+10,1+0.5+0.25);
\draw (5-0.3,1+0.5+0.25) to[out=180,in=180] (5-0.3,1+0.5-0.25);

%lower leg
\draw (5-0.3,1-0.5-0.25) -- (5+10,1-0.5-0.25);
\draw (5-0.3,1-0.5+0.25) -- (5+8.5,1-0.5+0.25);
\draw (5-0.3,1-0.5-0.25) to[out=180,in=180] (5-0.3,1-0.5+0.25);

\draw (10.5+3,1-0.5+0.25) to[out=0,in=-90] (10.5+3+2,1);
\draw (10.5+3,1+0.5-0.25) to[out=0,in=90] (10.5+3+2,1);

\draw (10.5+4.5,1-0.5-0.25) to[out=0,in=-90] (16+0.5,1);
\draw (10.5+4.5,1+0.5+0.25) to[out=0,in=90] (16+0.5,1);

\node at (16-0.2,1) {$u_i$};

%\draw (10.5+4.5,1+0.5+0.25) to[out=0,in=90] (16+0.5,1);

\end{scope}
%end scope

%broken
\node at (5,1) {};
%\draw[shift={(5,1)},rotate=0,shift={(0,0)}] (0,1) arc (90:270:1);
\draw[shift={(16,1)},rotate=180,shift={(0,0)}] (0,1) arc (90:270:1);
%\draw[shift={(0,-1)}] (5,1) -- (8,1);
%\draw[shift={(0,1)}] (5,1) -- (8,1);
%\draw[shift={(0,-1)},dashed] (8,1) -- (9.5,1);
%\draw[shift={(0,1)},dashed] (8,1) -- (9.5,1);
%\draw[shift={(0,1)},dotted,color=lightgray] (9.5,1) -- (10.5,1);
%\draw[shift={(0,-1)},dotted,color=lightgray] (9.5,1) -- (10.5,1);
\draw[shift={(0,-1)}] (10.5,1) -- (16,1);
\draw[shift={(0,1)}] (10.5,1) -- (16,1);
%\draw[shift={(0,-1)},dashed] (11.5,1) -- (13,1);
%\draw[shift={(0,1)},dashed] (11.5,1) -- (13,1);
%\draw[shift={(0,1)},dotted,color=lightgray] (10.5,1) -- (11.5,1);
%\draw[shift={(0,-1)},dotted,color=lightgray] (10.5,1) -- (11.5,1);
%\draw[color=lightgray] (5,1) ellipse (0.25 and 1);
\draw[color=lightgray] (16,1) ellipse (0.25 and 1);
\draw[] (10.5,1) ellipse (0.25 and 1);
\draw [
    thick,
    decoration={
        brace,
        mirror,
        raise=-0.1cm,
	amplitude=7pt
    },
    decorate
] (10.5,1-1.25) -- (16,1-1.25);
\node  at (13.25,1-1.75) {$[1,1+\epsilon) \times C$ };
\node  at (10.5-1,1+1) {$M _+ \cup C$ };
\node  at (10.5,1) {$C$ };

%Divisor E_+
\draw[line width=2pt] (16+0.5,1+0.5) -- (16+0.5,1-0.5);
\node at (16+0.8,1) {\tiny $E_+$};

%J holomorphic curve
\draw[] (10.5+0.1,1-0.5) ellipse (0.25*0.2 and 0.25);
\draw[] (10.5+0.1,1+0.5) ellipse (0.25*0.2 and 0.25);
\draw[shift={(1,0)},color=lightgray] (10.5+0.1,1-0.5) ellipse (0.25*0.2 and 0.25);
\draw[shift={(1,0)},color=lightgray] (10.5+0.1,1+0.5) ellipse (0.25*0.2 and 0.25);

\draw (10.5+0.1,1-0.5+0.25) -- (10.5+3,1-0.5+0.25);

\draw (10.5+0.1,1+0.5-0.25) -- (10.5+3,1+0.5-0.25);
\draw (10.5+0.1,1-0.5-0.25) -- (10.5+4.5,1-0.5-0.25);
\draw (10.5+0.1,1+0.5+0.25) -- (10.5+4.5,1+0.5+0.25);

\draw (10.5+3,1-0.5+0.25) to[out=0,in=-90] (10.5+3+2,1);
\draw (10.5+3,1+0.5-0.25) to[out=0,in=90] (10.5+3+2,1);

\draw (10.5+4.5,1-0.5-0.25) to[out=0,in=-90] (16+0.5,1);
\draw (10.5+4.5,1+0.5+0.25) to[out=0,in=90] (16+0.5,1);

\node at (16-0.2,1) {$\overline{u}_\infty$};

\draw[->] (8,1+0.2) -- (10.5,1+0.5);
\draw[->] (8,1-0.2) -- (10.5,1-0.5);
\node at (7,1) {Reeb orbits};

%\draw[] (10.5,1) ellipse (0.25 and 1);
%\draw [line width = 3pt,domain=25:55] plot ({10.5+0.25*cos(\x)},{1+sin(\x)});

\end{tikzpicture}

Let us return to our example
where $M = \C \P^1$ and $C = \R \P^1$.
We will look at the case where
$u_i : \C \P^1 \to M$ is a biholomorphism.
Then we can assume that $u_\infty : M_+ \to M_+$
is the identity map.
Here $\overline{u}_\infty : M_+ \cup C \to M_+ \cup C$
is the identity map.
In this case $\overline{u}_\infty$
has a negative end converging to the Reeb orbit
that wraps around $C = \R \P^1$ once.

\qed

\bigskip

We will now talk about {\it Gromov-Witten invariants}.
We have the following problem:
{\it How can we show that there is a $J$-holomorphic curve
representing a particular class $[A]$ for any almost complex structure
$J$ compatible with our symplectic form?}
In particular, how can we set things up so that
the conditions of Lemma \ref{label:neckstretchingdemonstration} hold?
The idea here is to use {\bf Gromov-Witten theory}.
Genus $0$ Gromov-Witten invariants for general symplectic manifolds
have now been defined in many different ways:
\cite{FukayaOno:Arnold},
\cite{CieliebakMohnke:symplectichypersurfaces},
\cite{HoferWysockiZehnder:polyfoldapplications1} and \cite{LiTian:sympGW}.
Earlier work for special symplectic manifolds such as symplectic manifolds
of real dimension $6$ or less are done in
\cite{Ruan:topologicalsigma}, \cite{Ruan:3folds} and
\cite{RuanTian:quantum}.
Because the proof of Corollary \ref{corollary:fillingissmooth}
only involves symplectic manifolds of dimension $6$,
one can use these earlier works in this special case.
These invariants can also be defined in a purely algebraic way
 \cite{LiTian:algGW}, \cite{BehrendFantechi:normalcone}
and \cite{Behrend:GW} but we will not use these theories here.
The definitions involved in Gromov-Witten theory are very technical,
and so we only state the properties that we need.

\begin{theorem} \label{theorem:gw}
 (\cite[Theorem 1.3]{FukayaOno:Arnold},
\cite[Theorem 1.12, and the following paragraph]{HoferWysockiZehnder:polyfoldapplications1}
or
\cite[Theorem 2.5]{LiTian:sympGW}).
Let $(M,\omega)$ be a compact symplectic manifold
with a class $[A] \in H_2(M;\Z)$
satisfying $c_1(M,\omega)([A]) + n-3 = 0$,
then we can assign an invariant
$\text{GW}_0(M,[A],\omega) \in \Q$ satisfying the following properties:
\begin{enumerate}
\item  \label{item:holomorphiccurveexistence}
If $\text{GW}_0(M,[A],\omega)\neq 0$ then
there exists a compact nodal $J$-holomorphic curve representing $[A]$
for any almost complex structure $J$ compatible with $\omega$.
\item \label{item:gwdeformationinvariance}
Given a smooth family of symplectic forms $(\omega_t)_{t \in [0,1]}$ on $M$
with $\omega_0 = \omega$, then
\[\text{GW}_0(M,[A],\omega_0) = \text{GW}_0(M,[A],\omega_1).\]
\item \label{item:actualcount}
\cite[Theorem 3.3.1 and Theorem 7.1.8]{McDuffSalamon:Jholomorphiccurves}.
Suppose that $M$ admits an almost complex structure $J$ compatible with $\omega$
so that $(M,J)$ is biholomorphic to a complex manifold and so that
for all arithmetic genus $0$ $J$-holomorphic curves $u : \Sigma \to M$, 
\begin{itemize}
\item $u$ is smooth and
\item
$u^*(TM)$ is a direct sum of complex line bundles of degree $\geq -1$,
\end{itemize}
then $\text{GW}_0(M,[A],\omega)$ is equal to the number of connected genus $0$ $J$-holomorphic curves representing $[A]$.
\end{enumerate}
\end{theorem}

Here $\text{GW}_0(M,[A],\omega) \in \Q$ is a `count' of genus zero holomorphic curves representing $[A]$
for a fixed almost complex structure $J$ compatible with $\omega$. In many cases this is not an actual count,
as the numbers can be negative or non-integer valued. Having said that in some cases
this is an actual count (such as in part (\ref{item:actualcount}) of Theorem \ref{theorem:gw}).
The formula $c_1(M,\omega)([A]) + n-3$ is the `dimension' of the
space of $J$-holomorphic curves representing $[A]$
(called the {\it virtual dimension}).
We want this to be zero so we can `count' the number of these curves.

%TODO:continue here.

We will now talk about {\it symplectic dilation}.
Another problem we will come across is that if $C \subset M$ is a contact hypersurface
with associated contact structure $\xi$ and $\beta$ is some contact form associated to $\xi$
(but not necessarily associated to the contact embedding),
then we would like to find a contact embedding $(C,d\beta) \hookrightarrow (M,\omega)$.
This is needed so that we can apply Lemma \ref{label:neckstretchingdemonstration}.
Quite often this is not possible and so we have to change the symplectic form.
The following definition tells us how to change the symplectic form and
Lemma \ref{lemma:stretchingembedding} below tells us how to construct
such a contact embedding.

\begin{defn} \label{defn:symplecticdilation}
Let $(M,\omega)$ be a symplectic manifold so that $C \subset M$ is a contact hypersurface
with associated contact form $\alpha$ and so that $M \setminus C$ has two connected components $M_-$
and $M_+$.
A {\bf symplectic dilation of $\omega$ along $C$ with respect to $\alpha$} is a family of symplectic forms $(\omega_t)_{t \in [0,\infty)}$
constructed as follows:
First we choose a neighborhood
$C \times (1-\epsilon,1+\epsilon)$ of $C$ in $M$
so that $\omega = d(r \alpha)$ as in Lemma \ref{lemma:contacthypersurfaceneighborhood}.
We will assume that $M_+$ contains $C \times (1,1+\epsilon)$.
Choose a smooth family of non-decreasing functions $\big(\rho_t : (1-\epsilon,1+\epsilon) \to (0,\infty)\big)_{t \in [0,\infty)}$
equal to $\frac{1}{t+1}$ near $1-\epsilon$ and equal to $1$ inside $(\epsilon,1+\epsilon)$.
We define
\[\omega_t = \left\{ \begin{array}{cc}
\omega & \text{inside} \quad M_+ \\
d(\rho_t(r) r \alpha) & \text{inside} \quad C \times (1-\epsilon,1+\epsilon) \\
\frac{1}{t+1} \omega & \text{inside} \quad M_- \setminus (C \times (1-\epsilon,1))
\end{array}
\right.
.
\]
Note that $\omega_t$ is a smooth family of symplectic forms with $\omega_0 = \omega$.
The {\bf support} of our dilation is the region $C \times (1-\epsilon,1+\epsilon) \subset M$.
The $1$-forms $\rho_t(r) r \alpha$ inside $C \times (1-\epsilon,1+\epsilon) \subset M$
will be called the {\bf stretching $1$-forms} of our symplectic dilation.
\end{defn}

We will now explain why this is called {\it symplectic dilation}.
The symplectic manifold $(M,\omega_t)$ as in Definition \ref{defn:symplecticdilation},
can be constructed in the following way (up to symplectomorphism):
Our manifold $M_t$ is the smooth manifold obtained by gluing
$M_-$, $(\frac{1}{1+t}(1-\epsilon), 1 + \epsilon) \times C$ and $M_+$ together
by identifying each $(x,r) \in (1-\epsilon,1) \times C \subset M_-$
with $(x,\frac{1}{1+t}r) \in (\frac{1}{1+t}(1-\epsilon), 1 + \epsilon) \times C$
and each $(x,r) \in (1,1+\epsilon) \times C \subset M_+$
with $(x,r)  \in (\frac{1}{1+t}(1-\epsilon), 1 + \epsilon) \times C$.
The symplectic form is $\frac{1}{t+1} \omega$ on $M_-$,
$d(r \alpha)$ on $C \times (\frac{1}{1+t}(1-\epsilon),1+\epsilon)$ and
$\omega$ on $M_+$.
Here is a picture:

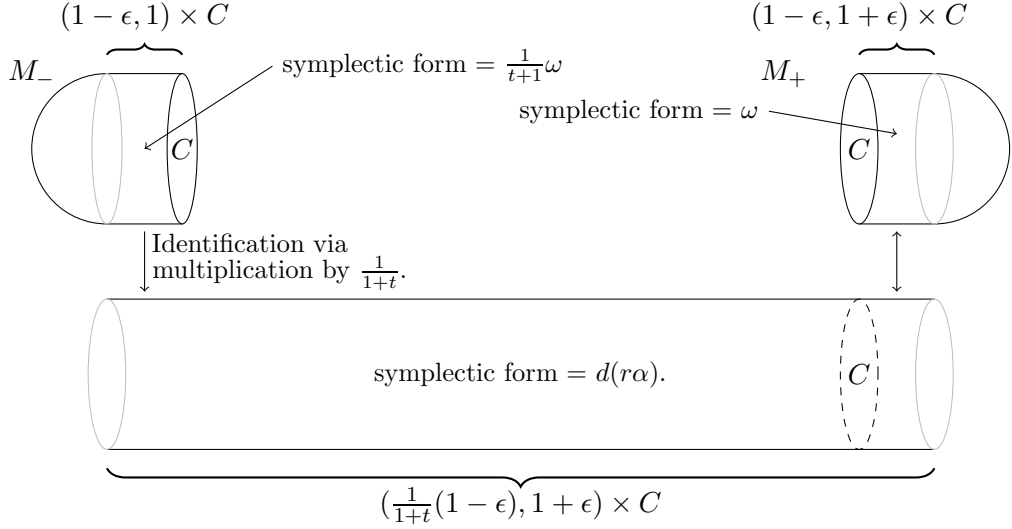
\scalebox{1.0}
{
\begin{tikzpicture} [scale=1.0]
% Unstretched
\draw[shift={(5,5)},rotate=0,shift={(0,0)}] (0,1) arc (90:270:1);
%\draw[shift={(6,5)},rotate=180,shift={(0,0)}] (0,1) arc (90:270:1);
\draw[shift={(0,-1)}] (5,5) -- (6,5);
\draw[shift={(0,1)}] (5,5) -- (6,5);
\draw[color=lightgray] (5,5) ellipse (0.2 and 1);
\draw[] (6,5) ellipse (0.2 and 1);
%\draw[] (5.5,5) ellipse (0.2 and 1);
\node  at (6,5) {$C$ };

\draw [
    thick,
    decoration={
        brace,
        raise=0.0cm,
	amplitude=4pt
    },
    decorate
] (5,5+1.25) -- (6,5+1.25);
\node  at (5.5,5+1.75) {$(1-\epsilon,1) \times C$ };
\node  at (5-1,5+1) {$M_-$ };

%Partially stretched
%\draw[shift={(9,5)},rotate=0,shift={(0,0)}] (0,1) arc (90:270:1);
\draw[shift={(16,5)},rotate=180,shift={(0,0)}] (0,1) arc (90:270:1);
\draw[shift={(0,-1)}] (15,5) -- (16,5);
\draw[shift={(0,1)}] (15,5) -- (16,5);
\draw[color=lightgray] (16,5) ellipse (0.25 and 1);
%\draw[color=lightgray] (16,5) ellipse (0.25 and 1);
\draw[] (15,5) ellipse (0.25 and 1);
\node  at (15,5) {$C$ };
\draw [
    thick,
    decoration={
        brace,
        raise=0.0cm,
	amplitude=4pt
    },
    decorate
] (15,5+1.25) -- (16,5+1.25);
\node  at (15,5+1.75) {$(1-\epsilon,1+\epsilon) \times C$ };
\node  at (15-1,5+1) {$M_+$ };

%broken
%\draw[shift={(5,1)},rotate=0,shift={(0,0)}] (0,1) arc (90:270:1);
%\draw[shift={(16,1)},rotate=180,shift={(0,0)}] (0,1) arc (90:270:1);
\draw[shift={(0,-1)}] (5,2) -- (16,2);
\draw[shift={(0,1)}] (5,2) -- (16,2);
%\draw[shift={(0,-1)},dashed] (8,1) -- (9.5,1);
%\draw[shift={(0,1)},dashed] (8,1) -- (9.5,1);
%\draw[shift={(0,1)},dotted,color=lightgray] (9.5,1) -- (10.5,1);
%\draw[shift={(0,-1)},dotted,color=lightgray] (9.5,1) -- (10.5,1);
%\draw[shift={(0,-1)}] (13,1) -- (16,1);
%\draw[shift={(0,1)}] (13,1) -- (16,1);
%\draw[shift={(0,-1)},dashed] (11.5,1) -- (13,1);
%\draw[shift={(0,1)},dashed] (11.5,1) -- (13,1);
%\draw[shift={(0,1)},dotted,color=lightgray] (10.5,1) -- (11.5,1);
%\draw[shift={(0,-1)},dotted,color=lightgray] (10.5,1) -- (11.5,1);
\draw[color=lightgray] (5,2) ellipse (0.25 and 1);
\draw[color=lightgray] (16,2) ellipse (0.25 and 1);
\draw[dashed] (15,2) ellipse (0.25 and 1);
\draw [
    thick,
    decoration={
        brace,
        mirror,
        raise=0.0cm,
	amplitude=7pt
    },
    decorate
] (5,2-1.25) -- (16,2-1.25);
\node  at (10.5,2-1.75) {$(\frac{1}{1+t}(1-\epsilon),1 + \epsilon) \times C$ };
%\node  at (5-1,1+1) {};
\node  at (15,2) {$C$ };

%identification arrows
\draw[->] (5+0.5,5 - 1.1) -- (5+0.5,2 + 1.1);
\draw[<->] (16-0.5,5 - 1.1) -- (16-0.5,2 + 1.1);
\node  at (5+1.9,3.5+0.25) {\small Identification via};
%\node  at (5+1.9,3.5) {\small via multiplication};
\node  at (5+2.3,3.5-0.15) {\small multiplication by $\frac{1}{1+t}$.};

\draw[<-] (5+0.5,5) -- (5+2.2,5 + 1.1) node[right] { \small symplectic form = $\frac{1}{t+1} \omega$};
\draw[<-] (16-0.5,5.2) -- (16-2.2,5 + 0.5) node[left] { \small symplectic form = $\omega$};

\node  at (10.5,2) {\small symplectic form = $d(r\alpha)$.};
\end{tikzpicture}
}

The following lemma explains how to re-embed contact hypersurfaces
into symplectic dilations once we change the contact form.

\begin{lemma} \label{lemma:stretchingembedding}
Let $(M,\omega)$ be a symplectic manifold and let $C \subset M$ be a contact hypersurface with  associated contact structure $\xi$
and $\beta$ some choice of contact form associated to $\xi$.
Let $(\omega_t)_{t \in [0,\infty)}$ be a symplectic dilation along $C$.
Then for all $t$ sufficiently large
there a contact embedding $\iota : (C,d(c\beta)) \hookrightarrow (M,\omega_t)$
for some constant $c>0$
with the property that $\iota(C)$ is isotopic to $C$ through contact hypersurfaces
all contained in the support of our symplectic dilation.
\end{lemma}

By using Lemma \ref{lemma:contacthypersurfaceneighborhood}
combined with the above construction of $M_t$,
one can prove this Lemma by embedding $C$ as contact hypersurface in
$(\frac{1}{t+1}(1-\epsilon),1+\epsilon) \times C$ with associated contact form $c\beta$ for $c>0$ small.
Here we give a direct proof.

\proof
Let $C \times (1-\epsilon,1+\epsilon)$ our neighborhood of $C$ where we perform our symplectic dilation
with associated function $\rho_t : (1-\epsilon,1+\epsilon) \to (0,\infty)$
exactly as in Definition \ref{defn:symplecticdilation}.
We have that $\beta = f \alpha$ for some $f > 0$ (because we are assuming $\xi$ is cooriented and $\beta$ respects coorientation).
Choose a constant $c>0$ so that $c f < 1$.
Choose $t$ large enough so that $\frac{1}{t+1}(1-\epsilon) < \text{inf}_{x \in C}(cf(x))$.
Define $q_t : (1-\epsilon,1) \to (\frac{1}{t+1}(1-\epsilon),1)$ by $q_t(x) = \rho_t(x) x$.
Here $q_t$ is invertible.
Define $\iota : C \to C \times (1-\epsilon,1)$ by
$\iota(x) = (x, q_t^{-1}(c f(x)))$.
Then $\iota$ is our contact embedding.
This contact embedding is smoothly deformation equivalent to $C \subset M$
through contact embeddings via a smooth family of maps of the form
$x \to (x,g(x)) \in C \times (1-\epsilon,1+\epsilon)$.
\qed

\bigskip

{\it Sketch of the proof of the harder argument.}
Let $\beta$ be a contact $1$-form associated to $\xi_A$.
We need to show that $\beta$ has a Reeb orbit
either with negative $\text{lSFT}$ index
or with $\text{lSFT}$ index $\leq 2 \text{md}(A,0)$.
First we compactify the line bundle $\widetilde{\pi} : \widetilde{A} \to X$
to a $\P^1$ bundle $\breve{S} := \P( \widetilde{A} \oplus \C )$.
This is a projective variety with a natural embedding in $\P^N$
coming from our embedding $X \subset \P^{N-1}$.
Let $\omega_{\breve{S}}$ be the symplectic form on $\breve{S}$
obtained by restricting the natural Fubini-Study symplectic form
on $\P^N$.
Let $\overline{\pi} : \breve{S} \to X$ be the natural projection map.
We now wish to embed $(C,\beta)$ as a contact hypersurface inside $\breve{S}$.
Recall, the contact hypersurface $L_A \subset \widetilde{A}$
is the the circle bundle of size $\epsilon$
inside $\widetilde{A} \subset \breve{S}$.
We have that $L_A$ divides $\breve{S}$ into two connected regions $\breve{S}_+$ and $\breve{S}_-$
which we will view as codimension $0$ submanifolds with boundary
$\partial \breve{S}_+ = \partial \breve{S}_- = L_A$.
We will assume $\breve{S}_-$ contains $X$.
We define $S$ to be the blowup of $\breve{S}$
at some point $x \in \breve{S}_+ \setminus L_A$.
One can put a symplectic structure $\omega_S$ on $S$
so that the complex structure is compatible $\omega_S$
and so that $\text{Bl}^* \omega_{\breve{S}} = \omega_S$
outside a small neighborhood (disjoint from $\breve{S}_-$) of the exceptional divisor
(see \cite[p. 182]{GriffithsHarris:algeraicgeometry}).
We define $S_\pm := \text{Bl}^{-1}(\breve{S}_\pm)$.
Define $[A] \in H_2(S;\Z)$ to be the proper transform of the fiber
through $x$
(i.e. the class represented by the closure of $\text{Bl}^{-1}(\overline{\pi}^{-1}(\overline{\pi}(x)) \setminus \{x\})$).
The reason why we blow up $\breve{S}_\pm$ is to ensure
we have a class $[A]$ satisfying $c_1(S,\omega_S)([A]) + n-3 = 0$ so that we can use 
Theorem \ref{theorem:gw}.

Let $(\omega_t)_{t \in [0,\infty)}$ be a symplectic dilation
of $S$ along $L_A$ as in Definition \ref{defn:symplecticdilation}.
By Lemma \ref{lemma:stretchingembedding}, there is some $t_{\text{max}} > 0$
and contact embedding $\iota : (L_A,d\beta) \hookrightarrow (S,\omega_{t_{\text{max}}})$
which is smoothly homotopic to $L_A$ through contact embeddings in $S$.

By part (\ref{item:actualcount}) of Theorem \ref{theorem:gw},
one can show that $\text{GW}_0(S,[A],\omega_S) = 1$.
Hence by part
(\ref{item:gwdeformationinvariance}) of
Theorem \ref{theorem:gw}, we have
$\text{GW}_0(S,[A],\omega_t) = 1$ for all $t \geq 0$.
Hence part (\ref{item:holomorphiccurveexistence}) of Theorem \ref{theorem:gw} combined with
Lemma \ref{label:neckstretchingdemonstration} implies that $\beta$
has a Reeb orbit $\gamma$ by using our contact embedding $\iota$.

How do we calculate $\text{lSFT}(\gamma)$?
In this paper, the calculation of $\text{lSFT}(\gamma)$ is contained in
the proof of Lemma \ref{lemma:reeborbits}.
Here are the main ideas of this calculation:
The proof of Lemma \ref{label:neckstretchingdemonstration}
tells us that there is a $J_\infty$-holomorphic curve
$u_\infty : \Sigma_\infty \to S_+ \setminus L_A$
with negative ends converging to some Reeb orbits
$\gamma_1,\cdots,\gamma_l$
for some appropriate almost complex structure $J_\infty$.
It turns out that this curve is somewhere injective
as it intersects the exceptional divisor exactly once
and a slightly harder argument involving a maximum principle
\cite[Lemma 7.2]{SeidelAbouzaid:viterbo} shows that it is irreducible.
%So for generic $J_\infty$, $u_\infty$ is {\it regular}
%(see the comment after Definition \ref{defn:pseudoholomorphiccurve}).
The result in \cite{Dragnev:transversality} then tells us that the
space of somewhere injective curves ${\mathcal M}$
with negative ends converging to
$\gamma_1,\cdots,\gamma_l$
and representing the same Borel-Moore homology class as
$u_\infty$ is a manifold for $C^\infty$ generic
choice of $J_\infty$ (which we can assume).
The dimension of ${\mathcal M}$ is some function of the Chern number
of $u_\infty^* \kappa^*_{M_+}$ and
the Conley-Zehnder indices of $\gamma_1,\cdots,\gamma_l$.
One can show that this dimension is
$2\text{md}(A,0) - \sum_i \text{lSFT}(\gamma_i)$.
%The details of this Chern class calculation are contained in
%the proof of Lemma \ref{lemma:reeborbits}.
Also because $u_\infty \in {\mathcal M}$,
we have that the dimension of ${\mathcal M}$ must be non-negative.
Hence $\sum_i \text{lSFT}(\gamma_i) \leq 2\text{md}(A,0)$.
This implies that if $\text{md}(A,0) < 0$, then $\text{lSFT}(\gamma_i) < 0$
for some $i$  and if $\text{md}(A,0) \geq 0$,
then $\text{lSFT}(\gamma_i) \leq 2\text{md}(A,0)$ for some
$i$.
\qed

\bigskip

We will now give a few short comments explaining how to extend the above proof to the case when $A$ is not a cone singularity.
\smallskip

{\it Very short comment on the easy part of the proof:}
In this case, we do not have such a nice Reeb flow for $\alpha_A$.
A rather involved construction (Theorem \ref{label:nicecontactneighbourhoodexistence})
gives us an appropriate $\alpha_A$.
Let $E_1,\cdots,E_l$ be the exceptional divisors of a resolution $\widetilde{A}$ of $A$.
For each formal sum $V := \sum_i d_i E_i$ where $\cap_{\{i | d_i \neq 0\}} E_i \neq \emptyset$,
we have a family of Reeb orbits $B_V$
which are near $\cap_{\{i | d_i \neq 0\}} E_i$
and which `wrap' around $E_j$ $d_j$ times for each $j$ with $d_j \neq 0$
and which have $\text{lSFT}$ index
$\sum_j d_j (a_j+1) - 2$.
Every Reeb orbit is contained in such a family.
Hence $\text{mi}(\alpha_A) = 2\text{md}(A,0)$.

{\it Very short comment on the hard part of the proof:}
Here we would like to construct some nice compactification $S$ of $\widetilde{A}$
so that we can calculate appropriate Gromov-Witten invariants.
The problem is that there is no nice way of compactifying $\widetilde{A}$.
Instead, we partially compactify $A$
(see Step 2 from the proof of Theorem \ref{theorem:reeborbitlowerboundarounddivisors}).
We do this by looking at a small neighborhood of a divisor $E_i$ with smallest discrepancy.
This is a fibration $\widetilde{\pi} : \text{neighborhood}(E_i) \twoheadrightarrow E_i$ whose fibers are symplectic disks.
We then take $Q := \widetilde{\pi}^{-1}(U)$ for some open subset $U \subset E_i \setminus \cup_{j \neq i} E_j$
and compactify the fibers of $\widetilde{\pi}|_Q$ to $\P^1$ fibers.
Actually this region is deformed so that it is a product (see Definition \ref{defn:productregion}).
If $\widetilde{\pi}$ is chosen appropriately then
this gives us a partial compactification $\breve{S}$.
We then blow up $\breve{S}$ at an appropriate point giving us a sympectic manifold $S$.

Now we need to show that we can continue proving the hard part of our theorem as we did earlier
using this non-compact symplectic manifold $S$.
In particular, we need to make sure that we can define Gromov-Witten invariants.
This is where we develop the notion of a GW triple to deal with this problem
(see Definition \ref{defn:gwtriple}).
The point here is that we try and ensure that all the $J$-holomorphic curves for
certain compatible almost complex structures $J$ stay inside some fixed compact subset of $S$,
and this will enable us to have a good theory of Gromov-Witten invariants
to complete our proof.

In fact, to show that we have a good theory of Gromov-Witten invariants for $S$,
we need to use
certain hypersurfaces called stable Hamiltonian hypersurfaces
(see Appendix A for a definition).
Stable Hamiltonian hypersurfaces are much like contact hypersurfaces.
For instance they have a Reeb flow, and one can neck stretch along them.
To ensure that all $J$-holomorphic curves inside $S$ stay inside a compact
set, one only considers almost complex structures $J$
which have been neck stretched along a specially constructed stable Hamiltonian hypersurface.
This stable Hamiltonian hypersurface is a modified version of the contact hypersurface with specific contact form $\alpha_A$
constructed in the proof of our easier statement.
The reason why this cannot be a contact hypersurface is that
for a particular technical reason, one needs to ensure that the map
$\widetilde{\pi}|_Q$ is holomorphic and this can only be done if we neck stretch along a stable Hamiltonian hypersurface
(see Step 2 from the proof of Theorem \ref{theorem:reeborbitlowerboundarounddivisors}).
Note that in our proof we perform {\it two} neck stretches.
We first neck stretch along an appropriate stable Hamiltonian hypersurface in $S$ to ensure Gromov-Witten invariants are well defined for this new stretched almost complex structure,
and {\it then} we neck stretch along $\iota(C)$ to find our Reeb orbit of the appropriate index.

\section{Minimal Discrepancy of Isolated Singularities} \label{section:minimaldiscrepancydefinition}

The main ideas in this section come from \cite{BFFU:valuation}.
Let $A \subset \C^N$ be a singularity which is isolated at $0$.
First of all, we will give two definitions of a numerically $\Q$-Gorenstein singularity.
One definition will be algebraic involving $\Q$-Cartier divisors
(see \cite{BFFU:valuation}), and the other will be topological involving
the first Chern class of our contact structure $\xi_A$
(See  \cite[Definition 1.2]{Durfee:signature}).

We will start with the algebraic definition and then we will give the topological one and prove they are equivalent.
Start with some resolution $\pi : \widetilde{A} \twoheadrightarrow A$
so that the preimage of $0$ is a union of smooth normal crossing divisors $E_i$ and so that
$\pi$ is an isomorphism away from these divisors
(If $A$ is smooth we blow up at least once, so $\pi$ is never an isomorphism).
Let $K_{\widetilde{A}}$ be the canonical bundle of $\widetilde{A}$ which we will view as
a $\Q$-Cartier divisor.
We say that $A$ is {\bf numerically $\Q$-Gorenstein} if there exists a $\Q$-Cartier divisor
$K^{\text{num}}_{\widetilde{A}/A} := \sum_j a_j E_j$ with the property that
$C \cdot (K^{\text{num}}_{\widetilde{A}/A}  - K_{\widetilde{A}}) = 0$
for any projective algebraic curve $C \subset \pi^{-1}(0)$.
By the negativity Lemma \cite[Theorem 4.39]{KollarMori:birational}
one can show that the coefficients $a_j$ are unique (see \cite[Proposition 5.3]{BFFU:valuation}).
Here $a_j \in \Q$ is called the {\bf discrepancy} of $E_j$.
In the literature, one usually calls the number $1+a_j$
the ``log-discrepancy'' of $E_j$. We won't use this denomination in the paper.

%Before we do this we need a definition:
%\begin{defn}
%Let $(C,\xi)$ be a contact manifold.
%The contact hyperplane distribution $\xi$ on $C$ has a natural symplectic structure given %by restricting $d\alpha$
%and so if we choose a compatible almost complex structure on this contact hyperplane %distribution then
%it has a natural $U(n-1)$ structure. 
%Because this bundle is a complex bundle we can take its highest exterior power.
%Such a bundle is called the {\bf anticanonical bundle} of $(C,\xi)$.
%We define $c_1(\xi)$ to be the first Chern class of the anticanonical bundle.
%The dual of such a bundle is called the {\bf canonical bundle} of $(C,\xi)$.
%\end{defn}
%
%We also have the following Lemma:

Before we give an alternative definition of being numerically $\Q$-Gorenstein, we need the following two Lemmas:

\begin{lemma} \label{lemma:equivalenceoflinebundles}
$c_1(\xi_A) = c_1(TA|_{L_A})$.
\end{lemma}
\proof
Let $\phi = \sum_j |z_j|^2|_A$ where $z_1,\cdots,z_N$ are coordinates for $\C^N$.
The contact structure $\xi_A$ is equal to:
\[\text{ker}(d^c(\phi)) \cap \text{ker}(d\phi)\]
restricted to $L_A$.
This is a complex subbundle of $TA|_{L_A}$.
The subbundle of $TA|_{L_A}$ spanned by the vector fields
$\nabla(\phi)$ and $i \nabla(\phi)$ (with respect to the induced metric on $A \subset \C^N$)
is orthgonal to $\xi_A$ and also a holomorphic subbundle.
It is also trivial.
Hence $TA|_{L_A}$ is equal to $\xi_A$ plus a trivial complex line bundle
and so $c_1(\xi_A) = c_1(TA|_{L_A})$.
\qed

\bigskip

Let $\widetilde{A}_\epsilon := \pi^{-1}(B^{2n}_\epsilon)$ where $B^{2n}_\epsilon \subset \C^n$ is the closed $\epsilon$ ball.
The boundary of $\widetilde{A}_\epsilon$ is equal to $L_A$.
\begin{lemma} \label{lemma:kernallemma}
The natural map $\overline{\mu} : H^2(\widetilde{A}_\epsilon,L_A;\Q) \to \text{ker}\left(H^2(\widetilde{A}_\epsilon;\Q) \to H^2(L_A;\Q)\right)$
is an isomorphism and
$\text{ker}\left(H^2(\widetilde{A}_\epsilon;\Q) \to H^2(L_A;\Q)\right)$
is freely generated by $c_1({\mathcal O}_{\widetilde{A}_\epsilon}(E_i);\Q)$.
\end{lemma}
\proof of Lemma \ref{lemma:kernallemma}.
The long exact sequence:
\[ H^2(\widetilde{A}_\epsilon,L_A;\Q) \to H^2(\widetilde{A}_\epsilon;\Q) \to H^2(L_A;\Q))\]
tells us that
$\overline{\mu}$ is surjective.

We will now show that $\overline{\mu}$ is injective.
We have the Lefschetz duality isomorphism given by $ H_{2n-2}(\widetilde{A}_\epsilon,\Q) \cong H^2(\widetilde{A}_\epsilon,L_A,\Q)$.
For a class $x$ in $H_{2n-2}(\widetilde{A}_\epsilon;\Q)$ we write $\text{LD}(x)$ for its  Lefschetz dual.
A Mayor-Vietoris argument tells us that $H_{2n-2}(\widetilde{A}_\epsilon,\Q)$ is freely
generated by classes $[E_i]$.
Also $\cup_j E_j \hookrightarrow \widetilde{A}_\epsilon$ is a homotopy equivalence.
Hence
$H^2(\widetilde{A}_\epsilon,L_A,\Q)$
is freely generated by classes $\text{LD}([E_i])$.

Now suppose $\overline{\mu}(\sum_i b_i \text{LD}(E_i)) = 0$ for for some $b_1,\cdots,b_l \in \Q$.
Then $\sum_i b_i c_1({\mathcal O}_{\widetilde{A}_\epsilon}(E_i);\Q) = 0$ (as $\overline{\mu}(\text{LD}([E_i])) = c_1({\mathcal O}_{\widetilde{A}_\epsilon}(E_i);\Q)$)
and so $C \cdot (\sum_i b_i E_i) = 0$ for all projective algebraic curves $C \subset \pi^{-1}(0)$.
By the negativity Lemma \cite[4.39]{KollarMori:birational} applied to $\sum_i b_i E_i$ and $\sum_i (-b_i) E_i$,
we get that $b_i \leq 0$ and $b_i \geq 0$ for all $i$ and hence $b_i = 0$ for all $i$.
Hence $\overline{\mu}$ is injective and hence an isomorpism.

Because $\overline{\mu}$ is an isomorphism and
$H^2(\widetilde{A}_\epsilon,L_A,\Q)$
is freely generated by classes $\text{LD}([E_i])$,
we get that
$\text{ker}(H^2(\widetilde{A}_\epsilon;\Q) \to H^2(L_A;\Q))$
is freely generated by $c_1({\mathcal O}_{\widetilde{A}_\epsilon}(E_i);\Q)$.
\qed

\bigskip

We will now give a topological characterization of being numerically $\Q$-Gorenstein for isolated singularities (which we regard as the alternate topological definition).

\begin{lemma} \label{lemma:topologicalnumericallygorenstein}
We have that $A$ is numerically $\Q$-Gorenstein if and only if \[c_1(\xi_A;\Q) = c_1(TA|_{L_A};\Q) = 0 \in H^2(L_A;\Q)\]
(i.e. $c_1(\xi_A)$ is torsion in $H^2(L_A;\Z)$).

Also if $A$ is numerically $\Q$-Gorenstein, then
$c_1(\widetilde{A}_\epsilon;\Q)$ lifts to a unique class in $H^2(\widetilde{A}_\epsilon,L_A;\Q)$
which is Lefschetz dual to
$\sum_j a_j [E_j] \in H_{2n-2}(M;\Q)$.
\end{lemma}

\proof of Lemma \ref{lemma:topologicalnumericallygorenstein}.
%Throughout this Lemma, we use the following fact:
%$c_1(T^* \widetilde{A}_\epsilon|_{L_A};\Q) =  -c_1(\xi_A;\Q)$.
%This follows from the identity
%$c_1(T^* \widetilde{A}_\epsilon|_{L_A};\Q) = -c_1(T \widetilde{A}_\epsilon|_{L_A};\Q)$
%and because
%the subbundle of $T \widetilde{A}_\epsilon|_{L_A}$
%consisting of vectors symplectically orthogonal to $\xi_A \subset T \widetilde{A}_\epsilon|_{L_A}$
%is trivial.
%
%
Suppose first that $A$ is numerically $\Q$-Gorenstein.
Then there exists a $\Q$-Cartier divisor $\sum_j a_j E_j$
with the property that $C \cdot (\sum_j a_j E_j  - K_{\widetilde{A}}) = 0$
for any projective algebraic curve $C \subset \pi^{-1}(0)$.
Because $\cup_i E_i \hookrightarrow \widetilde{A}_\epsilon$ is a homotopy equivalence, we have
$\sum_j a_j c_1({\mathcal O}_{\widetilde{A}_\epsilon}(E_j);\Q)  - c_1({\mathcal O}_{\widetilde{A}_\epsilon}(K_{\widetilde{A}_\epsilon});\Q) = 0 \in H^2(\widetilde{A}_\epsilon;\Q)$
by \cite[Lemma 5.13]{BFFU:valuation}.
If $\nu : H^2(\widetilde{A}_\epsilon,\Q) \to H^2(L_A,\Q)$ is the natural restriction map then
$\nu(c_1({\mathcal O}_{\widetilde{A}_\epsilon}(E_i);\Q) = 0$.
Hence $\nu(c_1({\mathcal O}_{\widetilde{A}_\epsilon}(K_{\widetilde{A}});\Q)) = \nu(c_1({\mathcal O}_{\widetilde{A}_\epsilon}(K_{\widetilde{A}} - \sum_j a_j E_j);\Q)) = 0$.
Now $\nu(c_1({\mathcal O}_{\widetilde{A}_\epsilon}(K_{\widetilde{A}});\Q)) = c_1(T^* \widetilde{A}_\epsilon|_{L_A};\Q) = -c_1(\xi_A;\Q)$
by Lemma \ref{lemma:equivalenceoflinebundles}
and so we get $c_1(\xi_A;\Q) = 0$ which implies $c_1(\xi_A;0) = 0 \in H^2(L_A,\Q)$.

Conversely suppose that $c_1(\xi_A;\Q)=0$.
Then $c_1(T^* \widetilde{A}_\epsilon|_{L_A};\Q) = -c_1(\xi_A;\Q) = 0 \in H^2(L_A;\Q)$ by Lemma \ref{lemma:equivalenceoflinebundles}.
Hence $c_1({\mathcal O}_{\widetilde{A}_\epsilon}(K_{\widetilde{A}});\Q) \in \text{ker}(H^2(\widetilde{A}_\epsilon;\Q) \to H^2(L_A;\Q))$.
Therefore by Lemma \ref{lemma:kernallemma}, there exists $a_1,\cdots,a_l \in \Q$ so that
$\sum_j a_j c_1({\mathcal O}_{\widetilde{A}_\epsilon}(E_j);\Q)  = c_1({\mathcal O}_{\widetilde{A}_\epsilon}(K_{\widetilde{A}_\epsilon});\Q)$.
Hence $C \cdot (\sum_j a_j E_j  - K_{\widetilde{A}}) = 0$
for any projective algebraic curve $C \subset \pi^{-1}(0)$ which implies that
$A$ is numerically $\Q$-Gorenstein.
This argument also implies that
$c_1(\widetilde{A}_\epsilon;\Q)$ lifts to a unique class in $H^2(\widetilde{A}_\epsilon,L_A;\Q)$
which is Lefschetz dual to
$\sum_j a_j [E_j] \in H_{2n-2}(M;\Q)$.
\qed
\bigskip

\begin{defn}
The {\bf minimal discrepancy} $\text{md}(A,0)$ of $A$ is the infimum of $a_j$ over all resolutions $\pi$.
\end{defn}
To calculate the minimal discrepancy, one only needs to look at any fixed resolution $\pi$ of
$A$. If $\pi$ is a fixed resolution and $\pi$ is not the identity map, then
\[ \text{md}(A,0) = \left\{ \begin{array}{cc}
\text{min}_j a_j & \text{if} \quad a_j \geq -1 \quad \forall j \in \{1,\cdots,l\} \\
-\infty & \text{otherwise}
\end{array}
 \right. \]
(see \cite[17.1.1]{kollar1992flips}).
In order to calculate the minimal discrepancy of $A$ when $A$ is smooth at $0$,
you have to blow up at least once, which gives a minimal discrepancy of $n-1$ where $n = \text{dim}_{\C} A$.

\section{The Conley-Zehnder Index} \label{section:conleyzehnderindexsection}

In this section we will define the Conley-Zehnder index of a Reeb orbit.
A restricted definition of this index was introduced in
\cite{arnold:characteristic} and in \cite{ConleyZehnder:Index}.
The definition we will use is contained in \cite{RobbinSalamon:maslov}
and also \cite{Gutt:GeneralizedConleyZehnder}.
We will
then define what a {\it pseudo Morse-Bott family of Reeb orbits} is and then prove some result
showing how the Conley-Zehnder indices of Reeb orbits change
when we perturb the contact form.

%--------------------
%TODO: define Maslov index, CZ index and properties of Maslov index.

%TODO: add axioms somewhere from Robbin Salamon

\subsection{Definition of Conley-Zehnder Index} \label{section:conleyzehnderindex}

We will first give a definition of the Conley-Zehnder index of a path of symplectic matrices in terms
of the Maslov index before we define it for Reeb orbits. This will be useful later on.
Let $\mathcal{L}$ be the set of Lagrangian vector subspaces of a symplectic vector space $Y$.
Fix some $L \in \mathcal{L}$.
To any smooth path $\Lambda : [a,b] \to \mathcal{L}$ we can assign a
{\bf Maslov index} $\text{Mas}_{Y,L}(\Lambda) \in \frac{1}{2} \Z$
(see \cite{RobbinSalamon:maslov}).
We do not need to know the exact definition of this index,
we will just need some properties of this index which we will cite when needed later on.
We will just write $\text{Mas}(\Lambda) = \text{Mas}_{Y,L}(\Lambda)$ when the context is clear.

Let $W$ be a symplectic vector space and $\text{Sp}(W)$ the space of linear symplectomorphisms of $W$.
We will write $\text{Sp}(2n)$ if it is clear which $2n$ dimensional symplectic vector space we are using.
Let $\overline{W} \times W$ be the product symplectic vector space with symplectic form $(-\omega_W,\omega_W)$
where $\overline{W}$ is equal to $W$
and let $\Delta$ be the diagonal Lagrangian.
We define the {\bf Conley-Zehnder index} $\text{CZ}(A) \in \frac{1}{2}\Z$
of a path of symplectic matrices $A : [a,b] \to \text{Sp}(2n)$ as follows
(see \cite[Rem. 5.4]{RobbinSalamon:maslov}):
$\text{CZ}(A) := \text{Mas}_{\overline{W} \times W, \Delta}(\Gamma(A))$
where $\Gamma(A)$ is the path of Lagrangians given by the graph of $A(t)$ in $\overline{W} \times W$
viewed as a map $A(t) : \overline{W} \to W$.
%This is an invariant of $A(t)$ up to homotopies fixing its endpoints
%(see \cite[Theorem 4.1]{RobbinSalamon:maslov}
%or \cite[Theorem 2]{Gutt:GeneralizedConleyZehnder}).
In order to define Conley-Zehnder index for a Reeb orbit,
it only turns out that we need paths which start at the identity.
Having said that we define it for general paths as one way of computing the Conley-Zehnder
index of a path
is to chop up a path into little pieces and sum up the Conley-Zehnder indices
of each piece separately (see property \ref{item:catenation} below).
This is used for instance in Lemma \ref{lemma:perturbedindexcalculation}
(it is also used elsewhere but the paths there all start from $\text{id}$).
Here are a few properties of Conley-Zehnder indices that we will need in this paper
(see
\cite[Theorem 2.3]{RobbinSalamon:maslov} or
\cite[Theorem 55]{Gutt:ConleyZehnderMatricies}):
\begin{CZ}
\item \label{item:homotopyrelendpoints}
If two paths $A_1,A_2$ are homotopic relative to their endpoints
then $\text{CZ}(A_1) = \text{CZ}(A_2)$
\item \label{item:zero}
Any constant path has Conley-Zehnder index $0$.
\item \label{item:catenation}
It is additive under catenation of paths.
\item  \label{item:additiveunderdirectsum}
Consider two paths $A_1 : [a,b] \to \text{Sp}(W_1)$,
$A_2 : [a,b] \to \text{Sp}(W_2)$,
then $\text{CZ}(A_1 \oplus A_2) = \text{CZ}(A_1) + \text{CZ}(A_2)$
where $A_1 \oplus A_2 : [a,b] \to \text{Sp}(W_1 \oplus W_2)$
is the direct sum.
\item \label{item:determinantproperty}
Suppose that $A : [a,b] \to U(n) \subset \text{Sp}(2n)$
with $A(a) = A(b)$ then
$\text{det}(A) : [a,b] \to S^1$ represents an integer $k \in \pi_1(S^1) = \Z$.
We have $\text{CZ}(A) = 2k$
\cite[Section 3]{SZ:morsetheory}.
\end{CZ}
Properties
\ref{item:homotopyrelendpoints},
\ref{item:zero} and
\ref{item:catenation}
 also tell us that the Conley-Zehnder index does not depend on the
parameterization of the path so long as such a parameterization respects orientation of the domain
(the idea here is to enlarge the domains of both the original and reparameterized path
by adding constant paths at each end
so that their domains become identical).
There three other properties of Conley-Zehnder indices that we will need.
One is stated in
Lemma \ref{lemma:homotopyofsymplecticpaths}
and is a direct consequence of a corresponding property for Maslov indices in
\cite[Theorem 2.4]{RobbinSalamon:maslov} below.
The second one is in Lemma
\ref{lemma:indexofsmallpath}
and uses
\cite[Theorem 2.3 (localization axiom)]{RobbinSalamon:maslov}.
The final property is called the
normalization property from \cite[Theorem 55]{Gutt:ConleyZehnderMatricies}
and is used in Lemma \ref{lemma:hamiltonianconleyzehnderindexcomparisonstandard}.

We need a Lemma enabling us to define a Conley-Zehnder index for Reeb orbits.
Let $S_k \subset \text{Sp}(2n)$ be the set of symplectic matrices $A$ with  $\text{dim ker}(A - \text{id}) = k$.
Let $o_1(t), o_2(t)$ be two paths in $\text{Sp}(2n)$.
We say that $o_1,o_2$ are {\bf stratum  homotopic} if there is a smooth family of paths
$\psi_s : [0,1] \to \text{Sp}(2n)$ where $\psi_s(0) \in S_{k_1}, \psi_s(1) \in S_{k_2}$
for all $s$ and some fixed $k_1,k_2$ and $\psi_0,\psi_1$ are homotopic to $o_1,o_2$ respectively relative to their endpoints.
\begin{lemma} \label{lemma:homotopyofsymplecticpaths}
If $o_1$,$o_2$ are stratum homotopic then they have the same Conley-Zehnder indices.
\end{lemma}
\proof %of Lemma \ref{lemma:homotopyofsymplecticpaths}.
This follows directly from \cite[Theorem 2.4]{RobbinSalamon:maslov}.
\qed

\bigskip

We will now define the Conley-Zehnder index for Reeb orbits.
Let $(C,\xi)$ be a cooriented contact manifold of dimension $2n-1$ and let $\alpha$ be a contact form
respecting the coorientation with $\text{ker}(\alpha) = \xi$.
Fix some complex structure on the bundle $\xi$ compatible with $d\alpha|_\xi$.
Let $\kappa$ be the canonical bundle of $(C,\xi)$ and $\kappa^*$ the anticanonical bundle.
%Recall that a closed Reeb orbit is a closed orbit of the Reeb vector field associated to $\alpha$.
%From now on we will call such an orbit a {\bf Reeb orbit} as opposed to a closed Reeb orbit
%because  we will not be dealing with non-compact flowlines of the Reeb vector field.
From now on we will assume that $c_1(\xi;\Q)=0 \in H^2(C;\Q)$.
This means there is some number $N \in \N$ so that $N c_1(\xi) = 0 \in H^2(C;\Z)$.
Therefore we can trivialize the $N$th power of the anticanonical bundle.
Let $\tau : (\kappa^*)^{\otimes N} \to C \times \C$ be a choice of such a trivialization.
Let $\gamma : \R / L \Z \to C$ be a Reeb orbit of $\alpha$.
Choose a trivialization of
$\tau_\gamma : \gamma^* \oplus_{j=1}^{N} \xi \to C \times \C^{(n-1)N}$
as a Hermitian vector bundle
so that its highest complex exterior
power coincides with our trivialization $\tau$.
Such a choice does not depend up to homotopy on our choice of
$d\alpha|_\xi$ compatible complex structure because the space
of such structures is connected.
Let $\phi_t : C \to C$ be the flow of the Reeb vector field of $C$.
Let $D\phi_t|_\xi : \xi \to \xi$ be the restriction of the linearization $D\phi_t : TC \to TC$
of the Reeb flow $\phi_t$ of $\alpha$.
Then using the above trivialization along $\gamma$ we have that $\oplus_{j=1}^{N} D\phi_t|_\xi$
gives us a smooth path of symplectic matrices $A_\gamma : [0,L] \to \text{Sp}(\R^{(2n-2) N})$ given by
\begin{equation} \label{equation:trivializationfamiliy}
A_\gamma(t) := \text{pr} \circ \tau_\gamma \circ \left(\oplus_{j=1}^{N} D\phi_t|_\xi \right) \circ \left(\text{pr} \circ \tau_\gamma|_{\gamma^* \oplus_{i=1}^N \xi_{\gamma(0)}}\right)^{-1}
\end{equation}
where $\text{pr} : C \times \C^{(n-1)N} \to \C^{(n-1)N} = \R^{(2n-2)N}$ is the natural projection map.
%
%We view this as a path parameterized by $[0,L]$ under the surjection $[0,L] \twoheadrightarrow \R / L \Z$.
%
\begin{defn}
Define the {\bf Conley-Zehnder index of $\gamma$}, $\text{CZ}_\tau(\gamma) \in \frac{1}{2N}\Z \subset \Q$,
to be $\frac{1}{N}\text{CZ}(A_\gamma)$.
If $H^1(C,\Q) = 0$ then we define $\text{CZ}(\gamma) := \text{CZ}_\tau(\gamma)$
(see the Lemma below).
\end{defn}

\begin{lemma}
If $H^1(C,\Q) = 0$ then $\text{CZ}_\tau(\gamma)$ does not depend on $\tau$ or $N$.
\end{lemma}
\proof
Suppose that we have some other choice of trivialization $\upsilon : (\kappa^*)^{\otimes N_1} \to C \times \C$ as above
where $N_1 \in \N$ then
$\tau$ and $\upsilon$ induce natural trivializations $\tau^{\otimes N_1}$
and $\upsilon^{\otimes N}$ of $(\kappa^*)^{\otimes N N_1}$.
By property
\ref{item:additiveunderdirectsum} we then get that
$\text{CZ}_{\tau}(\gamma) = \text{CZ}_{\tau^{\otimes N_1}}(\gamma)$
and similarly $\text{CZ}_{\upsilon}(\gamma) = \text{CZ}_{{\upsilon}^{\otimes N}}(\gamma)$.
Now $\tau^{\otimes N_1} \circ ({\upsilon}^{\otimes N})^{-1}$ is a bundle morphism between trivial bundles
and hence a section of the trivial bundle $C \times \C^*$ which is equivalent to a smooth map from $C \to \C^*$.
Because the pullback of $d\vartheta$ via $C \to \C^*$ is exact (as $H^1(C,\R) = 0$) we get that such a smooth map is homotopic to the constant
map and this means that $\tau^{\otimes N_1}$ is homotopic to ${\upsilon}^{\otimes N}$ through trivializations.
Let $(\tau_s)_{s \in [0,1]}$ be this family of trivializations
joining $\tau^{\otimes N_1}$ and ${\upsilon}^{\otimes N}$.
Then each $\tau_s$ gives us a family of symplectic matrices
$A_{\gamma,s}(t)$ as in Equation \ref{equation:trivializationfamiliy}.
Here $A_{\gamma,s}$ are stratum homotopic and so have the same
Conley-Zehnder index by Lemma \ref{lemma:homotopyofsymplecticpaths}.
Hence $\text{CZ}_{\tau^{\otimes N_1}}(\gamma) = \text{CZ}_{{\upsilon}^{\otimes N}}(\gamma)$
by \ref{item:homotopyrelendpoints}.
And so $\text{CZ}_{\tau}(\gamma) = \text{CZ}_{\upsilon}(\gamma)$.
\qed

\subsection{Pseudo Morse-Bott families}

Recall that a smooth function $f : X \to \R$
has a Morse-Bott family of critical points $B \subset X$
if $B$ is a submanifold
and the Hessian of $f$ at each $b \in B$ has kernal equal to $T_b B$.
We say $f$ is Morse-Bott if every critical point sits inside a Morse-Bott family.
The reason why the above definition is useful is because
many manifolds have natural Morse-Bott functions
 and it useful to exploit
such functions to tell us what the homology groups of the above manifold is
using Morse-Bott homology.

If we have a manifold $C$ with contact form $\alpha$,
then there is a natural function
sending loops $\gamma : S^1 \to C$
to $\int \gamma^* \alpha$.
The critical points of such a function are Reeb orbits.
Morally, contact homology is Morse homology of such a function
(in reality this is not true and there are many difficult problems to overcome),
and therefore in order to calculate contact homology one needs to know
when $\gamma \to \int \gamma^* \alpha$
is `Morse-Bott'.

\begin{defn} \label{defn:morsebott}
A {\bf Morse Bott family of Reeb orbits of $(C,\alpha)$ of period $T$}
is a closed submanifold $B \subset C$
where $B$ is a union of closed Reeb orbits of period $T$,
and if $\phi_t : C \to C$ is the Reeb flow then
$\text{Ker}(D\phi_T)|_B = TB$
(see \cite{Bourgeois:morsebott}).
\end{defn}

The families of critical points that we will be dealing with in Section
\ref{subsection:specificcontactform}
will be manifolds with corners.
It is not known how to define `Morse-Bott manifolds with corners'.
In this paper we do not need such a strong definition.
We are only interested in the indices of such Reeb orbits
and so we are content with the definition of a
`pseudo Morse-Bott submanifold'.

As motivation, here is a finite dimensional version:
\begin{defn}
Let $f : X \to \R$ be a smooth function.
Suppose that we have a connected set of critical points
$B \subset X$ which is isolated and contained in a level set of $f$.
Then $B$ is said to be {\bf pseudo Morse-Bott}
if the kernal of the Hessian has constant rank along $B$.
\end{defn}
This is a more general definition than just being a (connected) Morse-Bott family of critical points.
The good thing about this definition is that if one perturbs $f$
by a small amount so that it is a Morse function then one can estimate
the indices of the critical points near $B$
by looking at the Hessian.
We wish to do exactly the same thing for the function
$\gamma \to \int \gamma^* \alpha$.

\begin{defn}
Suppose $B_T$ is a subset of $C$ so that:
\begin{enumerate}
\item $B_T$ is a union Reeb orbits of period $T$.
\item There is a neighborhood $\mathcal{N}_{B_T}$ of $B_T$ and a constant $\delta>0$ so that
any Reeb orbit
with period in the interval $[T-\delta,T+\delta]$ meeting $\mathcal{N}_{B_T}$
is in fact contained in $B_T$ and has period $T$.
\end{enumerate}
Then we say $B_T$ {\bf is an isolated family of Reeb orbits of period} $T$.
We will say that $\mathcal{N}_{B_T}$ is {\bf an isolating neighborhood for} $B_T$.
\end{defn}
%Note that the union of two isolated families of Reeb orbits with the same period $T$ is an isolated family of Reeb orbits %of period $T$.

%An example of an isolated family of period $T$ is a Morse-Bott family of Reeb orbits of period $T$.
%We will be dealing with something more general than this.
%Let $\phi_t : C \to C$ be the time $t$ flow of $R_\alpha$ where $R_\alpha$ is the Reeb vector field.
\begin{defn} \label{defn:morsebottfamily}
A {\bf pseudo Morse-Bott family} is an isolated family of Reeb orbits $B_T$ of period $T$ with the additional property
that $B_T$ is path connected and for each point $p \in B_T$ we have
$\text{Size}_p(B_T) := \text{dim ker}(D_p\phi_T|_\xi - \text{id})$ is constant along $B_T$
(recall that $D_p\phi_t|_\xi : \text{ker}(\alpha)_p \to \text{ker}(\alpha)_{\phi_t(p)}$ is the restriction
of the linearization of $\phi_t$ to $\text{ker}(\alpha)_p$).
\end{defn}
An example of a pseudo Morse-Bott family of period $T$ is a Morse-Bott family of Reeb orbits of period $T$.

We are interested in indices of Reeb orbits and so
from now on we assume that we work with a fixed trivialization of a fixed power 
of the canonical bundle of $(C,\alpha)$.

By Lemma \ref{lemma:homotopyofsymplecticpaths} we have that the Conley-Zehnder index of the
period $T$ orbits starting in $B_T$ are all the same because $B_T$ is path connected.
Hence we define the {\bf Conley-Zehnder index of} $B_T$, $\text{CZ}(B_T)$, to be the Conley-Zehnder index of one of its period $T$ Reeb orbits.

\subsection{Perturbations of Contact Forms and Indices}

Suppose we have a contact form and we wish to perturb it,
then this new contact form has new Reeb orbits that are near the old Reeb orbits.
We wish to relate the indices of these new Reeb orbits with the indices of the old
ones.
The aim of this section is to give us such a relation.
The reason why we need such a relation is that
in order to bound the minimal discrepancy of an isolated singularity
from below, one needs to find a Reeb orbit of a particular Conley-Zehnder index for any
contact form associated to the contact distribution on our link.
To find such Reeb orbits we need to use $J$-holomorphic curves,
and such curves find Reeb orbits of the right index
when the contact form is $C^\infty$ {\it generic}.
So if we have a non-generic contact form, we need to show that
the orbits of a nearby generic contact form have related indices.
We need a linear algebra lemma and corollary first.

Let $Y$ be a symplectic vector space and $\mathcal{L}$ the set of linear Lagrangians inside $Y$.
Let $L \in \mathcal{L}$ be a fixed Lagrangian.
Define $\mathcal{L}_k$ to be the set of Lagrangians in $Y$ whose intersection with our fixed Lagrangian $L$
has dimension $k$.
\begin{lemma} \label{lemma:indexofsmallpath}
Fix $\Lambda_0 \in\mathcal{L}_k$. Then for a sufficiently small neighborhood $\mathcal{N}_{\Lambda_0}$
of $\Lambda_0$ we have that any path $\Lambda : [a,b] \to \mathcal{N}_{\Lambda_0}	$ with $\Lambda(a) = \Lambda_0$ has Maslov index
in $[-\frac{k}{2},\frac{k}{2}]$.
\end{lemma}
\proof of Lemma \ref{lemma:indexofsmallpath}.
First of all we can identify our symplectic vector space $Y$ with $T^* L$ via a linear symplectomorphism so that:
\begin{enumerate}
\item The zero section is identified with $L$.
\item By choosing a linear coordinate system $x_1,\cdots,x_n$ on $L$
we get associated coordinates $x_1,\cdots,x_n,y_1 dx_1,\cdots, y_n dx_n$ on $T^* L$
and this identifies $T^* L$ with $L \oplus L$ in a natural way.
We require that the Lagrangian $\Lambda_0$ is the graph of a symmetric
matrix $A$ with respect to this coordinate system.
Another way of seeing this is viewing $L$ as the graph of the differential of the
function $f_A(x) = \frac{1}{2} x^T A x$.
\end{enumerate}
Let $Z \subset L$ be a subspace of $L$ so that $(L \cap \Lambda_0) \oplus Z = L$.
We get that $A|_Z$ is non-degenerate as a quadratic form.
We choose $\mathcal{N}_{\Lambda_0}$ small enough so that every element of this set can also be expressed as the graph
of a symmetric matrix whose restriction to $Z$ is non-degenerate.
If $\Lambda : [a,b] \to \mathcal{N}_{\Lambda_0}$ is a smooth path starting at $\Lambda_0$ then this is represented by a smooth family of symmetric
matrices $(A_t)_{t \in [a,b]}$ with $A_a = A$ whose restriction to $Z$ is non-degenerate.
By the localization axiom from \cite[Theorem 2.3]{RobbinSalamon:maslov} we then get that the Maslov index is
$\frac{1}{2}(\text{sign}(A_b) - \text{sign}(A_a))=\frac{1}{2}(\text{sign}(A_b) - \text{sign}(A))$.
Here $\text{sign}$ means the sign of the symmetric matrix as a quadratic from.
Now $\text{sign}(A_t) = \text{sign}(A_t|_Z) + \text{sign}(A_t|_{L \cap \Lambda_0})$.
Because $A_t|_Z$ is a smooth family of non-degenerate quadratic forms, they have the same sign.
Also $\text{sign}(A)|_{L \cap \Lambda_0} = 0$.
Hence
\[\text{Mas}(\Lambda(t)) = \frac{1}{2}(\text{sign}(A_b) - \text{sign}(A)) = \frac{1}{2}\text{sign}(A_b)|_{{L \cap \Lambda_0}}.\]
Because the dimension of $L \cap \Lambda_0$ is $k$ this means that
$\text{Mas}(\Lambda(t)) \in [-\frac{k}{2},\frac{k}{2}]$.
\qed

\bigskip
As a result we have the following direct corollary:
\begin{corollary} \label{corollary:indexofsmallpathofsymplecticmatrices}
Let $A \in S_k$. Then for a sufficiently small neighborhood $\mathcal{N}_A$
of $A$ we have that any path $o : [a,b] \to \mathcal{N}_A$ with $o(a) = A$ has Conley-Zehnder index
in $[-\frac{k}{2},\frac{k}{2}]$.
\end{corollary}

The following lemma implies that if we perturb a pseudo Morse-Bott family then all the nearby Reeb orbits have a
bound on their Conley-Zehnder indices.
Recall that $D\phi_t$ is the linearization of the Reeb flow $\phi_t$ of $\alpha$.
Let $\xi := \text{ker}(\alpha)$ be the contact hyperplane distribution.
%Because $D\phi_t$ preserves $\xi$ we will write $D\phi_t|_\xi : \xi \to \xi$ to be the
%restriction of $D\phi_t$ to the hyperplane distribution.
We choose a trivialization $\tau$ of $\otimes_{j=1}^{N} \kappa^*$ where $\kappa^*$
is the anticanonical bundle.
%is the highest exterior power of $\xi$ viewed as a complex bundle with respect
%to some compatible almost complex structure.

\begin{lemma} \label{lemma:perturbedindexcalculation}
Let $\gamma$ be any Reeb orbit of $\alpha$ of period $T$ and define
$K := \text{dim ker}(D\phi_T|_\xi(\gamma(0)) - \text{id})$. Fix some riemannian metric on $C$.
There is a constant $\delta > 0$ and a neighborhood $N$ of $\gamma(0)$ so that
for any contact form $\alpha_1$ with $|\alpha - \alpha_1|_{C^2} < \delta$
and any Reeb orbit $\gamma_1$ of $\alpha_1$ starting in $N$ of period in $[T - \delta, T+ \delta]$
we have $\text{CZ}(\gamma_1) \in [\text{CZ}(\gamma)-\frac{1}{2}K,\text{CZ}(\gamma)+\frac{1}{2}K]$.
\end{lemma}
\proof of Lemma \ref{lemma:perturbedindexcalculation}.
Choose a sequence of contact forms $\alpha_i$ $C^2$ converging to $\alpha$,
and suppose there is a sequence of Reeb orbits $\gamma_i$ of $\alpha_i$ period $T_i$ where $T_i$  converges to $T$
and $\gamma_i(0)$ converges to $\gamma(0)$.
By Gray's stability theorem we can assume that $\text{ker}(\alpha_i) = \text{ker}(\alpha) = \xi$.
We wish to show that $\text{CZ}(\gamma_i) \in [\text{CZ}(\gamma)-\frac{1}{2}K,\text{CZ}(\gamma)+\frac{1}{2}K]$ for sufficiently large $i$.
Let $D_i(t) : T_{\gamma_i(0)} C \to T_{\gamma_i(t)} C$ be the linearization of the Reeb flow of $\alpha_i$
from $\gamma_i(0)$ to $\gamma_i(t)$.
Similarly define $D_\infty(t)$ for $\gamma(t)$.
Because these linearizations preserve the contact distribution we will actually view them as symplectic linear maps:
$D_i(t)|_\xi : \text{ker}(\alpha_i)_{\gamma_i(0)} \to \text{ker}(\alpha_i)_{\gamma_i(t)}$
and $D_\infty(t)|_\xi : \text{ker}(\alpha)_{\gamma_\infty(0)} \to \text{ker}(\alpha)_{\gamma_\infty(t)}$.
Choose trivializations of the Hermitian vector bundles $\tau_i : \gamma_i^* \oplus_{j=1}^{N} \xi \to (\R / T_i \Z) \times \R^{(2n-2) N}$
and $\tau_\infty : \gamma^* \oplus_{j=1}^{N} \xi \to (\R / T \Z) \times \R^{(2n-2) N}$ whose highest exterior power agrees with our chosen trivialization $\tau$. This means that
$\oplus_{j=1}^{N} D_i(t)|_\xi$ and $\oplus_{j=1}^{N} D_\infty(t)|_\xi$
are represented by paths of symplectic matrices in $\text{Sp}((2n-2)N)$ respectively, which we write as
$D^\oplus_i(t)$ and $D^\oplus_\infty(t)$
so that $D^\oplus_i(t)$,$t \in [0,T_i]$ $C^\infty$ converges (after linearly rescaling the $t$ coordinate)
to $D^\oplus_\infty(t)$, $t \in [0,T]$.
Note that $\text{dim ker}(D^\oplus_i(T_i) - \text{id}) = N \text{dim ker}(D_i(T_i)|_\xi - \text{id})$
and that $\text{CZ}(D^\oplus_i) = N \text{CZ}(\gamma_i)$
for all $i \in \N \cup \{\infty\}$.

By Corollary \ref{corollary:indexofsmallpathofsymplecticmatrices} there is a small neighborhood $\mathcal{N}$ of $D^\oplus_\infty(T)$ in
$\text{Sp}((2n-2)N)$
so that every path of symplectic matrices in $\mathcal{N}$ starting at $\gamma(T)$ has
Conley-Zehnder index in the interval $[-\frac{1}{2}NK,\frac{1}{2}NK]$.
We can also assume that $\mathcal{N}$ is contractible.
For $i$ large enough we have that $D^\oplus_i(T_i) \in \mathcal{N}$ and also that $D^\oplus_i(t)$
is homotopic to $D^\oplus_\infty(t)$ through paths whose starting point is fixed and whose endpoint is contained in $\mathcal{N}$.
Let $p_i(t) \in \mathcal{N}$ be a path starting at $D^\oplus_i(T_i)$ and ending at $D^\oplus_\infty(T)$.
The catenation of $D^\oplus_i(t)$ and $p_i(t)$ is a path homotopic to $D^\oplus_\infty(t)$ through paths with fixed endpoints equal to
$D^\oplus_\infty(0)$ and $D^\oplus_\infty(T)$
and so they have the same Conley-Zehnder index
by property \ref{item:homotopyrelendpoints}.
Using the fact that the Conley-Zehnder index of $p_i(t)$ is in $[-\frac{1}{2}NK,\frac{1}{2}NK]$ and
property \ref{item:catenation}
we get that the Conley-Zehnder index of $D^\oplus_i(t)$ is in $\text{CZ}([D^\oplus_\infty)-\frac{1}{2}NK,\text{CZ}(D^\oplus_\infty)+\frac{1}{2}NK]$.
Now $\text{CZ}(\gamma_i) = \frac{1}{N} \text{CZ}(D^\oplus_i)$
and $\text{CZ}(\gamma) = \frac{1}{N} \text{CZ}(D^\oplus_\infty)$
and hence $\text{CZ}(\gamma_i) \in  [\text{CZ}(\gamma)-\frac{1}{2}K,\text{CZ}(\gamma)+\frac{1}{2}K]$
for $i$ sufficiently large.
\qed

\begin{corollary} \label{corollary:lsftlowersemicontinuous}
Let $b,c \in \R$ be fixed constants.
Suppose that for all $\epsilon>0$, there is a contact form $\alpha_1$ admitting a Reeb orbit of $\text{lSFT}$ index $\leq b$ and period $\leq a$
satisfying $|\alpha - \alpha_1|_{C^2} < \epsilon$.
Then $\alpha$ admits a Reeb orbit with $\text{lSFT}$ index $\leq b$ and period $\leq a$.
\end{corollary}
\proof
Choose a sequence of such $\alpha_1$'s $C_2$ converging to $\alpha$
and let $\gamma_i$ be their respective orbits of $\text{lSFT}$ index $\leq b$ and period $\leq a$.
A compactness argument tells us that after passing to a subsequence,
$\gamma_i$ converges to a Reeb orbit $\gamma$ of $\alpha$ of period $\leq a$.
Lemma \ref{lemma:perturbedindexcalculation} then gives us our bound on
$\text{lSFT}(\gamma)$.
\qed

%We have the following direct Corollary of Lemma \ref{lemma:perturbedindexcalculation}:
%\begin{corollary} \label{corollary:lowersemicontinuouslsft}
%If $\alpha$ is a contact form on $C$, then
%$\text{mi}(\alpha_1) \geq \text{mi}(\alpha)$ for any contact form $\alpha_1$ sufficiently $C^2$ close to $\alpha$.
%In other words, $\text{mi}(\cdot)$ is lower semi-continuous in the $C^2$ topology.
%\end{corollary}

\section{Neighborhoods of Symplectic Submanifolds with Contact Boundary} \label{section:neighbourhoods}

We need a purely symplectic notion of divisor, and we need definitions
which make such a divisor look like the resolution of a singularity
which is numerically $\Q$-Gorenstein.

\begin{defn}
If we have some real codimension $2$ submanifolds $S_1,\cdots,S_l$ inside some manifold $M$
then they are said to be {\bf normal crossings}
if for each $p \in M$ there is a coordinate chart
$x_1,y_1,\cdots,x_k,y_k,w_1,\cdots,w_q$ centered at $p$
and distinct elements $i_1,\cdots,i_k \in \{1,\cdots,l\}$ so that
$S_{i_j} = \{x_j = y_j = 0\}$ near $p$ for all $j=1,\cdots,k$
and $S_i$ does not intersect this coordinate chart for all $i \in \{1,\cdots,l\} \setminus \{i_1,\cdots,i_k\}$.
\end{defn}

\begin{defn}
(\cite[Definition 5.1]{McLean:affinegrowth} or \cite[Definition 2.1]{McLeanTehraniZinger:normalcrossings}).
If $S_1,\cdots,S_l$ are submanifolds of a symplectic manifold $M$
then we say they are {\bf positively intersecting} if

\begin{enumerate}
\item
If $S_1,\cdots,S_l$ are connected and normal  crossings and
$S_I := \cap_{i \in I} S_i$ is a symplectic submanifold for each $I \subset \{1,\cdots,l\}$.
\item For each $I,J \subset \{1,\cdots,l\}$ with $I \cap J = \emptyset$
let $N_1$ be the symplectic normal bundle for $S_{I \cup J}$
in $S_I$ and $N_2$ the symplectic normal bundle for $S_{I \cup J}$ in $S_J$.
The condition $I \cap J = \emptyset$ ensures that $S_I$ and $S_J$ are transversally intersecting
and so $N_1 \oplus N_2 \oplus TS_{I \cup J} = TM|_{S_{I \cup J}}$.
Because $TS_{I \cup J}$, $N_1$, $N_2$ have natural orientations, we require
that the orientation of $N_1 \oplus N_2 \oplus TS_{I \cup J}$ agrees with the natural symplectic orientation on $TM|_{S_{I \cup J}}$ for all $I,J$
with $S_{I \cup J} \neq 0$.
\end{enumerate}
\end{defn}

\begin{defn}
A {\bf normal crossings exact divisor} consists of
a triple $(M,\cup_i S_i,\theta)$ where $M$ is a manifold with boundary, $S_i$ are real codimension $2$ submanifolds
and $\theta$ is a $1$-form on $M \setminus \cup_i S_i$ satisfying:
\begin{itemize}
\item $d\theta$ extends to a symplectic form $\omega$ on $M$.
\item $S_1,\cdots,S_l$ are positively intersecting symplectic submanifolds without boundary
and $\cup_i S_i \to M$ and $\partial M \to M \setminus \cup_i S_i$ are homotopy equivalances.
\end{itemize}
A normal crossings exact divisor $(M,\cup_i S_i,\theta)$ is {\bf compact} if $M$ is compact.
A normal crossings exact divisor $(M,\cup_i S_i,\theta)$ is {\bf strongly numerically $\Q$-Gorenstein}
if $H^1(M \setminus \cup_i S_i;\Q)=0$ and
$c_1(TM|_{M \setminus \cup_i S_i};\Q) = 0 \in H^2(M \setminus \cup_i S_i,\Q)$.
%A normal crossings Liouville divisor $(M,\cup_i S_i,\theta)$ is {\bf positively wrapped}
%if the wrapping number of $\theta$ around $S_i$ is positive for each $1 \leq i \leq l$.
%tA normal crossings exact divisor $(M,\cup_i S_i,\theta)$ is
It is {\bf orthogonal} 
if $S_1,\cdots,S_l$ are symplectically orthogonal submanifolds
(i.e. for each $i \neq j$ we have that the symplectic normal bundle of $S_i$ along $S_i \cap S_j$
is contained in $TS_j$).
\end{defn}

Let $(M,\cup_i S_i,\theta)$ be a compact normal crossings exact divisor.
Let $\theta_c$ be a $1$-form on $M$ equal to $\theta$ near $\partial M$ but equal to $0$ near $\cup_i S_i$.
Then $\omega - d\theta_c$ represents an element
$[\omega - d\theta_c]$ of $H^2(M,\partial M;\R)$
(the chain complex for this cohomology group consists of de Rham forms whose restriction to $\partial M$ is $0$).
Such a class only depends on $\theta$.
By Lefschetz duality we have $H^2(M,\partial M;\R) = H_{2n-2}(M;\R)$
and because $M$ is homotopic to $\cup_i S_i$ we get that $[\omega - d\theta_c]$ is Lefschetz
dual to $-\sum_i \lambda_i [S_i] \in H_{2n-2}(M;\R)$ for some $\lambda_1,\cdots,\lambda_l$.
\begin{defn}
The constant $\lambda_i$ is called the {\bf wrapping number of} $\theta$ {\bf around} $S_i$.
\end{defn}
We can also define wrapping numbers in the case where $M$ is non-compact, in which case
$H_{2n-2}(M;\R)$ is replaced by Borel-Moore homology.
The wrapping number was originally defined in \cite[Section 5.2]{McLean:affinegrowth} modulo rescaling by $2\pi$
and it is defined in the following way (which justifies the name `wrapping number'):
Choose a point $p \in S_i \setminus \cup_{j \neq i} S_j$
and an embedded symplectic $2$-dimensional disk $D \subset M$
of some small radius
intersecting $S_i$ at $p$ transversally at the origin with positive intersection number
and so that $(D \setminus \{0\}) \cap (\cup_i S_i) = \emptyset$.
Let $(r,\vartheta)$ be standard polar coordinates for $D \subset \C$
so that the symplectic form on $D$ is $\frac{1}{2} r^2 d\vartheta$.
We have that $[\theta|_{D \setminus \{0\}} - \frac{1}{2}r^2 d\vartheta] \in H^1(D \setminus \{0\};\R) = \R$
is equal to $(l_i/ 2\pi) d\vartheta$ for some $l_i \in \R$.
Here $l_i$ is the number of times $\theta - \frac{1}{2}r^2 d\vartheta$ `wraps' around $p$.
\begin{lemma} \label{lemma:wrappingequivalence}
$l_i = \lambda_i$.
\end{lemma}
\proof
By Stokes' Theorem one can add $df$ to $\theta$ without changing $\lambda_i$
or $l_i$ if the support of $f : M \to \R$ is contained in the interior of $M$.
Hence we can assume that $\theta|_{D\setminus \{0\}} = \frac{1}{2}r^2 d\vartheta + (l_i / 2\pi) d\vartheta$.
We can also assume that $\theta_c|_{D \setminus \{0\}} = \rho(r) \left(\frac{1}{2}r^2 d\vartheta + (l_i / 2\pi)\right) d\vartheta$
where $\rho(r) = 0$ near $p$ and $1$ near $\partial D$.

We have that $H^2(M,\partial M;\R) $ is
freely generated by the Lefschetz duals $\text{LD}([S_i])$ of $[S_i] \in H_{2n-2}(M;\R)$
by a Mayor-Vietoris argument.
Let $\rho : H^2(M,\partial M;\R) \to H^1(D \setminus \{0\};\R)$ be the natural map:
\[H^2(M,\partial M;\R) = H^2(M;M\setminus \cup_i S_i;\R) \to H^2(D;D\setminus \{0\};\R) = H^1(D \setminus \{0\};\R).\]
The identification $H^2(D;D\setminus \{0\};\R) = H^1(D \setminus \{0\})$
comes from the long exact sequence associated to the pair $(D,D \setminus \{0\})$.
The map $\rho$ sends $-\sum_j \lambda_j \text{LD}([S_j])$ to $-(\lambda_i / 2\pi) [d\vartheta]$.
Also $\rho$ sends
$[\omega-d\theta_c]$
% \big]
to $-(l_i / 2\pi) [d\vartheta]$
because
\[[(\omega-d\theta_c)|_D]= \big[d((1-\rho(r))\frac{1}{2}r^2 d\vartheta)  -  d(\rho(r) (l_i/2\pi) d\vartheta)\big]=\]
\[ -\big[d(\rho(r) (l_i/2\pi) d\vartheta)] \in H^2(D,\partial D;\R) \cong H^2(D,D\setminus \{0\};\R).\]
Hence $l_i = \lambda_i$.
\qed.

\bigskip

\begin{defn}
A normal crossings exact divisor $(M,\cup_i S_i,\theta)$
is called a {\bf positively wrapped divisor}
if the wrapping number of $\theta$ around $S_i$ is positive for each $1 \leq i \leq l$.
\end{defn}

Throughout this section, the following examples of
strongly numerically $\Q$-Gorenstein
positively wrapped divisors should be kept in mind:
If $A \subset \C^N$ is an affine variety with an isolated singularity at zero,
then we can intersect it with a small closed ball $B_\delta$.
By \cite{hironaka:resolution}
we resolve $A$ at $0$ by blowing up along smooth subvarieties
and take the preimage $\widetilde{A}_\delta$ of $B_\delta$ under this resolution map.
The preimage of zero under this resolution map is normal crossings,
and connected by Zariski's connectedness theorem
\cite{Zariski:connectedness}.
Later on (see Lemma \ref{lemma:symplecticformonresolution}) we will show that such a resolution has a compatible symplectic form and the exceptional divisors
will be positively intersecting symplectic submanifolds of $\widetilde{A}_\delta$
and we will also show that there is a $1$-form $\theta$ as above giving these submanifolds
positive wrapping numbers for an appropriate resolution.

\subsection{Boundary Construction and Uniqueness} \label{section:boundaryconstructionanduniqueness}

In this section we will construct a natural contact manifold associated to
positively wrapped divisors $(M,\cup_i S_i,\theta)$ which we will call the contact link.
This will be the boundary of an appropriate neighborhood of $\cup_i S_i$.
In the case where  $S_1,\cdots,S_l$
are exceptional divisors of a resolution of an isolated singularity, we will show later on that this contact
manifold is contactomorphic to the link
(see Lemma \ref{lemma:symplecticformonresolution}).
In order to do this we need a preliminary definition:
\begin{defn} \label{defn:radialcoordinate}
Suppose that $S$ is a smooth submanifold of a manifold $X$.
Choose a metric on $X$.
Let $TS^\perp \subset TM|_S$ be the subbundle of vectors that are
orthogonal to $TS$ with respect to this metric.
Let $TS^{\perp,\delta} \subset TS^\perp$ be the subset consisting of vectors
of length at most $\delta$.
Let $\delta_r>0$ be small enough so that the exponential map
$\text{exp} : TS^{\perp,\delta_r} \to X$
of our metric is well defined.

We define $r : \text{exp}(TS^{\perp,\delta_r}) \to [0,\delta_r)$
by $r(x) := |\text{exp}^{-1}(x)|$.
We say that $r$ is a {\bf radial coordinate associated to $S$}.
The metric $g$ is called a {\bf metric associated to $r$}.
\end{defn}

Let $(M,\cup_i S_i,\theta)$ be a compact positively wrapped divisor.
Let $r_i : N_i \twoheadrightarrow [0,\delta_r)$ be a radial coordinate
associated to $S_i$.
%along with a normal disk bundle $p_i : N_i \twoheadrightarrow S_i$ compatible with $r_i$.
Let $\rho : [0,\delta_r) \to [0,1]$ be a smooth function so that $\rho(x) = x^2$ near $0$
and $\rho(x) = 1$ near $\delta_r$ with $\rho' \geq 0$.
We define $\rho_{r_i}$ to be $\rho \circ r_i$ inside the domain of $r_i$
and $1$ elsewhere.
A smooth function $f : M \setminus \cup_i S_i \to \R$ is said to be {\bf compatible with } $\cup_i S_i$ if it is equal to
$\sum_i m_i \log(\rho_{r_i})) + \tau$ for some chosen radial coordinates $r_i$ associated to $S_i$, choice of function $\rho$ as above,
constants $m_i >0$,
and some smooth function $\tau : M \to \R$.
%For any $1$-form $\nu$, define $X_\nu$ to be the vector field which is $\omega$-dual to $\nu$.

The main proposition of this section is:
\begin{prop} \label{prop:boundaryexistence}
For any smooth function $f : M \setminus \cup_i S_i \to \R$ compatible with
$\cup_i S_i$,
there exists a smooth function $g : M \setminus \cup_i S_i \to \R$
so that $df(X_{\theta + dg}) > 0$ in some small neighborhood of 
$\cup_i S_i$.
\end{prop}

The above proposition gives us the following definition.
\begin{defn} \label{defn:linkofsi}
Let $(M,\cup_i S_i,\theta)$ be a compact positively wrapped divisor, let $f : M \setminus \cup_i S_i \to \R$ be compatible with $\cup_i S_i$
and let $g$ be as in the above proposition.
Let $y \in \R$ be sufficiently negative so that for any $y_1 \leq y$,
$df(X_{\theta + dg}) > 0$ along $f^{-1}(y_1)$
and so that $f^{-1}(y_1)$ is compact.
%(such a $y$ exists by the previous proposition).
The contact structure
$(f^{-1}(y),\text{ker}( (\theta+dg)|_{f^{-1}(y)}))$
is called the {\bf contact link} of $(M,\theta,\cup_i S_i)$.
Corollary \ref{corollary:canonicalboundary} below
tells us that this contact structure is independent of choice of radial coordinates, constants $m_i$ and
functions $\rho,\tau$ and $g$.
\end{defn}

We also have a parameterized version of Proposition \ref{prop:boundaryexistence}.
Suppose that $(M,\cup_i S_i,\theta_t)_{t \in [0,1]}$ are compact positively wrapped divisors
where $\theta_t$ smoothly depends on $t$.
Let $\omega_t$ be the symplectic form on $M$ equal to $d\theta_t$ on $M\setminus \cup_i S_i$.
%For any $1$-form $\eta$ defined on an open subset of $M$ we define $X^t_\eta$ to be the $\omega_t$-dual of
%$\eta$. 
\begin{prop} \label{proposition:parameterizedboundaryexistence}
Let $f_t: M \setminus \cup_i S_i \to \R$ by a smooth family of functions compatible with
$\cup_i S_i$ so that the associated radial coordinates $r_i$, constants $m_i$ and
functions $\rho$ and $\tau$ smoothly vary as $t$ varies.
Then there is a smooth family of functions $g_t : M \setminus \cup_i S_i \to \R$
so that $df_t(X^{\omega_t}_{\theta_t + dg_t}) > 0$ near $\cup_i S_i$. 
\end{prop}
The proof of this proposition is exactly the same as the proof of
Proposition  \ref{prop:boundaryexistence} except that we now
have everything parameterized by $t$.

We have the following corollary:
\begin{corollary} \label{corollary:canonicalboundary}
Let $f_0,f_1 : M \setminus \cup_i S_i \to \R$ be compatible with $\cup_i S_i$
and let $g_0,g_1 : M \setminus \cup_i S_i \to \R$ be such that for $j=0,1$, $df_j(X^{\omega_j}_{\theta_j + dg_j}) > 0$ near $\cup_i S_i$. 
Then for all sufficiently negative $l$ we have that $\left(f_0^{-1}(l),\text{ker}\big((\theta_0 + dg_0)|_{f_0^{-1}(l)}\big)\right)$ is contactomorphic to $\left(f_1^{-1}(l),\text{ker}\big((\theta_1 + dg_1)|_{f_1^{-1}(l)}\big)\right)$.
\end{corollary}

This corollary tells us that the contact link only depends on the deformation class of $(M,\theta,\cup_i S_i)$.

\proof of Corollary \ref{corollary:canonicalboundary}.
By definition we have radial coordinates $r_i^j : N^j_i \twoheadrightarrow [0,\delta_r)$
associated to $S_i$ defined using metrics $g_i^j$ along with
constants $m^j_i > 0$
and smooth functions $\tau_j$ defined on a neighborhood of $\cup_k S_k$
so that $f_j = \sum_i m^j_i \log(\rho_{r^j_i}) + \tau_j$ for both $j = 0,1$.

Now choose a smooth family of metrics $g_i^t$, $t \in [0,1]$
joining $g_i^0$ with $g_i^1$.
Let $S^\perp_{i,t} \subset TM|_{S_i}$ be the set of vectors $g_i^t$
orthogonal to $TS_i \subset TM|_{S_i}$.
Let $S^\perp_{i,t,\delta_r} \subset N_{i,t}$ be the set of vectors of length less than $\delta_r$
and $\text{exp}_{g_i^t} : S^\perp_{i,t,\delta_r} \to M$ the associated exponential map
(which is well defined for $\delta_r$ small enough).
Let $N_i^t$ be the image of this exponential map.
Define $r_i^t : N_i^t \twoheadrightarrow [0,\delta_r)$
by $r_i^t(x) := |\text{exp}_{g_i^t}^{-1}(x)|_t$
where $|\cdot|_t$ is the norm induced by $g_i^t$ on $S^\perp_{i,t}$.
%Define 
%$p^t_i : N^t_i \twoheadrightarrow S_i$ by
%$p^t_i(x) := P(\text{exp}_{g_i^t}^{-1}(x))$.
Then $r_i^t$ is a smooth family of radial coordinates associated to $S_i$
joining $r_i^0$ and $r_i^1$.
% and
%$p^t_i$ is a smooth family of normal disk bundles associated to $r_i^t$
%joining $p^0_i$ and $p^1_i$.

Let $\tau_t$ be a smooth family of functions on $M$ joining $\tau_0$ and $\tau_1$
and $m^t_i > 0$ a smooth family of constants joining $m^0_i$ with $m^1_i$.
Then $f_t := \sum_i m^t_i \log(\rho_{r^t_i}) + \tau_t$
is a smooth family of functions compatible with $\cup_k S_k$
joining $f_0$ and $f_1$.
We can now apply Proposition \ref{proposition:parameterizedboundaryexistence}
to the functions $f_t$ and symplectic forms $\omega_t$
combined with Gray's stability theorem to conclude that
$(f_0^{-1}(l),\theta_0 + dg^1_0|_{f_0^{-1}(l)})$ is contactomorphic to $(f_1^{-1}(l),\theta_1 + dg^1_1|_{f_1^{-1}(l)})$
for some smooth functions
$g^1_0,g^1_1 : M \setminus \cup_i S_i \to \R$.
By linearly interpolating between $g^1_0$ and $g_0$ and $g^1_1$ and $g_1$
we get that

$\left(f_0^{-1}(l),\text{ker}\big((\theta_0 + dg^1_0)|_{f_0^{-1}(l)}\big)\right)$ is contactomorphic to
$\left(f_0^{-1}(l),\text{ker}\big((\theta_0 + dg_0)|_{f_0^{-1}(l)}\big)\right)$ and
$\left(f_1^{-1}(l),\text{ker}\big((\theta_1 + dg^1_1)|_{f_0^{-1}(l)}\big)\right)$ is contactomorphic to
$\left(f_1^{-1}(l),\text{ker}\big((\theta_1 + dg_1)|_{f_0^{-1}(l)}\big)\right)$.
Hence
$\left(f_0^{-1}(l),\text{ker}\big((\theta_0 + dg_0)|_{f_0^{-1}(l)}\big)\right)$ is contactomorphic to $\left(f_1^{-1}(l),\text{ker}\big((\theta_1 + dg_1)|_{f_1^{-1}(l)}\big)\right)$.
\qed

\bigskip
To prove Proposition \ref{prop:boundaryexistence} we need the following Lemma:
\begin{lemma} \label{lemma:localboundaryexistence}
Let $(U,\cup_i \Sigma_i,\theta_U)$ be a non-compact positively wrapped divisor with positive wrapping numbers.
Let $\pi : U \twoheadrightarrow S$ be a smooth fibration
whose fibers are all symplectomorphic to  a ball of small radius in $\R^{2n}$ and so that the natural symplectic connection (i.e. the Ehresmann connection given by the subbundle of vectors symplectically orthogonal to the fibers)
has parallel transport maps in $U(n)$.
Hence we can view $S$ as a submanifold of $U$ given by the $0$ section.
We will assume that $S$ is diffeomorphic to a ball and $S = \cap_i \Sigma_i$.
Let $\| \cdot\|$  be any metric on $U$ so that $\|\omega_U\|$ is bounded.

Then
(after shrinking the fibers of $\pi$) there is a function
$g_U : M \setminus \cup_i S_i \to \R$ so that
for all smooth functions $f : U \setminus \cup_i \Sigma_i \to \R$
compatible with $\cup_i S_i$ we have that
$df(X_{\theta_U+dg_U}) > c_f \|\theta_U+dg_U\|\|df\|$ near $\cup_i \Sigma_i$ where $c_f>0$ is a constant depending on $f$.
Also $a_1 \|db\| < \|\theta_U + dg_U\| < a_2\|db\|$
near $\cup_i S_i$ for some smooth function $b$ compatible with $\cup_i \Sigma_i$
and constants $a_1,a_2>0$.
%We also have a parameterized version of this theorem where $S'_1,\cdots,S'_n$, $\theta'$ and $g' : M \setminus \cup_i S_i \to \R$
%smoothly depend on a parameter $t$ lying inside a compact manifold. 
\end{lemma}

\proof of Lemma \ref{lemma:localboundaryexistence}.
Before we get in to the details of this technical Lemma, we will first give a sketch of the proof
in the case when $S$ is a point and $U = \C^2$ with the standard symplectic form
and $\Sigma_1,\Sigma_2$ are linear (the more general case is similar).
First of all we pass to the universal cover $\widetilde{U}$ of $\C^2 \setminus \Sigma_1$
which is symplectomorphic to the open subset
$\widetilde{x}_1 > 0$ in $\C^2$ where $\widetilde{x}_1 + i \widetilde{y}_1, \widetilde{x}_2 + i \widetilde{y}_2$
are complex coordinates for $\C^2$.
Here $\widetilde{x}_1$ is equal to the pullback of the square of some radial coordinate associated to $\Sigma_1$, $\widetilde{y}_1$ is an angle coordinate and $\widetilde{x}_2,\widetilde{y}_2$
are the pullbacks of linear coordinates.
The closure of the preimage of $\Sigma_2$ in this cover is a manifold with boundary $\widetilde{\Sigma}_2$
containing the line $\nu$ given by $\widetilde{x}_1 = \widetilde{x}_2 = \widetilde{y}_2 = 0$.
We have that $\widetilde{\Sigma}_2$ as an oriented real blowup of $\Sigma_2$.
The key technical step in the proof is to construct a
closed $1$-form $\widetilde{q}_1$ in a neighborhood of $\nu$
so that $X_{\widetilde{q}_1}$ is non-zero and tangent to $\widetilde{\Sigma}_2$.
One can also ensure that $\widetilde{q}_1$ descends to $q_1$ on $\C^2 \setminus S_1$
so that it represents $\lambda_1 \in \R = H^1(\C^2 \setminus \Sigma_1;\R)$.
The positive wrapping condition ensures that
$X_{q_1}$ points away from $S_1$.
We construct $q_2$ in a similar way. 
Because our symplectic form $\omega_U$ on $U$ is closed,
it is exact on a neighborhood of $0$ in $\C^2$ and hence is equal to $d\theta_1$.
Then
$\theta_1 + q_1 + q_2$ is cohomologous to $\theta$
and the tangency condition ensures that
$X_{\theta_1 + q_1 + q_2}$ has the right properties.
The point being that $X_{q_1+q_2}$ can be computed,
and $X_{\theta_1}$ is tiny compared with
$X_{q_1 + q_2}$.

We will now proceed with the full proof (it is good to work through this proof
in the case where $(U,\cup_i \Sigma_i,\theta_U)$ is as in the above example).
Because the structure group of $\pi$ is $U(n)$,
the fibers have a natural metric.
Hence we can shrink the fibers of $\pi$
by reducing their radius so that each submanifold $\Sigma_i$ intersects each fiber of $\pi$ transversally.
Let $\lambda_1,\cdots,\lambda_n >0$ be the wrapping numbers of $\theta_U$ around $\Sigma_1,\cdots,\Sigma_n$ respectively.
For the moment we will fix some $1 \leq i \leq n$.
Again after shrinking the fibers of $U$ we have functions $x_1,y_1,\cdots,x_n,y_n \in C^\infty(U)$
so that $\Sigma_i = \{x_1,y_1=0\}$ and so that the restriction of these functions to each fiber of $\pi$
gives us a {\it symplectic} coordinate chart centered at $0$.
Note that these coordinates are dependent on $i$.
We also have polar coordinates $r,\vartheta$ depending on $x_1,y_1$ only so that $x_1 = r\cos(\vartheta)$ and $y_1 = r\sin(\theta)$.
The universal cover $\widetilde{U}_i$ of $U \setminus \Sigma_i$
admits a fibration $\widetilde{\pi}_i : \widetilde{U}_i \twoheadrightarrow S$ equal to the covering map composed with $\pi$.
Also $\widetilde{U}_i$ has functions
$\widetilde{x}_1,\widetilde{y}_1,\widetilde{x}_2,\widetilde{y}_2,\cdots,\widetilde{x}_n,\widetilde{y}_n$
which are pullbacks of
$\frac{1}{2}r^2,\vartheta,x_2,y_2\cdots,x_n,y_n$ respectively.
Each fiber of $\widetilde{\pi}_i$ is the universal cover of each fiber of $\pi|_{U \setminus \Sigma_i}$.
The coordinates $(\widetilde{x}_1,\widetilde{y}_1,\widetilde{x}_2,\widetilde{y}_2,\cdots,\widetilde{x}_n,\widetilde{y}_n)$
naturally identify the fibers of $\widetilde{\pi}_i$ with an open subset of $\R^{2n}$ hence $\widetilde{\pi}_i : \widetilde{U}_i \twoheadrightarrow S$ enlarges to a smooth fibration
$\widehat{\pi}_i : \widehat{U}_i \twoheadrightarrow S$ whose fibers are diffeomorphic to
$\R^{2n}$ with standard coordinates  $(\widetilde{x}_1,\widetilde{y}_1,\widetilde{x}_2,\widetilde{y}_2,\cdots,\widetilde{x}_n,\widetilde{y}_n)$
but the symplectic form does not necessarily extend.
Having said that the restriction of the symplectic form to each fiber of $\widetilde{\pi}_i$ extends to
$\sum_i d\widetilde{x}_i \wedge d\widetilde{y}_i$ on the fibers of $\widehat{\pi}_i$.
The group of deck transformations $\Z$ acts on each fiber of $\widetilde{\pi}_i$ by sending the coordinate
$(\widetilde{x}_1,\widetilde{y}_1,\widetilde{x}_2,\widetilde{y}_2,\cdots,\widetilde{y}_n)$ to
$(\widetilde{x}_1,\widetilde{y}_1+ 2\pi k,\widetilde{x}_2,\widetilde{y}_2,\cdots,\widetilde{y}_n)$
for each $k \in \Z$. This group action extends to one on $\widehat{U}_i$.
%
%

%Note that if we have submanifold $\widetilde{V} \subset \widehat{U}_i$
%which intersects each fiber of $\widehat{\pi}_i$ transversally and whose intersection
%with each fiber is an affine linear subspace of $\R^{2n}$
%in the coordinates $(\widetilde{x}_1,\widetilde{y}_1,\widetilde{x}_2,\widetilde{y}_2,\cdots,\widetilde{x}_n,%\widetilde{y}_n)$ 
%where the restriction of $\widetilde{y}_1$ to this submanifold is constant and invariant under translations in the %$\widetilde{x}_1$ direction,
%then it maps via the covering map to a fiberwise linear submanifold $V$ of $U$ in the coordinates
%$x_1,y_1,\cdots,x_n,y_n$.
%We will call $V$ the {\it fiberwise linear image} of $\widetilde{L}$.
%We can view the tangent space at $0$ of $\pi^{-1}(s) \cap \cap_{j  \neq i} \Sigma_j$
%as a subspace of $\pi^{-1}(s)$ linear in the coordinates $(x_1,y_1,\cdots,x_n,y_n)$.
Let $L_{\vartheta} \subset \widehat{U}_i$ be the unique submanifold of $\widehat{U}_i$ so that:
\begin{enumerate}
\item it intersects each fiber transversally.
\item its intersection with each fiber $F = (\widehat{\pi}_i)^{-1}(p)$ is a straight line $L_{\vartheta}^F$
contained in the plane $F \cap \{\widetilde{y}_1 = \vartheta\}$.
%passing through the point
%$(0,\vartheta,0,0,\cdots,0)$ in the coordinates
%$(\widetilde{x}_1,\widetilde{y}_1,\widetilde{x}_2,\widetilde{y}_2,\cdots,\widetilde{x}_n,\widetilde{y}_n)$
%(and not necessarily passing though $0$).
%\item  $\widetilde{y}_1$ restricted to $L_{\vartheta}^F$ is constant for all $F$.
\item The projection of $L_{\vartheta}^F \cap \widetilde{U}_i$ to $\pi^{-1}(p)$ extends to a
unique straight line, which is contained
in the tangent space at zero of
$\pi^{-1}(p) \cap \cap_{j \neq i} \Sigma_j$.
\end{enumerate}
% under the quotient map $\widetilde{U}_i \twoheadrightarrow U_i$.
%Define $L := \cup_{\vartheta'} L_{\vartheta'}$.
%This is a submanifold of $\R^{2n}$ which is invariant under our $\Z$ action of deck transformations.
If $q$ is a $1$-form in $\widehat{U}_i$ then define $X^v_q$
to be the unique vector field tangent to the fibers $(\widehat{\pi}_i)^{-1}(p)$ of $\widehat{\pi}_i$
so that $q|_{(\widehat{\pi}_i)^{-1}(p)} = i_{X^v_q}(\sum_j d\widetilde{x}_j \wedge d\widetilde{y}_j)$.

There is a $1$-form $\widetilde{q}_i$ in $\widehat{U}_i$ with the following properties:
\begin{enumerate}
\item The restriction of $\widetilde{q}_i$ to each fiber is closed.
\item Inside each fiber $F$, $X^v_{\widetilde{q}_i}$ is tangent to the line $L^F_{\vartheta}$ at the point $(0,\vartheta,0,0,\cdots,0)$
and pointing in the direction in which $\widetilde{x}_1$ is strictly increasing.
\item \label{item:wrappingproperty}
Define $\nu := \cap_j \{\widetilde{x}_j = 0\} \cap \cap_{j \neq 1} \{\widetilde{y}_j = 0\}$.
The integral from $\widetilde{y}_1=0$ to $\widetilde{y}_1=2\pi$
of $\widetilde{q}_i$ along each line $\nu \cap (\widehat{\pi}_i)^{-1}(p)$ is
$\lambda_i$.
\item $\widetilde{q}_i$ is invariant under the $\Z$ action near $\nu$.
%\item The kernel of $\widetilde{q}_i$ at each point $p \in \nu$
%contains the tangent space to $\cap_j  \{\widetilde{x}_j = \widetilde{x}_j(p)\} \cap \cap_j \{\widetilde{y}_j = \widetilde{x}_j(p)\}$.
\end{enumerate}
Because the fibers of $\widehat{\pi}_i$ are contractible
we have that $\widetilde{q}_i$ is fiberwise exact and so there is a smooth function $g_i : \widehat{U}_i \to \R$ whose differential restricted to each fiber
is $\widetilde{q}_i$.
Let $\widehat{q}_i = dg_i$. This is a closed $1$-form with exactly the same properties as $\widetilde{q}_i$.
After shrinking the fibers of $\pi$, we have that $\widehat{q}_i$ descends to a
closed $1$-form $q_i$ on $U \setminus \Sigma_i$.

Because $dq_i = 0$ we have that the symplectic form $\omega_U$ is equal to the exterior derivative of
$\Theta := \theta_1 + \sum_i q_i$
where $\theta_1$ is a $1$-form on $U$ whose norm is bounded.
This is equal to $\theta_U + dg_U$ near the zero section for some smooth function $g_U : U \setminus \cup_i \Sigma_i \to \R$
because both $\theta_U$ and $\Theta$ have the same wrapping numbers
by property (\ref{item:wrappingproperty}) stated for $\widetilde{q}$ earlier.

We now wish to show that $df(X_{\theta_U + dg_U}) > c_f \|\theta_U+dg_U\|\|df\|$ near $S$ for each $f$ compatible with $\cup_i \Sigma_i$.
Because the norm of $\theta_1$ is bounded
it is sufficient to show that $df(X_{\sum_i q_i}) > c\|\sum_i q_i\|\|df\|$ for some constant $c$ near $S$.
Because $f$ is compatible with $\cup_i \Sigma_i$,
we have radial coordinates $r_i$ associated to $S_i$ for each $i$,
constants $m_i>0$ and a smooth function $\tau \in C^\infty(U)$  so that
$f = \sum_i m_i \log(r_i) + \tau$ near $0$.
For any sequence of points $p_k \in U \setminus \cup_j \Sigma_j$ tending to $p \in S$ we have that (after passing to a subsequence)
$X_{q_i}/\|q_i\|$ at $p_k$ tends to a vector $v \in TM|_S \setminus T\Sigma_i$.
Because each such vector $v$ is transverse to $\Sigma_i$
we have that $d\log(r_i)(X_{q_i}) > c_i \|q_i\| \|d\log(r_i)\|$  near $0$ for some constant $c_i$.
Also for any sequence of points $p_k \in U \setminus \cup_j \Sigma_j$ tending to $0$ and any $i_1\neq i_2$ we have that (after passing to a subsequence)
$X_{q_{i_1}}/\|q_{i_1}\|$ at $p_k$ tends to a vector $v$ tangent to $S_{i_2}$.
This implies that $dr_i(X_{q_{i_1}}/\|q_{i_2}\|)$
tends to $0$ and hence $\frac{d\log(r_i)(X_{q_{i_1}})} {\|q_{i_1}\| \|d\log(r_i)\|}$ tends to $0$.
Putting everything together we get that
\[d\log(r_i)(\sum_j X_{q_j}) > c^f_i \|q_i\| \|d\log(r_i)\|\] for some constant $c^f_i$.
Hence
\[d\Big(\sum_i m_i\log(r_i)\Big)\Big(\sum_j X_{q_j}\Big) > c \sum_i \| q_i\| \|d\log(r_i)\|\] for some constant $c>0$.
Now there are constants $c_1,c_2$ so that $c_1\|q_i\| < \|d\log(r_i)\| < c_2\|q_i\|$.
Hence $\sum_i \| q_i\| \|d\log(r_i)\|$ 
is greater than some constant times $\sum_i m_i \|d\log(r_i)\|^2$
which in turn is greater than a constant times $(\sum_i  \|d\log(r_i)\| )\|\sum_i m_i d\log(r_i)\|$
is greater than some constant times $\|\sum_i q_i\|\|\sum_i m_i d\log(r_i)\|$.
This implies that there is a constant $c>0$ so that:
\[d\Big(\sum_im_i \log(r_i)\Big)\Big(\sum_j X_{q_j}\Big) > c \|\sum_i q_i\| \|\sum_i m_i d\log(r_i)\|.\]

Now  because $\tau$ is smooth at $S$ and near $S$, we have
$d\tau(\sum_i q_i) < \gamma \|\sum_i q_i\|$ for some constant $\gamma>0$.
Because $\|d\log(r_i)\|$ and $\|q_i\|$
tend to infinity as we approach $0$ for each $i$, we get our bound
$df(X_{\sum_i q_i}) > c \|\sum_i q_i\|\|df\|$ for some $c>0$.
Hence $df(X_{\theta_U+dg_U}) > c_f \|\theta_U+dg_U\|\|df\|$ near $S$ for some $c_f>0$.
Because $\theta_U +dg_U = \theta_1 + \sum_i q_i$ near $S$
and that $\|\theta_1\|$ is bounded near $0$, there are constants $a_1,a_2>0$ so that
$a_1 \|db\| < \|\theta_U+dg_U\| < a_2\|db\|$
where $b = \sum_i \log(r_i)$.

\qed

\bigskip

\proof of Proposition \ref{prop:boundaryexistence}.
Fix some metric $\|\cdot\|$ on $M$.
For each $I \subset \{1,\cdots,l\}$
define $S_I :=  \cap_{i \in I} S_i$ and
also choose open sets $V^I_1,\cdots,V^I_{k_I},~I \subset \{1,\cdots,l\}$ of $M$
so that:
\begin{enumerate}
\item $\cup_{I, j} V^I_j$ contains  $\cup_{i \in I} S_i$.
\item $W^I_j := V^I_j \cap S_I$ is diffeomorphic to an open ball.
\item There are smooth fibrations $\pi^I_j : V^I_j \twoheadrightarrow W^I_j$
whose fibers are symplectomorphic to small open balls in $\C^{|I|}$
and so that the natural symplectic connection has parallel transport maps in $U(|I|)$.
\end{enumerate}

Choose some total order $\preceq$ on the set of pairs $(I,j)$ where $I \subset \{1,\cdots,l\}$
and $1 \leq j \leq k_I$.
We write $(I_1,j_1) \prec (I_2,j_2)$ if $(I_1,j_1) \neq (I_2,j_2)$ and $(I_1,j_1) \preceq (I_2,j_2)$.
We will also choose slightly smaller subsets ${\breve{V}}^I_j$ whose closure is contained in $V^I_j$
but which still cover $\cup_i S_i$ and we define $\breve{W}^I_j := \breve{V}^I_j \cap S_I$.
Suppose (inductively) we have constructed some function $g^\prec$
so that $df(X_{\theta + dg^\prec}) > c \|\theta + dg^\prec\|\|df\|$
and $\|\theta + dg^\prec\| < \|db^\prec\|$
 in some small open set ${\mathcal N}$ containing
$\cup_{(I_1,j_1) \prec (I,j)} {\breve{W}}^{I_1}_{j_1}$  where $b^\prec$ is some function compatible with $\cup_i S_i$.
We wish to prove the same thing
in some open set containing 
$\cup_{(I_1,j_1) \preceq (I,j)} {\breve{W}}^{I_1}_{j_1}$.
By Lemma \ref{lemma:localboundaryexistence} there is a smooth function $g^= : M \setminus \cup_{i \in I} S_i \to \R$
so that $df(X_{\theta + dg^\prec + dg^=}) > c_f \|\theta + dg^\prec + dg^=\| \|df\|$ in some neighborhood of $W^I_j$.
Also there is a smooth function $b$ compatible with $\cup_i S_i$
so that $\|\theta + dg^\prec + dg^=\| < \|db\|$.
Let $\rho : M \to \R$ be a bump function equal to $0$ outside $V^I_j$
and equal to $1$ in ${\breve{V}}^I_j$.
Define $g^\preceq := g^\prec + \rho g^=$.
Inside ${\mathcal N} \cap V^I_j$ we have
$\|\theta + dg^\prec\| < \|db^\prec\|$ and 
$\|\theta + dg^\prec +dg^=\| < \|db\|$ and so
$\|dg^=\| < \|d\beta\|$ inside
${\mathcal N} \cap V^I_j$
for some function $\beta$
compatible with $\cup_i S_i$.
This means that $|g^=|$
is bounded above by some function $\nu$ so that $-\nu$ is
compatible with $\cup_i S_i$
inside ${\mathcal N} \cap V^I_j$.
Now \[X_{\theta + dg^\preceq} = X_{(1-\rho)(\theta +dg^\prec)} +X_{\rho(\theta +dg^\prec + dg^=)} +
g^=X_{\rho}.\] 
Because $\frac{\|g^=X_{\rho}\|}{\|d\beta\|}$
is bounded inside ${\mathcal N} \cap V^I_j$
we get
$df(X_{\theta + dg^\preceq}) > c \|df\| \|\theta + dg^\preceq\|$ for some $c>0$ in some open
set containing 
$\cup_{(I_1,j_1) \preceq (I,j)} {\breve{W}}^{I_1}_{j_1}$.
Also by construction we have $\|\theta + dg^\preceq\|$ is bounded above by $d\beta_2$
in this same open set
where $\beta_2$ is compatible with $\cup_i S_i$.
Hence by induction we have shown
$df(X_{\theta + dg}) > 0$ near
$\cup_i S_i$ for some $g : M \setminus \cup_i S_i$. 
\qed

The proof of Proposition \ref{proposition:parameterizedboundaryexistence}
is identical to the proof of Proposition
\ref{prop:boundaryexistence}
but now $\theta$, 
$\pi^I_j : V^I_j \twoheadrightarrow W^I_j$,
$f$,$g$,$g^\prec$,$g^=$,$g^\preceq$,$\rho$ and ${\mathcal N}$ all smoothly depend on the parameter $t$.

\bigskip

\subsection{Constructing a Specific Contact Form on Our Boundary}
\label{subsection:specificcontactform}

We will now define minimal discrepancy of
a strongly $\Q$-Gorenstein exact divisor
$(M,\cup_i S_i,\theta)$
by using the ideas from Lemma \ref{lemma:topologicalnumericallygorenstein}.
There is a long exact sequence
\[0=H^1(M \setminus \cup_i S_i,\Q)
\hookrightarrow
H^2(M,M \setminus \cup_i S_i,\Q)
\hookrightarrow H^2(M,\Q) \to
H^2(M \setminus \cup_i S_i,\Q).
\]
Now $c_1(TM) \in  H^2(M,\Q)$
maps to zero in this long exact sequence and so there is a unique element
$c_1(M,M \setminus \cup_i S_i)
\in H^2(M,M \setminus \cup_i S_i,\Q)$
called the relative first Chern class.
Because
 $\cup_i S_i \hookrightarrow M$ is a homotopy equivalence, we get that 
$H^2(M,M \setminus \cup_i S_i,\Q) = H_{2n-2}(M,\Q)$
is generated by fundamental classes $[S_i]$
and so 
$c_1(M,M \setminus \cup_i S_i)
= -\sum_i a_i [S_i]$
for some $a_1,\cdots,a_l \in \Q$
(we chose to write minus signs in order to indicate that
we are dealing with Chern classes of tangent bundles
instead of cotangent bundles).
\begin{defn} \label{defn:symplecticminimaldiscrepancy}
We define the {\bf discrepancy} of $S_i$ in a strongly $\Q$-Gorenstein
normal crossings exact divisor $(M,\cup_i S_i,\theta)$
to be $a_i$ as defined above.
We define the {\bf minimal discrepancy} $\text{md}(M,\cup_i S_i,\theta)$
of $(M,\cup_i S_i,\theta)$
to be $\text{min}_i a_i$ if $\text{min}_i a_i \geq -1$
and $-\infty$ otherwise.
\end{defn}

\begin{defn} \label{defn:linkregion}
Let $(M,\cup_i S_i,\theta)$ be a positively wrapped divisor.
A pair $(M_1,g : M \setminus \cup_i S_i \to \R)$ is called
is called a {\bf link region} if
$M_1 \subset M$ is a codimension $0$ submanifold containing $\cup_i S_i$ so that
there is a function $f$ compatible with $\cup_i S_i$ satisfying
\begin{itemize}
\item  $df(X_{\theta + dg}) > 0$ inside $M_1 \setminus \cup_i S_i$  and
\item $M_1 = \{f^{-1}(-\infty,c] \cup \cup_i S_i$.
\end{itemize}
This implies that the boundary of $M_1$ is contact submanifold of $M$ contactomorphic to the link of $\cup_i S_i$.
The function $f$ will be called a {\bf function compatible with $(M_1,g : M \setminus \cup_i S_i \to \R)$}.
%Here $c$ is negative enough so that $(\partial M_1,\theta + dg|_{\partial M_1})$ is contactomorphic to the link of %$(M,\cup_i S_i,\theta)$.
\end{defn}

%TODO:continue here.
We now need a specific technical definition which will be used in the proof of
Theorem \ref{theorem:reeborbitlowerboundarounddivisors}.

\begin{defn} \label{defn:productregion}
Let $(M,\cup_i S_i,\theta)$ be a positively wrapped divisor.
Let $(M_1,g : M \setminus \cup_i S_i \to \R)$ be a link region
and let $B \subset \partial M_1$ be a pseudo Morse-Bott family of Reeb orbits
of $\theta + dg|_{\partial M_1}$.
Let $B^{2n}_\nu \subset \C^{2n}$ be the symplectic ball of radius $\nu$.

An $(\epsilon_1,\epsilon_2)$-{\bf ball product associated to $B$ and $S_k$ inside $M_1$}
consists of a subset $W \subset M_1 \setminus \cup_{j \neq k} S_j$ symplectomorphic to $B^{2n-2}_{\epsilon_1} \times B^2_{\epsilon_2}$ so that:
\begin{enumerate}
\item $W \cap S_k = B^{2n-2}_{\epsilon_1} \times \{0\} \subset S_k \setminus \cup_{j \neq i} S_j$,
\item the restriction of $\theta + dg$ to $B^{2n-2}_{\epsilon_1} \times B^2_{\epsilon_2}$ is $\text{pr}_1^* \beta + \text{pr}_2^* ((\frac{1}{2}r_k^2 + \frac{1}{2\pi}\lambda_k)d\vartheta_k)$
where $\beta$ is a $1$-form, $\lambda_k$ is the wrapping number of $\theta$ around $S_k$,$(r_k,\vartheta_k)$ are polar coordinates on $B^2_{\epsilon_2}$ and $\text{pr}_1$ and $\text{pr}_2$ are the projections to  $B^{2n-2}_{\epsilon_1}$ and $B^2_{\epsilon_2}$ respectively,
\item $\partial M_1 \cap W = B^{2n-2}_{\epsilon_1} \times \partial B^2_{\epsilon_2}$
and the Reeb orbits inside $\partial M_1 \cap W_k$ wrapping around $S_1$ once
are contained inside the family $B$ (in other words, the Reeb orbits
$\gamma_x : \R / (\pi \epsilon_2^2 + \lambda_k) \Z \to B^{2n-2}_{\epsilon_1} \times \partial B^2_{\epsilon_2}$ given by $\gamma_x(t) = (x,\epsilon_2 e^{(2i\pi/(\pi \epsilon_2^2 + \lambda_k))t})$
for each $x \in B_{\epsilon_1}^{2n-2}$).
\end{enumerate}
\end{defn}

The reason why we need such a product disk is that in the proof of
Theorem \ref{theorem:reeborbitlowerboundarounddivisors},
we will need to partially compactify $M_1$ by replacing $B^{2n-1}_{\epsilon_1} \times B^2_{\epsilon_2}$
with $B^{2n-1}_{\epsilon_1} \times S^2_{\epsilon_3}$ where $S^2_{\epsilon_3}$
is a symplectic sphere of symplectic area $\epsilon_3 > \epsilon_2$
(see Step $2$ of the proof of Theorem \ref{theorem:reeborbitlowerboundarounddivisors}).
This partial compactification is needed so that we can create a GW triple which in turn will be used to find Reeb orbits.

\begin{theorem} \label{label:nicecontactneighbourhoodexistence}
Let $(M,\cup_i S_i,\theta)$ be an orthogonal strongly numerically $\Q$-Gorenstein
positively wrapped divisor with discrepancies $a_1,\cdots,a_l$
and wrapping numbers $\lambda_1,\cdots,\lambda_l$ associated to $S_1,\cdots,S_l$ respectively.
Let $\epsilon>0$ be small.
Then there exists a link region $(M_1,g : M \setminus \cup_i S_i \to \R)$ so that:
\begin{enumerate}
\item \label{item:propertyreeborbitlocationsandindices}
 Let $V_S =\N \big< S_1, \cdots, S_l \big>$ be the free commutative monoid  generated by $S_1,\cdots, S_l$.
For each element of $V_S$ given by $V := \sum_{j \in I_V} d_j S_j$
where $d_j \neq 0$ for all $j \in I_V$
and $\cap_{i \in I_V}S_j \neq \emptyset$
we have a non-empty
pseudo Morse-Bott family $B_V$ of Reeb orbits of $(\theta + dg)|_{\partial M_1}$
whose Conley-Zehnder index is given by $2 \sum_{i \in I_V} (a_i+1) d_i + \frac{1-|I_V|}{2}$ where
$\text{Size}(B_V) = 2n - |I_V|-1$.
The period of the respective Reeb orbits minus $\sum_{i \in I_V} d_i (\pi \epsilon^2 + \lambda_i)$
has absolute value less than $\epsilon^3 \Big(\sum_{i \in I_V} d_i \Big)$.
Also there exists a disk bounding each Reeb orbit in $B_V$ whose intersection with $S_i$
is $d_i$.
 We require that every Reeb orbit sits inside some family $B_V$.
\item \label{item:productballexistence}
For each $i$ there is an $(2\epsilon,\epsilon_M)$-ball product associated to $B_{S_i}$ and $S_i$ inside $M_1$
for some $\epsilon > \epsilon_M$.
\end{enumerate}
\end{theorem}

%TODO:continue here
We have the following Corollary needed in Theorem \ref{theorem:reeborbitlowerboundarounddivisors}.
\begin{corollary} \label{corollary:specificformaintheorem}
 Let $(M,\cup_i S_i,\theta)$ be an orthogonal strongly numerically $\Q$-Gorenstein
positively wrapped divisor with discrepancies $a_1,\cdots,a_l$
and wrapping numbers $\lambda_1,\cdots,\lambda_l$ associated to $S_1,\cdots,S_l$ respectively.
%Let $N \in \N$ be such that $N c_1(TM|_{M \setminus \cup_i S_i}) = 0 \in H_2(M \setminus \cup_i S_i;\Z)$.
Let $\epsilon>0$ be small.
% and let $a_i \in \Q$ be the discrepancy of $S_i$ and $\lambda_i$ the wrapping number.
Choose $j \in \{1,\cdots,l\}$ so that
\begin{itemize}
\item if $\text{md}(M,\cup_i S_i,\theta) < 0$, then $a_j < 0$ and 
$\lambda_i \geq \lambda_j$ for any $i$ with $a_i < 0$,
\item
if $\text{md}(M,\cup_i S_i,\theta) \geq 0$, then $a_j = \text{md}(M,\cup_i S_i,\theta)$ and
$\lambda_i \geq \lambda_j$ for any $i$ with $a_i = a_j$.
\end{itemize}
Then there exists a link region $(M_1,g : M \setminus \cup_i S_i \to \R)$,
a pseudo Morse-Bott family of Reeb orbits $B$ of $(\theta + dg)|_{\partial M_1}$
and a $(2\epsilon,\epsilon_M)$-ball product associated to $B$ and $S_j$ inside $M_1$ for some $\epsilon_M < \epsilon$
so that:
\begin{enumerate}
%\item $\epsilon > \epsilon_M$.
\item  $\text{lSFT}(B) = 2\text{md}(M,\cup_i S_i,\theta)$
and the period of $B$ is $\pi \epsilon_M^2 + \lambda_j$.
\item Any other Reeb orbit $\gamma$ not contained in $B$ of period strictly less than the period of $B$
minus $\epsilon^2$
satisfies $\text{lSFT}(\gamma) \geq 0$,  $\text{lSFT}(\gamma) > \text{lSFT}(B)$ and
is non-degenerate.
\end{enumerate}
\end{corollary}
\proof of Corollary \ref{corollary:specificformaintheorem}.
Define $\alpha := (\theta + dg)|_{\partial M_1}$.
Without loss of generality, we will assume that $j=1$.
First of all for $\epsilon>0$ small we choose a link region $(M_1,g : M \setminus \cup_i S_i \to \R)$ and a constant $\epsilon_M$
as in Theorem \ref{label:nicecontactneighbourhoodexistence} above.
We will perturb this link region slightly later on.
We let $B$ be equal to $B_{S_1}= B_{S_j}$ from this theorem.
By property (\ref{item:productballexistence}) of Theorem \ref{label:nicecontactneighbourhoodexistence}
we have an $(2\epsilon,\epsilon_M)$-ball product associated to $B$ and $S_i$ inside $M_1$.
%where and $\epsilon_M$ can be made arbitrarily close to $\epsilon$.
This implies in particular that the period of $B$ is $\pi \epsilon_M^2 + \lambda_1$.

Property (\ref{item:propertyreeborbitlocationsandindices}) of Theorem \ref{label:nicecontactneighbourhoodexistence}
combined with the fact that the period of $B$ is $\pi \epsilon_M^2 + \lambda_1$
tells us that $\pi(\epsilon^2 - \epsilon^2_M) < \epsilon^3$.
Hence we can make sure that  
%$\epsilon>0$ from Theorem \ref{label:nicecontactneighbourhoodexistence}
$\epsilon$ is small enough so that for any positive integers $d_i$ labelled by $I_V \subset \{1,\cdots,l\}$
that satisfy 
$\sum_{i \in I_V} d_i(\pi \epsilon^2 + \lambda_i) + \epsilon^3(\sum_{i \in I_V} d_i) < \pi \epsilon_M^2 + \lambda_1 - \epsilon^2$
also satisfy
$\sum_{i \in I_V} d_i \lambda_i < \lambda_1$.
Hence Property (\ref{item:propertyreeborbitlocationsandindices}) of Theorem \ref{label:nicecontactneighbourhoodexistence}
tells us that if the period of a pseudo Morse Bott family of Reeb orbits $B_V$ represented by $V = \sum_{i \in I_V} d_i$
has period less than the period of $B$ minus $\epsilon^2$ then $\sum_{i \in I_V} d_i \lambda_i < \lambda_1$.

This means every pseudo Morse-Bott family $B_V$ of Reeb orbits described as in  Theorem \ref{label:nicecontactneighbourhoodexistence}
of period less than the period of $B$ minus $\epsilon^2$
satisfies \[\quad \text{CZ}(B_V) - \frac{1-|I_V|}{2} \geq 2,  \quad \text{and} \quad \text{CZ}(B_V) - \frac{1-|I_V|}{2} > 2(a_1+1).\]
Hence
$\text{lSFT}(B_V) = \text{CZ}(B_V) - \frac{1}{2}\text{Size}(B_V) + (n-3)$
satisfies:
\[\text{lSFT}(B_V) \geq 0 \quad \text{and} \quad \text{lSFT}(B_V) > 2a_1 = \text{lSFT}(B).\]

Now the only problem is that the orbits $B_V$ might be degenerate.
To fix this we use the fact that a generic perturbation of our contact form
$\alpha$ has non-degenerate Reeb orbits.
By Gray's stability theorem we can ensure this perturbation is of the form
$q \alpha$ for some function $q : \partial M \to (0,\infty)$.
Also by rescaling $q$ slightly we can ensure that if we have any Reeb orbit $\gamma$ of $q\alpha$
that is $C^\infty$ close to a Reeb orbit of $\alpha$
of period $\geq 2 \pi \epsilon_M^2 + \lambda_1$
then the period of $\gamma$ also satisfies the same inequality.

Let ${\mathcal N}$ be a small neighborhood of $B$ in $\partial M_1$
with the property that
all the pseudo Morse-Bott families $B_V$ of period $< \pi \epsilon_M^2 + \lambda_1$
do not intersect ${\mathcal N}$.
Choose a bump function $\beta$ on $\partial M_1$ equal to $0$ near $B$
and $1$ outside ${\mathcal N}$.
Then $\alpha_1 := (\beta q + (1-\beta) )\alpha$ (for $q$ sufficiently $C^\infty$ close to $1$) has the property that
all its Reeb orbits of period $< \pi \epsilon_M^2 + \lambda_1$
are non-degenerate and $C^\infty$ close to Reeb orbits of $\alpha$
of period $< \pi \epsilon_M^2 + \lambda_1$.

Lemma \ref{lemma:perturbedindexcalculation} then tells us that for a small enough perturbation,
all Reeb orbits of $\alpha_1$ of period less than $\pi \epsilon_M^2 + \lambda_1$ minus $\epsilon^2$
have lSFT index greater than or equal to $0$
and also strictly greater than $2a_1 = \text{lSFT}(B)$.
A standard Moser argument then tells us that instead of perturbing $\alpha$
we can just perturb $M_1$ slightly giving us a new
link region with the properties we want.

%Because $\text{CZ}(B_V) - \frac{1-|I_V|}{2}$ is an even multiple of $\frac{1}{N}$
%we get that:
% \[\text{CZ}(B_V) - \frac{1-|I_V|}{2} > \text{max}\Big(2(-\frac{1}{2N}+1)+\frac{1}{2N},2(a_1+1)\Big)
%= \text{max}(2-\frac{1}{2N},2a_1+2)\]
%Hence:
%\begin{align}
%\begin{split} \label{indexinequality_eqn}
%\text{CZ}(B_V) - \frac{1}{2}\text{Size}(B_V) = \text{CZ}(B_V) - \frac{1}{2}( 2n -2 + (1-|I_V|)) \\
%> \text{max}(3-n-\frac{1}{2N},2a_1 + 3-n).
%\end{split}
%\end{align}

\qed

\bigskip

We will now start proving Theorem \ref{label:nicecontactneighbourhoodexistence}.
From now on we fix an orthogonal strongly $\Q$-Gorenstein positively intersecting divisor $(M,\cup_i S_i,\theta)$ with discrepancies $a_1,\cdots,a_l$
and wrapping numbers $\lambda_1,\cdots,\lambda_l$ associated to $S_1,\cdots,S_l$ respectively.
We need some preliminary lemmas
(Lemmas \ref{lemma:choiceoftheta}, \ref{lemma:hamiltonianconleyzehnderindexcomparisonstandard} and \ref{lemma:hamiltonianconleyzehnderindexcomparison})
before we prove Theorem \ref{label:nicecontactneighbourhoodexistence}.
Before we state these Lemmas, we construct nice fibrations around each $S_i$.
Lemma \ref{lemma:choiceoftheta} will give us a function $g : M \setminus \cup_i S_i \to \R$ so that
$\theta_g := \theta + dg$ behaves well with respect to these fibrations.
In the proof of Theorem \ref{label:nicecontactneighbourhoodexistence}, we construct some nice Hamiltonian $H : M \to \R$
so that $C := H^{-1}(\delta)$ for $\delta > 0$ small will be a contact manifold with contact form $\theta + dg$.
Note that there is a one-to-one correspondence between Hamiltonian orbits of $H$ inside $H^{-1}(\delta)$
and Reeb orbits of $C$.
In order to calculate the indices of the Reeb orbits of $C$, we need a relationship between these indices and the
Conley-Zehnder indices of the respective Hamiltonian orbits of $H$ inside $H^{-1}(\delta)$.
This is the reason why we have Lemmas \ref{lemma:hamiltonianconleyzehnderindexcomparisonstandard} and \ref{lemma:hamiltonianconleyzehnderindexcomparison}.
Lemma \ref{lemma:hamiltonianconleyzehnderindexcomparisonstandard}  is an easier case of Lemma
\ref{lemma:hamiltonianconleyzehnderindexcomparison} and is used to prove Lemma
\ref{lemma:hamiltonianconleyzehnderindexcomparison}.

\bigskip

For $I \subset \{1,\cdots,l\}$ we write $S_I := \cap_{i \in I} S_i$.
By \cite[Lemma 5.14]{McLean:affinegrowth}
(or \cite[Theorem 2.12]{McLeanTehraniZinger:normalcrossings})
 we have that for each $I \subset \{1,\cdots,l\}$
there are open neighborhoods $U_I$ of $S_I$,
and smooth fibrations $\pi_I : U_I \twoheadrightarrow S_I$ satisfying the following properties:
\begin{enumerate}
\item \label{item:containsa}
$U_{I \cup J} = U_I \cap U_J$.
\item
Each fiber of $\pi_I$ is a symplectic submanifold symplectomorphic to $\prod_{i \in I} \dot{\D}_\epsilon^i$
where $\dot{\D}^i_\epsilon \subset \C$ is the open disk of radius $\epsilon$ labeled by $i \in I$. Here $\epsilon$ is a fixed constant
that can be made as small as we like.
For $J \subset I$, the fibers of $\pi_J|_{U_I}$ are contained in the fibers $\prod_{i \in I} \dot{\D}_\epsilon^i$ of $\pi_I$ and are of the form
$\prod_{i \in J} \dot{\D}_\epsilon^i \times \prod_{i \in I \setminus J} \{z_i\}$ for points $z_i \in \dot{\D}_\epsilon^i$.
Also for all $I \subset \{1,\cdots,l\}$ we have $\pi_I(U_I \setminus \cup_{i \notin I} U_{\{i\}})$ is a compact subset of $S_I$.
\item The set of vectors in $U_I$ which are symplectically orthogonal to the fibers of $\pi_I$
induce an Ehresmann connection whose parallel transport maps are in $U(1)^{|I|}$ where each $U(1)$ in this product
rotates the disk $\dot{\D}_\epsilon^i$ in the product $\prod_{i \in I} \dot{\D}_\epsilon^i$.
The zero section is equal to $S_I$.
\end{enumerate}
We define $\pi_i,U_i$ as $\pi_{\{i\}}$ and $U_{\{i\}}$ respectively.
Because the structure group of $\pi_i$ is $U(1)$ we have functions $r_i$ given by the radial coordinate in the fibers
which generates the $U(1)$ action rotating the fibers.
Here $r_i$ is a radial coordinate associated to $S_i$ for each $1 \leq i \leq l$.

We define $\rho : [0,\epsilon^2) \to \R$
to be equal to $1$ near $0$ and equal to $0$ near $\epsilon^2$
and define $\rho_{r_i^2} : M \setminus S_i \to \R$ to be $\rho(r_i^2)$ inside $U_i$ and $0$ outside $U_i$.
Now on each fiber we have an angle coordinate $\vartheta_i$ which is only well defined up to adding a constant,
and so $d\vartheta_i$ is well defined on each fiber.
By abuse of notation we will let $d\vartheta_i$ be a $1$-form on $U_i \setminus S_i$
whose restriction to each fiber is $d\vartheta_i$ with the additional property that $d\vartheta_i(X) = 0$ for any vector $X$
symplectically orthogonal to the fibers of $\pi_i$. Note this $1$-form may not be closed.

\begin{lemma} \label{lemma:choiceoftheta}
After making the radius $\epsilon$ of the disks $\dot{\D}^i_\epsilon$ slightly smaller,
we can find a function $g : M \to \R$ with the property that
$\theta + dg$ restricted to any fiber $\prod_{i \in I} \dot{\D}^i_\epsilon$ of $\pi_I$
is equal to $\prod_{i \in I} (\frac{1}{2}r_i^2 + \frac{1}{2\pi}\lambda_i) d\vartheta_i$.
%(here $r_i,\vartheta_i$ are polar coordinates on $\D^i_\epsilon$).
We also have that the norm with respect to some metric on $M$ of $\theta+dg - \sum_i \rho_{r_i^2}\frac{1}{2\pi}\lambda_id\vartheta_i$
is bounded.
\end{lemma}
\proof of Lemma \ref{lemma:choiceoftheta}.
To give an idea of the proof, we will give a sketch of the proof in the case where we only have one divisor $S_1$.
In this case we choose some $1$-form whose restriction to each fiber of
$\pi_i$ is $\frac{1}{2}r_1^2 + \frac{1}{2\pi} \lambda_1 d\vartheta_1$.
The difference between this form and $\theta$ is fiberwise exact by the definition of wrapping number
and so there is a function $g$ so that $\theta + dg$ restricted to each fiber is
$\frac{1}{2}r_1^2 + \frac{1}{2\pi} \lambda_1 d\vartheta_1$.
When we have more than one divisor, one proceeds by induction on the strata of $\cup_i S_i$
and so that each inductive step is similar to the case of one divisor as above.

Choose a total ordering $\preceq$ on the set of subsets $I \subset \{1,\cdots,l\}$ so that if $|I| > |J|$
then $J \preceq I$.
We write $J \prec I$ if $J \preceq I$ and $J \neq I$.
Suppose that on some neighborhood of $\cup_{J \prec I} S_J$
there is a smooth function $g^\prec : M \setminus \cup_i S_i$ so that for each $J \prec I$,
$\theta + dg^\prec$ restricted to a fiber $\prod_{i \in J} \D^i_\epsilon$ of $\pi_J$
is equal to $\prod_{i \in J} (\frac{1}{2}r_i^2 + \frac{1}{2\pi}\lambda_i) d\vartheta_i$ near
$\cup_{K \prec I} S_K$.

For each $p \in S_I \setminus \cup_{j \notin I} S_i$ there is a function $f_p \in C^\infty(\pi_I^{-1}(p))$ so that $(\theta + dg^\prec + df_p)|_{\pi_I^{-1}(p)} = \prod_{i \in I} (\frac{1}{2}{r_i^2 + \frac{1}{2\pi}} \lambda_i) d\vartheta_i$ near $S_I$
using the definition of wrapping number.
We can ensure that $f_p$ smoothly depends on $p$ and that $f_p = 0$ for $p$ near $\cup_{j \notin I} S_j$.
We can construct a smooth function $f$ on $M \setminus \cup_i S_i$
so that $f|_{\pi^{-1}(p)} = f_p$ for all $p \in S_I \setminus \cup_{j \notin I} S_i$
and so that $f = 0$ near $\cup_{j \notin I} S_i$.
Hence $g^\preceq := g^\prec + f$ has the property that
$\theta + dg^\preceq$ restricted to a fiber $\prod_{i \in J} \D^i_\epsilon$ of $\pi_J$
is equal to $\prod_{i \in J} (\frac{1}{2}r_i^2 + \frac{1}{2\pi}\lambda_i) d\vartheta_i$ near $\cup_{K \preceq I} S_K$
for each $J \preceq I$.
Hence our induction step is done and we have proved the first part of the lemma after
shrinking the radius of our fibers.
Therefore we define $g$ to be $g^\preceq$ in the case where $I$
is the last element of our total ordering.

We now need to show that the norm of $\theta+dg - \sum_i \rho_{r_i^2}\frac{1}{2\pi} \lambda_i d\vartheta_i$ is bounded.
Let $p \in \cup_i S_i$. Let $I \subset \{1,\cdots,l\}$ be the set of $i$ with $S_i$ containing $p$.
Choose a $1$-form $\beta$ near $p$ with the property that $d\beta = \omega$
and with the property that the restriction of $\beta$ to each fiber of $\pi_I$ is $\sum_{i \in I} \frac{1}{2}r_i^2 d\vartheta_i$.
Near $p$ we have $\theta + dg - \sum_{i \in I} \rho_{r_i^2} \frac{1}{2\pi}\lambda_i d\vartheta_i - \beta$ is a closed $1$-form
because $\rho_{r_i^2} = 1$ near $p$ for $i \in I$.
Hence it is equal to $dg^p$ near $p$ for some function $g^p : M \setminus \cup_i S_i \to \R$.
Also by construction its restriction to the fibers of $\pi_I$ vanish
which means that $dg^p$ restricted to the fibers of $\pi_I$ vanish near $p$ which means that $g^p$ is constant along the fibers of $\pi_I$
and hence $g^p$ smoothly extends over $\cup_I S_i$ near $p$.
Putting all of this together we get that the norm of $\theta + dg - \sum_{i \in I} \rho_{r_i^2}\frac{1}{2\pi} \lambda_i d\vartheta_i - \beta$
is bounded near $p$.
Becuse the norm of $\beta$ and $\rho_{r_i^2}$ is also bounded
we then get that the norm of $\theta+dg - \sum_i \rho_{r_i^2} \frac{1}{2\pi} \lambda_i d\vartheta_i$ is bounded near $p$.
Because $\cup_i S_i$ is compact we then get that the norm of $\theta+dg - \sum_i \rho_{r_i^2} \frac{1}{2\pi} \lambda_i d\vartheta_i$ is bounded.
\qed

\bigskip
%Recall that for a Hamiltonian $H$ on a symplectic manifold $(M,\omega)$ we define $X_H$ to be the unique	vector
%field satisfying $i_{X_H} \omega = dH$.
The following definition shows that the Conley-Zehnder index of a Hamiltonian orbit is defined in a very similar way to the Conley-Zehnder
index Reeb orbit.

\begin{defn}
Let $(X,\omega_X)$ be a symplectic manifold with a choice of compatible almost complex structure and
a choice trivialization of
the canonical bundle $\kappa$ of $M$ and let $H \in C^\infty(M)$ be a Hamiltonian.
Then for any $T$-periodic orbit $\gamma$ of $H$ we can define its {\bf Conley-Zehnder index} $\text{CZ}(\gamma)$
as follows:
First choose a trivialization of the Hermitian bundle $\gamma^* \oplus_{j=1}^{N} TM$
so that its highest exterior power coincides with $(\kappa^*)^{\otimes N}$.
Then define $\text{CZ}(\gamma)$ to be $\frac{1}{N}$ multiplied
by the Conley-Zehnder index of the family of symplectic matrices
induced by the flow of $X_H$ along $\gamma$ in our trivialization.
\end{defn}

\begin{lemma} \label{lemma:hamiltonianconleyzehnderindexcomparisonstandard}
Let $(C,\xi)$ be a contact manifold with associated contact form $\alpha$
and let $h : \R \to \R$ be a function with $h' < 0,h''>0$ and $h'(0) = -1$.
Let $\widehat{C} := C \times \R$ be the symplectization of $C$ with symplectic form $d(e^r \alpha)$
where $r$ parameterizes $\R$.
Let $\gamma(t)$ be a Reeb orbit of $\alpha$ of period $L$ with a choice of  trivialization of the symplectic vector bundle $\oplus_{j=1}^{N}TM$ along this orbit.
This choice of trivialization induces a choice of trivialization of $\gamma^* \oplus_{j=1}^{N} \xi$ in a natural way.
Then the Hamiltonian $L h(e^r)$ has a $1$ periodic orbit equal to $\gamma(Lt)$ inside $C \times \{0\} = C$ and its Conley-Zehnder index is equal to
$\text{CZ}(\gamma)-\frac{1}{2}$.

\end{lemma}
\proof of Lemma \ref{lemma:hamiltonianconleyzehnderindexcomparisonstandard}.
The key idea of the proof here is that
the family of matrices associated to the flow
of $L h(e^r)$ along our Reeb orbit is equal to
the linearization of the Reeb flow restricted to the contact structure
plus a very specific family of $2 \times 2$ matrices.
The Conley-Zehnder index of such a family of $2 \times 2$ matrices is $-\frac{1}{2}$
and hence we get our result by \ref{item:additiveunderdirectsum}.

We identify $C$ with $C \times \{0\}$ in the natural way.
%Let $R$ be the Reeb vector field of $\alpha$.
Let $K := L h(e^r)$ and
let $\gamma_L(t) := \gamma(Lt)$ be the Hamiltonian $1$-periodic orbit of $K$.
We will define $\text{CZ}(\gamma_L)$ to be the Conley-Zehnder index of the orbit $\gamma_L$.
By abuse of notation we define $\alpha$ to be a $1$-form on $\R \times C$ by pulling back $\alpha$ along the natural projection
$\R \times C \twoheadrightarrow C$.
Let $R$ be the vector field tangent to $C \times \{r\}$ and equal to the Reeb flow of $\alpha|_{C \times \{r\}}$ for each $r$.
Along $\gamma$ we have that the symplectic vector bundle $T\widehat{C}$ splits into two symplectically orthogonal subspaces
$\xi \oplus \xi^\perp$ where $\xi$ is the contact distribution and $\xi^\perp$
is the symplectic vector space spanned by $X_{\alpha}$ and $R$.
The vectors $(X_{\alpha},R)$ form a symplectic basis and so give us a trivialization of the symplectic vector bundle $\xi^\perp$ along $\gamma$ and hence because
$\oplus_{j=1}^{N} T\widehat{C}$ also has a trivialization along $\gamma$ we also get a trivialization of the symplectic vector bundle $\oplus_{j=1}^{N}  \xi$
along $\gamma$.
The linearization of the flow of $R$ along $\gamma$ gives us a sequence of linear maps $\xi_{\gamma(0)} \to \xi_{\gamma(t)}$.
Using the trivialization of $\gamma^* \oplus_{j=1}^{N} \xi$ we get that this linearization can be viewed as
a family of symplectic matrices $S_t$ in $\C^{(n-1)N}$.
The Conley-Zehnder index of $S_t$ is $N$ times the Conley-Zehnder index of $\gamma$.
The linearization of the flow of $X_{K}$ along $\gamma$ preserves our splitting $\xi \oplus \xi^\perp$
and its restriction to $\xi$ is equal to the linearization of the flow of $R$ after a linear re-parameterization of time.
Also the linearization of the flow of $X_{K}$ preserves $R$ and sends $X_{\alpha}$
to $X_{\alpha} - a t R$ where $a > 0$ and $t$ is time.
Hence using our trivializations of $\gamma^* \oplus_{j=1}^{N} \xi$ and $\gamma^* \xi^\perp$
we get a family of symplectic matrices $S_{Lt} \oplus \oplus_{j=1}^{N} Q_t$ on $\C^{(n-1)N} \oplus
\oplus_{j=1}^{N}\C$, $t \in [0,1]$
where the first factor comes from a trivialization of $\gamma^* \oplus_{j=1}^{N} \xi$ and the
remaining $N$ factors come from our trivialization of $\xi^\perp$.
The family of symplectic matrices $Q_t$
is equal to $\Big(\begin{array}{cc} 1 & 0 \\ -t a & 1 \end{array}\Big)$ and has Conley-Zehnder index equal to $-\frac{1}{2}$
(this is calculated using the normalization property from \cite[Theorem 55]{Gutt:ConleyZehnderMatricies}).
Also the Conley-Zehnder index of $ S_{Lt}$, $t \in [0,1]$ is the Conley-Zehnder index of $S_t$, $t \in [0,L]$.
So the Conley-Zehnder index of $S_{Lt} \oplus \oplus_{j=1}^{N} Q_t$
is $\text{CZ}(S_{Lt}) - \frac{1}{2} N$ by \ref{item:additiveunderdirectsum}.
Hence $\text{CZ}(\gamma_L) = \text{CZ}(\gamma) - \frac{1}{2}$.
\qed

\bigskip

\begin{defn}
Let $K$ be a Hamiltonian on a symplectic manifold $(X,\omega_X)$ and $B$ is a set of fixed points
of its time $T$ flow.
We say that $B$ is {\bf isolated} if any such fixed point near $B$ is contained in $B$.
If $B$ is a path connected topological space and the canonical bundle of our symplectic manifold has a choice of trivialization then
by Lemma \ref{lemma:homotopyofsymplecticpaths} we have that every such Hamiltonian orbit has the same Conley-Zehnder index
and we will write $\text{CZ}(B,K)$ for the Conley-Zehnder index.
The set $B$ is said to be {\bf Morse-Bott} if $B$ is a submanifold and
$\text{ker}(D\phi_{K}^T-\text{id}) = TB$ along $B$ where $\phi_K^T : X \to X$ is the time $T$ Hamiltonian flow
of $K$.
Sometimes we write $B$ {\bf is Morse-Bott for $K$}.
If we wish to emphasize which symplectic form we are using then if $\omega$ is our symplectic form
we will say that $(K,B,\omega)$ is Morse-Bott.
\end{defn}
%If we have some Hamiltonian $K$ and we have some set of fixed points $B$ of $\phi^1_{K}$ (the time $1$-flow of %$X_{K}$)
%then we say that such a set of fixed points is {\it isolated} if any sufficiently close fixed point is contained in $B$.
%Sometimes we write $\text{CZ}(B,K,\omega)$ if we wish to stress the fact that we are using the symplectic form %$\omega$.
%A Hamiltonian $K$ with a family of fixed points $B$ is said to be {\it Morse-Bott} if $B$ is a submanifold
%and 
%We say a family of Reeb orbits $B$ is {\it Morse-Bott} of period $L$ if $B$ is a submanifold $M$ and every
%Reeb orbit of period $L$ starting in $B$ has linearized return map $D$ satisfying $\text{ker}(D-\text{id}) = TB$.
%Sometimes we will write $\text{CZ}(B,\alpha)$ for the Conley-Zehnder index of our family of Reeb orbits $B$ using our %contact form $\alpha$.

\begin{lemma} \label{lemma:hamiltonianconleyzehnderindexcomparison}
Let $(W,\omega_W)$ be a symplectic manifold
with a choice of trivialization of
the $N$th tensor power of its canonical bundle
and let $\theta_W$ be a $1$-form satisfying $d\theta_W = \omega_W$.
Let $K$ be a Hamiltonian with the property that $b := -i_{X_{\theta_W}} dK > 0$.
This means $C_r := {K}^{-1}(r)$ is a contact manifold with contact form $\alpha_r := \theta_W|_{C_r}$.
Let $B \subset W$ be a connected submanifold transverse to $C_r$ for each $r$ so that $B_r := C_r \cap B$
is a Morse-Bott submanifold of the contact manifold $(C_r,\alpha_r)$ of period $L_r$ where $L_r$ smoothly depends on $r$.
Suppose that $b = L_0$ along $B_0$
and that $db(V) < \frac{d(L_r)}{dr}|_{r = 0}$ along $B_0$ where $V$ is a vector field tangent to $B$
satisfying $dK(V) = -1$.
Then $B_0$ is Morse-Bott for $K$ and $\text{CZ}(B_0,K)=\text{CZ}(B_0,\alpha_0)-\frac{1}{2}$.
\end{lemma}
\proof Lemma \ref{lemma:hamiltonianconleyzehnderindexcomparison}.
The key idea of this Lemma is to deform $K$ through appropriate Hamiltonians
so that $K$ looks like the Hamiltonian
$Lh(e^r)$ from Lemma \ref{lemma:hamiltonianconleyzehnderindexcomparisonstandard}.
We will deform $\theta_W$ and $K$ through appropriate forms, and then use Lemmas \ref{lemma:homotopyofsymplecticpaths} and \ref{lemma:hamiltonianconleyzehnderindexcomparisonstandard} for our Conley-Zehnder index calculations. This will be done in $3$ steps.
In the first step we deform $\theta_W$ so that a neighborhood of $C_0$
looks like a symplectization. We do this by `stretching' $\theta_W$
in the direction of the gradient of $K$.
In the second step we deform our Hamiltonian $K$
through functions of the form $q \circ K$ where $q$ has very large second derivative at $0$
(ensuring that the families of orbits in $C_0$ are still isolated).
The third step just applies the previous Lemma to our new Hamiltonian giving us our result.

{\it Step $1$}. In this step we will deform $\theta_W$ through certain primitives of symplectic forms.
Because $b= -i_{X_{\theta_W}} dK = i_{X_{K}} \theta_W$ and $b = L_0$ along $B_0$
we have that for each Reeb orbit $\gamma(t)$,$t \in [0,L_0]$ of $C_0$ starting in $B_0$, there is a $1$-periodic Hamiltonian orbit
$\gamma(L_0 t)$ of $K$.
By considering the flow of $\frac{1}{b} X_{\theta_W}$ we have that a neighborhood of $C_0$
is identified with $C_0 \times (-\epsilon,\epsilon)$ for some $\epsilon>0$ satisfying:
\begin{enumerate}
\item If $r : C_0 \times (-\epsilon,\epsilon) \twoheadrightarrow (-\epsilon,\epsilon)$ is the natural projection map
then $K = -r$.
\item $X_{\theta_W} = b \frac{\partial}{\partial r}$.
\end{enumerate}
Let $\Phi : [0,1] \times C_0 \times (-\epsilon,\epsilon) \to [0,1] \times C_0 \times (-\epsilon,\epsilon)$
be a smooth map sending $(w,y,r)$ to $(w,y,(1-w)r)$.
Let $\pi_{C_0 \times (-\epsilon,\epsilon)} :  [0,1] \times C_0 \times (-\epsilon,\epsilon) \twoheadrightarrow C_0 \times (-\epsilon,\epsilon)$ be the natural projection.
Define $\widetilde{\theta}_W := \Phi^* \pi_{C_0 \times (-\epsilon,\epsilon)}^*\theta_W$.
Let $\iota_w : C_0 \times (-\epsilon,\epsilon) \hookrightarrow [0,1] \times C_0 \times (-\epsilon,\epsilon)$
send $(y,r)$ to $(w,y,r)$.
Let $\pi_{C_0} : C_0 \times (-\epsilon,\epsilon) \twoheadrightarrow C_0$
be the natural projection.
Let $\eta : [0,1] \times (-\epsilon,\epsilon) \to \R$
be a smooth function satisfying:
\begin{enumerate}
\item $\frac{\partial}{\partial r}\eta(w,r) \geq 0$.
\item $\eta(0,r) = 1$, $\eta(1,r) = e^r$ and $\eta(w,0) = 1$ for all $w,r$.
\end{enumerate}
Define $\theta_{W,w} := \eta(w,r) \iota_w^* \widetilde{\theta}_W.$

For $r$ near $0$ we have that the restriction of $\theta_{W,w}$ to $C_r$ is a contact form.
So after shrinking $\epsilon$ we can assume that the restriction of $\theta_{W,w}$ to $C_r$ is a contact form for all $r$.
We will write $b$ as a function $b(y,r)$ of two variables $(y,r) \in C_0 \times (-\epsilon,\epsilon)$. 
Let $R_w$ be a vector field tangent to $C_r$ and equal to the Reeb vector field of $\theta_{W,w}|_{C_r}$  for each $r$.
Note that $\zeta_w := \pi_{C_0 \times (-\epsilon,\epsilon)} \circ \Phi \circ \iota_w$ sends $(y,r)$ to $(y,(1-w)r)$.
Because $b \frac{\partial}{\partial r} =  X_{\beta}$, $\zeta_w^* b(y,r) = b(y,(1-w)r)$
and $\zeta_w^* \frac{\partial}{\partial r} = \frac{1}{(1-w)} \frac{\partial}{\partial r}$
we have $\zeta_w^* X_{\beta} = \frac{b(y,(1-w)r)}{(1-w)} \frac{\partial}{\partial r}$.
Hence $i_{\frac{\partial}{\partial r}}d\theta_{W,w}
= \Big(\frac{\eta(w,r)(1-w)}{b(y,(1-w)r)} + \frac{\partial}{\partial r}\eta(w,r)\Big) \widetilde{\theta}_W$.
So near $C_0$ we then get that $i_{R_w} i_{\frac{\partial}{\partial r}}d\theta_{W,w} > 0$.
This means that $d\theta_{W,w}$ is a symplectic form after shrinking $\epsilon$.

\bigskip

{\it Step $2$}. In this step we construct a smooth family of Hamiltonians depending on $w$ so that $B_0$ is a Morse-Bott submanifold using the symplectic form $d\theta_{W,w}$.
Let $B^w := \zeta_w^{-1}(B)$.
We have that $B^w|_{C_r} \subset C_R$ is Morse-Bott for the Reeb flow of $\theta_{W,w}$ whose period is a smooth function $L(w,r)$ of $w$ and $r$.

Let $q_\lambda : \R \to \R$ be a smooth family of functions with the property that $q_\lambda(0) = 0$,
$q_\lambda'(0) = 1$ and $q_\lambda''(0) = \lambda$. We will also assume that $q_0 = \text{id}$.
For $0<\eta<1$ let $\rho_\eta : [0,1] \to [0,1]$ be a bump function with the property that $\rho_\eta(x) = 0$ near $0$ 
and $\rho_\eta(x) = 1$ inside $[\eta,1]$.
Define $K_{w,\eta,\lambda} :=  q_{\lambda \rho_\eta(w)} \circ K$.
Here $\frac{1}{\lambda} \ll \eta \ll 1$.
For any $1$-form $\alpha$ on $C_0 \times (-\epsilon,\epsilon)$ define
$X^w_\alpha$ to be the $d\theta_{W,w}$
dual of $\alpha$.
We define $s_w := i_{X^w_{dK_{w,\eta,\lambda}}} \theta_{W,w}$. Here $s_w$ depends on $\lambda$ and $\eta$
but we suppress this notation for simplicity.
We have that $s_w = (q'_{\lambda \rho_\eta(w)} \circ K) b_w$ where $b_w := i_{X^w_{dK}} \theta_{W,w}$.
Let $V_w$ be a smooth family of vector fields parameterized by $w \in [0,1]$
with the property that $V_w$ is tangent to $B^w$ with $dK(V_w) = -1$
which means that $dr(V_w) = 1$. We will also assume that $V_0|_B = V$.
Now
\[ds_w(V_w) = (q'_{\lambda \rho_\eta(w)} \circ K) db_w(V_w) + (q''_{\lambda \rho_\eta(w)} \circ K) dK(V_w) b_w = \]
\[(q'_{\lambda \rho_\eta(w)} \circ K) db_w(V_w) - (q''_{\lambda \rho_\eta(w)} \circ K) b_w.\]

Because $db_0(V_0) < \frac{\partial L(0,r)}{\partial r}|_{r=0}$ and $q''_{\lambda \rho_\eta(w)}>0$
we get for $\eta>0$ small enough that $ds_w(V_w) < \frac{\partial L(w,r)}{\partial r}|_{r=0}$ along $C_0$
for $w < \eta$. Such a bound is true for any $\lambda$ as long as $\eta \ll 1$.
Now we fix such a small $\eta$ and increase $\lambda$.
Because $q'_{\lambda \rho_\eta(w)} > 0$ and $q''_{\lambda \rho_\eta(w)}(0) = \lambda$ for $w \in [\eta,1]$
we get for $\lambda$ sufficiently large that $ds_w(V_w) < \frac{\partial L(w,r)}{\partial r}|_{r=0}$ along $C_0$
for $w \in [\eta,1]$.
Hence $ds_w(V_w) < \frac{\partial  L(w,r)}{\partial r}|_{r=0}$ along $C_0$ for $\frac{1}{\lambda} \ll \eta \ll 1$.

Let $p$ be a point in $B_0$ and let $D_p^w : TM_p \to TM_p$
be the linearization of the time $1$ flow of $X^w_{K_{w,\eta,\lambda}}$.
Because $X^w_{K_{w,\eta,\lambda}} = s_w R_w$ and $R_w = R$ along $C_0$
we get that any vector $W$ at $p$ tangent to $C_0$ but not tangent to $B_0$
satisfies $D^w_p(W) \neq W + l R_w$ for any $l \in \R$.
The inequality $ds_w(V_w) < \frac{\partial  L(w,r)}{\partial r}|_{r=0}$ combined with the fact that the Reeb
vector field $R_w$ of $\theta_{W,w}$ is tangent to $B^w$ ensures that
$D^w_p(V_w) = V_w - \kappa(p) R_w$ at $p \in B_0$ for some positive smooth function $\kappa : B_0 \to \R$.
Any vector $X$ in $M$ at $p \in B_0$ is equal to a unique sum $W + k V_w $ for some $k \in \R$ where $W$ is tangent to $C_0$
which means $D^w_p(X) = D^w_p(W) + k V_w - k\kappa(p) R_w$.
Now suppose $X$ is not tangent to $B_0$. This means $k \neq 0$ or $W$ is not tangent to $B_0$.
If $W$ is not tangent to $B_0$ then $D^w_p(W) \neq W + l R_w$ for any $l \in \R$
and so $W + k V_w - D^w_p(X) = W -D^w_p(W) + k \kappa(p) R_w$ is non-zero.
If $W$ is tangent to $B_0$ then $k \neq 0$ and $W = D^w_p(W)$ which implies that
$W + k V_w - D^w_p(X) = k \kappa(p) V_w \neq 0$.
Putting everything together we get $\text{ker}(D^w_p - \text{id}) = TB_0$ along $B_0$ for all $w$.

\bigskip

{\it Step $3$}. In this step we finish the proof by calculating Conley-Zehnder indices.
Lemma \ref{lemma:homotopyofsymplecticpaths} combined with the fact that
$K = K_{0,\eta,\lambda}$ implies that
\[\text{CZ}(B_0,K,d\theta_W) = \text{CZ}(B_0,K_{0,\eta,\lambda},d\theta_{W,0}) = \text{CZ}(B_0,K_{1,\eta,\lambda},d\theta_{W,1})\]
By Lemma \ref{lemma:hamiltonianconleyzehnderindexcomparisonstandard}, the Conley-Zehnder index of $B_0$ with respect to the Hamiltonian $K_{1,\eta,\lambda}$
and the symplectic form $d\theta_{W,1}$ is
$\text{CZ}(B_0,\alpha_0)-\frac{1}{2}$.
Putting everything together tells us that $\text{CZ}(B_0,K) = \text{CZ}(B_0,\alpha_0)-\frac{1}{2}$ which gives us our result.
\qed

\bigskip

\proof of Theorem \ref{label:nicecontactneighbourhoodexistence}.
The contact manifold we want will be a smoothing of the hypersurface with corners
given by $\cup_i \left(\{ r_i = \epsilon \} \setminus \cup_{j  \neq i} \{r_j < \epsilon \}\right)$.
The smoothing is realized by $H^{-1}(\delta)$ for $\delta > 0$ small
where $H  = \sum_i q(r_i^2)$ and where $q$ is the function pictured below.
The Hamiltonian flow of $H$ is the Reeb flow after a
time reparameterization.
The key point is that $H$ looks like a product Hamiltonian and so the
Reeb orbits should be fairly easy to find and it should also be fairly
straightforward to calculate their Conley-Zehnder indices using the previous Lemma.
The difficult part of this Lemma is showing that the associated Reeb orbits
are pseudo-Morse Bott and also setting things up so that we can calculate their indices
using Lemma \ref{lemma:hamiltonianconleyzehnderindexcomparison}.

We will prove this Theorem in $6$ steps. In the first step, we will construct our contact manifold and functions $f$ and $g$.
In the second step,  we will find our family of Reeb orbits $B_V$, and then we will show that they are pseudo Morse-Bott. In the third Step we will calculate their indices using Lemma \ref{lemma:hamiltonianconleyzehnderindexcomparison}.
The second and third steps are the longest and most technical steps.
In Step $4$ we estimate the period of these families of Reeb orbits.
In Step $5$ we construct the $(2\epsilon,\epsilon_M)$-ball product associated to $B_{S_i}$ and $S_i$ inside $M_1$ for each $i$.
Finally in Step $6$ we show that $df(X_{\theta_g}) > 0$ in $M_1 \setminus \cup_i S_i$.
Steps $4,5,6$ are very short.

\bigskip
{\it Step 1}: In this Step we will construct our link region $(M_1,g : M \setminus \cup_i S_i \to \R)$
along with its function $f$ compatible with $\cup_i S_i$.
We define $q : [0,\epsilon^2) \to \R$ to be a smooth function so that:
\begin{enumerate}
\item There is some $\epsilon_q \in (\epsilon-\frac{1}{2}\epsilon^3,\epsilon)$ with $q(x) = 0$ if and only if $x \in [\epsilon_q^2,\epsilon^2)$.
Also $q(x) = 1 - x^2$ near $x = 0$.
\item We also assume that the derivative of $q$ is non-positive
and that it is strictly negative when $q(x)$ is positive and $x \neq 0$.
\item There is a unique point $x$ with $q''(x) = 0$ and $q(x) \neq 0$.
\end{enumerate}

%TODO modify this picture
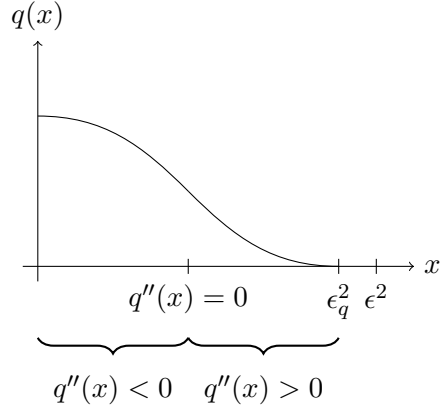
\begin{center}
\begin{tikzpicture}[domain=0:4]
%   \draw[very thin,color=gray] (-0.1,-1.1) grid (3.9,3.9);
   \draw[->] (-0.2,0) -- (5,0) node[right] {$x$};
   \draw[->] (0,-0.2) -- (0,3) node[above] {$q(x)$};
%region with positive 2nd derivative
    \draw (0,2) to[out=0,in=135] (2,1);
%region with negative 2nd derivative
    \draw (2,1) to[out=-45,in=180] (4,0);
%labels
    \draw (2,0.1) -- (2,-0.1) node[below,scale=1] {$q''(x)=0$};
    \draw (4,0.1) -- (4,-0.1) node[below,scale=1] {$\epsilon_q^2$};
    \draw (4.5,0.1) -- (4.5,-0.1) node[below,scale=1] {$\epsilon^2$};
\draw [
    thick,
    decoration={
        brace,
        mirror,
        raise=0.2cm,
	amplitude=6pt
    },
    decorate
] (0.0,-0.75) -- (2,-0.75) 
node [pos=0.5,anchor=north,yshift=-0.55cm] {$q''(x)<0$}; 
\draw [
    thick,
    decoration={
        brace,
        mirror,
        raise=0.2cm,
	amplitude=6pt
    },
    decorate
] (2,-0.75) -- (4,-0.75) 
node [pos=0.5,anchor=north,yshift=-0.55cm] {$q''(x)>0$}; 
\end{tikzpicture}
\end{center}
%Fix $0<\epsilon_1,\cdots,\epsilon_l < \epsilon$.
%We have a function $q\Big(\frac{\epsilon^2 r_i^2}{\epsilon_i^2} \Big)$ which we define to be equal to $0$
%outside the region $\{r_i \leq \epsilon_i\}$.
We define $H := \sum_i q(r_i^2)$.

Define $g : M \setminus \cup_i S_i \to \R$ as in Lemma \ref{lemma:choiceoftheta} and define $\theta_g := \theta + dg$.
%We let $q$ be the function described above.
Define $M_1 := H^{-1}([0,\delta))$.
Let $\nu : (0,\epsilon] \to \R$ be a smooth function equal to $0$ near $\epsilon^2$ and satisfying $\nu(x) = \log(x)$ for $x \leq \epsilon_q^2$. Define $\nu_{r_i^2}$ to be $\nu(r_i^2)$ inside $U_i$ and $0$ outside $U_i$. Our function $f : M \setminus \cup_i S_i \to \R$ will be defined as $\sum_i \nu_{r_i^2}$ is our function compatible with $(M_1, g : M \setminus \cup_i S_i \to \R)$.
Let $C_\delta := \partial M_1$ and $\alpha_\delta := \theta_g|_{C_\delta}$.

\bigskip
{\it Step 2}:
In this step we will find our family of Reeb orbits $B_V$, and then we will show that they are pseudo Morse-Bott.
For each $i \in \{1,\cdots,l\}$ define $h_i : [0,\epsilon^2) \to \R$ by $h_i(x) := -2q'(x) (\frac{1}{2}x + \frac{1}{2\pi}\lambda_i)$.
We can extend the function $h_i(r_i^2) : U_i \to \R$ to a function $M \setminus \cup_i S_i \to \R$ by defining it to be $0$ outside $U_i$.
Because the derivative	of $\log(h_i(x))$ tends to $-\infty$ as we approach $\epsilon_q$ from below we
get that $\frac{h_i'}{h_i}$ tends to $-\infty$ as we approach $\epsilon_q^2$ from below.
Similarly by looking at $\log(-q')$ we have that $\frac{q''}{q'}$ tends to $-\infty$ as we approach $\epsilon_q^2$ from below.
We choose $\delta$ small enough so that:
\begin{enumerate}
\item \label{item:deltaprop1}
For each $i \in I$ we have $q'(r_i^2)<0, q''(r_i^2) > 0$, $h_i'(r_i^2)<0$,
$-\frac{q''(r_i^2)}{q'(r_i^2)} > n$ and $r_i > \frac{\epsilon_q}{2}$
along $H^{-1}(\delta) \cap U_i$.
\item \label{item:deltaprop2}
For each $I \subset \{1,\cdots,l\}$
we require $\sum_{j \in I} h'_j(r_j^2) > -1$ 
along $H^{-1}(\delta) \cap U_I$.
%\item \label{item:deltaprop3}
%For each $I \subset \{1,\cdots,l\}$ and each $j \in \{1,\cdots,l\}$ we have \[\frac{1}{\sum_{i \in I} h_i(r_i^2)} \sum_{i \in I} h'_i(r_i^2) < -\frac{2\pi}{\lambda_j}\] along $H^{-1}(\delta) \cap U_I$.
\end{enumerate}

The Reeb flow $R_\delta$ of $\alpha_\delta$ is equal to $\frac{1}{b}X_H$ where $b = i_{X_H}\theta_g$.
Now let's look inside the region $U_I$ for some $I \subset \{1,\cdots,l\}$.
Define $H_I : S_I \setminus \cup_{j \notin I} S_j \to \R$
by $\sum_{i \notin I} q(r_i^2)$.
Inside $(H_I \circ \pi_I)^{-1}(0)$
we have that $H$ is a function of $(r_i^2)_{i \in I}$
given by $\sum_{i \in I} q(r_i^2)$ and so inside $C_\delta$
the coordinates $r_i$ are subject to the conditions $\sum_{i \in I} q(r_i^2) = \delta$.
Note that if $(r_i,\vartheta_i)$
are polar coordinates for $\dot{\D_i}$ then $\theta_g$ restricted to a fiber $\prod_{i \in I} \dot{\D_i}$
is $\sum_{i \in I} (\frac{1}{2}r_i^2 + \frac{1}{2\pi}\lambda_i) d\vartheta_i$.
Because $X_H$ is tangent to the fibers of $\pi_I$ inside $(H_I \circ \pi_I)^{-1}(0)$
we then get that $X_H = -\sum_{i \in I} 2q'_i(r_i^2) \frac{\partial}{\partial \vartheta_i}$.
Hence $b = i_{X_H} \theta_g = \sum_{i \in I} h_i(r_i^2)$ inside $(H_I \circ \pi_I)^{-1}(0)$.

Let $F_I$ be the set of $|I|$ tuples $(d_i)_{i \in I}$ of positive integers.
Let ${{d}} := (d_i)_{i \in I} \in F_I$.
Define $O_{a,{{d}}} := (H_I \circ \pi_I)^{-1}(0) \cap \cap_{i \in I}\{ -2q'(r_i^2)  = a d_i \}$
and $O_{{d}} := \cup_{a > 0} O_{a,{{d}}}$.
Near $C_\delta$ we have that $q'(r_i^2)$ can be expressed as a smooth function of $q(r_i^2)$ with negative derivative and so
 $O_{{d}}$ is a submanifold near $C_\delta$ which intersects $C_\delta$ transversally.
We define $O^\delta_{{d}} := O_{{d}} \cap C_\delta$.
Now every point on $O^\delta_{{d}}$ is the starting point of an orbit of $H$ contained inside a fiber $\pi_I^{-1}(p)$
representing the homology class $(d_i)_{i \in I} \in \prod_{i \in I} \Z \cong H_1(\pi_I^{-1}(p) \setminus \cup_i S_i)$.
Let $l^\delta_{{d}}$ be the period of such a family of Reeb orbits. 
To connect these computations back to the statement of the Theorem we should note that
if $V = \sum_{i \in I} d_i E_i$ is an element of our monoid $V_S$ then $B_V := O^\delta_{{d}}$.

We will now show that this is a pseudo Morse-Bott family of Reeb orbits.
Now the set $O^\delta_{{d}}$ is a manifold with corners diffeomorphic to a torus bundle over
the manifold with corners $H_I^{-1}(0)$
where the tori have dimension $|I|$.
Let $\phi^{R_\delta}_t : C_\delta \to C_\delta$ be the flow of the Reeb vector field $R_\delta$ of $\alpha_\delta$.
Let $D_t : T(C_\delta)_p \to T(C_\delta)_{\phi^{R_\delta}_t(p)}$ be the linearization of $\phi^{R_\delta}_t$
and similarly let $D^H_t : TM_p \to TM_{\phi^H_t(p)}$ the linearization of $\phi^H_t$.
We will show that $\text{ker}(D_{l_{\delta,{{d}}}}(p) - \text{id}) = TO^\delta_{{d}}$.
One can think of this as a Morse-Bott manifold with corners.
Because $X_H = b R_\delta$ and $b$ is constant along $O^\delta_{{d}}$
we have for a vector $V$ tangent to $C_\delta$ at $p \in O^\delta_{{d}}$ that $D_{t}(V) = D^H_{\frac{t}{b(p)}}(V) - \int_0^t ((\phi^{R_\delta}_s)^* (db|_{C_\delta})))(V)ds R_\delta$ for $V$.
We have
\[D^H_t(\frac{\partial}{\partial r_j^2}) = \frac{\partial}{\partial r_j^2} + 2tq''(r_j^2)\frac{\partial}{\partial \vartheta_j}\]
and $D^H_t$ is the identity on any other vector symplectically orthogonal to the fibers of $\pi_I$ and also on $\frac{\partial}{\partial \vartheta_j}$ for each $j$.

Let $V$ be a vector in $TC_\delta$ at a point $p \in O^\delta_{{d}}$.
We have
\[D_{l_{\delta,{{d}}}}(V)  = D^H_{\frac{l_{\delta,{{d}}}}{b(p)}}(V) -
\int_0^{l_{\delta,{{d}}}} ((\phi^R_s)^* (db|_{C_\delta}))(V)ds R_\delta\]
\[ =  V  
+ \frac{l_{\delta,{{d}}}}{b(p)} \sum_{j \in I} 2q''(r_j) (i_V dr_j^2) \frac{\partial}{\partial \theta_j}\]
\[ -\int_0^{l_{\delta,{{d}}}} ((\phi^R_s)^* (db|_{C_\delta}))(V)ds R_\delta.\]
Now $V = W + \sum_{j \in I} \alpha_j \frac{\partial}{\partial \vartheta_j} + \beta_j \frac{\partial}{\partial r_j^2}$ for some $\alpha_j,\beta_j \in \R$ where $W$ is symplectically orthogonal to the fibers of $\pi_I$.
Because $db|_{C_\delta}$ does not change as we flow along $O^\delta_{{d}}$, we have that $\int_0^{l_{\delta,{{d}}}} ((\phi^R_s)^* (db|_{C_\delta}))(V)ds = l_{\delta,{{d}}} \sum_{j \in I} \beta_j h_j'(r_j^2)$.
Also, $R_\delta = \frac{1}{b} \sum_{i  \in I} q'(r_i^2) \frac{\partial}{\partial \vartheta_i}$.
Hence 
\[ D_{l_{\delta,{{d}}}}(V)=  W + \sum_{j \in I} a_j \frac{\partial}{\partial \vartheta_j} + \beta_j \frac{\partial}{\partial r_j^2} \]\[ 
+ \frac{l_{\delta,{{d}}}}{b(p)}   \sum_{i \in I} \Big(2q''(r_i^2) \beta_i   - \Big( \sum_{j \in I} \beta_j h'(r_j^2)\Big)q'(r_i)  \Big)  \frac{\partial}{\partial \vartheta_i}.\]

If $V$ is not tangent to  $O^\delta_{{d}}$ then $dr_k(V) = \beta_k \neq 0$ for some $k$ where $\beta_k$ satisfies $|\beta_k| \geq |\beta_i|$ for all
$i$.
Conditions (\ref{item:deltaprop1}) and
(\ref{item:deltaprop2}) above then imply $\Big|\Big( \sum_{j \in I} \beta_j h'(r_j^2)\Big)q'(r_k^2)\Big|  < \Big|\beta_k q''(r_k^2)\Big|$.
This implies $D_{l_{\delta,{{d}}}}(V) \neq V$.
Hence $O^\delta_{{d}}$ is pseudo Morse-Bott.
Let $\partial O^\delta_d$ be the union of the boundary and corners of $O^\delta_d$.
We have also shown that $O^\delta_d \setminus \partial O^\delta_d$ is Morse-Bott.
We will use this fact in the next step.

\bigskip
{\it Step 3}:
We now need to calculate the Conley-Zehnder index of $O^\delta_{{d}}$.
Recall $b = i_{X_H} \theta_g = \sum_{i \in I} h_i(r_i^2)$ inside $(H_I \circ \pi_I)^{-1}(0)$.
Let $X$ be a vector in $TM|_{O^\delta_{{d}}}$ tangent to $O_{{d}}$ but transverse to $O^\delta_{{d}}$
and satisfying $dH(X)=-1$.
We have that $db(X) = \sum_{i \in I} h'_i(r_i^2) d(r_i^2)(X)$.
Now $(r_i)_{i \in I}$ restricted to $O^\delta_{{d}}$ is a constant $(c_i)_{i \in I}$.
Hence $b$ is a constant $b_0 = \sum_{i \in I} h_i(c_i^2) \in \R$ along $O^\delta_{{d}}$.

Recall that each orbit of $O^\delta_{{d}}$ is contained in some fiber $\pi_I^{-1}(q)$ and the restriction of $\theta$ to this fiber is $\sum_{i \in I} (\frac{1}{2}r_i^2 + \frac{1}{2\pi}\lambda_i)d\vartheta_i$.
So we have that the period $l^\delta_{{d}}$ of our family of Reeb orbits $O^\delta_{{d}}$ is
$\sum_{j \in I} d_j (\pi r_j^2  + \lambda_j)$.
Hence we have \[\frac{d l^s_{{d}}}{ds}\Big|_{s=\delta} = \sum_{j \in I} 2 d_j \pi dr_j^2(X). \]

%Choose $k$ so that $dr_k^2(W') \leq dr_j^2(W')$ for all $j$.
Let \[K := \frac{1}{b_0}l^\delta_{{d}}H = \frac{\sum_{j \in I} d_j (\pi c_j^2  + \lambda_j)}{\sum_{i \in I} h_i(c_i^2)}H.\]
%where $c$ is some fixed point in $O^\delta_d$.
We have rescaled $K$ here so that $O^\delta_d$ becomes an isolated family of time $1$ orbits.
We have by condition (\ref{item:deltaprop1}) that \[d(i_K\theta)(X) = \frac{\sum_{j \in I} d_j (\pi c_j^2  + \lambda_j)}{\sum_{i \in I} h_i(c_i^2)} \sum_{i \in I} h'_i(c_i^2) d(r_i^2)(X) \]\[ < 0 < \sum_{j \in I} 2 d_j \pi dr_j^2(X)=\frac{d l^s_{{d}}}{ds}\Big|_{s=\delta}.\]

Using this fact combined with the fact that
$O^\delta_d \setminus \partial O^\delta_d$ is Morse-Bott we have by Lemma \ref{lemma:hamiltonianconleyzehnderindexcomparison} that
$\text{CZ}(O^\delta_{{d}},K) = \text{CZ}(O^\delta_{{d}}) - \frac{1}{2}$.
Hence 
$\text{CZ}(O^\delta_{{d}}) = \text{CZ}(O^\delta_{{d}},K) + \frac{1}{2}$.
We will now calculate $\text{CZ}(O^\delta_{{d}},K)$. Because $K$ is a constant multiple of $H$, this is the same as
$\text{CZ}(O^\delta_{{d}},H)$ where we are looking at orbits of $X_H$ of period $\frac{l^\delta_d}{b_0}$.
Let $p \in O^\delta_{{d}}$.
The tangent space $TM$ at some point $p \in U_I$ splits as $A \oplus \oplus_{i \in I} B_i$
where $B_i$ consists of the tangent space to \[\dot{\D}^i_\epsilon = \dot{\D}^i_\epsilon \times \prod_{j \in I \setminus \{i\}} \{r_i(p)\} \subset  \prod_{j \in I}  \dot{\D}^j_\epsilon = \pi_I^{-1}(\pi_I(p))\]
and $A$ is the subspace symplectically orthogonal to the fiber of $\pi_I$ through $p$.
The Hamiltonian flow of $H$ preserves this splitting and is the identity on $A$
and is equal to the differential of the flow of $q(r_i^2)$ on $\D^i_\epsilon$ on $B_i$.
So we have $\text{CZ}(O^\delta_{{d}},H,d\theta_g ) = \sum_{i \in I} \Big(2(a_i + 1) d_i - \frac{1}{2}\Big)$.
Hence the index of our Reeb orbit is:
$2 \big(\sum_{i \in I} (a_i+1) d_i\big) + \frac{1-|I|}{2}$.
To connect these computations back to the statement of the Theorem we have
$B_V = O^\delta_{{d}}$ for  $V = \sum_{i \in I} d_i E_i$.
The size of $B_V$ is $2n - |I| - 1$ and its index is the index of $O^\delta_{{d}}$
which is $2 \big(\sum_{i \in I} (a_i+1) d_i\big) + \frac{1-|I|}{2}$.

\bigskip
{\it Step 4}:
Now we need to estimate the period of each pseudo Morse-Bott family of Reeb orbits.
In our case we need to calculate the period of $O^\delta_d$ for $d = (d_i)_{i \in I}$.
This pseudo Morse-Bott family corresponds to some element $V = \sum_{i \in I} d_i S_i$ in our monoid which is our family of Reeb orbits.
The period of $O^\delta_d$ is $\sum_{i \in I} d_i (\pi r_i^2 + \lambda_i)$.
Because $|\epsilon_q - \epsilon| < \frac{1}{2}\epsilon^3$, we have
for $0 < \delta \ll \epsilon$ small enough, that the period of
$O^\delta_d$ minus $\sum_{i \in I} d_i (\pi \epsilon^2 + \lambda_i)$ is less than
$\epsilon^3 \sum_{i \in I} d_i$.

\bigskip
{\it Step 5}:
We now need to construct our submanifolds $(2\epsilon,\epsilon_M)$-ball product associated to $B_{S_i}$ and $S_i$ inside $M_1$ for each $i$.
Let $p \in S_i$ be a point in $H_i^{-1}(0) \subset S_i$.
We let $\epsilon_M := q^{-1}(\delta)$.
For $\epsilon>0$ small enough there is some so that $B^{2n-2}_{3\epsilon}$ symplectically embeds in $H_i^{-1}(0)$
centered at $p$. Now $\pi_i^{-1}(B^{2n-2}_{3\epsilon})$ is a symplectic fibration with fibers symplectomorphic to $B^2_{\epsilon}$.
Inside this region, we can deform our fibration $\pi_i$ (along with our Hamiltonian $H$) near $\pi_i^{-1}(B^{2n-2}_{2\epsilon})$
so that $\pi_i^{-1}(B^{2n-2}_{2\epsilon}) \cap \{r_i \leq \epsilon_M\}$
becomes symplectomorphic
to $B^{2n-2}_{2\epsilon} \times B^2_{\epsilon_M}$ and so that $\theta_g$ becomes a product $\pi_i^* \beta + (\frac{1}{2}r_i^2 + \lambda_i)d\vartheta_i$ in this region.
We can also assume that $r_i$ and $d\vartheta_i$ become polar coordinates for $B^2_{\epsilon_M}$.
For $\delta$ small enough, this does not change the properties of $(H^{-1}(\delta),\theta_g|_{H^{-1}(\delta)})$.
Our region $W_i$ is now $B^{2n-2}_{2\epsilon} \times B^2_{\epsilon_M}$.

\bigskip
{\it Step 6}.
Here we show that $df(X_{\theta_g}) > 0$ in $M_1$.
We have inside $(H_I \circ \pi_I)^{-1}(0)$ that $X_{\theta_g} = X_1 + X_2$ where $X_1$ is tangent to the fibers
$\prod_{i \in I} \dot{\D_\epsilon}$ of $\pi_I$ and $X_2$ is symplectically orthogonal to these fibers.
Because the restriction of $\theta_g$ to the fiber is $\sum_{i \in I} (\frac{1}{2}r_i^2 + \frac{1}{2\pi}\lambda_i) d\vartheta_i$
we get that $X_1 = 2\sum_{i \in I}(\frac{1}{2}r_i^2 + \frac{1}{2\pi}\lambda_i) \frac{\partial}{\partial r_i^2}$.
Because $f = \sum_{i \in I} \nu_{r_i^2}$, and $df(\frac{\partial}{\partial r_j^2}) > 0$
in the region $\cap_{i \in I} \{r_j^2 < \epsilon_q\} \cap (H_I \circ \pi_I)^{-1}(0)$ for each $I \subset \{1,\cdots,l\}$, we have that $df(X_{\theta_g}) > 0$ in $M_1$.
This completes the proof of the theorem.
\qed

%TODO: continue here
\subsection{Bounding the Minimal Discrepancy of a Singularity from Above}

In this subsection we will use the results from the last two subsections to give an upper bound for the  minimal discrepancy of an isolated singularity in terms of Conley-Zehnder indices of Reeb orbits of its link.
The main theorem of this section is the following: 

\begin{theorem} \label{theorem:boundingdiscrepancyfromabove}
Let $A \subset \C^N$ be an affine variety which has an isolated singularity at zero.
Suppose that its link $(L_A,\xi_A)$ satisfies $H^1(L_A;\Q) = 0$ and $c_1(\xi_A;\Q) = 0 \in H^2(L_A;\Q)$.
Then $\text{hmi}(L_A,\xi_A) \geq 2\text{md}(A,0)$.
\end{theorem}

We will need some preliminary lemmas before we prove the above theorem. 
Lemma \ref{lemma:symplecticformonresolution} will also be useful in Section \ref{section:boundingminimaldiscrepancyfrombelow} where we bound minimal discrepancy from below.

\begin{lemma} \label{lemma:radialcoordinate}
Let $X$ be complex manifold, and let $L$ be a line bundle on $X$ with Hermitian metri $\|\cdot\|$.
Let $E$ be a smooth complex hypersurface.
Suppose that $L$ has a meromorphic section $s$ with a simple pole along $E$
and no other poles. We will also assume that $s^{-1}(0) = \emptyset$.
Define $r \in C^0(X)$ so that $r = \frac{1}{\|s\|}$ outside $E$
and $r = 0$ on $E$.
Then $r$ is a radial coordinate for $E$ as in Definition \ref{defn:radialcoordinate}.
% for some choice of
%Hermitian metric $\|\cdot\|$.
\end{lemma}
\proof
Fix some metric on $X$ and let
$\pi_E : TX^\perp \to E$
be the bundle of vectors orthogonal to $TX$ with respect to our metric.
Let $TX^{\perp,<\delta} \subset TX^\perp$ be the subset of vectors of length $< \delta$.
We have the exponential map $\text{exp} : TX^{\perp,<\delta} \to X$
which is a diffeomorphism onto its image for $\delta$ small. Let $U$ be the image of this map.
Let $\pi_U : U \to E$ be defined by
$\pi_U(x) = \pi_E \circ \text{\text{exp}}^{-1}(x)$.
Let $U_1 := r^{-1}([0,\delta_1))$
so that $U_1 \subset U$.
The fibers of $\pi_U|_{U_1}$ are disk fibers for $\delta_1 >0$ small because
in each local holomorphic cart $z_1,\cdots,z_n$
centered at each point of $E$
with $E = \{z_1 = 0\}$,
we have $r = e^{-\mu} |z_1|$ for a smooth function $\mu$.
This fiber bundle has a $U(1)$ structure group rotating these disk fibers
where the natural radial coordinate on each fiber is $r$.
Choose an Ehresmann connection $L \subset TU_1$ on $\pi_U|_{U_1}$
compatible with this $U(1)$ structure group (i.e. so that the parallel
transport maps rotate the disk fibers).
Let $q$ be a $1$-form on $U_1$ which vanishes on $L$
and so that $q$ restricted to each fiber is $d\vartheta$ where $\vartheta$
is the angle coordinate on some trivialization of this fiber.
Choose a metric on $U_1$ so that
it makes the $1$-forms $dr$ and $rq$ orthogonal and of length $1$
and so that $L$ is orthogonal to the fibers of $\pi_U|_{U_1}$.
We also make sure that the metric restricted to $L$ is the pullback of the induced metric on $TE$ via $D\pi_U|_L$.
This metric on $U_1$ makes the fibers of $\pi$ totally geodesic by
%\cite[Theorem 9.59]{besseeinstein} or
\cite[Theorem 3.5]{Vilmsgeodesics}.
Hence an exponential map calculation inside each of the fibers ensures that
$r$ is a radial coordinate as in Definition
\ref{defn:radialcoordinate}.
\qed

%We can find an $S^1$ action on a neighborhood of $E$
%tangent to

%Fix some metric on $X$. We have a bundle isomorphism $\iota : L^*|_E \to TX^\perp$
%where $TX^\perp$ is the bundle of vectors orthogonal to $TX$ with respect to our metric.
%Let $TX^{\perp,<\delta} \subset TX^\perp$ be the subset of vectors of length $< \delta$.
%We have the exponential map $\text{exp} : TX^{\perp,<\delta} \to X$
%which is a diffeomorphism onto its image, denoted by $\text{Im}(\text{exp})$.
%Let $\widetilde{s}$, $\|\cdot\|^v$ be the pullback of $s^*$
%and $\|\cdot\|^*$ respectively to $(\iota \circ \exp)^* L^*$.
%Choose a metric on the total sp
%\qed

\bigskip

Let $A$ be an affine variety with an isolated singularity at $0$ and
let $A_\delta$ be its intersection with a small ball of radius $\delta$.
We resolve $A$ at $0$ by blowing up along smooth loci by \cite{hironaka:resolution} and take the preimage $\widetilde{A}_\delta$ of $A_\delta$ under this resolution map.
We suppose that the exceptional divisors $E_1,\cdots,E_l$ are smooth normal crossing
and that $E_I := \cap_{i \in I} S_i$ is connected for each $I \subset \{1,\cdots,l\}$.
Note that the preimage of $0$ is connected by  \cite{Zariski:connectedness}.

\begin{lemma} \label{lemma:symplecticformonresolution}
There exists a $1$-form $\theta_A$ on $\widetilde{A}_\delta \setminus \cup_i E_i$
making $(\widetilde{A}_\delta,\cup_i E_i,\theta_A)$ into a positively wrapped divisor
whose associated symplectic form $\omega_A$ is a K\"{a}hler form
and whose contact link is contactomorphic to the link $L_A$ of $A$.
If $A$ is numerically $\Q$-Gorenstein and $H^1(L_A;\Q) = 0$
then $(A_\delta,\cup_i E_i,\theta_A)$ is strongly numerically $\Q$-Gorenstein
and the discrepancy of $E_i$ and minimal discrepancy 
defined in Section \ref{section:minimaldiscrepancydefinition} coincides with
the associated discrepancy and minimal discrepancy from Definition
\ref{defn:symplecticminimaldiscrepancy}.
\end{lemma}
\proof of Lemma \ref{lemma:symplecticformonresolution}.
Because our resolution is obtained by blowing up along smooth loci starting from $A \subset \C^N$, we can also blow up $\C^N$ along the same loci giving us $\widetilde{\C^N}$.
Let $E^\C_1,\cdots,E^\C_m \subset \widetilde{\C^N}$ be the corresponding exceptional divisors.
We can reorder them so that the proper transform of the exceptional divisor of the $i$th blowup of $\C^N$ is $E^\C_i$.
Hence for positive integers $\widetilde{\nu}_1 \gg \cdots \gg \widetilde{\nu}_l >0$
we have $-\sum_{i \in I} \widetilde{\nu}_i E^\C_i$
is ample in some neighborhood of $\cup_i E^\C_i$.
The restriction of $-\sum_{i=1}^m  \widetilde{\nu}_i E^\C_i$ to $\widetilde{A}_\delta$
is then an ample divisor $-\sum_{i=1}^l \nu_i E_i$ for some integers $\nu_1,\cdots,\nu_l > 0$.
%So for $\delta > 0$ small enough we have $-\sum_{i \in I} \nu_i E_i$ is ample in  $\widetilde{A}_\delta$.
This means that if $L \to \widetilde{A}_\delta$ is the line bundle associated to $-\sum_{i \in I} \nu_i E_i$
then it has a Hermitian metric $\|\cdot\|$ with positive curvature $F$ and a meromorphic section $s$ with divisor $-\sum_{i = 1}^l \nu_i E_i$.
Our symplectic form is $\omega_A := -dd^c \log(\|s\|)$ and $\theta_A = -d^c \log(\|s\|)$.
The wrapping numbers of $E_1,\cdots,E_l$ with respect to $\theta_A$ are the positive numbers $2\pi \nu_1,\cdots,2\pi \nu_l$ respectively.
Hence $(\widetilde{A}_\delta,\cup_i S_i,\theta_A)$ is a positively wrapped divisor.

We will now show that $f_A:= -\log(\|s\|)$ is compatible with $\cup_i E_i$.
We have $L = \prod_i L_i^{\otimes \nu_i}$ where $L_i$ is the line bundle coming from $-E_i$
and $s = u . \otimes_i s_i^{\otimes \nu_i}$ for meromorphic sections $s_i$ of $L_i$  with divisor $-E_i$
and where $u$ is a non-zero holomorphic function.
The metric $\|\cdot\|$ is equal to $e^{-\widetilde{\mu}} \prod_i \|\cdot\|_i^{\nu_i}$
where $\|\cdot\|_i$ is a metric on $L_i$ for some smooth function $\widetilde{\mu} \in C^\infty(\widetilde{A}_\delta)$.
Let $\rho$ be the function described in Section \ref{section:boundaryconstructionanduniqueness}
(in other words, for some small $\delta_r>0$, $\rho : [0,\delta_r) \to [0,1]$ is a smooth function so that $\rho(x) = x^2$ near $0$
and $\rho(x) = 1$ near $\delta_r$ with $\rho' \geq 0$).
We define $\rho_{1/\|s_i\|}$ to be $\rho(1/\|s_i\|_i)$ inside the domain of $r_i$
and $1$ elsewhere.
Then $f_A = \widehat{\mu} +
\sum_i \frac{1}{2}\nu_i (\log(\rho_{1/\|s_i\|_i}))$
where $\widehat{\mu} = \widetilde{\mu}
+ \sum_i \nu_i (\frac{1}{2}\log(1/\|s_i\|)-\log(\rho_{1/\|s_i\|}))$
is bounded.
We have that $\frac{1}{\|s_i\|_i}$ is a radial coordinate by Lemma \ref{lemma:radialcoordinate} and so $f_A$ is compatible with $\cup_i E_i$.

We will now show that the link of $A$ is contactomorphic to the contact link of $\cup_i E_i$
using our function $f_A$.
The proof of this fact is very similar to the proofs of Lemmas 4.3 and 4.4 in
\cite{Seidel:biasedview}.
Let $z_1,\cdots,z_N$ be the natural coordinates on $\C^N$. Define $\phi : \widetilde{A}_\delta \setminus \cup_i E_i \to \R$ to be the pullback of $\sum_{i \in I} |z_i|^2$.
For all small enough $\eta>0$ we have that $\left(\phi^{-1}(\eta),\text{ker}(-d^c\phi(\eta)|_{\phi^{-1}(\eta)})\right)$ is contactomorphic to $(L_A,\xi_A)$
by definition.

For $t \in [0,1]$ define $\phi_t := (1-t) \phi - t\log(\|s\|)$.
For each $p \in \cup_i E_i$ choose local holomorphic  coordinates $w_1,\cdots,w_n$ so that $\cup_i E_i = \{\prod_{i = 1}^k w_i = 0\}$.
Let $V$ be the radial vector field emanating from $0$ with respect to this coordinate system
(in other words the vector field $\sum_j x_j \frac{\partial}{\partial x_j} + y_j \frac{\partial}{\partial y_j}$
where $w_j = x_j + i y_j$).
Now $\phi = |b|^2 \prod_{i = 1}^k |w_i|^{2q_i}$ for some integers $q_1,\cdots,q_k > 0$ and some non-zero holomorphic function $b$ defined near $p$
and so $d(\log(\phi))(V) = d(\log(|b|^2))(V) + \sum_{i=1}^k q_i d(\log(|w_i|^2)(V)$.
The term $d(\log(|b|^2))(V)$ is bounded near $p$ and the term
$\sum_{i=1}^k q_i d(\log(|w_i|^2)(V)$ tends to $\infty$
as we approach $p$ inside $\widetilde{A}_\delta \setminus \cup_i E_i$.
Hence $d\phi(V) > 0$ inside $\widetilde{A}_\delta \setminus \cup_i E_i$ if we are near $p$.
Also with respect to some trivialization of $L$ near $p$ we have $\|\cdot\| = e^\mu |\cdot|$ for some smooth function $\mu$ defined near $p$. Hence $-\log(\|s\|) = \mu + \sum_{i=1}^k e_i\log(|w_i|)$ for positive integers $e_1,\cdots,e_k$.
Again this implies that $d(-\log(\|s\|))(V) > 0$
near $p$ and inside $\widetilde{A}_\delta \setminus \cup_i E_i$.
Hence $d\phi_t(V) > 0$ near $p$ and inside $\widetilde{A}_\delta \setminus \cup_i E_i$.
Because $\cup_i E_i$ is compact we then get $d\phi_t \neq 0$ near $\cup_i E_i$
inside $A_\delta \setminus \cup_i E_i$.
Hence there is a smooth function $\lambda : [0,1] \to \R$ satisfying:
\begin{enumerate}
\item For any $l \leq L := \lambda(1)$ we have that $l$ is a regular value of $\phi_1 = \log(\|s\|)$.
\item We have $\lambda(0) > 0$ and
for any $l \in (0,\lambda(0))$, $l$ is a regular value of $\phi_0 = \phi$.
\item $\lambda(t)$ is a regular value of $\phi_t$ for all $t \in [0,1]$.
\end{enumerate}
Hence by Gray's stability theorem we have that $\left(\phi_0^{-1}(\lambda(0)),\text{ker}(-d^c\phi_0|_{\phi_0^{-1}(\lambda(0))})\right)$
is contactomorphic to $\left(\phi_1^{-1}(\lambda(1)),\text{ker}(-d^c\phi_1|_{\phi_1^{-1}(\lambda(1)})\right)$.
Hence  $\left(\phi_1^{-1}(\lambda(1)),\text{ker}(-d^c\phi_1|_{\phi_1^{-1}(\lambda(1))})\right)$
is contactomorphic to $(L_A,\xi_A)$
and
$\left(f_A^{-1}(l),\text{ker}(\theta_A|_{f_A^{-1}(l)})\right)$ for $l \leq L$.
Hence $\left(f_A^{-1}(l),\text{ker}(\theta_A|_{f_A^{-1}(l)})\right)$
is contactomorphic to $(L_A,\xi_A)$ for all $l \leq L$.

The statement in this theorem about discrepancy and minimal discrepancy
comes from the fact that $\omega_A$ is K\"{a}hler
which means that the canonical bundle of $\omega_A$
is identical to the canonical bundle of $A_\delta$
and also from Lemma \ref{lemma:topologicalnumericallygorenstein}.
\qed

\proof of Theorem \ref{theorem:boundingdiscrepancyfromabove}.
%
%
%Let $\omega_A$, $\widetilde{A}_\delta$
%,$E_1,\cdots,E_l$, $\theta_A$ , $f_A : \widetilde{A}_\delta \setminus \cup_i E_i \to %\R$, $L\geq 0$ and $(L_A,\alpha_A)$
%be from the proof of Lemma \ref{lemma:symplecticformonresolution}.
By Lemma \ref{lemma:symplecticformonresolution},
there exists a $1$-form $\theta_A$ on $\widetilde{A}_\delta \setminus \cup_i E_i$
making $(\widetilde{A}_\delta,\cup_i E_i,\theta_A)$ into a
strongly numerically $\Q$-Gorenstein
positively wrapped divisor
%whose associated symplectic form $\omega_A$ is a K\"{a}hler form
%and
whose contact link is contactomorphic to the link $(L_A,\xi_A)$ of $A$.
Also $\text{md}(\widetilde{A}_\delta,\cup_i E_i,\theta_A) = \text{md}(A,0)$. 
By \cite[Theorem 5.3]{McLean:affinegrowth} we can deform $S_1,\cdots,S_l$ through positively intersecting submanifolds so that they become orthogonal.
This induces a deformation of 
strongly numerically $\Q$-Gorenstein
positively wrapped divisors, which all have the same
minimal discrepancy and which have contactomorphic links by
Corollary \ref{corollary:canonicalboundary}.
So from now on we can assume that $S_1,\cdots,S_l$ are symplectically orthogonal.
By Corollary \ref{corollary:specificformaintheorem},
we have a contact form associated to $\xi_A$
and a family of Reeb orbits $B$ satisfying:
$\text{lSFT}(B) = \text{md}(\widetilde{A}_\delta,\cup_i E_i,\theta_A)$.
Hence:
$\text{hmi}(L_A,\xi_A) \geq \text{md}(A,0)$.
 
%Let $a_1,\cdots,a_l \in \Q$ be the discrepancies of $E_1,\cdots,E_l$ respectively.
%Now $f_A$ is compatible with $\cup_i E_i$
%and so by Corollary \ref{corollary:canonicalboundary} and Theorem  %\ref{label:nicecontactneighbourhoodexistence}
%we have that $(f_A^{-1}(l),\theta_A|_{f_A^{-1}(l)})$ admits a contact form $\beta$
%with $\text{ker}(\beta) = \text{ker}(\theta_A|_{f_A^{-1}(l)})$ and
%so that every Reeb orbit sits inside a pseudo Morse-Bott family of size $2n - |I_V| -1$
%and Conley-Zehnder index $2 \sum_{i \in I_V} (a_i+1) d_i + \frac{1-|I_V|}{2}$ where %$I_V \subset \{1,\cdots,l\}$ is a subset of cardinality at most $n=\text{dim}_{\C}%(A)$
%and $d_i$ are positive integers labeled by $i \in I_V$.
%Then $\text{mi}(\beta)$ is equal to the infimum over all $I_V$
%and $d_i$ of 
%\[2 \sum_{i \in I_V} (a_i+1) d_i + \frac{1-|I_V|}{2}-\frac{1}{2}(2n - |I_V| -1) + %n-3\]
%\[ = 2\sum_{i \in I_V} (a_i + 1)d_i  -2.\]
%Now if $a_ i <-1$ then both
%$\text{md}(A,0)$ and $\text{mi}(\beta)$
%are $-\infty$ and so
% $2~\text{md}(A,0) \leq \text{hmi}(L_A,\xi_A)$ in this case.
%Otherwise the infimum over all $I_V$
%and $d_i$ of 
%$2\sum_{i \in I_V} (a_i + 1)d_i  -2$ is $2~\text{inf}(a_i)$.
%So \[2~\text{md}(A,0) =2~\text{inf}(a_i) = \text{mi}(\beta) \leq \text{hmi}(L_A,%\xi_A).\]
\qed

\section{Gromov-Witten Invariants on Open Symplectic Manifolds} \label{section:gromovwittenopen}
%TODO define Gromov-Witten invariants and state properties as in log  Kodaira dimension paper

%We will use the genus $0$  Gromov-Witten invariants defined for general symplectic manifolds.
In Section \ref{subsection:sketchofproof}
we gave a short survey of Gromov-Witten invariants for closed symplectic manifolds.
We wish to do this for certain open symplectic manifolds. The problem is that we cannot have
such general properties such as
properties (\ref{item:holomorphiccurveexistence})
and (\ref{item:gwdeformationinvariance}) of Theorem \ref{theorem:gw}.
Instead we have to consider subfamilies ${\mathcal J}$ of compatible almost complex structures
satisfying certain special properties and then define Gromov-Witten invariants for these families only.

%The symplectic manifolds in our case are open but all the holomorphic curves stay inside a fixed compact subset %ensuring that Gromov-Witten invariants are still well defined.
%not a problem (with the possible exception of \cite{CieliebakMohnke:symplectichypersurfaces}
%which is reliant on Donaldson hypersurface techniques).

\begin{defn} \label{defn:gwtriple}
The triple $(S,[A],\mathcal{J})$ is {\bf a GW triple} if
 $(S,\omega_S)$ is a (possibly non-compact) symplectic manifold,
$\mathcal{J}$ is a family of compatible almost complex structures in $S$ and $[A] \in H_2(S;\Z)$ so that:
\begin{enumerate}
 \item $\mathcal{J}$ is non-empty and path connected.
\item \label{item:compact}
There is a relatively compact open subset $U_S$ of $S$
so that for every $J \in \mathcal{J}$, every compact genus $0$ nodal $J$-holomorphic curve representing $[A]$ is contained in $U_S$.
\item \label{item:dimensioncondition}
$c_1(S,\omega_S)([A]) + (n-3) = 0$.
\end{enumerate}
\end{defn}

The definition of a GW triple is new to this article, although many Gromov-Witten
invariants have in the past been calculated for open symplectic manifolds
(see for instance \cite{Ritternegativelinebundles}).
Condition (\ref{item:dimensioncondition}) tells us that the space of $J$-holomorphic curves
representing $[A]$ has `dimension' $0$ for generic $J$ (also known as virtual dimension,
see \cite[Theorem 1.3]{FukayaOno:Arnold} for a dimension formula for instance).
Combining this with condition (\ref{item:compact})
which tells us that the space of such curves stay inside a compact set
means that one can `count' the number of these curves representing $[A]$.
This count will be denoted by $\text{GW}_0(S,[A],J) \in \Q$.
In some cases, the dimension $0$ condition means that these curves form a discrete family,
and the compactness condition tells us that there are only finitely many such curves.
Sometimes $\text{GW}_0(S,[A],J) \in \Q$ is the number of these curves.
In general though, this counting process is quite complicated and is not just the number of
$J$-holomorphic curves representing $[A]$.

\begin{theorem} \label{theorem:gromovwitteninvariants}
(\cite[Theorem 1.3]{FukayaOno:Arnold},
\cite[Theorem 1.12, and the following paragraph]{HoferWysockiZehnder:polyfoldapplications1}
or
\cite[Theorem 2.5]{LiTian:sympGW}).
We can assign an invariant
$\text{GW}_0(S,[A],J) \in \Q$ satisfying the following properties:
\begin{enumerate}
\item If $\text{GW}_0(S,[A],J)\neq 0$ then
there exists a compact nodal $J$-holomorphic curve representing $[A]$.
\item Given a smooth path $(J_t)_{t \in [0,1]}$ in $\mathcal{J}$
we have $\text{GW}_0(S,[A],J_0) = \text{GW}_0(S,[A],J_1)$.
\item
\cite[Theorem 3.3.1 and Theorem 7.1.8]{McDuffSalamon:Jholomorphiccurves}.

Suppose that every connected genus $0$ $J$-holomorphic curve $u : \Sigma \to S$ representing $[A]$ satisfies:
\begin{itemize}
\item $J$ is integrable near $u(\Sigma)$,
\item $u^*(TS)$ is a direct sum of complex line bundles of degree $\geq -1$,
\item \label{item:actualcounting}
$u$ is smooth,
\end{itemize}
then $\text{GW}_0(S,[A],J)$ is equal to the number of connected genus $0$ $J$-holomorphic curves representing $[A]$.
\end{enumerate}
\end{theorem}

The `count' $\text{GW}_0(S,[A],J) \in \Q$ will be called the {\bf genus $0$ Gromov-Witten invariant in the class [A]}.
If $J_1 \in \mathcal{J}$ then there is a smooth path of almost complex structures in $\mathcal{J}$ joining $J$
and $J_1$.
This implies that $\text{GW}_0(S,[A],J) = \text{GW}_0(S,[A],J_1)$.
So from now on we will define $\text{GW}_0(S,[A],\mathcal{J}) := \text{GW}_0(S,[A],J)$ for some $J \in \mathcal{J}$.

Sometimes we need to perturb $J$ by a small amount. For instance when defining the invariant $\text{GW}_0(S,[A],J)$
one quite often needs to perturb $J$ to some $C^\infty$ generic $J$ before counting the holomorphic curves.
In this case we do the following.
By a Gromov compactness argument
(see \cite{fish:compactness})
we have that for any $J_1$ sufficiently $C^\infty$ close to $J \in \mathcal{J}$
we either have that every $J_1$-holomorphic curve is contained inside some fixed
relatively compact subset of $U_S$
 or some part of the curve maps outside $U_S$.
This means that $(U_S,[A],J_1)$ is a GW triple and so we define  $\text{GW}_0(S,[A],J_1) :=  \text{GW}_0(U_S,[A],J_1)$.
Note that this count is independent of the choice of open set $U_S$ because if we had some other
relatively compact set $V_S$
containing all $J$-holomorphic curves then we can make
$J_1$ sufficiently close to $J$ to ensure that all $J_1$-holomorphic curves are inside $U_S \cap V_S$.
We also have the following similar lemma.
\begin{lemma} \label{lemma:openness}
Let $(S,[A],\mathcal{J})$ be a GW triple.
Suppose that $S_1$ is a large relatively compact open subset of $S$ containing the closure of
our relatively compact open subset $U_S$ with the additional
property that there is a unique homology class $[A_1]$ mapping to $[A]$ under the inclusion map.
Then there is an open subset $\mathcal{J}^0$ of the space of compatible almost complex structures with respect
to the $C^\infty$ topology so that $\mathcal{J}^0$ contains $\mathcal{J}$ and so that:
$(S_1,[A_1],\mathcal{J}^0)$ is a GW triple with $\text{GW}_0(S,[A],\mathcal{J}) = \text{GW}_0(S_1,[A_1],\mathcal{J}^0)$.
We can ensure that any genus $0$ nodal $J_1$-holomorphic curve for $J_1 \in \mathcal{J}^0$ whose image is contained in $S_1$
has its image contained in $U_S$.
\end{lemma}
\proof of lemma \ref{lemma:openness}.
Suppose that the Theorem above is false. This means that there is a $J \in \mathcal{J}$
and a sequence of compatible almost complex structures $J_i$ $C^\infty$ converging to $J$
and a sequence of $J_i$-holomorphic curves $u_i : \Sigma_i \to S_1$ representing $[A_1]$ and not contained in $U_S$.
By a Gromov compactness argument, a subsequence then converges to a $J$-holomorphic curve
mapping to the closure of $S_1$. Because $(S,[A],\mathcal{J})$ is a GW triple we then get that such a limit curve is contained in $U_S$.
But this means that $u_i$ maps to $U_S$ for $i$ large enough which is a contradiction.
Hence $(S_1,[A_1],\mathcal{J}^0)$ is a GW triple.
The reason why $\text{GW}_0(S,[A],\mathcal{J}) = \text{GW}_0(S_1,[A_1],\mathcal{J}^0)$
is because we can count our curves with respect to some almost complex structure
$J \in \mathcal{J} \subset \mathcal{J}^0$.
\qed

\bigskip

\begin{defn} \label{defn:gwtripledeformation}
Let $(S,[A],\mathcal{J}_0)$ and $(S,[A],\mathcal{J}_1)$ be GW triples with associated symplectic forms
$\omega_{S,0}$ and $\omega_{S,1}$.
A {\bf smooth deformation of GW triples joining} $(S,[A],\mathcal{J}_0)$ and $(S,[A],\mathcal{J}_1)$
consists of a family of GW triples $((S,[A],\mathcal{J}_t))_{t \in [0,1]}$ such that:
\begin{enumerate}
\item For each $t \in [0,1]$ the associated symplectic form $\omega_{S,t}$ for $(S,[A],\mathcal{J}_t)$ smoothly varies with $t$.
\item There is a smooth family of almost complex structures $J_t,t\in [0,1]$ such that $J_t \in \mathcal{J}_t$ for all $t \in [0,1]$.
\item \label{item:deformationmaximumprinciple}
There is a relatively compact open subset $U_S$ with the property that any $J_t$-holomorphic curve representing $[A]$ is contained in
$U_S$.
\end{enumerate}
We will say that $(S,[A],\mathcal{J}_0)$ is {\bf deformation equivalent} to $(S,[A],\mathcal{J}_1)$
if there exists a smooth deformation of GW triples joining
$(S,[A],\mathcal{J}_0)$ and $(S,[A],\mathcal{J}_1)$.
\end{defn}

\begin{lemma} \label{lemma:deformationinvariance}
Let $(S,[A],\mathcal{J}_0)$ and $(S,[A],\mathcal{J}_1)$ be GW triples which are deformation equivalent to each other.
Then $GW_0(S,[A],\mathcal{J}_0) = GW_0(S,[A],\mathcal{J}_1)$.
\end{lemma}
The proof of this lemma is basically the same as the proof that Gromov-Witten invariants
do not change when deforming the symplectic and almost complex structure.
Our deformation is $(\omega_{S,t},J_t)$.
The only difference in our argument is that we are in an open symplectic manifold but this is OK as all the $J_t$
holomorphic curves stay inside a fixed relatively compact open subset.
Again sometimes  we would like to perturb $J_t$ slightly in which case we fix a relatively compact open subset containing
$U_S$ and count curves inside this subset using the same ideas from Lemma \ref{lemma:openness}.
%EXPLAIN THIS MORE - AXIOMATICALLY?

Lemmas \ref{lemma:gwtripletest} and \ref{lemma:reeborbits} below are important tools which will be used later to give a lower bound on minimal discrepancy.
Lemma \ref{lemma:gwtripletest} gives us some sufficient conditions for a symplectic manifold,
a homology class and a family of almost complex structures to be a GW triple
and Lemma \ref{lemma:reeborbits} helps us find Reeb orbits whose Conley-Zehnder index has an explicit bound.
Both Lemmas involve neck stretching, but for two different purposes.
Lemma \ref{lemma:gwtripletest} uses such an argument to show that
property (\ref{item:deformationmaximumprinciple}) of Definition \ref{defn:gwtripledeformation}
 is satisfied
(i.e. all holomorphic curves representing our chosen homology class stay inside a compact set).
On the other hand Lemma \ref{lemma:reeborbits} uses a neck stretching argument to find Reeb orbits with an appropriate upper bound
on their Conley-Zehnder index.
In particular the stable Hamiltonian hypersurface from Lemma \ref{lemma:gwtripletest} is different and has a different purpose
to the contact hypersurface in Lemma \ref{lemma:reeborbits}.
These Lemmas are used in the proof of Theorem \ref{theorem:reeborbitlowerboundarounddivisors}.
We refer the reader to Appendix A (Section \ref{section:appendixA}) for definitions concerning neck stretching along stable Hamiltonian
structures.
%For the non-expert, one can initially pretend that all the stable Hamiltonian hypersurfaces are contact hypersurfaces %because all the properties we use in this section are also properties contact hypersurfaces (also one can pretend %stable Hamiltonian cobordisms are symplectic manifolds whose boundary is a contact hypersurface, where such a %boundary is negative (resp. positive) if the contact form induces an inward (resp. outward) orientation).

\begin{lemma} \label{lemma:gwtripletest}
Suppose that $(S,\omega_S)$ is a symplectic manifold,
$C \subset S$ a compact stable Hamiltonian hypersurface,
$J_i$ a sequence of almost complex structures on $S$ compatible with the symplectic form $\omega_S$ and
$[A] \in H_2(S;\Z)$ satisfying:
\begin{enumerate}
\item 
We have $S \setminus C$ is a disjoint union $\dot{S}_+$, $\dot{S}_-$
where the closure of $\dot{S}_\pm$ inside $S$, which is equal to $\dot{S}_\pm \cup C \subset S$,
is a stable Hamiltonian cobordism $S_\pm$.
Here $C$ has a standard neighborhood $(-\epsilon_h,\epsilon_h) \times C$
% symplectomorphic to its symplectization
and we require that $[0,\epsilon_h) \times C \subset S_+$ and $(-\epsilon_h,0] \times C \subset S_-$
which means that $S_+$ has a negative boundary and no positive boundary and $S_-$ has a positive boundary but no negative boundary.
We will assume that $S_-$ is compact but $S_+$ may be non-compact.
\item $J_i$ is a sequence of almost complex structures stretching the neck along $C$ with respect to $(-\epsilon_h,\epsilon_h) \times C$
so that $J_i|_{\dot{S}_+}$ converges in $C^\infty_{\text{loc}}$ inside $\dot{S}_+$ to  an almost complex structure $J_+$ compatible with the completion of $S_+$.
% after neck stretching along $C$.
\item \label{item:stayinginsidecompactset}
There exists a properly embedded
%$J_+$-holomorphic
real codimension $2$ submanifold
$E$ in $\dot{S}_+$ satisfying $[A].E = 1$ and whose closure in $S$ does not intersect $C$.
Also there is a relatively compact open subset $U_S \subset S$ so that
any finite energy proper $J_+$-holomorphic curve in $\dot{S}_+$ intersecting $E$ with multiplicity $1$
is contained inside $U_S \cap \dot{S}_+$.
\item $c_1(S,\omega_S)([A]) + n-3 = 0$.
% and $S$ admits a Morse function with finitely many critical points.
\end{enumerate}

Let $\omega_t$ in $S$ ($t \in [0,1]$) be a smooth family of symplectic forms agreeing with
$\omega_S$ inside $(-\epsilon_h,\epsilon_h) \times C \cup \dot{S}_+$ and define
$\mathcal{J}_{i,\omega_t}$ to be the set of compatible almost complex structures $J_{i,t}$
agreeing with
$J_i$ inside $(-\epsilon_h,\epsilon_h) \times C \cup \dot{S}_+$.
Then if $S_1$ is some relatively compact open subset containing $U_S \cup S_- \cup ([0,\epsilon_h) \times C)$ and homotopic to $S$ then
there exists an $i_0$ so that $(S_1,[A],\mathcal{J}_{i,\omega_t})$ is a smooth deformation of GW triples for all $i \geq i_0$.
\end{lemma}

%See the appendix for a definition of an almost complex structure being compatible with the completion.
%By Lemma \ref{lemma:openness} we can assume that
%$J_-$ is generic among almost complex structures equal to $J_-$ near $C$.
%The point is that we an enlarge $K_S$ very slightly and get the same result and then by a compactness
%result ensure that all of the $J_-$-holomorphic curves of energy at most $[A].\omega_S + 1$
%sit inside $K_S$.

Roughly, the above Lemma says that if we have a symplectic manifold $S$ with $[A] \in H_2(S;\Z)$
and families of complex structures ${\mathcal J}_i$
\begin{itemize}
\item which `stretch' $S$ near infinity so that $S$ breaks into two pieces $S_+,S_-$
\item so that these almost complex structures `converge' to $J_+$ in $S_+$,
\item and so that all $J_+$-holomorphic curves stay inside a compact subset of $S$
\end{itemize}
then $(S,[A],\mathcal{J}_i)$ is a GW triple for $i$ large enough.
The key point is that we can pretend that $(S,[A],J_+)$ is sort of like a GW triple
(even though $J_+$ is not well defined on $S_-$) and we have a sequence of triples
`converging' to this triple and hence for $i$ large enough they are in fact GW triples.

To provide some context for this lemma, we will apply it in the proof of Theorem \ref{theorem:reeborbitlowerboundarounddivisors} later on
in the following way:
Let $\widetilde{A}$ be a resolution of our singularity $A$ and let $\widetilde{A}_\delta$ be the preimage of a small
ball under the resolution map. By Lemma \ref{lemma:symplecticformonresolution},
$\widetilde{A}_\delta$ along with the exceptional divisors has the structure of
a strongly numerically $\Q$-Gorenstein positively wrapped divisor.
One can show by
Corollary \ref{theorem:reeborbitlowerboundarounddivisors} that there is some link
region $M_1$ whose boundary is contactomorphic to $(L_A,\xi_A)$
and an $(2\epsilon,\epsilon_M)$-ball product associated to some family of Reeb orbits $B$
and some divisor $D$ of lowest discrepancy.
Our manifold $S$ will be obtained from $M_1$ by first partially compactifying $M_1$
by replacing our  $(2\epsilon,\epsilon_M)$-ball product with
a product $B^{2n-2}_{2\epsilon} \times S^2_{\epsilon}$ of a ball and a sphere. Then $S$ is obtained from
this product by blowing up this partial compactification at a single point inside
$\{0\} \times S^2_{\epsilon}$.
The stable Hamiltonian hypersurface in $S$ will be some modification of $\partial M_1$
inside our product region. For technical reasons we need this to be stable Hamiltonian
(and not just a contact hypersurface) so that part (\ref{item:stayinginsidecompactset}) of Lemma \ref{lemma:gwtripletest} above holds
(see Step 3 from the proof of Theorem \ref{theorem:reeborbitlowerboundarounddivisors}).
This is also the reason why we chose the divisor $D$ of smallest discrepancy.
The homology class $[A]$ is represented by the proper transform of $\{0\} \times S^2$.
The submanifold $E$ is of the form $B^{2n-2}_{2\epsilon} \times \{x\}$ for some appropriate $x$
and $U_S$ is some relatively compact neighborhood of $M_1$ in $S$.
The reason why we need a well defined GW triple is so that we can apply a neck stretching argument
along contact hypersurfaces isotopic to $\partial M_1$ so that we can find a Reeb orbit of an appropriate $\text{lSFT}$ index in order to bound the minimal discrepancy of $A$ from below.

%The reason why we need a stable Hamiltonian hypersurface and not just a contact hypersurface is that
%%we need a hypersurface which behaves well with respect to the $\P^1$ bundle $\pi_D$,
%and we need to stretch the neck so that this $\P^1$ bundle has a particular form.
%For more details about this see Step 3 in the proof of
%Theorem \ref{theorem:reeborbitlowerboundarounddivisors}.

\bigskip

\proof of Lemma \ref{lemma:gwtripletest}.
The hard part of showing that $(S_1,[A],\mathcal{J}_{i,\omega_t})$  is a smooth deformation of GW triples
is showing that they all satisfy property
(\ref{item:deformationmaximumprinciple}) of Definition \ref{defn:gwtripledeformation}.
%Here is the idea of how this is proven:
%One should think of a sequence of almost complex structures

Let $U^1_S := \dot{S}_- \cup \big((-\epsilon_h,\frac{\epsilon_h}{2}) \times C\big) \cup U_S$.
We will show that for $i$ large enough, every genus $0$ nodal $J$-holomorphic curve in $S_1$
representing $[A]$ for $J \in \mathcal{J}_{i,\omega_t}$ is contained inside the relatively compact open set $U^1_S$.
Suppose for a contradiction that we have a sequence $t_i \in [0,1]$,
a sequence  $J^1_{i,t} \in \mathcal{J}_{i,\omega_t}$
and a sequence of genus $0$ $J^1_{i,t_i}$-holomorphic curves
$u_i : \Sigma_i \to S_1$ representing $[A]$ not contained in $U^1_S$ but all contained in some much larger relatively compact set $S$.
After passing to a subsequence we can assume that $t_i$ converges to $t_\infty$ for some $t_\infty \in [0,1]$
and hence $\omega_{t_i}$ $C^\infty$ converges to $\omega_{t_\infty}$.
Because $J^1_{i,t_i}$ stretch the neck along $C$ and converge in $C^\infty_{\text{loc}}$ to $J_+$
inside $\dot{S}_+$ as $i \to \infty$ we have by Proposition \ref{proposition:compactnessresult}
a finite energy $J_+$-holomorphic curve in $\dot{S}_+$ intersecting $Q$ with multiplicity $1$.
%with $u_i|_{u_i^{-1}(\dot{S}_+)}$ Gromov converging in .
Such a curve is contained inside $U_S$ by assumption.
Proposition \ref{proposition:compactnessresult} then tells us that for some
$i$ large enough, the image of $u_i$
is contained in $U^1_S$ which gives us a contradiction.
Hence for $i$ large enough we have that the image of every $J^1_{i,t}$-holomorphic curve
is contained in the relatively compact open set $U^1_S$.
Hence $(S_1,[A],J^1_{i,t})$ is a smooth deformation of GW triples.
\qed

\bigskip

We are interested in indices of Reeb orbits and the following lemma gives us an upper bound for the Conley-Zehnder index of a Reeb orbit.
In order to state the Lemma we need another definition of discrepancy.
This definition is more general than
Definition \ref{defn:symplecticminimaldiscrepancy}.
Here we will be dealing with discrepancies of other symplectic submanifolds
which do not correspond to exceptional divisors of some resolution of our singularity but instead
correspond to `divisors' which partially compactify our singularity.
\begin{defn} \label{defn:moregeneraldiscrepancy}
Let $(M,\omega)$ be a not-necessarily compact symplectic manifold of dimension $2n$ with compact boundary
and with a choice of trivialization of the $N$th power of the canonical bundle $K_M$ along $\partial M$ for some $N>0$.
Let $S_1,\cdots,S_l$ be properly embedded (not-necessarily compact) connected codimension $2$ normal crossings
submanifolds which are symplectic
so that $\partial M \hookrightarrow M \setminus \cup_i S_i$ is a homotopy equivalence
and so that $\cup_i S_i$ does not intersect $\partial M$.
The Borel-Moore homology group $H^{BM}_{2n-2}(M,\Q)$ is isomorphic to
$H^2(M,\partial M;\Q) = H^2(M,M \setminus \cup_i S_i;\Q)$.
We also have that $H^2(M,M \setminus \cup_i S_i;\Q)$ is naturally isomorphic to
$\oplus_i H^2(M;M \setminus S_i) \cong \oplus_i H^0(S_i;\Q) \cong \Q^l$ by a repeated use of Mayor-Vietoris and the Thom isomorphism theorem.
Hence $H^{BM}_{2n-2}(M;\Q)$ is freely generated
by the fundamental classes of $S_i$.
Choose a smooth section of $K_M$ which is
transverse to $0$ and equal  to a non-zero constant with respect to our trivialization along $\partial M$.
The zero set could be non-compact. Nonetheless it represents a locally-finite cocycle
in Borel-Moore homology $H^{BM}_{2n-2}(M;\Q)$.
The homology class of its zero set
is homologous to $\sum_i a_i [S_i]$ for some $a_1,\cdots,a_l \in \Q$.
The coefficient $a_i \in \Q$ is called the {\bf discrepancy} of $S_i$
with respect to our trivialization.
\end{defn}

%We will assume that $\partial M$ is compact but $M$ may not be compact.
%Now let $K$ be the canonical bundle of $M$ (i.e. the highest exterior power of $T^* M$ viewed as a complex vector bundle)
%and let $\tau$ be a trivialization of $K^{\otimes N}$ along $\partial M$ for some $N \in \N$.
%Using the trivialization $\tau$ we have a complex line bundle $\overline{K}^{\otimes N}$
%in the quotient $M / \partial M$.
%We define the relative first Chern class $c_1(K;\Q) \in H^2(M,\partial M;\Q)$ of the pair $(M,\partial M)$ with respect to $\tau$
%to be $\frac{1}{N} c_1(\overline{K};\Q) \in H^2(M / \partial M;\Q) = H^2(M,\partial M;\Q)$.

Unlike Lemma \ref{lemma:gwtripletest}, the following Lemma involves a contact hypersurface instead of a stable Hamiltonian hypersurface.
This is because the stable Hamiltonian hypersurface is used to establish a maximum principle whereas the contact hypersurface
corresponds to our link which is a contact manifold. In particular these are {\it different} hypersurfaces with different purposes.

\begin{lemma} \label{lemma:reeborbits}
Suppose that $(S,[A],\mathcal{J})$ is a GW triple and $C \subset S$ a contact hypersurface
with associated contact form $\alpha$ satisfying the following properties:
\begin{enumerate}
\item The hypersurface $C$ splits $S$ into two stable Hamiltonian cobordisms $S_+$, $S_-$ where $\partial_- S_+ = C$ and
$\partial_+ S_- = C$.
We will also assume that there is a compact codimension $2$
(not necessarily symplectic) submanifold $Q$ of $\dot{S}_-$ so that $[A].Q \neq 0$.
\item
%The contact form $\alpha$ on $C$ only has non-degenerate Reeb orbits.
We will suppose that $C$ has a natural trivialization of the $N$th tensor power of its canonical bundle so that we can define Conley-Zehnder
indices for Reeb orbits, and hence also the $\text{lSFT}$ index.
This also implies that the $N$th tensor power of the canonical bundle of $S_+$ restricted to $\partial_- S_+$ also has a chosen trivialization.
\item \label{item:submanifoldproperties}
We have two properly embedded codimension $2$ symplectic submanifolds $D_\infty$, $E$ of $\dot{S}_+$ which intersect transversally
and are holomorphic with respect to some $J \in \mathcal{J}$ and disjoint from a standard neighborhood $C \times (-\epsilon_h,\epsilon_h)$
of $C$ (see Appendix A for a definition of standard neighborhood).
We will assume that $E$ is compact.
We require $[A] \cdot E = 1$, $[A] \cdot D_\infty = 0$
and $S_+ \setminus (D_\infty \cup E)$ deformation retracts onto $C = \partial_- S_+$.
\item \label{item:almostcomplexstructureproperties}
For any compatible almost complex structure $J$ on $S$ equal to some $J_1 \in \mathcal{J}$ outside a small (fixed) neighborhood of the closure
of $S_- \cup ([0,\epsilon_h) \times C)$, we have $J \in \mathcal{J}$.
We will assume that $\mathcal{J}$ is open in the $C^\infty$ topology.
Finally we assume $\text{GW}(S,[A],\mathcal{J}) \neq 0$.
\end{enumerate}
%Let $a,b$ the unique rational numbers with the property that $a E + b D_\infty$ is
%Lefschetz dual to the relative first Chern class of $S_+$.
Then there exists a Reeb orbit $\gamma$ of $(C,\alpha)$
satisfying $\text{lSFT}(\gamma) \leq 2(n-3 - a)$
%of Conley-Zehnder index less than or equal to $n-3 - 2a$
if $n-3  - a \geq 0$ and satisfying $\text{lSFT}(\gamma) < 0$ otherwise
where $n = \frac{1}{2}\text{dim}_\R (S)$ and $a$ is the discrepancy of $E$ inside $S_+$.

The period of  $\gamma$ is $\leq -w$
where $w$ is the wrapping number of
$\theta$ around $E$
for any $1$-form $\theta$ on $S_+ \setminus (D_\infty \cup E)$
satisfying $d\theta = \omega_S$ and $\theta|_{\partial S_+} = \alpha$.
\end{lemma}

\proof of Lemma \ref{lemma:reeborbits}.
%A small neighborhood of $C \subset S$ can be identified with
%In the standard neighborhood
%$C \times (-\epsilon_h,\epsilon_h)$, we have $\omega = d(e^{r_C} \alpha)$
%where $r_C$ parameterizes $(-\epsilon_h,\epsilon_h)$.
%For any function $\phi \in C^\infty(C)$ sufficiently small,
%we have a contact embedding $C \hookrightarrow C \times (-\epsilon_h,\epsilon_h) \subset S$ of $(C,e^{\phi} \alpha)$
%sending $c$ to $(e^{\phi(c)},c)$.
Because $C^\infty$ generic perturbations of $\alpha$ only have non-degenerate Reeb orbits,
we have in turn that generic perturbations of $C$ inside $S$ only have non-degenerate Reeb orbits.
Corollary \ref{corollary:lsftlowersemicontinuous} then tells us that if our result is true for any generic perturbation of
$C$ inside $S$, then our result is true for $C$ itself.
So from now on we will assume that $C$ only has non-degenerate Reeb orbits.

%We can assume that $\alpha$
%only has non-degenerate Reeb orbits
%This is because $C$ admits a $C^\infty$ small perturbation inside $S$ along with $\alpha$
%so that all the Reeb orbits of $\alpha$ are non-degenerate,
%and as soon as we have our result for generic perturbations of $\alpha$
%we have it for $\alpha$ as well by Lemma \ref{corollary:lsftlowersemicontinuous}.

We stretch the neck along $C$ using almost complex structures $J_i$ from $\mathcal{J}$.
We will also assume that $J_i$ is cylindrical near $C$ (see Appendix A for a definition).
This can be done by property (\ref{item:almostcomplexstructureproperties}) from above.
We can also ensure that $J_i|_{\dot{S}_+}$ converges in $C^\infty_{\text{loc}}$ to
an almost complex structure $J_+$ compatible with the completion of $S_+$.
We will also assume that $J_+$ makes $E$ and $D_\infty$ holomorphic.
%Also $D_\infty$ and $E$ are $J_+$-holomorphic submanifolds by property \ref{item:submanifoldproperties}.
%Note that any $J_+$-holomorphic curve of finite energy with a smooth domain and which intersects
%$E$ with multiplicity $1$ is somewhere injective.

%Because $\text{GW}(S,[A],\mathcal{J}) \neq 0$ we have after passing to a subsequence a $J_i$-holomorphic curve
%$u_i : \Sigma_i \to S$ representing $[A]$ where $\Sigma_i$ is a genus $0$ nodal curve.
%Each curve intersects $E$ once and does not intersect $D_\infty$ because $J_+$ is
%$C^\infty$ close to an almost complex structure making $D_\infty$ and $E$ holomorphic.
%The point here is that the compactness result \cite{BEHWZ:compactnessfieldtheory} tells us that if $J_+$ was actually equal to $J$
%near $D_\infty$ and $E$ then the curves $u_i$ would converge to some curve intersecting $E$ once and $D_\infty$ $0$ times by positivity of intersection. The problem is that $J_+$ is supposed to only be $C^\infty$ close to $J$ near $E$ and $D_\infty$.
%Having said that perturbing $J_+$ slightly near $E$ and $D_\infty$ does not change the property by a similar compactness argument.

Because $\text{GW}(S,[A],\mathcal{J}) \neq 0$ we have a  $J_i$-holomorphic curve from a genus $0$ connected
nodal curve representing $[A]$.
By SFT compactness \cite[Theorem 10.3]{BEHWZ:compactnessfieldtheory}
we get a connected $J_+$-holomorphic curve $u_\infty : \Sigma_\infty \to \dot{S}_+$
which intersects $E$ once and does not intersect
$D_\infty$ by positivity of intersection,
which is contained in $\kappa \cap \dot{S}_+$ where $\kappa \subset S$ is a compact set and
where $u_\infty$ is a proper map.
%This is the top level of our holomorphic building.
%Now the condition $[A].E = 1$ ensures that $\Sigma_\infty$ is non-empty.
The holomorphic curve $\Sigma_\infty$ has irreducible components
$\Sigma_\infty^1,\Sigma_\infty^2,\cdots$.
%We assume that $J_+$ is sufficiently $C^\infty$ close to $J'_+$ so that any $J_+$
%each irreducible component has non negative intersection number with $E$ and $D_\infty$.
Exactly one of these components intersects $E$
with multiplicity $1$ after mapping them to $\dot{S}_+$ under $u_\infty$.
After relabeling such components we can assume that this component is $\Sigma_\infty^1$.
Also by positivity of intersection we have that no irreducible component intersects $D_\infty$
after mapping them to $\dot{S}_+$ by $u_\infty$.

Because the image of $u_\infty|_{\Sigma^2_\infty}$ is contained in the region $\dot{S}_+ \setminus (D_\infty \cup E)$
where $\omega_S$ is exact, we have by the maximum principle in \cite[Lemma 7.2]{SeidelAbouzaid:viterbo}
that such a curve does not exist.
Hence $\Sigma_\infty^2$ does not exist and so $\Sigma_\infty$ has exactly one irreducible component $\Sigma^1_\infty$.

Now we need to show that $u_\infty|_{\Sigma_\infty}$ converges to one or more Reeb orbits.
In other words, $\Sigma_\infty$ compactifies to a manifold $\overline{\Sigma}_\infty$ with non-empty boundary
and the map $u_\infty|_{\Sigma_\infty}$ continuously extends to
$\overline{u}_\infty : \overline{\Sigma}_\infty \to S_+$ where
$\overline{u}_\infty$ restricted to each connected component of
$\partial \overline{\Sigma}_\infty$ is a Reeb orbit of $(C,\alpha)$ after parameterizing the boundary appropriately
(see \cite[Proposition 5.6]{BEHWZ:compactnessfieldtheory}).
%The results in \cite{BEHWZ:compactnessfieldtheory} tell us that in fact $u_i$ converges
%to a holomorphic building.
%The lowest level of such a holomorphic building is contained in $S_-$.
%Such a lowest level is non-empty because $[A].Q \neq 0$ where $Q$ is a compact submanifold of $\dot{S}_-$.
%Now because our holomorphic building is connected we get that the top level of our holomorphic building
%$u_\infty$ must converge to a Reeb orbit.
Theorem 10.3 from \cite{BEHWZ:compactnessfieldtheory} also tells us that
because $[A].Q \neq 0$ and all of our $J_i$-holomorphic curves are connected,
one has that $\partial \overline{\Sigma}_\infty \neq \emptyset$.
Hence $u_\infty$ is an irreducible $J_+$-holomorphic curve intersecting $E$ once
and $D_\infty$ $0$ times and converging to at least one Reeb orbit in $C = \partial_- S_+$.

By property (\ref{item:almostcomplexstructureproperties}),
we can assume that $J_+$ is $C^\infty$ generic among almost complex structures equal to $J_+$
outside a compact subset of $\dot{S}_+$ and making $D_\infty$ and $E$ holomorphic.
In particular we can assume by \cite{Dragnev:transversality} that any irreducible $J_+$
holomorphic curve intersecting $E$ with multiplicity $1$ and with finite energy is regular
by perturbing $J_+$ near $E$ (but not along $E$ or $D_\infty$). This is because it is somewhere injective near $E$ as its intersection multiplicity
with $E$ is $1$.
In particular we can assume that $u_\infty$ is regular
(see the comment after Definition \ref{defn:pseudoholomorphiccurve}).
Now \cite[Proposition 5.6]{BEHWZ:compactnessfieldtheory} tells us that $\Sigma_\infty$ compactifies to a surface with boundary
$\overline{\Sigma}_\infty$ and $u_\infty$ extends continuously to a map $\overline{u}_\infty : \overline{\Sigma}_\infty \to S_+$
so that $\overline{u}_\infty(\partial \overline{\Sigma}_\infty)$ is a union of Reeb orbits $\gamma_1,\cdots,\gamma_l$.

We now need to compute the sum of the Conley-Zehnder indices of these orbits using \cite{Dragnev:transversality}.
The orientation of the boundary of $\overline{\Sigma}_\infty$ coming from the {\it inward} normal is equal to
the natural orientation of the Reeb orbits $\gamma_1,\cdots, \gamma_l$ (i.e. these orbits are negative ends).
Now $\overline{u}_\infty$ intersects $E$ once and $D_\infty$ $0$ times
and so let $p \in \overline{\Sigma}_\infty$ be the unique point satisfying $\{p\} = \overline{u}_\infty^{-1}(E \cup D_\infty)$.
Let $K_+$ be the canonical bundle of $\dot{S}_+$.
Because $S_+ \setminus (E \cup D_\infty)$ deformation retracts onto $C = \partial_- S_+$ we have
a trivialization of $K_+^{\otimes N}|_{S_+ \setminus (E \cup D_\infty)}$ coming from our trivialization of $K_+^{\otimes N}|_{\partial_- S_+}$ and this gives us a trivialization
\[\tau : \C \times (\overline{\Sigma}_\infty \setminus \{p\})  \to
\overline{u}_\infty^* K_+^{\otimes N}|_{\overline{\Sigma}_\infty \setminus \{p\}}.\]
Let  $p_\C : \C \times (\overline{\Sigma}_\infty \setminus \{p\}) \twoheadrightarrow \C$ be the natural projection.
Using our trivialization $\tau$, any smooth section $\sigma$ of $\overline{u}_\infty^* K_+^{\otimes N}|_{\overline{\Sigma}_\infty \setminus \{p\}}$ is given by a function $p_\C \circ \tau^{-1} \circ \sigma : \overline{\Sigma}_\infty \setminus \{p\} \to \C$.
By Definition \ref{defn:moregeneraldiscrepancy},
we can choose such a section $\sigma$ so that $p_\C \circ \tau^{-1} \circ \sigma$ is equal to $z^{Na}$ near $p$ where $z$ is a local holomorphic coordinate chart around $p$
and is non-zero away from $p$ and constant near $\partial \overline{\Sigma}_\infty$.
This means we can view $p_\C \circ \tau^{-1} \circ \sigma$ as a map from $\overline{\Sigma}_\infty \setminus \{p\}$ to $\C^*$.
Because $\overline{\Sigma}_\infty$ is homotopic to a wedge of circles we can
choose some trivialization $\upsilon$ of $u_\infty^* K_+$ (i.e. $\upsilon : \C \times \Sigma_\infty \to u_\infty^* K_+$
is a complex bundle isomorphism). Unlike $\tau$ our choice is not unique up to homotopy.
Using our trivializations $\tau$ and ${\upsilon}^{\otimes N}$
respectively we can get trivializations of the $N$th tensor power of the canonical bundle
along $\gamma_1,\cdots,\gamma_l$ and this means we get two Conley-Zehnder indices $\text{CZ}_\tau(\gamma_j)$ and
$\text{CZ}_{\upsilon}(\gamma_j)=\text{CZ}_{{\upsilon}^{\otimes N}}(\gamma_j)$.
Now $T := (({\upsilon}^{\otimes N})^{-1} \circ \tau)|_{\partial \overline{\Sigma}_\infty}$
is an automorphism of trivial bundles and so we can view $T$ as a map $T : \partial \overline{\Sigma}_\infty \to \C^*$.
Because $\partial \overline{\Sigma}_\infty$ is a union of $l$ oriented circles, this means that $T$
represents $l$ elements of $\pi_1(\C^*) = \Z$ given by $q_1,\cdots,q_l \in \Z$ (from now on we will give $\partial \overline{\Sigma}_\infty$
the orientation coming from the inward normal, which is the same orientation as the Reeb orbits).
Properties \ref{item:catenation} and \ref{item:determinantproperty} from Section \ref{section:conleyzehnderindex}
imply
$-2 q_j + N\text{CZ}_{\upsilon}(\gamma_j) =  N\text{CZ}_{\tau}(\gamma_j)$
(remember that the Conley-Zehnder index is calculated from a trivialization of some multiple of the tangent bundle and not the cotangent bundle
and so the dual of our morphism $T$ sends the trivialization of the $N$th tensor power of the {\it anti}canonical bundle
induced by the dual of $\upsilon$ to the trivialization induced by the dual of the $N$th tensor power of $\tau$ which is represented by $-2q_j$).
Also because our section $\sigma$
is constant along $\partial \overline{\Sigma}_\infty$ with respect to our trivialization $\tau$ we get that
$({\upsilon}^{\otimes N})^{-1} \circ \sigma|_{\partial \overline{\Sigma}_\infty} : \partial \overline{\Sigma}_\infty \to \C^*$ is homotopic to $T$.
This implies that $\sum_j q_j = Na$.
Hence $-2 a + \sum_{j=1}^l \text{CZ}_{\upsilon}(\gamma_j) =  \sum_{j=1}^l \text{CZ}_{\tau}(\gamma_j)$.

By the remark after Corollary 2 in \cite{Dragnev:transversality},
we get that $(n-3)(2 - l) - \sum_{j=1}^l \text{CZ}_{\upsilon}(\gamma_j) \geq 0$.
Hence $(n-3)(2 - l) - 2a - \sum_{j=1}^l \text{CZ}_{\tau}(\gamma_j) \geq 0$.
This implies:
$\sum_{j=1}^l (\text{CZ}_{\tau}(\gamma_j) + (n-3)) \leq 2(n-3) - 2a$.
Hence
$\sum_{j=1}^l \text{lSFT}(\gamma_j) \leq 2(n-3 - a)$.
If $n-3 -a \geq 0$ then there exists a $j$ so that
$\text{lSFT}(\gamma_j) \leq 2(n-3 - a)$.
If $n-3 - a  < 0$ then there exists a $j$ so that
$\text{lSFT}(\gamma_j) < 0$.
And so $\text{lSFT}(\gamma_j) < 0$.
%Hence by choosing $\gamma = \gamma_j$ depending on
%whether $n-3 - a \geq 0$ or not, we have
%$\text{lSFT}(\gamma) \leq 2(n-3 - a)$ if $2(n-3-a) \geq 0$
%and $< 0$ otherwise.
The bound on the period of $\gamma_j$ comes from the fact that $\int_{\Sigma_\infty} u_\infty^* \omega_S \geq 0$
and that $-w - \sum_i \text{period}(\gamma_i) = \int_{\Sigma_\infty} u_\infty^* \omega_S$ by Stokes' theorem.
This gives us our result with $\gamma := \gamma_j$.
\qed

\bigskip

\section{Bounding Minimal Discrepancy from Below} \label{section:discrepancyfrombelow}

\subsection{Bounding Minimal Discrepancy of Positively Intersecting Submanifolds from Below}

The main theorem of this section is:
\begin{theorem} \label{theorem:reeborbitlowerboundarounddivisors}
Let $(C,\xi)$ be contactomorphic to the contact link of
a strongly $\Q$-Gorenstein positively intersecting divisor $(M,\cup_i S_i,\theta)$.
Then
\begin{itemize}
\item If $\text{md}(M,\cup_i S_i,\theta) \geq 0$ then $\text{hmi}(C,\xi) \leq 2\text{min}_i (a_i)$.
\item If $\text{md}(M,\cup_i S_i,\theta) < 0$ then $\text{hmi}(C,\xi) < 0$.
\end{itemize}
%where $D\gamma : \text{ker}(\alpha)|_{\gamma(0)} \to \text{ker}(\alpha)|_{\gamma(0)}$
%is the linearized first return map of the Reeb flow along $\gamma$, restricted to the contact hyperplane distribution.
\end{theorem}

In the next subsection we will apply this to isolated singularities. The symplectic manifold $M$ should be thought of as a neighborhood
of the exceptional divisors of a resolution of an isolated singularity and $S_i$ should be thought of as an irreducible component of the exceptional divisor.

\proof of Theorem \ref{theorem:reeborbitlowerboundarounddivisors}.
It is sufficient to prove that for any contact form $\alpha$ associated to $\xi$,
there is a Reeb orbit $\gamma$ of $\alpha$ which either has negative $\text{lSFT}$ index,
or which has $\text{lSFT}$ index $\leq 2\text{md}(M,\cup_i S_i,\theta)$.
So from now on we fix such an $\alpha$.
We will prove this theorem in 5 steps.
\begin{enumerate}
\item
In the first step we first deform the submanifolds $S_i$ so that they are orthogonal.
Next we use use Corollary \ref{corollary:specificformaintheorem} to construct a link region $(M_1,g : M \setminus \cup_i S_i \to \R)$,
a chosen `minimal' pseudo Morse Bott family of Reeb orbits $B$ and an
$(2\epsilon,\epsilon_M)$-ball product associated to  $B$ and $S_1$ inside $M_1$
for some $\epsilon > \epsilon_M$ after relabeling the $S_i$'s.
\item
In the next step we partially compactify $M_1$ to a symplectic manifold $\breve{S}$
by replacing the $(2\epsilon,\epsilon_M)$ ball product $B^{2n-2}_{2\epsilon} \times B^2_{\epsilon_M}$ with $B^{2n-2}_{2\epsilon} \times S^2_\epsilon$ where $S^2_\epsilon$ is the $2$-sphere with area $\pi \epsilon^2$.
We then blow up $\breve{S}$ at a point outside $M_1$ and inside $\{0\} \times S^2_\epsilon$ giving us a partial compactification $S$.
Now the problem is that we want $S$ to be part of a GW triple by using Lemma \ref{lemma:gwtripletest}
which in turn uses a neck stretching argument.
In Step 3,
it turns out that we need a sequence of almost complex structures which is are a product
in the region $B^{2n-2}_{2\epsilon} \times S^2_\epsilon \subset \breve{S}$ and which stretch the neck along $\partial M_1$.
This is needed for us to apply Lemma  \ref{lemma:gwtripletest} in Step 4.
But the problem is that if we do this then we cannot have a product almost complex structure in $B^{2n-2}_{2\epsilon} \times S^2_\epsilon$
if $\partial M_1$ is a contact manifold.
Instead we have to deform $\partial M_1$ inside $S$ near $B^{2n-2}_{2\epsilon} \times S^2_\epsilon$
so that it is {\it stable Hamiltonian} inside $B^{2n-2}_{2\epsilon} \times S^2_{\epsilon}$, so that its properties from Step 1 remain the same
and so that one can neck stretch while retaining the product complex structure in $B^{2n-2}_{2\epsilon} \times S^2_\epsilon$.
\item
The point of this step is to construct an appropriate almost complex structure on $S$
so that the last statement in part (\ref{item:stayinginsidecompactset}) of Lemma \ref{lemma:gwtripletest} is satisfied.
We now choose a sequence of almost complex structures on $\breve{S}$ which are a product inside $B^{2n-2}_{2\epsilon} \times S^2_\epsilon$
and ones on $S$ for which the blowdown map $\text{Bl} : S \to \breve{S}$ is holomorphic
and which stretch the neck along our stable Hamiltonian hypersurface.
Let $[A]$ be the homology class represented by the proper transform of the curve $\{0\} \times S^2$.
We show that holomorphic curves representing $[A]$
stay inside the union of a compact subset of $\text{Bl}^{-1}(B^{2n-2}_{2\epsilon} \times S^2_\epsilon)$
with the compact set bounded by our stable Hamiltonian hypersurface.
This is done using an index calculation combined with an energy estimate.
To calculate the energy estimate one needs
the almost complex structure to be a product inside $B^{2n-2}_{2\epsilon} \times S^2_\epsilon$
and hence why one needs a stable Hamiltonian hypersurface and not just a contact hypersurface.
\item
In the fourth step we apply Lemma  \ref{lemma:gwtripletest} to make $S$ along with the almost complex structure from Step 3
and the homology class $[A]$ into a GW triple.
It turns out that every genus $0$ holomorphic curve representing our homology class inside $S$ also stays inside
$B^{2n-2}_{2\epsilon} \times S^2_\epsilon$ and hence one can show that the associated genus $0$
Gromov-Witten invariant is non-zero.

Next we symplectically dilate the symplectic structure inside our GW triple so that $(C,\alpha)$
can be embedded as a contact hypersurface inside $S$ ready for step 5.
%This process is called a symplectic dilation.
This stretching induces a deformation of GW triples.
All have non-zero GW invariant by Lemma \ref{lemma:deformationinvariance}.
\item
Here we use Lemma \ref{lemma:reeborbits}
with our `stretched' GW triple with $(C,\alpha)$ embedded inside it
to show that
$\alpha$ has a Reeb orbit whose lower SFT index
is at most the minimal discrepancy.
\end{enumerate}
\bigskip

{\it Step 1}. In this step we will construct a symplectic manifold and a codimension $0$ submanifold
whose boundary is contactomorphic to $(C,\alpha)$ and whose Reeb flow has nice properties.
By \cite[Theorem 5.3]{McLean:affinegrowth} we can deform $S_1,\cdots,S_l$ through positively intersecting submanifolds
so that $S_1,\cdots,S_l$ becomes symplectically orthogonal.
This induces a deformation of 
strongly numerically $\Q$-Gorenstein
positively wrapped divisors, which all have the same
minimal discrepancy and which have contactomorphic links by
Corollary \ref{corollary:canonicalboundary}.
%So from now on we can assume that $S_1,\cdots,S_l$ are symplectically orthogonal.
Hence we can assume that $(M,\cup_i S_i,\theta)$ is an orthogonal $\Q$-Gorenstein positively intersecting divisor.
Let $\lambda_1,\cdots,\lambda_l$ be the wrapping numbers
and let $a_1,\cdots,a_l$ be the discrepancies of
$S_1,\cdots,S_l$ respectively.
% as defined in Definition \ref{defn:symplecticminimaldiscrepancy}.
After relabeling the divisors $S_1,\cdots,S_l$, we can assume that:
\begin{itemize}
\item if $\text{md}(M,\cup_i S_i,\theta) < 0$, then $a_1 < 0$ and 
$\lambda_i \geq \lambda_1$ for any $i$ with $a_1 < 0$,
\item
if $\text{md}(M,\cup_i S_i,\theta) \geq 0$, then $a_1 = \text{md}(M,\cup_i S_i,\theta)$ and
$\lambda_i \geq \lambda_1$ for any $i$ with $a_i = a_1$.
\end{itemize}

Hence by Corollary \ref{corollary:specificformaintheorem}  we have for $\epsilon>0$ small,
 a link region $(M_1,g : M \setminus \cup_i S_i \to \R)$,
a pseudo Morse-Bott family of Reeb orbits $B$ of $(\theta + dg)|_{\partial M_1}$
and an $(2\epsilon,\epsilon_M)$-ball product associated to $B$ and $S_1$ inside $M_1$ for $\epsilon_M < \epsilon$
so that:
\begin{enumerate}
%\item $\epsilon > \epsilon_M$.
\item  $\text{lSFT}(B) = 2\text{md}(M,\cup_i S_i,\theta)$
and the period of $B$ is $\pi \epsilon_M^2 + \lambda_1$.
\item \label{item:indexupperbounds}
Any other Reeb orbit $\gamma$ not contained in $B$ of period strictly less than the period of $B$ minus $\epsilon^2$
satisfies $\text{lSFT}(\gamma) \geq 0$,  $\text{lSFT}(\gamma) > \text{lSFT}(B)$ and
is non-degenerate.
\end{enumerate}

Define $C_1 := \partial M_1$, $\theta_g := \theta + dg$.
We have that $(C_1,\theta_g|_{C_1})$ is contactomorphic to the link of $\cup_i S_i$ by Definition \ref{defn:linkofsi}.
We write $\partial_+ M_1 = C$
and $\partial_- M_1 = \emptyset$.

\bigskip
{\it Step 2}. We will now construct an appropriate symplectic manifold containing $M_1$
along with some natural compatible almost complex structures on this manifold.
Define $\omega_{g,+} := d\theta_g|_{C_1}$.
By flowing backwards along $X_{\theta_g}$ we will assume that a neighborhood of $\partial_+ M_1 = C_1$
is diffeomorphic to $(-\epsilon,0] \times C_1$ with $\theta_g = e^r \alpha_1$ where $r$ parameterizes $(-\epsilon,0]$.
%Define $C_2 := \{r = -\frac{\epsilon}{2}\}$.

 Now let $S^2_{\epsilon}$ be the two dimensional symplectic sphere of area $\pi \epsilon^2$ and let
$B^2_{\epsilon} \subset S^2_{\epsilon}$ be some symplectic embedding of the disk of radius $\epsilon$ into the sphere
and let $q_\infty \in S^2_{\epsilon} \setminus B^2_{\epsilon}$
be the unique point in the complement of this embedding.
We define $(\breve{S},\omega_{\breve{S}})$ to be the symplectic manifold given by the interior of the gluing of $M_1$ with $B^{2n-2}_{2\epsilon} \times S^2_{\epsilon}$
along $W_1 = B^{2n-2}_{2\epsilon} \times B^2_{\epsilon_M}$ where $W_1 \subset M_1$ is the submanifold described earlier.
We define $(S,\omega_S)$ to be the symplectic blowup of $(\breve{S},\omega_{\breve{S}})$ along $\{0\} \times \{q_\infty\} \in B^{2n-2}_{2\epsilon} \times S^2_{\epsilon}$. We let $\text{Bl} : S \twoheadrightarrow S^1$ be the blowdown map which is a diffeomorphism away from some
symplectic submanifold $E$ (the exceptional divisor) and a symplectomorphism outside some very small open subset containing $E$. Hence because $q_\infty$ is disjoint from $M_1$, we have that $M_1$ is naturally
a submanifold of $S$.
We require that this blow up should be small enough so that the restriction of $\omega_S$ to $M_1$ is $\omega_{\breve{S}}$.

We will now put a stable Hamiltonian structure on $C_1$.
Define $\upsilon_M= \frac{3}{2}\epsilon$.
Define $\alpha_+ \in \Omega^1(C_1)$ to be $\alpha_1$ outside $B^{2n-2}_{2\epsilon} \times \partial B^2_{\epsilon_M} \subset C_1$
and equal to $(\frac{1}{2}r_1^2 + \frac{1}{2\pi}\lambda_1) d\vartheta$ inside $B^{2n-2}_{\upsilon_M} \times \partial B^2_{\epsilon_M} \subset C_1$.
We do this by choosing some bump function $\rho$ on $B^{2n-2}_{2\epsilon}$ equal to $0$ near its boundary and equal to $1$ along
$B^{2n-2}_{\upsilon_M}$ and then defining $\alpha_+$ to be $(\text{pr}_1^* (1-\rho)) \alpha_1 + (\text{pr}_1^*\rho) (\frac{1}{2}r_1^2 + \frac{1}{2\pi}\lambda_1) d\vartheta$ where $\text{pr}_1 : B^{2n-2}_{2\epsilon} \times B^2_{\epsilon_M} \to B^{2n-2}_{2\epsilon}$ is the natural projection map.
Then $(\omega_{g,+},\alpha_+)$ is a Stable Hamiltonian structure on $C_1$ with exactly the same Reeb flow as $\alpha_1$.
In particular the Reeb orbits of $(\omega_{g,+},\alpha_+)$ are exactly the same as the Reeb orbits of
$\theta_g|_{C_1}$ and they have exactly the same periods and indices.
The restriction of $\alpha_+$ to $\{a\} \times \partial B^2_{\epsilon_M}$ is  $(\frac{1}{2}r_1^2 + \frac{1}{2\pi}\lambda_1) d\vartheta$
for each $a \in B^{2n-2}_{2\epsilon}$.
%Because $(\omega_{g,+},\alpha_+)$ is a contact structure away from our pseudo Morse-Bott family $B$
%we can perturb $\alpha_+$ outside $B$ and hence also $\omega_{g,+} = d\alpha_+$ by a $C^\infty$ small amount so that all Reeb orbits of $\alpha_+$
%whose period is less than $2\pi \epsilon_M^2 + \lambda_1$
%are non-degenerate (note that these orbits are disjoint from $R_{S_1}$).
%By Lemma \ref{lemma:perturbedindexcalculation} we have that the
%Conley-Zehnder indices of any of these non-degenerate Reeb orbits of period less than
%$2\pi \epsilon_M^2 + \lambda_1$
%is strictly greater than $\text{max}(3-n-\frac{1}{2N},2a_1 + 3-n)$ from the inequality
%(\ref{indexinequality_eqn}) in Step $1$.

A small neighborhood of $C_1$ inside $S$ is symplectomorphic to $(-\eta,\eta) \times C_1$
with $\omega_S = \omega_{g,+}  + d(r_1 \alpha_+)$ where $\eta > 0$ is small and $r_1$ parameterizes $(-\eta,\eta)$.
Let $J_i$ be a sequence of almost complex structures on $S$ compatible with $\omega_S$
which stretch the neck along $C_1$ with respect to $(-\eta,\eta) \times C_1$.
Now $C_1$ splits $S$ into two regions $S_{+}$ and $S_{-}$ where $S_{+},S_{-}$
are cobordisms of stable Hamiltonian structures with $\cup_i S_i \subset S_{-}$,
$\partial_+ S_{+} = \partial_- S_{-} = \emptyset$ and $\partial_- S_{+} = \partial_+ S_{-} = C_1$.
Let $\breve{S}_{+} := \text{Bl}(S_{+}) \subset \breve{S}$
where $\text{Bl} : S \to \breve{S}$ is our blowdown map. 
This is also a cobordism of stable Hamiltonian structures.
We will assume that $J_i|_{\dot{S}_{+}}$ converges in $C^{\infty}_{\text{loc}}$ to some almost complex structure $J_{+}$
compatible with the completion of $S_{+}$.
We will choose $J_{+}$ so that the blowdown map $\text{Bl}|_{S_+}$ is a $(J_{+},\breve{J}_{+})$-holomorphic map
where $\breve{J}_{+}$ is compatible with the completion of $\breve{S}_{+}$.
We can assume that $\breve{D}_\infty := B^{2n-2}_{2\epsilon} \times \{q_\infty\}$ is $\breve{J}_{+}$-holomorphic
and that the preimage of $\breve{D}_\infty$ under the blowdown map is a union of
transversely intersecting $J_{+}$-holomorphic submanifolds
$E \cup D_\infty$ where $D_\infty$ the proper transform of $\breve{D}_\infty$
%(proper transform here means we take the closure of the preimage of $\breve{D}_\infty \setminus \{q_\infty\}$ under %the blowdown map).
Because $(\omega_{g,+},\alpha_+)$
is a specific stable Hamiltonian structure
in the region $B^{2n-2}_{\upsilon_M} \times S^2_{\epsilon}$ we can assume (maybe after deforming our neighborhood $(-\eta,\eta) \times C_1$) that $\breve{J}_{+}$
is a product almost complex structure $(J_{B^{2n-2}},J_{\infty,S^2})$ on
$B^{2n-2}_{\upsilon_M} \times (S^2_{\epsilon} \setminus B^2_{\epsilon_M})$ where $J_{B^{2n-2}}$
is the standard complex structure on $B^{2n-2}_{\upsilon_M} \subset \C^{n-1}$ and
where $J_{\infty,S^2}$
is some compatible complex
structure on
$S^2_{\epsilon} \setminus B^2_{\epsilon_M}$.
For each $i$ we will assume that $\breve{J}_i$ restricted to
$B^{2n-2}_{\upsilon_M} \times S^2_{\epsilon}$
is the product complex structure
$(J_{B^{2n-2}},J_{i,S^2})$
where $J_{i,S^2}$ is a complex structure on the sphere $S^2_{\epsilon}$.
Finally we can assume that $J_i$ and $J_+$
in the regions $\text{Bl}^{-1}(B^{2n-2}_{\upsilon_M} \times S^2_{\epsilon})$
and
$\text{Bl}^{-1}(B^{2n-2}_{\upsilon_M} \times (S^2_{\epsilon} \setminus B^2_{\epsilon_M}))$
respectively are integrable
(this can be done my ensuring that $\text{Bl}$ corresponds to a K\"{a}hler blowup
with respect to $J_i$ for all $i$ and $J_+$).

\bigskip
{\it Step 3}. In this step we will show that any finite energy proper $J_{+}$-holomorphic curve intersecting $E$ once and $D_\infty$ $0$ times is contained inside a fixed compact subset of $S_{+}$.
In fact we will show that it is contained inside the interior of $\text{Bl}^{-1}(B^{2n-2}_{\upsilon_M} \times S^2_{\epsilon})$
where $\text{Bl} : S \to \breve{S}$ is the blowdown map.
The reason why we want to do this is that in Step 4 we create an appropriate family of GW triples using Lemma \ref{lemma:gwtripletest} combined with results from this step.

Let $u :  \Sigma \to \dot{S}_{+}$ be such a curve.
Let $\breve{u} : \Sigma \to \dot{\breve{S}}_{+}$ be the composition of $u$ with the blowdown map.
We will first show $\Sigma$ is irreducible. Let $\Sigma_1,\Sigma_2,\cdots$
be the irreducible components of $\Sigma$.
Exactly one of these components
intersects $\breve{D}_{\infty}$
after mapping them with $\breve{u}$
by positivity of intersection.
We will assume after reordering
these components that $\Sigma_1$
intersects $\breve{D}_{\infty}$.
If $\Sigma_2$ exists then it must
map to the stable Hamiltonian
cobordism $\dot{\breve{S}}_{+} \setminus \breve{D}_{\infty}$.
Because $\theta_g$  evaluated on
the Reeb vector field along $\partial_- \dot{\breve{S}}_{+}$
is $1$
and because $\theta_g$ extends to
a $1$-form
on  $\dot{\breve{S}}_{+} \setminus \breve{D}_{\infty}$
whose exterior derivative is the symplectic form
 we then get that
$\Sigma_2$ cannot exist
by Proposition \ref{proposition:maximumprinciple}.
Hence
$\Sigma$ is irreducible.

Now $(\text{pr}_1 \circ \breve{u})|_{\text{pr}_1^{-1}(B^{2n-2}_{\upsilon_M})} : \text{pr}_1^{-1}(B^{2n-2}_{\upsilon_M}) \to B^{2n-2}_{\upsilon_M}$ is a
$J_{B^{2n-2}}$-holomorphic map, whose domain is a punctured Riemann surface,
and so by the removable singularity theorem this extends to a proper map.
Hence this
is either the constant map, or a map whose energy is greater than or equal to $\pi (\upsilon_M)^2$ by the monotonicity formula.
Let $\sigma_i \in \Sigma$ be a sequence of points so that $\breve{u}(\sigma_i)$ gets
arbitrarily close to
$\partial_- \dot{\breve{S}}_{+}$.
Now Lemma \ref{lemma:convergence} tells us that after passing to a subsequence, $\breve{u}(q_i)$
converges to a point on a Reeb orbit of period at most $\pi \epsilon^2+ \lambda_1$ because the intersection of $\breve{u}$ with $\breve{D}_\infty$ is $1$
and the wrapping number of $\theta_g$ around $\breve{D}_\infty$ is $-2 \pi (\epsilon)^2 - \lambda_1$ after extending $\theta_g$ in some way
to $\dot{\breve{S}}_+ \setminus \breve{D}_\infty$
so that $d\theta_g$ is the symplectic form
on $\breve{S}_+ \setminus \breve{D}_\infty$.
%Here is a more detailed construction of this extension of $\theta_g$:
%Using the $1$-form $\theta_g$, we can create a new $1$-form $\breve{\theta}_{+}$ on $\breve{S}_{+} \setminus \cup_i S_i$  with the following properties:
%\begin{enumerate}
%\item $d\breve{\theta}_{+} = \omega_{\breve{S}}$.
%\item Inside $B^{2n-2}_{\upsilon_M} \times B^2_{\epsilon}$ it is equal to
%$\text{pr}_1^* \beta + \text{pr}_2^* ((r_i^2 + \frac{1}{2\pi}\kappa_1) d\vartheta_i)$. Here $B^2_{\epsilon} = S^2_{\epsilon} \setminus \{q_\infty\}$.
%\end{enumerate}
%Hence the wrapping number of $\breve{\theta}_{+}$ around $\breve{D}_\infty$ is $-(2\pi (\epsilon)^2 + \kappa_1)$.

We have two cases: {\it Case 1}: This Reeb orbit has period $\geq \pi (\epsilon_M)^2 + \lambda_1 - \epsilon^2$
for some sequence of points $\sigma_i$ as described above.
Recall that here $\pi (\epsilon_M)^2 + \lambda_1$ is the period of the pseudo Morse-Bott family $B$ constructed in Step $1$.
{\it Case 2}: Every such Reeb orbit has period $< \pi (\epsilon_M)^2 + \lambda_1 -\epsilon^2$
for any sequence of points. We will in fact show that this case cannot occur
using a proof by contradiction.

\smallskip
{\it Case 1}:
Suppose this Reeb orbit has period $\geq \pi (\epsilon_M)^2 + \lambda_1 - \epsilon^2$.
The second half of Lemma \ref{lemma:convergence}
tells us $\int_{\Sigma} (\breve{u})^* \omega_{\breve{S}} \leq \pi (\epsilon)^2 + \lambda_1 - (\pi (\epsilon_M)^2 + \lambda_1) < \pi (\upsilon_M)^2$.
The energy argument above then tells us that the image of $\breve{u}$ under $\text{pr}_1$ is contained inside $B^{2n-2}_{\upsilon_M} \subset B^{2n-2}_{2\epsilon}$. Hence the image of $\breve{u}$ is contained inside $B^{2n-2}_{\upsilon_M} \times S^2_{\epsilon}$
in this case.

\smallskip
{\it Case 2}: This case involves a lot more work. 
Here we suppose (for a contradiction) that every such sequence $\sigma_i$ converges to a Reeb orbit of length $< \pi (\epsilon_M)^2 + \lambda_1-\epsilon^2$.
All such Reeb orbits are non-degenerate by Corollary \ref{corollary:specificformaintheorem}.
% as they cannot come from the pseudo Morse-Bott family $B$ or any multiple of such Reeb orbits.
Lemma \ref{lemma:convergence} tells us that $\Sigma$ is biholomorphic to $\P^1 \setminus \{w_1,\cdots,w_l\}$
where $w_1,\cdots,w_l$ are $l$ distinct points.
By \cite{BEHWZ:compactnessfieldtheory} these punctures converge to non-degenerate Reeb orbits $R_1,\cdots,R_l$.
We have that $u$ is somewhere injective because it intersects our holomorphic submanifold $E$ with multiplicity $1$
and because $\Sigma$ is irreducible.
%Now by a compactness argument it is %sufficent to show that every $J_{+}$ %holomorphic curve is contained inside %the interior of
%$B^{2n-2}_{\upsilon_M} \times %S^2_{\epsilon}$ even if we perturb %$J_{+}$ by a $C^\infty$ small amount.
By perturbing $J_+$ appropriately
we can ensure that $J_+$
has the property that
every somewhere injective
$J_+$-holomorphic curve not contained
in $E \cup D_\infty \cup (B^{2n-2}_{\upsilon_M} \times (S^2_{\epsilon} \setminus B^2_{\epsilon_M}))$
and whose domain is a punctured
sphere is regular by \cite{Dragnev:transversality}
(see the comment about regularity after Definition \ref{defn:pseudoholomorphiccurve}).
We can still ensure that all the properties of $J_{+}$ hold.
The point here is that
to achieve regularity we only
need to perturb $J_+$
outside $E \cup D_\infty \cup (B^{2n-2}_{\upsilon_M} \times (S^2_{\epsilon} \setminus B^2_{\epsilon_M}))$
as every somewhere injective
$J_+$-holomorphic curve not contained
in $E \cup D_\infty \cup (B^{2n-2}_{\upsilon_M} \times (S^2_{\epsilon} \setminus B^2_{\epsilon_M}))$ has a somewhere
injective point outside this region.

Let $K$ be the canonical bundle of $S$ and $\breve{K}$ the canonical bundle of $\breve{S}$.
Because $H_1(C_1;\Q) = 0$ and
$N c_1(TS|_{C_1}) = 0$ we have a canonical trivialization $\tau : C_1 \times \C \to K^{\otimes N}|_{C_1} = {\breve{K}}^{\otimes N}|_{C_1}$ of $K^{\otimes N}$ and ${\breve{K}}^{\otimes N}$
along $C_1 = \partial_- S_{+} = \partial_- \breve{S}_{+}$.
Also because $Q := \breve{S}_{+} \cap (B^{2n-2}_{2\epsilon} \times S^2_{\epsilon})$ is contractible we also have
a canonical trivialization $\breve{\tau}_+$ of ${\breve{K}}^{\otimes N}$ in this region. Similarly we have a trivialization $\breve{\tau}_-$ of ${\breve{K}}^{\otimes N}$ inside
$M_1 \cap (B^{2n-2}_{2\epsilon} \times S^2_{\epsilon})$.
Now $\partial Q$ is homotopic to a circle (oriented to coincide with the boundary of $\{0\} \times B^2_{\epsilon_M}$)
and so ${\breve{\tau}_-}^{-1} \circ {\breve{\tau}}_+$ is represented by a map from $S^1 \simeq \partial Q$ to $S^1 = U(1)$.
The degree of this map is $-2 N$ because the Chern number of
the canonical bundle of $S^2_{\epsilon}$ is $-2$ and the normal bundle of $\{0\} \times S^2_{\epsilon}$
inside $\breve{S}$ is trivial as a complex vector bundle.
Hence the degree of ${\breve{\tau}_+}^{-1} \circ {\breve{\tau}}_-$ is $2 N$.
Also ${\breve{\tau}_-}^{-1} \circ \tau$ has degree $a_1 N$ by the definition of discrepancy.
Hence ${\breve{\tau}_+}^{-1} \circ \tau$ has degree $(2 + a_1) N$.
The trivialization $\breve{\tau}_+$ also gives a trivialization $\tau_+$ of $K^{\otimes N}$ along
$\text{Bl}^{-1}(Q) \setminus E$ because $\text{Bl}$ is a biholomorphism onto its image
away from $E$.
Because $E$ is the exceptional divisor of a blowup we can construct a section
$\sigma$ of $K^{\otimes N}$ along
$\text{Bl}^{-1}(Q)$ which is constant with respect to $\tau_+$ near $\partial_- S_{+}$
and non-zero away from a small neighborhood of $E$ and whose zero set is homologous to $(n-1) N [E]$ near $E$.
All of this implies that there is a section $\sigma_1$ of $K^{\otimes N}|_{S_+}$
which is constant with respect to the trivialization $\tau$ along $\partial S_+$
and whose zero set is homologous to $(-2-a_1) N D_\infty + ((n-1) - 2 - a_1) N [E] = (-2 - a_1) N D_\infty + (-a_1 + (n-3)) N [E]$
in Borel-Moore homology
(note that $2+a_1$ now becomes $-2-a_1$ in our calculations because if we have a disk inside $(B^{2n-2}_{2\epsilon} \times S^2_{\epsilon}) \cap S_+$ intersecting $D_\infty$ positively once
with boundary $\{0\} \times \partial B^2_{\epsilon_M}$ then its boundary orientation is the opposite of $\{0\} \times \partial B^2_{\epsilon_M}$).
This means that the relative Chern number of $u^* K^{\otimes N}$ with respect to the pullback of the trivialization $\tau$
is $(-a_1 + (n-3)) N$.
Hence the relative Chern number of the pullback via $u$ of the $N$th tensor power of
the {\it anti}-canonical bundle is $(a_1 - (n-3)) N$.

Now choose some trivialization $\widehat{\tau}$ of the pullback of the anticanonical bundle along $u$.
%Note that near the Reeb orbits $R_i$, we have
By the statement after Corollary 2 in \cite{Dragnev:transversality} we get that
$(n-3)(2 - l) - \sum_{i=1}^l \text{CZ}_{\widehat{\tau}}(R_i) \geq 0$.
Now $N\sum_{i=1}^l \text{CZ}_{\widehat{\tau}}(R_i) = N \sum_{i=1}^l \text{CZ}(R_i) - 2(a_1 - (n-3)) N$
by the above Chern number calculation.
Hence
$(n-3)(2 - l) + 2a_1 - 2(n-3)  - \sum_{i=1}^l \text{CZ}(R_i) \geq 0$.
Hence $\sum_i (\text{CZ}(R_i)+ (n-3)) \leq 2a_1$ and so
$\sum_i \text{lSFT}(R_i) \leq 2a_1$.
Now if $a_1 \geq 0$ then $\text{lSFT}(R_i) \leq 2a_1$ for some $i$
contradicting property (\ref{item:indexupperbounds}) from Step $1$
(here we have used the fact that the stable Hamiltonian structure
$(\omega_{g,+},\alpha_+)$ has exactly the same Reeb orbits with the same periods
and $\text{lSFT}$ indices as the Reeb orbits of $\theta_g|_{\partial M_1}$).
If $a_1 < 0$ then $\text{lSFT}(R_i) < 0$ for some $i$
contradicting property (\ref{item:indexupperbounds}) from Step $1$.
Therefore we have a contradiction and so case $2$ doesn't happen.
Hence $u$ is contained inside $B^{2n-2}_{\upsilon_M} \times B^2_{2\epsilon}$.

\bigskip
{\it Step 4}. In this step we construct an appropriate family of symplectic structures and almost complex structures on $S$
giving us a family of GW triples.
We will use Lemma \ref{lemma:gwtripletest} combined with Step 3 to do this.

We now have two regions containing $C_1$. One is $(-\epsilon,0] \times C_1$ with $\theta_g = e^r\alpha_1$
and the other is $(-\eta,\eta) \times C_1$
with $\omega_S = \omega_{g,+}  + d(r_1 \alpha_+)$.
To avoid confusion we will write the first neighborhood $(-\epsilon,0] \times C_1$  as
$(-\frac{\epsilon}{2},\frac{\epsilon}{2}] \times C_2$ with $C_2 := \{r=-\frac{\epsilon}{2}\}$ and $\theta_g = e^{r_2} \alpha_2$ where $\alpha_2 = \theta_g|_{C_2}$
and $r_2$ parameterizes the interval $(-\frac{\epsilon}{2},\frac{\epsilon}{2}]$.
%We will assume that $\eta$ is small enough so that $(-\eta,\eta) \times C_1$
%is disjoint from $[-\frac{\epsilon}{4},\frac{\epsilon}{4}] \times C_2$.
We will assume that $J_i$ is cylindrical inside $[-\frac{\epsilon}{4},\frac{\epsilon}{4}] \times C_2$ (see Appendix A for a definition)
and also translation invariant in this region with respect to the coordinate $r_2$.
We can also assume that $J_i = J_j$ in this region for all $i,j$.

The hypersurface $C_2$ splits $S$ into two manifolds with boundary $S_{2,+}$ and $S_{2,-}$
where $S_{2,-}$ is the region containing our exceptional divisors.
Let $(\omega_t)_{t \in [0,\infty)}$ be a symplectic dilation of
$\omega_S$ along $C_2$ with respect to the contact form
$\theta_g|_{C_2}$ as in Definition \ref{defn:symplecticdilation}.
We can perform this symplectic dilation in such a way so that
the associated family of stretching $1$-forms
$\rho_t(r) r (\theta_g|_{C_2})$
as in Definition \ref{defn:symplecticdilation}
is equal to $\theta_g$ near $C_2$ inside $S_{2,+}$.
We will assume that the support of our symplectic dilation is disjoint
from $S_+ \cup (-\eta,\eta) \times C_1$.
Hence we can define $\theta_t$ to be a smooth family of $1$-forms
on $S \setminus \cup_i S_i$ as:
\[ \theta_t := \left\{ 
\begin{array}{cc}
\theta_g & \text{inside} \quad S_{2,+} \\
\rho_t(r) r (\theta_g|_{C_2}) & \text{inside the support of our symplectic dilation}  \\
\frac{1}{1+t} \theta_g & \quad \text{elsewhere.} \end{array}
 \right. \]

We let ${\mathcal J}_{i,t}$ be the set of almost complex structures compatible with $\omega_t$
and equal to $J_i$ in the region $S_+ \cup ((-\eta,\eta) \times C_1)$.
By Lemma \ref{lemma:stretchingembedding} there is some $t_{\text{max}} > 0$
and some embedding $\iota : C \hookrightarrow S$ 
with the property that $\iota^* \theta_{t_{\text{max}}} = c\alpha$
for some constant $c>0$
and so that $\iota$
 isotopic to $C_2$
through contact embeddings in the support of our
symplectic dilation.
Hence $\iota$ is isotopic to $C_1$ through contact embeddings inside $S_-$.
We define our homology class $[A]$ to be represented by the proper transform of $\{0\} \times S^2_{\epsilon}$ (i.e. represented by the closure of the preimage of
$\{0\} \times (S^2_{\epsilon} \setminus q_\infty)$ inside $S$).
We have $c_1(S,\omega_S)([A]) +n-3 = 0$.
So by using Lemma \ref{lemma:gwtripletest} combined with Step 3
(where the objects
$C,J_i,J_+,t,\omega_t,{\mathcal J}_{i,\omega_t},[A],\epsilon_h$
in the statement of Lemma \ref{lemma:gwtripletest}
correspond to
$C_1,J_i,J_+,t,\omega_t, {\mathcal J}_{i,t},[A],\eta$
here respectively)
we get for $i$ large enough that $(S,[A],{\mathcal J}_{i,t})$
with symplectic structures $\omega_t$
is a smooth family of GW triples parameterized by $t \in [0,t_{\text{max}}]$.
%Hence  $(S,[A],{\mathcal J}_{i,t})$
%with the original smooth family of symplectic structures $\omega_R$ is a
%smooth family of GW triples parameterized by
%$R \in [\frac{\epsilon}{4},R_{\text{max}}]$.
%Now for $R = \frac{\epsilon}{4}$ we can assume that $\phi_R$ is the identity map with $J_{i,R} = J_i$
%and $\theta_R = \theta_g$.
%So for $R = \frac{\epsilon}{4}$, we have that every

We now need to calculate $\text{GW}_0(S,[A],{\mathcal J}_{i,t})$.
Now $\theta_0 = \theta_g$.
We have $J_i \in {\mathcal J}_{i,0}$.
We have that any
$J_i$-holomorphic curve representing $[A]$ has energy less than $\upsilon_M$ with respect to $\theta_g$.
If $u$ is such a curve then $\text{pr}_1 \circ \text{Bl} \circ u|_{u^{-1}(\text{Bl}^{-1}(B^{2n-2}_{\upsilon_M} \times S^2_{\epsilon}))}$
must have image equal to a point 
by the monotonicity formula.
Hence the image of every such curve $u$ is contained inside $\text{Bl}^{-1}(B^{2n-2}_{\upsilon_M} \times S^2_{\epsilon})$.
Hence in fact the image of $u$ must be contained inside
the proper transform of
$\{0\} \times S^2_{\epsilon}$ and so there is a unique such curve up to re-parameterization.
The almost complex structure here is integrable near the image of $u$,
the normal bundle of $u$ is
$\prod_{i = 1}^{n-1} \O(-1)$
and $u$ is injective and so the $u$ is regular.
Hence $\text{GW}_0(S,[A],{\mathcal J}_{i,0}) = \text{GW}_0(S,[A],J_i) = 1$
by part (\ref{item:actualcounting}) of Theorem \ref{theorem:gromovwitteninvariants}.
By Lemma \ref{lemma:deformationinvariance} we then get
$\text{GW}_0(S,[A],{\mathcal J}_{i,t}) = 1$ for each $t \in [0,t_{\text{max}}]$.
%Here we are using the $1$-form %$\theta_R$ and symplectic form %$\omega_R$.

\bigskip
{\it Step 5}.
The contact embedding $\iota : C \hookrightarrow S$ splits $S$ into two stable Hamiltonian cobordisms
$S_{+,\iota}$, $S_{-,\iota}$ where $S_{-,\iota}$ contains $\cup_i S_i$.
Hence by Lemma \ref{lemma:reeborbits}, using
our contact embedding $\iota : (C,d(c\alpha)) \hookrightarrow (S,\omega_{t_{\text{max}}})$, cobordisms $S_{+,\iota},S_{-,\iota}$ and $D_\infty$, $E$ as above,
we get that $(C,\alpha)$ admits a Reeb orbit $\gamma$ with $\text{lSFT}$ index $\leq 2(n-3-(-a_1 + (n-3)))=2a_1$
if $n-3 - (-a_1 + (n-3)) = a_1 \geq 0$ and less than $0$ otherwise.
Here we used the Borel-Moore homology class calculation of the zero set of $\sigma_1$ from Step $3$.
Hence $\text{mi}(\alpha) \leq 2a_1 = 2\text{md}(A,0)$ if $a_1 \geq 0$ or $\text{mi}(\alpha) < 0$ if $\text{mi}(A,0)=a_1 < 0$.
This proves our Theorem.

\qed

\subsection{Bounding Minimal Discrepancy of a Singularity from Below} \label{section:boundingminimaldiscrepancyfrombelow}

Let $A \subset \C^N$ be an affine variety of dimension $n$ with an isolated singularity at $0$.
Recall that the link of $A$ is a contact manifold $(L_A,\xi_A)$.

\begin{theorem} \label{theorem:lowerbound}
Suppose $L_A$ satisfies $H^1(L_A;\Q) = 0$ and $c_1(\xi_A;\Q) = 0\in H^2(L_A;\Q)$.
Then
\begin{itemize}
\item
If $\text{md}(A,0) \geq 0$ then
$\text{hmi}(L_A,\xi_A) \leq 2 \text{md}(A,0)$ and
\item
if $\text{md}(A,0) < 0$ then $\text{hmi}(L_A,\xi_A) < 0$.
\end{itemize}
\end{theorem}
Proof of Theorem \ref{theorem:lowerbound}.
Let $A_\delta$ be the intersection of $A$ with a small ball of radius $\delta$.
We resolve $A$ at $0$ by blowing up along smooth subvarieties by \cite{hironaka:resolution} and take the preimage $\widetilde{A}_\delta$ of $A_\delta$ under this resolution map.
We suppose that the exceptional divisors $E_1,\cdots,E_l$ are smooth normal crossing.

By Lemma \ref{lemma:symplecticformonresolution},
there exists a $1$-form $\theta_A$ on $A_\delta \setminus \cup_i E_i$
making $(A_\delta,\cup_i E_i,\theta_A)$ into a strongly numerically $\Q$-Gorenstein
positively wrapped divisor
where the discrepancy of $E_i$
defined in Section \ref{section:minimaldiscrepancydefinition} coincides with
the associated discrepancy from Definition
\ref{defn:symplecticminimaldiscrepancy}.
Our result then follows from Theorem \ref{theorem:reeborbitlowerboundarounddivisors}.
\qed

\section{Appendix A: Stable Hamiltonian Cobordisms} \label{section:appendixA}

Most material from this section is taken from \cite{BEHWZ:compactnessfieldtheory} and from \cite{cieliebak:stableHamiltonian}.	
We need to cover this material as we have to deal with certain compactness results  involving  holomorphic curves
in stable Hamiltonian structures whose associated Reeb flow is not necessarily Morse-Bott.
Having said that the structures we are interested in are pseudo Morse-Bott submanifolds although we will not need this condition here.
\begin{defn} A {\bf stable Hamiltonian structure} on a manifold $C$ of dimension $2n-1$
is a pair $(\omega_h,\alpha_h)$ where $\omega_h$ is a closed $2$-form
and $\alpha_h$ is a $1$-form with the property that $\alpha_h \wedge \omega_h^{n-1}$ is a volume form.
We also require that $\text{ker}(\omega_h) \subset \text{ker}(d\alpha_h$).
\end{defn}
Here for any differential form $\gamma$, $\text{ker}(\gamma)$ means the set of vectors $V$ so that $i_V \gamma = 0$.
Stable Hamiltonian structures are called `stable Hamiltonian' for the following reason (\cite[Lemma 2.3]{cieliebakmohnke:compactness}):
Let $(M,\omega)$ be a symplectic manifold and suppose we have a Hamiltonian $H : M \to \R$
so that there is a smooth family of diffeomorphisms $(\Psi_s)_{s \in (-\epsilon,\epsilon)}$
on $M$ sending $(H^{-1}(0),X_H|_{H^{-1}(s)})$ to $(H^{-1}(s),X_H|_{H^{-1}(s)})$,
then $H^{-1}(0)$ has a stable Hamiltonian structure with $2$-form given by $\omega|_{H^{-1}(0)}$.
In other words, $H^{-1}(0)$ has a stable Hamiltonian structure with $2$-form $\omega|_{H^{-1}(0)}$
if the dynamics of $H^{-1}(s)$ do not change (or are stable) for all $s$ near $0$.
Conversely any manifold with a stable Hamiltonian structure can be embedded in such a way into a symplectic manifold.
Such a definition was initially motivated by the Weinstein conjecture (which states that every contact manifold has at least one Reeb orbit),
and also by the pseudo-holomorphic curve theory used to prove specific instances of this conjecture (see \cite{hoferzehnder:symplecticinvariantsdynamics}).

\begin{defn}
The {\bf Reeb vector field} $R$ of a stable Hamiltonian structure $(\omega_h,\alpha_h)$,is the unique vector field $R$ on $C$ satisfying
$i_R \omega_h = 0$ and $i_R \alpha_h = 1$.
\end{defn}
The condition $R \in \text{ker}(\omega_h) \subset \text{ker}(\alpha_h)$ and $i_R \alpha_h = 1$ ensure that the flow of $R$
preserves $\omega_h$ and $\alpha_h$. The vector field $X_H$ inside $H^{-1}(0)$ in the above description is a scalar multiple of the Reeb vector field.
Hence one should imagine a Reeb vector field as a Hamiltonian vector field constrained to a `stable'  hypersurface inside a symplectic manifold.
Any contact structure $\lambda$ gives us a stable Hamiltonian structure $(\omega_h,\alpha_h)$
where $\alpha_h = \lambda$ and $\omega_h = d\lambda$.
An almost complex structure $J_C$ on the hyperplane bundle $\text{ker}(\alpha_h)$ is {\bf compatible} with $(\omega_h,\alpha_h)$
if it is compatible with the symplectic structure $\omega_h|_{\text{ker}(\alpha_h)}$.
%(i.e. for all non-zero $V,W \in \text{ker}(\alpha_h)$ we have $\omega_h(V,JV)>0$ and $\omega_h(JV,JW)=\omega_h(V,W)$).
We can define an almost complex structure $\widehat{J}_C$ on $\R \times C$ in the following way:
For vectors of the form $(0,V)$ where $V \in \text{ker}(\alpha_h)$ it is defined by $\widehat{J}_C(V) := J_C(V)$.
Also $\widehat{J}_C(\frac{\partial}{\partial r_C}) = R$ and $\widehat{J}_C(R) = -\frac{\partial}{\partial r_C}$
where $r_C$ parameterizes $\R$. We say that $\widehat{J}_C$ is a {\bf cylindrical almost complex structure}
associated to $J_C$.
A {\bf symplectization} of $(\omega_h,\alpha_h)$ is the product $(-\epsilon_h,\epsilon_h) \times C$ for $\epsilon_h>0$ small
with symplectic form $\widetilde{\omega}_h := \omega_h + r_C d\alpha_h + dr_C \wedge \alpha_h$.
Here by abuse of notation we have identified $\alpha_h$ with $\pi_C^* \alpha_h$ where $\pi_C : (-\epsilon_h,\epsilon_h) \times C \twoheadrightarrow C$
is the natural projection.
Also $\epsilon_h$ has to be sufficiently small to ensure that the closed $2$-form $\widetilde{\omega}_h$ is symplectic.
If $C \subset M$ is a compact and connected hypersurface of a symplectic manifold $(M,\omega)$ then we say that it is a {\bf stable Hamiltonian hypersurface} if
$\omega|_C = \omega_h$. A neighborhood of $C$ is symplectomorphic to its symplectization
$(-\epsilon_h,\epsilon_h) \times C$. We will call such a neighborhood a {\bf standard neighborhood}.

A complex structure $\widetilde{J}_C$ on $(-\epsilon_h,\epsilon_h) \times C$ is said to be {\bf compatible with the symplectization}
$(-\epsilon_h,\epsilon_h) \times C$ if
\begin{enumerate}
 \item it is compatible with the symplectic form $\widetilde{\omega}_h$,
 \item $\widetilde{J}_C(\text{ker}(\alpha_h)) = \text{ker}(\alpha_h)$,
 \item and $\widetilde{J}_C\Big(\frac{\partial}{\partial r_C}\Big) = f(r_C) R$ and $\widetilde{J}_C(R) = -\frac{1}{f(r_C)} \frac{\partial}{\partial r_C}$
 for some smooth positive function $f : (-\epsilon_h,\epsilon_h) \to (0,\infty)$.
\end{enumerate}
If $\widetilde{J}_C$ is only defined on $I \times C$ where $I \subset (-\epsilon_h,\epsilon_h)$
then we also say it is {\bf compatible with the partial symplectization} $I \times C$ if it satisfies the above properties.
The smooth positive function $f$ may only have domain given by our subset $I \subset (-\epsilon_h,\epsilon_h)$.

Let $(M,\omega)$ be a not-necessarily compact symplectic manifold whose boundary is a disjoint union
of compact submanifolds $\partial_- M \sqcup \partial_+ M$.
Suppose that we have a stable Hamiltonian structure $(\omega^\pm_h,\alpha^\pm_h)$ on $\partial_\pm M$.
We say that $(M,\omega)$ is a {\bf stable Hamiltonian  cobordism} from $\partial_- M$ to $\partial_+ M$ if there is a neighborhood of $\partial_+ M$
diffeomorphic to $(-\epsilon_h,0] \times \partial_+ M$ with $\omega = \omega^+_h + r_+ d\alpha^+_h + dr_+ \wedge \alpha_h^+$
and a neighborhood of $\partial_- M$
diffeomorphic to $[0,\epsilon_h) \times \partial_- M$ with $\omega = \omega^-_h + r_- d\alpha^-_h + dr_- \wedge \alpha_h^-$.
Here $\epsilon_h>0$ is a small constant and $r_+$ parameterizes $(-\epsilon_h,0]$ and $r_-$ parameterizes $[0,\epsilon_h)$.
We will call the neighborhoods  $(-\epsilon_h,0] \times \partial_+ M$ and  $[0,\epsilon_h) \times \partial_- M$
the positive and negative {\bf cylindrical ends} of $M$.
In the literature, $M$ is quite often compact as well, but we need a slightly more general definition for our purposes.

An almost complex structure $J$ is {\bf compatible with the stable Hamiltonian cobordism structure on} $M$ if it is compatible with the symplectic form and it is equal to almost complex structures
$\widetilde{J}_-,\widetilde{J}_+$ compatible with the partial symplectizations $[0,\epsilon_h) \times \partial_- M$ and $(-\epsilon_h,0] \times \partial_+ M$
on the cylindrical ends respectively.
An almost complex structure $J$ defined on the interior of $M$ is said to be {\bf compatible with the completion of } $M$
if it is compatible with the symplectic form and there are smooth maps
$\phi_+ : (-\epsilon_h,0) \to \R$, $\phi_- : (0,\epsilon_h) \to \R$ so that:
\begin{enumerate}
 \item On the cylindrical ends $J$ is equal to almost complex structures $\widetilde{J}_\pm$
 compatible with the partial symplectizations $(0,\epsilon_h) \times \partial_- M$ and $(-\epsilon_h,0) \times \partial_+ M$.
 \item \label{item:diffeomorphismontoimage}
 $\phi_+$ is an orientation preserving diffeomorphism onto its image $(0,\infty)$ and $\phi_-$ is an orientation preserving diffeomorphism onto its image $(-\infty,0)$.
 \item \label{item:cylindricalcondition} There is a compatible almost complex structure $J_+$ (resp. $J_-$) on $\text{ker}(\alpha_+)$ (resp. $\text{ker}(\alpha_-$)) with the property that
 $(s_r \circ \phi_+,\text{id})_* \widetilde{J}_+$
 (resp. $(s_r \circ \phi_-,\text{id})_* \widetilde{J}_-$) converges in $C^\infty_{\text{loc}}$ to $\widehat{J}_+$ in $\R \times \partial_+ M$ as $r \to -\infty$
 (resp. $r \to \infty$)
 where $s_r : \R \to \R$ sends $x$ to $x+r$.
 \end{enumerate}

Let $u : \Sigma \to \R \times C$ be any smooth map where $\Sigma$ is some surface.
We define $E_{\omega_h}(u)$ to be $\int_\Sigma (\pi_C \circ u)^* \omega_h$ where $\pi_C$
is the projection $\R \times C \twoheadrightarrow C$.
Let $r$ be the coordinate parameterizing $\R$ in $\R \times C$.
We define $E_{\alpha_h}(u) := \text{sup}_{\phi \in \mathcal{C}} \int_\Sigma (\phi \circ r \circ u) u^* (dr \wedge d\alpha_h)$
where $\mathcal{C}$ is the set of compactly supported smooth maps $\phi : \R \to \R$ whose integral is $1$.
We define $E_C(u) := E_{\omega_h}(u) + E_{\alpha_h}(u)$.
Now suppose $J_M$ is an almost complex structure compatible with the completion of $M$
and $\breve{u} : \breve{\Sigma} \to \dot{M}$ a $J_M$-holomorphic curve where $\dot{M}$ is the interior of $M$.
This means that there are smooth  maps
$\phi_+ : (-\epsilon_h,0) \to \R$, $\phi_- : (0,\epsilon_h) \to \R$ satisfying the conditions (\ref{item:diffeomorphismontoimage}), (\ref{item:cylindricalcondition}) above.
Define $N_+ M := (-\epsilon_h,0) \times \partial_+ M$ and $N_- M := (0,\epsilon_h) \times \partial_- M$.
We define $E_{\text{int}}(\breve{u}) := \int_{(\breve{u})^{-1}(\dot{M} \setminus (N_- M \cup N_+ M))} (\breve{u})^*\omega_M$.
Let $\breve{u}_- := \breve{u}|_{(\breve{u})^{-1}(N_- M)}$.
We define $E_-(\breve{u}) := E_{\partial_- M}((\phi_-,\text{id}_{\partial_- M}) \circ \breve{u}_-)$.
%Here $\phi_-^{-1} : (-\infty,0) \to (-%\epsilon_h,0)$ is the inverse of the %diffeomorphism $\phi_- : (-\epsilon_h,0) %\to (-\infty,0)$
%and so $(\phi_-^{-1},\text{id}_{\partial_- %M}) : (-\infty,0) \times \partial_- M \to (-%\epsilon_h,0) \times \partial_- M$
%is the inverse of the map $(\phi_-,%\text{id}_{\partial_- M})$.
We have a similar definition of $E_+(\breve{u})$.
We define $E_M(\breve{u}) := E_{\text{int}(\breve{u})} + E_-(\breve{u}) + E_+(\breve{u})$ and we will call this the {\bf energy} of $\breve{u}$.

 Suppose we have a symplectic manifold $(S,\omega_S)$.
Let $C \subset S$ be a stable Hamiltonian hypersurface with
stable Hamiltonian structure $(\omega_h,\alpha_h)$.
A small neighborhood of $C$ is symplectomorphic to the symplectization $(-\epsilon_h,\epsilon_h) \times C$
for small enough $\epsilon_h$.
We define $S_{\text{split}}$ to be the manifold with boundary obtained from $S$ by gluing $S \setminus C$ 
with the disjoint union $[0,\epsilon_h) \times C \sqcup (-\epsilon_h,0] \times C$ along
$(0,\epsilon_h) \times C \sqcup (-\epsilon_h,0) \times C \subset S \setminus C$.
Here $S_{\text{split}}$ is a cobordism of stable Hamiltonian structures with $\partial_- S_{\text{split}} = \partial_+ S_{\text{split}} = C$.
Note that the interior $\dot{S}_{\text{split}}$ of $S_{\text{split}}$ is diffeomorphic in a canonical way to $S \setminus C$.
\begin{defn}
A sequence of almost complex structures $J_i$ compatible with $\omega_S$
{\bf stretch the neck} along $C$ if there is
a compatible almost complex structure $J_C$ on $\text{ker}(\alpha_h)$ and a sequence of orientation preserving diffeomorphisms
$\phi_i : (-\epsilon_h,\epsilon_h) \to (-D_i,D_i)$ so that:
\begin{enumerate}
\item $\phi_i(0) = 0$ and
$\phi_i' = 1$ near $-D_i$ and $D_i$. Also $D_i$ is strictly increasing and $D_i \to \infty$ as $i \to \infty$.
\item $J_i|_{(-\epsilon_h,\epsilon_h) \times C}$ is compatible with the symplectization $(-\epsilon_h,\epsilon_h) \times C$.
\item $(\phi_i,\text{id}_C)_* J_i$ viewed as an almost complex structure in
$(-D_i,D_i) \times C \subset \R \times C$ converges in $C^\infty_{\text{loc}}$ to $\widehat{J}_C$ as $i \to \infty$. 
In the region $(-\frac{\epsilon_h}{2},\frac{\epsilon_h}{2}) \times C$ we will assume that the restriction of $J_i$
to $\text{ker}(\pi_C^* \alpha_h) \cap \text{ker}(dr_C)$ where $\pi_C$ is the projection to $C$ is uniformly convergent.
%Also if $\xi$ is the vector subbundle of $T\R \times C$ consisting of vectors in the kernel of $\alpha_h$
%and $dr_C$ then we require that $(\phi_i,\text{id}_C)_* J_i|_{\xi}$ uniformly converges to $\widehat{J}_C|_{\xi}$.
\item There is an almost complex structure $J_{\text{split}}$ compatible with the completion of $S_{\text{split}}$ and
an open neighborhood $NC$ of $C$ so that
$J_i|_{NC \setminus C}$ converges in $C^\infty_{\text{loc}}$ to
$J_{\text{split}}|_{NC\setminus C}$.
\end{enumerate}
\end{defn}
Note that we do not have any constraints on $J_i$ away from $(-\epsilon_h,\epsilon_h) \times C$
although when we prove a compactness result later we will assume that $J_i$ converges $C^\infty_{\text{loc}}$
to an almost complex structure in one region inside $S$ but there will be no constraints in other regions of $S$.
We need to do this in order to prove certain symplectic manifolds are in fact GW triples.

\begin{prop}
We can construct a sequence of almost complex structures $J^i_S$
stretching the neck along $C$ so that $J^i_S|_{S \setminus C}$
converges in $C^\infty_{\text{loc}}$ to an almost complex structure compatible
with the completion of $S_{\text{split}}$.
\end{prop}
\proof
Let $J_S$ be an almost complex structure on $S$ compatible with $\omega_S$ whose restriction to
$(-\epsilon_h,\epsilon_h) \times C$ is compatible with the symplectization
$(-\epsilon_h,\epsilon_h)  \times C$.
Let $J_C$ be an almost complex structure on $\text{ker}(\alpha_h)$ given by $J_S|_{TC \cap \text{ker}(\alpha_h)}$
and let $D_i>0$ be a sequence tending to $\infty$.
% $D_i > \text{max}(\epsilon_h^2,\epsilon_h,1)$ be a sequence satisfying $D_{i+1} > 3D_i$.
Define $\phi_\infty : (-\epsilon_h,\epsilon_h) \setminus \{0\} \to \R$ by
\[
\phi_\infty(x) := \left\{
\begin{array}{cc}
-1 - \frac{\epsilon_h}{x} & \text{when} \quad x < 0 \\
1 - \frac{\epsilon_h}{x} & \text{when} \quad x > 0.
\end{array}
\right.
\]
We choose a sequence of orientation preserving diffeomorphisms $\phi_i : (-\epsilon_h,\epsilon_h) \to (-D_i,D_i)$
so that for each $j>0$, $\phi'_i|_{(-D_j,D_j) \setminus \{0\}}$ converges in $C^\infty_{\text{loc}}$ to $\phi'_\infty|_{(-D_j,D_j) \setminus \{0\}}$,
and so that $\phi'_i = \phi'_\infty$ near $\pm \epsilon_h$.
%\begin{enumerate}
%\item $\frac{d\phi_i(x)}{dx} = 1$ near $\pm D_i$,
%\item $\phi_i(x) = \frac{D_i^2 x}{\epsilon_h}$ for $x \in \Big[-\frac{\epsilon_h}{3D_i},\frac{\epsilon_h}{3D_i}\Big]$.
%\item $\phi_i(x) = \phi_{i-1}(x)+D_i - D_{i-1}$ for $x \geq \frac{\epsilon_h}{2D_i}$
%and $\phi_i(x) = \phi_{i-1}(x)-D_i + D_{i-1}$ for $x \leq -\frac{\epsilon_h}{2D_i}$.
%\end{enumerate}
%
%By the definition of $J_S$ we have a function $f_S : (-\epsilon_h,\epsilon_h) \to (0,\infty)$ so that if
%$R$ is the Reeb vector field of $(\omega_h,\alpha_h)$ which we extend in the natural way to $(-\epsilon_h,\epsilon_h) \times C$
%then $J_S\Big(\frac{\partial}{\partial r_C}\Big) = f_S(r_C) R$ and $J_S(R) = -\frac{1}{f_S(r_C)} \frac{\partial}{\partial r_C}$ where $r_C$ is the natural coordinate parameterizing $(-\epsilon_h,\epsilon_h)$.
Let $R$ be the Reeb vector field of $(\omega_h,\alpha_h)$ which we extend
in the natural way to $\R \times C$.
Let $\pi_C : (-\epsilon_h,\epsilon_h) \times C \twoheadrightarrow C$ be the natural projection map.
%Let $\rho : (-\epsilon_h,\epsilon_h) \to (0,\infty)$ be a smooth function with $\rho(x) = f_S(x)$ for $|x| \leq \frac{\epsilon_h}{2}$ and $\rho(x) = 1$ for $x$ near $\pm \epsilon_h$.
For $i \in \N \cup \{\infty\}$,
we define $J_i(\frac{\partial}{\partial r_C}) := \phi'_i(r_C) R$
and $J_i(R) := -\frac{1}{\phi'_i(r_C)} \frac{\partial}{\partial r_C}$.
Also for $i \in \N \cup \{\infty\}$ and inside $(-\epsilon_h,\epsilon_h) \times C$
(and if $i = \infty$, outside $C$)
we define $J_i(V) := J_S(V)$ for $V \in \text{ker}(\pi_C^* \alpha_h) \cap \text{ker}(dr_C)$.
We also define $J_i$ to be any fixed compatible almost complex structure on $S \setminus (-\epsilon_h,\epsilon_h) \times C$ so that $J_i|_{S \setminus (-\epsilon_h,\epsilon_h) \times C} = J_j|_{S \setminus (-\epsilon_h,\epsilon_h) \times C}$ for all $i,j \in \N \cup \{\infty\}$.
We define a diffeomorphism $\phi_+ : (0,\epsilon_h) \to (-\infty,0)$ by $\phi_+(x) = \phi_\infty(x)$.
Similarly we define a diffeomorphism $\phi_- : (-\epsilon_h,0) \to (0,\infty)$  by $\phi_-(x) = \phi_\infty(x)$.
We define
$J_{\text{split}} := J_\infty$.
% on $S \setminus (-\epsilon_h,\epsilon_h) \times C$ and
%in the region:
% $S \setminus (-\frac{\epsilon_h}{2D_i},\frac{\epsilon_h}{2D_i})$ it is equal to $J_i$.
%This definition makes sense because $J_i = J_k$ in this region for all $k \geq i$.
Under the identification $S \setminus C = \dot{S}_{\text{split}}$ we get that $J_{\text{split}}$ is compatible
with the completion of $S_{\text{split}}$ using the functions $\phi_\pm$.
The sequence of almost complex structures $J_i$ is a sequence stretching the neck
along $C$ using the functions $\phi_i$
so that $J_i$ $C^\infty_{\text{loc}}$ converges in $S \setminus C$ to
$J_{\text{split}}$.
%Also by construction the pushforward $(\phi_i,\text{id}_C)_* J_i$ converges in $C^\infty_{\text{loc}}$
%to $\widehat{J}_C$ inside $(-\frac{\epsilon_h}{2},\frac{\epsilon_h}{2})$.
\qed

\bigskip

The main goal of this section is to prove a compactness result coming from neck stretching.
First we need a compactness result from \cite{fish:compactness} which we state here as a Theorem
(the compactness result from \cite{fish:compactness} is much stronger but we do not need the full force
of the theorem here).
Let $\Omega$ be a non-degenerate (not necessarily closed) $2$-form on a manifold $B$ and $J$ an almost complex structure.
We say that $(\Omega,J)$ is an {\bf almost Hermitian structure}  if $J$ is compatible with $\Omega$.
A nodal curve with boundary is given by a closed subset in the Euclidean topology
of a complex algebraic curve with nodal singularities with the property that the boundary
(i.e. the closure minus the interior) is a $1$-dimensional submanifold of the smooth part of this complex curve.
A $J$-holomorphic map from a nodal curve $\Sigma$ with boundary is a continuous map $\Sigma \to B$
which is smooth and $J$-holomorphic away from these singularities.

\begin{defn}
Suppose we have a sequence of almost Hermitian structures $(\Omega_i,J_i)$ $C^\infty$ converging to $(\Omega_\infty,J_\infty)$ in a manifold $B$.
Suppose that we have:
\begin{enumerate}
 \item a sequence of $J_i$-holomorphic curves $u_i : \Sigma_i \to B$ and a nodal $J_\infty$-holomorphic curve
 $u_\infty : \Sigma_\infty \to B$.
 \item a smooth surface $\widetilde{\Sigma}$ and smooth embeddings $v_i : \widetilde{\Sigma} \to \Sigma_i$,
 a continuous surjection $v_\infty : \widetilde{\Sigma} \twoheadrightarrow \Sigma_\infty$,
 \item a compact set $K \subset B$.
\end{enumerate}
so that
\begin{enumerate}
 \item If $\widetilde{\Sigma}^{\text{ns}} := \widetilde{\Sigma} \setminus v_\infty^{-1}(\Sigma_\infty^{\text{sing}})$
where $\Sigma_\infty^{\text{sing}}$ is the set of nodal points of $\Sigma_\infty$
then $v_\infty^{\text{nonsing}} := v_\infty|_{\widetilde{\Sigma}^{\text{ns}}}$ is a diffeomorphism onto its image
$\Sigma_\infty \setminus \Sigma_\infty^{\text{sing}}$.
\item  $u_i \circ v_i$ $C^0$ converges to $u_\infty \circ v_\infty$
and $u_i \circ v_i|_{\widetilde{\Sigma}^{\text{ns}}}$ $C^\infty_{\text{loc}}$ converges to $u_\infty \circ v_\infty|_{\widetilde{\Sigma}^{\text{ns}}}$.
\item $u_i \circ v_i$ and $u_\infty \circ v_\infty$ map the boundary of $\widetilde{\Sigma}$
to $B \setminus K$
\end{enumerate}
then we say that $u_i$ {\bf converges in the Gromov sense to $u_\infty$ near $K$}.
\end{defn}

The following theorem is proven by Fish in \cite{fish:compactness}.
\begin{theorem} \label{theorem:fishcompactnesscorollary}
Let $B$ be a compact manifold with boundary, $(\Omega_i,J_i)$ a sequence of almost Hermitian structures
$C^\infty$ converging to a
almost Hermitian structure $(\Omega_\infty,J_\infty)$.
Let $u_i : \Sigma_i \to B$ be a sequence of smooth genus $0$ $J_i$-holomorphic curves whose boundary maps to $\partial B$
with $\int_S u_i^* \Omega_i$ bounded above by a fixed constant independent of $i$
and let $K$ be a compact subset of the interior of $B$.
Then after passing to a subsequence there is a genus $0$ compact nodal $J_\infty$-holomorphic curve $u_\infty : \Sigma_\infty \to B$
so that $u_i$ converges in the Gromov sense to $u_\infty$ near $K$.
\end{theorem}

The following result is a much weaker version of the main result of
\cite{BEHWZ:compactnessfieldtheory} except that we allow our stable Hamiltonian structure to have Reeb orbits that
are not necessarily Morse-Bott (here Morse-Bott can be defined for stable Hamiltonian structures
in exactly the same way as Definition \ref{defn:morsebott}).
%Basically this compactness result only remembers the top level of the holomorphic building.
\begin{prop} \label{proposition:compactnessresult}
Let $(S,\omega_S)$ be a connected symplectic manifold and $C$ a stable Hamiltonian hypersurface in $S$ with a standard neighborhood
 $(-\epsilon_h,\epsilon_h) \times C \subset S$.
We suppose that $S \setminus C$ is a disjoint union $\dot{S}_+ \sqcup \dot{S}_-$ where $\dot{S}_+$ contains $(0,\epsilon_h) \times C$.
We define $S_+ := \dot{S}_+ \cup C$
which is a stable Hamiltonian cobordism from $\partial_- S_+ := C$ to $\partial_+ S_+ := \emptyset$.
Let $\breve{J}_i$ be a sequence of almost complex structures stretching the neck along $C$ using our symplectization $(-\epsilon_h,\epsilon_h) \times C$
with the additional property that $\breve{J}_i|_{\dot{S}_+}$ converges in $C^\infty_{\text{loc}}$
to an almost complex structure $J_+$ compatible with the completion of $S_+$.
Let $\omega_i$ be a sequence of symplectic structures in the same cohomology class equal to $\omega_S$ inside $(-\epsilon_h,\epsilon_h) \times C \cup \dot{S}_+$,
$C^\infty$ converging to $\omega_S$ and $J_i$ a sequence of $\omega_i$ compatible almost complex structures equal
to $\breve{J}_i$ inside $(-\epsilon_h,\epsilon_h) \times C \cup \dot{S}_+$.
Let $u_i : \P^1 \to S$ be a sequence of genus zero $J_i$-holomorphic curves so that:
\begin{enumerate}
\item their images stay inside a fixed compact set $\kappa \subset S$.
\item $\int_{\P^1} u_i^* \omega_i$ is bounded above by a fixed constant $E$.
\end{enumerate}
Then there exists a proper $J_+$-holomorphic curve $u_\infty : \Sigma_\infty \to \dot{S}_+$ so that:
\begin{enumerate}
\item $\Sigma_\infty$ is a genus $0$ nodal $J_+$-holomorphic curve without boundary mapping to $\kappa \cap \dot{S}_+$.
\item If $N \subset S$ is a neighborhood of the image of $u_\infty$ inside $S$ and $\epsilon_h > \eta>0$
then for all $i$, the image of $u_j$ is contained in $\dot{S}_- \cup ([0,\eta) \times C) \cup N$ for some $j \geq i$.
\item Its energy satisfies $E_{S_+}(u_\infty) < \infty$ and $\int_{\Sigma_\infty} u_\infty^* \omega_S \leq E$.
\end{enumerate}
Suppose $Q$  is a properly embedded codimension $2$ submanifold $Q \subset \dot{S}_+$ whose closure in $S$ is disjoint from $C$.
Then each $u_i$ for $i \gg 1$ has a positive intersection number with $Q$ if and only if $u_\infty$ also has a positive intersection number.
If $(u_i)_*([\P^1]) \cdot [Q] = 1$ for all $i$ then the intersection number of $u_\infty$ with $Q$ is also $1$.
\end{prop}

\proof of Proposition \ref{proposition:compactnessresult}.
We proceed in $3$ steps. In the first step we use the compactness result Theorem \ref{theorem:fishcompactnesscorollary}
to construct our holomorphic curve $u_\infty$.
In the second step we show that the energy of $u_\infty$ is bounded.
In the third step we prove the remaining parts of the Proposition.

\bigskip
{\it Step 1}:
Define $S^k_+ := \dot{S}_+ \setminus (0,\frac{\epsilon_h}{k}) \times C \subset S$.
We define $\Sigma_{k,+}^i := u_i^{-1}(S^k_+ )$. Maybe after perturbing $\epsilon_h$ slightly,
we can assume that $\Sigma_{k,+}^i$ is a manifold with boundary.
For each $k$ we have by Theorem \ref{theorem:fishcompactnesscorollary} that
a subsequence $u_{i_{k,j}}$ of $u_i|_{\Sigma_{k+1,+}^i} : \Sigma_{k+1,+}^i \to S^{k+1}_+$
 converges in the Gromov sense as $j \to \infty$ to a $J_+$-holomorphic nodal curve $u_{k,\infty} : \Sigma_{k,\infty} \to S^{k+1}_+$ near $S^k_+$.
By an inductive argument we can assume $\{ i_{k_1,j} | j \in \N\} \subset \{i_{k_2,j} | j \in \N\}$ for all $k_1 \leq k_2$.
We then get that the diagonal sequence $u_{i_{j,j}}$ has the property that for each $k \in \N$,
$u_{i_{j,j}}|_{\Sigma_{k+1,+}^{i_{j,j}}} : \Sigma_{k+1,+}^{i_{j,j}} \to S^{k+1}_+$ converges in the Gromov sense to $u_{k,\infty}$ near $S^k_+$.
After replacing $u_i$ with a subsequence we will assume that $u_j = u_{i_{j,j}}$ for all $j$.
The curve $u_{k,\infty}$ might have multiply covered components and so let $\overline{u}_{k,\infty} : \overline{\Sigma}_{k,\infty} \to S^{k+1}_+$ be the respective holomorphic curve without multiply covered components whose image is the same. Let $\overline{\Sigma}^{\text{nonsing}}_{k,\infty} \subset \overline{\Sigma}_{k,\infty}$ be the maximal open subset of the set of nonsingular points
with the property that 
$v_k := \overline{u}_{k,\infty}|_{\overline{\Sigma}^{\text{nonsing}}_{k,\infty}}$ is a diffeomorphism onto its image. Such an open set has a finite complement. For $k_1 \leq k_2$ we have that the image of $v_{k_1}$ is contained in the image of $v_{k_2}$. Hence we can construct a new Riemann surface $\Sigma^{\text{nonsing}}_\infty$ without boundary by gluing $\overline{\Sigma}^{\text{nonsing}}_{k_1,\infty}$ to $\overline{\Sigma}^{\text{nonsing}}_{k_2,\infty}$ using the maps $v_{k_2}^{-1} \circ v_{k_1}$ for all $k_1 \leq k_2$. The maps $v_k$ also glue and give us a $J_+$-holomorphic map $v : \Sigma^{\text{nonsing}}_\infty \to \dot{S}_+$.
By using the removable singularity theorem we then get $\Sigma_\infty^{\text{nonsing}}$ extends
to a nodal Riemann surface $\Sigma_\infty$ containing $\Sigma_\infty^{\text{nonsing}}$
and a $J_+$-holomorphic map $u_\infty : \Sigma_\infty \to \dot{S}_+$
extending $v$. Such a map is proper.

\bigskip
{\it Step 2}:
We now need to prove our energy bounds.
The bound $\int_{\Sigma_\infty} u_\infty^* \omega_S \leq E$ follows from the fact that $J_i$ is compatible with $\omega_i$ and the fact that $\int_{\P^1} u_i^* \omega_i \leq E$ for all $i \in \N$.
We now need to show $E_{S_+}(u_\infty) < \infty$.
Let $C_{(a,b)} := (a,b) \times C \subset (0,\epsilon_h) \times C \subset S$
where $0<a<b<\epsilon_h$ and let $C_a := \{a\} \times C$.
Let $\Sigma^i_{(a,b)} := u_i^{-1}(C_{(a,b)})$ where $0 < a <b < \epsilon_h$ and $i \in \N \cup \{\infty\}$.
We define $\Sigma^i_a := u_i^{-1}(C_a)$. 

From now on $i \in \N$ (i.e. $i \neq \infty$).
We let $r_C$ be the variable parameterizing $(-\epsilon_h,\epsilon_h)$ in $(-\epsilon_h,\epsilon_h) \times C$.
Let $\delta>0$ be small.
Let $\psi_\delta : (-\epsilon_h,\epsilon_h) \to \R$ be a smooth function whose derivative is non-negative and so that $\psi(x) = x$ for $x \notin (-\delta,\delta)$. Then we define $\omega_{i,\psi_\delta}$ to be equal to $\omega_i$ outside $(-\epsilon_h,\epsilon_h) \times C$
and equal to $\omega_h + d(\psi_\delta(r_C)\alpha_h)$ inside $(-\epsilon_h,\epsilon_h) \times C$.
This is in the same cohomology class as $\omega_i$.
For $i$ large enough and for $\delta>0$ small enough we then have that $\omega_{i,\psi_\delta}(V,J_iV) \geq 0$ for all vectors $V$.
The reason why this is true is because the restriction of $\omega_i$ and $J_i$ to the bundle $\text{ker}(dr_C) \cap \text{ker}(\alpha_h)$ is uniformly convergent to some symplectic form with compatible almost complex structure on this bundle in the region
$(-\frac{\epsilon_h}{2},\frac{\epsilon_h}{2}) \times C$.

From now on we fix $\delta >0$ small enough and pass to a subsequence of $u_i$'s to ensure the semi-positivity condition above holds for all $i \in\N$.
Then $\int_{\Sigma^i_{(a,b)}} u_i^* \omega_{i,\psi_\delta} \geq 0$ and $\int_{\P^1 \setminus \Sigma^i_{(a,b)}} u_i^* \omega_{i,\psi_\delta
} \geq 0$ for all $0 \leq a \leq b \leq \delta$, and so combined with $0 \leq \int_{\P^1} u_i^* \omega_{i,\psi_\delta} \leq E$
we get $0 \leq \int_{\Sigma^i_{(a,b)}} u_i^* \omega_{i,\psi_\delta} \leq E$.

Now let $\Psi_1 : (-\epsilon_h,\epsilon_h) \to \R$ be a smooth function with $\Psi_1(x) = x$ for $x \notin (-\delta,\delta)$, $\Psi_1' \geq 0$
and $\Psi_1(x) = 0$ for $x \in (\frac{\delta}{2},\frac{\delta}{2})$.
We then have $0 \leq \int_{\Sigma^i_{(a,\delta/2)}} u_i^* \omega_{i,\Psi_1} \leq E$ for all $0 \leq a \leq \delta/2$ and so $0 \leq \int_{\Sigma^i_{(a,\delta/2)}} u_i^* \omega_h \leq E$.
In a similar way, by considering the function $\Psi_2 : (-\epsilon_h,\epsilon_h) \to \R$  with $\Psi_2(x) = x$ for $x \notin (-\delta,\delta)$, $\Psi_2' \geq 0$
and $\Psi_2(x) = \frac{\delta}{2}$ for $x \in (\frac{\delta}{2},\frac{\delta}{2})$ we can show $0 \leq \int_{\Sigma^i_{(a,\delta/2)}} u_i^* (\omega_h + \frac{\delta}{2} d\alpha_h) \leq E$ for all $0 \leq a \leq \delta/2$.
The inequalities $0 \leq \int_{\Sigma^i_{(a,\delta/2)}} u_i^* \omega_h \leq E$ and $0 \leq \int_{\Sigma^i_{(a,\delta/2)}} u_i^* (\omega_h + \frac{\delta}{2} d\alpha_h) \leq E$ then tell us that $\Big| \int_{\Sigma^i_{(a,\delta/2)}} d\alpha_h)\Big| \leq \frac{2E}{\delta}$ for all $a$.

Because the curves $u_i$ converge to $u_\infty$ uniformly near $C_{\delta/2}$ we can assume that
$\Big|\int_{\Sigma^i_{\delta/2 }} u_i^* \alpha_h\Big| \leq K_1$ for some $K_1>0$ independent of $i$.
By Stokes' theorem we have
\[\int_{\Sigma^i_{(a,\delta/2)}} u_i^*d(\alpha_h) =
 \int_{\Sigma^i_{(\delta/2)}} u_i^*\alpha_h - \int_{\Sigma^i_a} u_i^*\alpha_h
\]
Hence we get $\Big|\int_{\Sigma^i_a} u_i^*\alpha_h\Big| \leq \frac{E}{\delta} + K_1$ for all $i \in \N$.
This implies that $\Big|\int_{\Sigma^\infty_a} u_\infty^*\alpha_h\Big| \leq \frac{E}{\delta} + K_1$ which in turn implies that
$E_{S_+}(u_\infty)$ is finite.

\bigskip
{\it Step 3}: We will now prove the remaining parts of this Proposition.
The set $\kappa \setminus (\dot{S}_- \cup ([0,\delta) \times C))$ is a compact subset of $\dot{S}_+$.
Now the $u_i$'s converge in the Gromov sense to some holomorphic curve near this compact subset whose image is contained in the image of $u_\infty$. 
So for $i$ large enough, we have that the image of $u_i$ is contained in $\dot{S}_- \cup ((-\epsilon_h,\delta) \times C) \cup N$.

Suppose $Q$  is a properly embedded codimension $2$ submanifold $Q \subset \dot{S}_+$ whose closure in $S$ disjoint from $C$.
Then each $u_i$ for $i \gg 1$ has a positive intersection number with $Q$ if and only if $u_\infty$ also has a positive intersection number
due to the fact that $u_i$ converges in the Gromov sense to a curve whose image is equal to the image of $u_\infty$ near $\kappa \cap Q$.

Now suppose $(u_i)_*([\P^1]) \cdot [Q] = 1$ for all $i$.
Then if (for a contradiction), $u_\infty$ has intersection greater than $1$ with $Q$ then
$u_i$ also has an intersection number greater than $1$ because $u_i$ converges in the Gromov sense to a $J_+$
holomorphic curve whose image is some branched cover of $u_\infty$ near $\kappa \cap Q$.
Hence the intersection number of $u_\infty$ with $Q$ is also $1$.
\qed
 
\bigskip

Now we need a result telling us how a holomorphic curve inside the interior of
a cobordism of stable Hamiltonian structures behaves near the boundary.
Again this result is weaker than the result in \cite{BEHWZ:compactnessfieldtheory}
except that we allow possibly degenerate stable Hamiltonian structures.
%Let $\mathcal{F}$ be a set of Reeb orbits in a stable Hamiltonian manifold of period $l$.
%We say such a set is an {\bf isolated family} if the set of Reeb orbits $\mathcal{F}$ is path connected
%and there is an $\epsilon_{\mathcal{F}}>0$ so that any orbit starting near
%$\mathcal{F}$ with period in $(l-\epsilon_{\mathcal{F}}, l+\epsilon_{\mathcal{F}})$
%is in $\mathcal{F}$.

\begin{lemma} \label{lemma:convergence}
Suppose we have a symplectic cobordism of stable Hamiltonian structures $(M,\omega_M)$
with $\partial_+ M= \emptyset$ and let $J$ be an almost complex structure
compatible with the completion of $M$.
%Suppose that every Reeb orbit $R$ is contained inside some isolated family
%$\mathcal{F}_R$ and 
Let $u : \Sigma \to M$ be a nodal, non-compact, proper $J$-holomorphic curve with finite energy.
Then for any sequence $\sigma_i \in \Sigma$ with $u(\sigma_i)$ getting
arbitrarily close to
$\partial_- M$ we have that
after passing to a subsequence, $u(\sigma_{i})$ converges to a point on some Reeb orbit $\gamma$ of $(\partial_- M,\omega_h^-,\alpha^-_h)$.
%If every Reeb orbit of $\partial_+ M$ is isolated 
We also show that if $\Sigma$ is smooth and connected of genus $0$  then $\Sigma$ is biholomorphic to $\P^1$ minus a finite number of points.

Let $Q$ be a properly embedded $J$-holomorphic hypersurface $Q$ of $M$ not intersecting $\partial_- M$ and let $\theta_M$ be a $1$-form on $M \setminus Q$  satisfying $d\theta_M =  \omega_M|_{M \setminus Q}$ and $i_X (\theta_M|_{\partial_- M}) = 1$ where $X$ is the Reeb vector field on $\partial_- M$.
Let $\eta$ be the wrapping number of $\theta_M$ around $Q$
and let $[Q] \cdot u$ be the intersection number between $u$ and $Q$.
Then the period of $\gamma$ is bounded above by $-([Q] \cdot u)\eta - \int_{\Sigma} u^* \omega_M$.

\end{lemma}
There are similar results when $\partial_+ M$ is non-empty but we decided to omit this so the statement
and the proof become less cluttered.
Note that if the Reeb orbit $\gamma$ from this Lemma is non-degenerate then by
\cite[Lemma 5.1]{BEHWZ:compactnessfieldtheory} there is some
holomorphic subset $(-\infty,1] \times S^1 \subset \Sigma$
with coordinates $(s,t)$
so that this region extends continuously to $[-\infty,1] \times S^1$
with $u(-\infty,t) = \gamma(Tt)$ where $T$ is the period of our Reeb orbit.

\proof of Lemma \ref{lemma:convergence}.
We prove this in $3$ steps. In Step $1$ we find our Reeb orbit.
In Step $2$ we construct our bound for the period of our Reeb orbit $\gamma$.
In Step $3$ we show that the domain of our curve is $\P^1$ minus a finite number of points.

\bigskip
{\it Step $1$}:
Because $\partial_- M$ is compact we can, after passing to a subsequence, ensure that $u(\sigma_{i})$ converges to some
point $x \in \partial_- M$.
Because $J$ is compatible with the completion of $M$ we have smooth embeddings
$\Phi_r : (0,\epsilon_h) \times \partial_- M \to \R \times \partial_- M$ defined by a pair
$(\tau_r \circ \phi_-, \text{id}_{\partial_- M})$.
Here $\phi_- : (0,\epsilon_h) \to (-\infty,0)$ is a diffeomorphism,
$\tau_r$ is the translation map sending $x \in \R$ to $x+r$.
These have the property that there is a cylindrical almost complex structure
$\widehat{J}_-$ on $\R \times \partial_- M$ so that
$(\Phi_r)_* J$ converges in $C^\infty_{\text{loc}}$ to $\widehat{J}_-$
as $r$ tends to $\infty$ for some $J_-$ compatible with the stable Hamiltonian structure on $\partial_-M$.
After passing to a subsequence again we will assume that $u(\sigma_{i}) = (\phi_-^{-1}(x_i),y_i) \in (0,\epsilon_h) \times \partial_- M$ for a sequence
$x_i \in (-\infty,0)$, $y_i \in \partial_- M$ where $x_i \to -\infty$ and $y_i \to y_\infty \in \partial_- M$.

Fix a small constant $\Delta>0$ and consider the interval $I_i := [x_i-\Delta, x_i+\Delta] \times \partial_- M$.
Let $T_i : I_i \to [-\Delta,\Delta] \times \partial_- M$ be the natural translation map
and define $J_i$ to be the pushforward $((T_i)_*(\Phi_0)_* J)|_{[-\Delta,\Delta] \times \partial_- M}$.
These almost complex structures $C^\infty$ converge to $\widehat{J}_-$.
If $\pi_{\partial_- M} : \R \times \partial_- M \to \partial_- M$ is the natural projection map
then by abuse of notation we write $\alpha_- = \pi_{\partial_- M}^* \alpha_-$
and $\omega_h^- = \pi_{\partial_- M}^* \omega_h^-$ where $(\omega_h^-,\alpha_-)$ is the stable Hamiltonian structure on $\partial_- M$.
We define $\Omega_i := ((T_i)_* (\Phi_0)_* \omega_S) + e^r dr \wedge \alpha_-$ where $r$ is the coordinate parameterizing the $\R$ factor in $\R \times \partial_- M$.
Here $\Omega_i$ may not be closed but it is compatible with $J_i$.
We have that $(\Omega_i,J_i)$ $C^\infty$ converges to an almost Hermitian structure
$(\Omega_\infty,\widehat{J}_-)$ where $\Omega_\infty = \omega_h^- + e^r dr \wedge \alpha_-$.
Define $\Sigma_i := u^{-1}(\phi_-^{-1}(I_i))$.
Because $E_{\dot{M}}(u)$ is bounded we then get that $\int_{\Sigma_i} u^*\Phi_0^*T_i^* \Omega_i$ is bounded
by a constant independent of $i$.
Hence by Theorem \ref{theorem:fishcompactnesscorollary}
we have after passing to a subsequence that $T_{i} \circ \Phi_0 \circ u_{i}$ 
converges in the Gromov sense to a $\widehat{J}_-$-holomorphic curve
$u_\infty : \Sigma_\infty \to [-\Delta,\Delta] \times \partial_- M$
near
$[-\frac{\Delta}{2},\frac{\Delta}{2}] \times \partial_- M$.

Let $r_C$ be the coordinate parameterizing $(0,\epsilon_h)$ in
$(0,\epsilon_h) \times \partial_- M$.
%Let $\Sigma' := u^{-1}((0,\epsilon_h) %\times \partial_- M)$.
%Because $\int_{\Sigma_\infty} u^* %\omega_M$ is finite and because
%$dr_C \wedge \alpha_-(V,JV) \geq 0$ %for all vectors $V$
%we have that $\int_{\Sigma_\infty} %u^*_\infty (\omega_h^- + r_Cd\alpha_-)$ %is bounded.
For $r_C$ small enough we get that
$r_C d\alpha_- (V,JV) \leq B_1 r_C \|V\|^2$ for any non-zero vector tangent to $\text{ker}(\alpha_h)$
where $\|\cdot\|$ is the natural metric for some constant $B_1>0$.
%
%But for any such vector $V$, $\omega_h^-(V,JV) \geq B_2 \|V\|^2$ for some constant $B_2>0$
%if $r_C$ is small enough.
%
%These facts imply that $\int_{\Sigma_\infty} u_\infty^* \omega_h^-$ is finite.
Using this fact combined with
the fact that $J$ is compatible
with our stable Hamiltonian
cobordism
and the fact that
$\int_\Sigma u^* \omega$ is finite,
we have that $\int_{\Sigma_i} u^*\Phi_0^*T_i^* \omega_h^-$ tends to $0$.
Hence because  $u_{i}$ Gromov converges to $u_\infty$
we then get that $\int_{\Sigma_\infty} u_\infty^* \omega_h^- = 0$.
Because $\widehat{J}_-$ is cylindrical this means that the projection of the image of $u_\infty$
under $\pi_{\partial_- M}$ is contained in a Reeb flowline of $(\omega_h^-,\alpha_-)$.
For generic small $\delta>0$, the intersection $u_\infty(\Sigma_\infty) \cap \{\delta\} \times \partial_- M$
is transverse and hence
is a one dimensional manifold $A$ with the property that $\pi_{\partial_- M}(A)$ is contained in a Reeb flowline.
Also the tangent space of $A$ maps to a non-trivial $1$-dimensional subspace of $T\partial_- M$
under the map $\pi_{\partial_- M}$.
Hence $\pi_{\partial_- M}(A)$ is a union of Reeb orbits.
Hence $u_{i}(\sigma_{i})$ converges to a point in $\pi_{\partial_- M}(A)$ which contains some Reeb orbit $\gamma$.

\bigskip
{\it Step $2$}:
Now let $Q$ be as in the statement of this Lemma. We will now find an upper bound for the period of our Reeb orbit $\gamma$.
Let $A_i := (T_i \circ \Phi_0 \circ u_i)^{-1}(\{\delta\} \times \partial_- M)$ for some generic small $\delta>0$ and for all $i \in \N$.
Define $A_\infty := u_\infty^{-1}(\{\delta\} \times \partial_- M)$.
We have that
$\int_{A_i} u_i^* \Phi_0^* T_i^* \alpha_-$
converges to $\int_{A_i} u_i^* \theta_M$ as $i$ tends to infinity.
Let $q_1,\cdots,q_l$ be the set of points in $\Sigma$ intersecting $Q$
and let $l_1,\cdots,l_i$ be small loops around these points oriented positively.
Now if the intersection multiplicity between $u$ and $Q$ at $q_j$
is $\iota_j$ then $\int_{l_j} u^* \theta_M$ converges to $\iota_j \eta$
as $l_j$ gets smaller and smaller.
Hence $\sum_j \int_{l_j} u^* \theta_M$ converges to $([Q] . u) \eta$ as the loops $l_j$ shrink.
Stokes' theorem tells us that
\[-\sum_j \int_{l_j} u^* \theta_M - \int_{A_i} u_i^* \theta_M = \int_{\breve{\Sigma}} \omega_M\]
where $\breve{\Sigma} \subset \Sigma$ is the surface which is bounded by the loops
$l_j$ and $A_i$.
As $i$ gets large and the the loops $l_j$ get smaller,
the $\omega_M$ area of $\breve{\Sigma}$ converges to the $\omega_M$ area of $\Sigma$.
Hence as $i$ tends to infinity and as the loops $l_j$ get smaller in length around $q_i$,
we get $\int_{A_i}  u_i^* \theta_M$ converges to some limit $L \leq -([Q] . u) \eta - \int_{\Sigma} \omega_M$.
Also $\int_{A_i} u_i^* \theta_M$ converges to
$\int_{A_\infty} u_\infty^* \pi_{\partial_- M}^* (\theta_M|_{\partial_- M})$.
Here $\int_{A_\infty} u_\infty^* \pi_{\partial_- M}^* (\theta_M|_{\partial_- M}) = \int_{A_\infty} u_\infty^* \alpha_-$ because $i_R \alpha_- = i_R \theta_M$
(here $R$ is the Reeb vector field)
and so putting everything together we get $\int_{A_\infty} u_\infty^* \alpha_- \leq -([Q]. u) \eta - \int_{\Sigma} \omega_M$.
Because $\pi_{\partial_- M}(A_\infty)$ is a union of Reeb orbits we then get that the sum of their periods
is bounded above by $-(([Q] \cdot u)\eta - \int_{\Sigma} \omega_M$ and in particular this is an upper bound for the period of $\gamma$.

\bigskip
{\it Step $3$}:
We now need to show that  $\Sigma$ is biholomorphic to $\P^1$ minus a finite number of points.
We will first show that there exists some $\delta_\Sigma >0$ and an increasing sequence of compact codimension $0$ submanifolds $\kappa_i$
whose union is $\Sigma$ so that
we can holomorphically embed the annulus $[-\delta_\Sigma,\delta_\Sigma] \times S^1$ into each connected component of
$\Sigma \setminus \kappa_i$ and so that $\partial \kappa_i$ has a bounded number of connected components independent of $i$.
Standard Riemann surface theory then implies that $\Sigma$ is biholomorphic to $\P^1$ minus finitely many points.
The point here is that the bounded number of connected components condition combined with the fact that we are in genus $0$
implies that the topology of $\Sigma$ is bounded and hence can be embedded as an open subset of $\P^1$.
The annulus embedding condition then implies that the complement of $\Sigma$ in $\P^1$ has no accumulation points.

Let $\breve{I}_i := [-\Delta-i,\Delta-i] \times \partial_- M \subset \R \times \partial_- M$
and let $\breve{T}_i : [-\Delta-i,\Delta-i] \times \partial_- M \to [-\Delta,\Delta] \times \partial_- M$ be the natural translation map.
Let $\breve{\Sigma}_i := ((\breve{T}_i \circ \Phi_0 \circ u)^{-1})([-\Delta,\Delta] \times \partial_- M)$.
Now a similar Gromov compactness argument as above tells us that after passing to a subsequence we have
that $(\breve{T}_{i} \circ \Phi_0 \circ u|_{\breve{\Sigma}_i})$ converges in the Gromov
sense to some $\widehat{J}_-$-holomorphic curve $v : \breve{\Sigma}_\infty \to [-\Delta,\Delta] \times \partial_- M$ near
$[-\frac{\Delta}{2},\frac{\Delta}{2}] \times \partial_- M$.
Now the annulus $[-\delta_\Sigma,\delta_\Sigma] \times \partial_- M$ holomorphically embeds into
$v^{-1}((-\frac{\Delta}{4},\frac{\Delta}{4}) \times \partial_- M)$ for some $\delta_\Sigma > 0$.
By the nature of Gromov convergence this means such an annulus holomorphically embeds into the interior of
$\breve{\Sigma}_{i}$ for $i$ large enough.
This means that such an annulus is embedded into $\Sigma$ minus the compact subset
$\kappa_i:= (\Phi_0 \circ u)^{-1}([\frac{\Delta}{2} - i,0)) \cup u^{-1}(M \setminus ((0,\epsilon_h) \times \partial_- M))$.
We can choose $\Delta$ generically so that $(\breve{T}_{i} \circ \Phi_0 \circ u)^{-1}(\{\frac{\Delta}{2}\} \times \partial_- M)$ is a submanifold of $\Sigma$.
Such a submanifold is equal $\partial \kappa_i$.
Again by Gromov convergence this sequence of submanifolds converges to some compact submanifold of $\breve{\Sigma}_\infty$
hence the number of connected components of $\partial \kappa_i$ is bounded.
Hence $\Sigma$ is biholomorphic to $\P^1$ minus finitely many points.

\qed

\section{Appendix B: A Maximum Principle}

Let $(M,\omega_M)$ be a stable Hamiltonian cobordism where $\partial_- M$ is compact and has a stable Hamiltonian structure
$(\omega_h,\alpha_h)$ and where $\partial_+ M = \emptyset$. Suppose we have a $1$-form $\theta$ on $M$ satisfying:
\begin{enumerate}
\item $d\theta = \omega_M$. In particular $d\theta|_{\partial_- M} = \omega_h$.
\item We have $X_\theta$ points inwards along $\partial_- M$.
\item $\theta(R) = \alpha_h(R)$ where $R$ is the Reeb vector field of $(\omega_h,\alpha_h)$.
\end{enumerate}
Let $\dot{M} = M \setminus \partial_- M$.

\begin{prop} \label{proposition:maximumprinciple}
Let $J$ be an almost complex structure  on $\dot{M}$ compatible with the completion of $M$.
There are no finite energy properly embedded nodal $J$-holomorphic curves $u : \Sigma \to \dot{M}$
such that the closure of the image of $u$ in $M$ is compact.
\end{prop}

The proof of this proposition borrows its main idea from \cite[Lemma 7.2]{SeidelAbouzaid:viterbo}.

\proof of Proposition \ref{proposition:maximumprinciple}.
%Let $u : \Sigma \to M$ be a properly embedded $J$-holomorphic curve with no positive ends.
We will create a contradiction by showing that $\int_{\Sigma} u^* \omega_M = 0$.
A neighborhood of $\partial_- M$ is symplectomorphic to $[0,\epsilon_h) \times \partial_- M$
with $\omega_M = \omega_h + r_M d\alpha_h + dr_M \wedge \alpha_h$
where $r_M$ parameterizes $[0,\epsilon_h)$.
We have a diffeomorphism $\phi : (0,\epsilon_h) \to (-\infty,0)$ and a cylindrical almost complex structure
$\widehat{J}_h$ on $\R \times \partial_- M$ so that
 $(\tau_r \circ \phi, \text{id})_* J$ converges in $C^\infty_{\text{loc}}$ to $\widehat{J}_h$
in $\R \times \partial_- M$ as $r \to +\infty$ where $\tau_r : \R \to \R$ sends $x$ to $x+r$.

Now consider a sequence $r_i <-1$ tending to $-\infty$
and let $I_i$ be the region $\phi^{-1}( [r_i-1,r_i+1] ) \times \partial_- M$ near $\partial_- M$.
Let $\Sigma_i := u^{-1}(I_i)$ where we choose $r_i$ generically so $\Sigma_i$ is a manifold with boundary.
We let $u_i : \Sigma_i \to I_0 \times \partial_- M$ be the sequence of maps
$u|_{\Sigma_i}$ composed with the natural translation map from
$I_i$ to $I_0$.
Let $\breve{I}_i := \phi^{-1}( [r_i-\frac{1}{2},r_i+\frac{1}{2}] ) \times \partial_- M$.

By Theorem \ref{theorem:fishcompactnesscorollary} we get after passing to a subsequence
a genus $0$ nodal $\widehat{J}_h$-holomorphic curve $u : \Sigma_\infty \to I_0 \times \partial_- M$
and genus zero compact curve  $\widetilde{\Sigma}$ with boundary, smooth embeddings
$v_i : \widetilde{\Sigma} \hookrightarrow \Sigma_i$ and a continuous surjection
$v_\infty : \widetilde{\Sigma} \twoheadrightarrow \Sigma_\infty$ such that:
\begin{enumerate}
\item Let $\widetilde{\Sigma}^{\text{ns}} := \Sigma_\infty \setminus \Sigma_\infty^{\text{sing}}$
where $\Sigma_\infty^{\text{sing}}$ is the preimage under $v_\infty$ of the nodal points of $\Sigma_\infty$.
Then $v_\infty^{\text{nonsing}} := v_\infty|_{\widetilde{\Sigma}^{\text{ns}}}$ is a diffeomorphism onto its image.
\item We have $u_i \circ v_i$ $C^0$ converges to $u_\infty \circ v_\infty$
and $u_i \circ v_i|_{\widetilde{\Sigma}^{\text{ns}}}$ $C^\infty$ converges to $u_\infty \circ v_\infty|_{\widetilde{\Sigma}^{\text{ns}}}$.
\item The $u_i \circ v_i$ and $u_\infty \circ v_\infty$ map the boundary of $\widetilde{\Sigma}$
outside $\breve{I}_0$.
\end{enumerate}

Let $B_\infty := (u_\infty \circ v_\infty)^{-1}(\{\delta\} \times \partial_- M)$ where $\delta$ is small and generic making $B_\infty$
into a manifold. This has a natural orientation coming from the outward normal of the
complex surface $ (u_\infty \circ v_\infty)^{-1}([\delta,\delta + \delta_1))$
for some even smaller $\delta_1>0$.
By Stokes' theorem we have
$\int_{\Sigma} u^* \omega_M = \int_{B_\infty} (u_\infty \circ v_\infty)^* \pi_{\partial_M}^*( \theta|_{\partial_- M})$
where $\pi_{\partial_M}$ is the natural projection to $\partial_- M$.
Because $\int_{\Sigma_\infty} u_\infty^* \omega_h  = 0$
we have that the tangent space $T\Sigma_\infty$ at each non-singular point gets mapped via
$(u_\infty \circ v_\infty)$ to the kernel of $\omega_h$.
Because $\theta$ restricted to the kernel of $\omega_h$ is equal to $\alpha_h$
restricted to the kernel of $\omega_h$ we then get:
\[\int_{\Sigma} u^* \omega_M = \int_{B_\infty} (u_\infty \circ v_\infty)^* \alpha_h.\]
This is equal to:
\[\int_{B_\infty} \alpha_h \circ d(u_\infty \circ v_\infty) =
\int_{B_\infty} \alpha_h \circ \widehat{J}_h \circ d(u_\infty \circ v_\infty) \circ j\]
where $j$ is the pullback of the complex structure on $\Sigma_\infty$ via $v_\infty$.
Now $\alpha_h \circ \widehat{J}_h = dr_M$ and $j(\zeta)$ points outwards along the surface
$(u_\infty \circ v_\infty)^{-1}([\delta,\delta+\delta_1) \times \partial_- M)$ where $\zeta$ is tangent to $B_\infty$ and respecting its orientation.
This implies that $\int_{\Sigma} u^* \omega_M$ is non-positive.
This gives us our contradiction.
\qed

\bibliography{references}

\end{document}